\documentclass[12pt,oneside]{article}

\usepackage{epsfig,lscape}
\usepackage{amssymb,amsfonts,amsmath}
\usepackage{rotating}

\usepackage{amsthm}

\usepackage{color}
\usepackage[dvipsnames]{xcolor}

\usepackage{float}
\usepackage[hypertexnames=false]{hyperref}

\def\dif{\mathrm d}

\def\N{\mathrm{N}}

\def\bbR{{\mathbb{R}}}
\def\argmin{\mathrm{argmin}}

\newcommand\mytext[1]{\text{\scriptsize{#1}}}

\def\one{\mathbf{I}}
\def\two{\mathbf{II}}
\def\three{\mathbf{III}}
\def\four{\mathbf{IV}}
\def\five{\mathbf{V}}

\def\six{\mathbf{VI}}
\def\seven{\mathbf{VII}}

\def\overL{\overline{L}}
\def\overA{\overline{A}}
\def\overl{\overline{l}}
\def\overa{\overline{a}}
\def\over1{\overline{1}}
\def\overlam{\overline{\lambda}}
\def\overq{\overline{q}}

\def\overLam{\overline{\Lambda}}

\def\underL{\underline{L}}

\def\underl{\underline{l}}
\def\underlam{\underline{\lambda}}
\def\underq{\underline{q}}

\def\calE{{\mathcal{E}}}

\def\myZ{{_Z}}

\newenvironment{prf}
{\noindent \textbf{Proof.}}{\hfill $\Box$ \vspace{.1in}}

\newtheorem{lem}{Lemma}
\newtheorem{pro}{Proposition}
\newtheorem{cor}{Corollary}

\theoremstyle{definition}

\theoremstyle{definition}

\allowdisplaybreaks

\usepackage[medium,compact]{titlesec}

\titlespacing*{\section} {0pt}{1.5ex}{1ex}
\titlespacing*{\subsection} {0pt}{1.5ex}{1ex}
\titlespacing*{\subsubsection} {0pt}{1ex}{1ex}

\begin{document}

\setlength{\abovedisplayskip}{5pt}
\setlength{\belowdisplayskip}{5pt}


\begin{titlepage}

\begin{center}
{\Large Sensitivity models and bounds under sequential unmeasured confounding in longitudinal studies}

\vspace{.1in} Zhiqiang Tan\footnotemark[1]

\vspace{.1in}
\today
\end{center}

\footnotetext[1]{Department of Statistics, Rutgers University. Address: 110 Frelinghuysen Road,
Piscataway, NJ 08854. E-mail: ztan@stat.rutgers.edu.
The author thanks Matteo Bonvini and James Robins for helpful discussions.}

\paragraph{Abstract.}
Consider sensitivity analysis to assess the worst-case possible values of counterfactual outcome means and
average treatment effects under sequential unmeasured confounding in a longitudinal study with time-varying treatments and covariates.
We formulate several multi-period sensitivity models to relax the corresponding versions of the assumption of sequential unconfounding.
The primary sensitivity model involves only counterfactual outcomes, whereas the joint and product sensitivity models
involve both counterfactual covariates and outcomes. We establish and compare explicit representations
for the sharp and conservative bounds 
at the population level through convex optimization, depending only on the observed data.
These results provide for the first time a satisfactory generalization from the marginal sensitivity model in the cross-sectional setting.

\paragraph{Key words and phrases.}
Marginal sensitivity model; Sensitivity analysis; Sequential unconfounding; Sequential unmeasured confounding; Time-varying treatments.

\end{titlepage}

\section{Introduction} \label{sec:intro}

Drawing inferences about average treatment effects (ATEs) from longitudinal studies is interesting but challenging
from both methodological and practical perspectives.
Typically, the first line of attack is to invoke the assumption of sequential unconfounding or exchangeability
such that the ATEs can be point identified at the population level
and hence consistently estimated from sample data under additional modeling assumptions (Hernan \& Robins 2022).
However, the assumption of sequential unconfounding is untestable from observed data and
may potentially violated. It is important to examine how the estimates of ATEs might change under
sequential unmeasured confounding, which is known as sensitivity analysis.

The main contribution of this paper is to formulate several multi-period sensitivity models
and derive explicit representations for the sharp and conservative bounds on the counterfactual outcomes means and ATEs,
in the longitudinal setting with time-varying treatments and covariates.
These results generalize those for sensitivity analysis
based on the marginal sensitivity model in the cross-sectional setting (Tan 2006, 2022; Dorn \& Guo 2022).

We adopt the existing counterfactual framework for longitudinal causal inference (Hernan \& Robins 2022).
The observed covariates and outcome and the treatment variables, referred to as the observed data,
are determined from the counterfactual covariates and outcome and the treatment variables, referred to as the full data,
through a consistency assumption.
Each of our sensitivity models restricts the distribution of the full data to allow
sequential unmeasured confounding up to a certain degree controlled by the sensitivity parameters, which can be varied over a range of
pre-specified values.
The sharp upper (or lower) bound on the mean of a counterfactual outcome is defined as
the maximum (or minimum) value over all possible distributions on the full data
which satisfy the sensitivity model and are fully compatible with the distribution of the observed data.
As the main technical results, we establish and compare explicit representations for the sharp and conservative bounds
at the population level through convex optimization, depending only on the observed-data distribution even though the sharp bounds are
originally defined from the full-data distribution.
These sensitivity bounds can be seen to capture some intrinsic property of the observed data
in terms of the worst-case values of the counterfactual outcome mean.

\textbf{Related work.}\; There has been an extensive literature on sensitivity analysis for cross-sectional studies.
Examples include the latent sensitivity model (Rosenbaum 2002) and recent work (Yadlowsky et al.~2022),
the selection odds model (Robins et al.~2000) and recent work (Franks et al.~2020; Scharfstein et al.~2021),
and the marginal sensitivity model (Tan 2006) and recent work (Zhao et al.~2019; Dorn \& Guo 2022; Dorn et al.~2021; Tan 2022).
Additional examples can also be found from these references.

Sensitivity analysis for longitudinal causal inference is more complex and
the related literature is more limited than for cross-sectional studies.
The selection odds model (Robins et al.~2000) requires parametric specifications
of the sensitivity ratios, for example, defined in (\ref{eq:lam-def-odds}).
By comparison, only the range of the sensitivity ratios are constrained in our sensitivity models,
and then sensitivity bounds are derived such that the implied sensitivity ratios can be interpreted (see for example Proposition~\ref{pro:upper-2period}).
Among recent examples, Kallus \& Zhou (2020) considered in a longitudinal setting
a sensitivity model similar to the marginal sensitivity model, but explicitly excluded the existence of time-varying unmeasured confounders.
Bovini et al.~(2022) also used a propensity sensitivity model related to the marginal sensitivity model, but
did not account for all inherent constraints, similar to for example
the normalization constraints (\ref{eq:constr-2period}) in our sensitivity analysis.

\textbf{Organization.}\;
The remainder of the paper is organized as follows. In Section~\ref{sec:setup}, we describe the existing
counterfactual framework and the assumptions of sequential unconfounding.
In Sections~\ref{sec:2period} and \ref{sec:Kperiod}, we develop our sensitivity models and bounds with a single treatment strategy
for $2$ periods and general $K \ge 2$ periods respectively.
Then we present a numerical study in Section~\ref{sec:numerical} and further discussion about multiple treatment strategies
in Section~\ref{sec:further-discussion}.
All technical details are provided in the Supplement Material.

\section{Setup and sequential unconfounding}  \label{sec:setup}

Consider a longitudinal setup with time-dependent treatments and covariates over $K$ periods (Hernan \& Robins 2022, Chapter 19).
For $k=0,\ldots,K-1$, let $A_k$ be a binary treatment at time $k$, and
$L_k$ be a covarite vector measured prior to treatment $A_k$ at time $k$.
Let $Y$ be a terminal outcome measured after treatment $A_{K-1}$.
We use an overline to indicate the history, for example, $\overline{L}_k = (L_0, \ldots, L_k)$ and
$\overline{A}_k = (A_0, \ldots, A_k)$,
and an underline to indicate the follow-up, for example, $\underline{L}_k = (L_k,\ldots, L_{K-1})$ and
$\underline{A}_k = (A_k,\ldots, A_{K-1})$.

For a static treatment strategy $\overline{a} = (a_0, \ldots, a_{K-1})$ with $a_0,\ldots,a_{K-1}\in\{0,1\}$,
the counterfactual (or potential) outcome $Y ^{\overline{a}}$ is defined as the outcome that would be observed
under the treatment  strategy $\overline{a}$.
There are in general $2^K$ possible static treatment strategies, leading to $2^K$ counterfactual outcomes.
The average treatment effects are defined as the contrasts between the means of the associated counterfactual outcomes.
For example, the average treatment effect (ATE) of ``always treat'' strategy $\overline{1}_{K-1}=(1,\ldots,1)$ and
``never treat'' strategy $\overline{0}_{K-1} = (0,\ldots,0)$ is $E( Y^{\overline{1}_{K-1}} ) - E( Y^{\overline{0}_{K-1}} )$.
In the following, we discuss point identification of the mean $\mu^{\overline{a}} = E ( Y ^{\overline{a}})$,
while assuming certain conditional expectations are known.
Given sample data, the empirical version of an identification formula can be used for estimating  $\mu^{\overline{a}}$
after the conditional expectations are estimated.

For a static treatment strategy $\overline{a}$, the mean $\mu^{\overline{a}} = E( Y^{\overline{a}} ) $
can be point identified under the following assumptions (Robins 1986). \vspace{-.05in}
\begin{itemize} \addtolength{\itemsep}{-.1in}
\item[(A1)] Consistency:
$ Y = Y^{\overline{a}}$ if $\overline A_{K-1} = \overline{a}$.

\item[(A2)] Sequential unconfounding (or exchangeability): for $k=0,1,\ldots,K-1$,
\begin{align*}
 Y^{\overline{a}} \perp A_k \,|\,  \overline{A}_{k-1} = \overline{a}_{k-1}, \overline{L}_k,   
\end{align*}
where $\perp$ denotes conditional independence, and $\overline{A}_{-1}$ is set to the null.

\item[(A3)] Positivity:  for $k=0,1,\ldots,K-1$, it holds that almost surely,
$$ P( A_k =a_k \,|\,  \overline{A}_{k-1} = \overline{a}_{k-1}, \overline{L}_k ) >0 .$$
\end{itemize}\vspace{-.05in}
Assumptions A1--A3 are contingent on the specific treatment sequence $\overline{a}$, but this dependency is suppressed in the notation.
These assumptions may not be satisfied by all $2^K$ static treatment strategies, each being a sequence of $K$ binary treatments.

Assumptions A1--A2 are concerned with the terminal counterfactual outcome $Y^{\overline{a}}$ only,
but not any other counterfactual variables prior to time $K$.
There is an extended version of Assumptions A1--A2 as follows, where counterfactual covariates are introduced.
For a static treatment strategy $\overline{a}$, the counterfactual covariates
$\overline{L}_{K-1}^{\overline{a}} = (L_0^{\overline{a}}, L_1^{\overline{a}}, \ldots, L_{K-1}^{\overline{a}})$
are defined as time-dependent covariates that would be observed under strategy $\overline{a}$ from time $0$ to $K-1$,
where $L^{\overline{a}}_0$ is set to $L_0$, independently of $\overline{a}$.\vspace{-.05in}
\begin{itemize} \addtolength{\itemsep}{-.1in}
\item[(A1$^\dag$)] Consistency:
$ Y = Y^{\overline{a}}$ if $\overline A_{K-1} = \overline{a}$;
$\overline{L}_k = \overline{L}_k^{\overline{a}}$ if $\overline{A}_{k-1} = \overline{a}_{k-1}$ for $k=1,\ldots,K-1$.
\item[(A2$^\dag$)] Sequential unconfounding (or exchangeability): for $k=0,1,\ldots,K-1$,
\begin{align*}
 (\underline{L}_{k+1}^{\overline{a}},  Y^{\overline{a}} ) \perp A_k \,|\,  \overline{A}_{k-1} = \overline{a}_{k-1}, \overline{L}_k,   
\end{align*}
where $\underline{L}_{k+1}^{\overline{a}} = (L_{k+1}^{\overline{a}}, \ldots, L_{K-1}^{\overline{a}})$,
and $\overline{A}_{-1}$ and $\underline{L}_K^{\overline{a}}$ are set to the null.
\end{itemize}\vspace{-.05in}
For distinction, we refer to Assumption A2 as \textit{primary} sequential unconfounding (being relevant to $Y^{\overline{a}}$ only)
and A2$^\dag$ as \textit{joint} sequential unconfounding.
On one hand, Assumptions A1$^\dag$--A2$^\dag$ are more restrictive than A1--A2 on the counterfactual data,
even though both sets of assumptions place no restriction on the observed data $(\overline{A}_K, \overline{L}_K, Y)$.
On the other hand,
Assumptions A1$^\dag$, A2$^\dag$, and A3 are then sufficient for point identification of not only
the mean of $Y^{\overline{a}}$ (i.e., $\mu^{\overline{a}}$), but also the mean of any function of
$(\overline{L}_{K-1}^{\overline{a}}, Y^{\overline{a}})$.
Point identification of such general expectations is required for that of the counterfactual outcome mean
from a dynamic treatment strategy $\overline{a}= (a_0, \ldots, a_{K-1})$, where $a_k$
can be randomly chosen depending on $(A_{k-1}, \overline{L}_k )$ for $k=0,1,\ldots,K-1$.
While only static treatment strategies will be dealt with in this paper,
our interest in Assumptions A1$^\dag$--A2$^\dag$ is motivated by another consideration:
relaxing the sequential unconfounding represented by either Assumption A2 or A2$^\dag$ may lead to different sensitivity models
as discussed in Section \ref{sec:2period}.


There are various approaches in which $\mu^{\overline{a}} = E( Y^{\overline{a}} ) $ can be point identified under Assumptions A1--A3
or Assumptions A1$^\dag$, A2$^\dag$, and A3.
For illustration, we describe the identification formulas for $\mu^{1,1} = E( Y^{1,1})$ in the case of $K=2$ and $\overline{a} = (1,1)$.
The primary sequential unconfounding Assumption A2 becomes
\begin{align}
Y^{1,1} \perp A_1 \,|\, A_0=1, \overline{L}_1, \quad
 Y^{1,1}  \perp A_0 \,|\, L_0 .  \label{eq:unconfounded-2period}
\end{align}
The joint sequential unconfounding Assumption A2$^\dag$ becomes
\begin{align}
Y^{1,1} \perp A_1 \,|\, A_0=1, \overline{L}_1, \quad
  (L_1^{1,1}, Y^{1,1} ) \perp A_0 \,|\, L_0 ,  \label{eq:unconfounded-2period-Jnt}
\end{align}
which, for the convenience of later discussion, can be equivalently stated as
\begin{align}
Y^{1,1} \perp A_1 \,|\, A_0=1, \overline{L}_1, \quad
  L_1^{1,1} \perp A_0 \,|\, L_0 , \quad Y^{1,1} \perp A_0 \,|\, L_0, L_1^{1,1} , \label{eq:unconfounded-2period-Jnt2}
\end{align}
by a decomposition of the second condition in (\ref{eq:unconfounded-2period-Jnt}).
The iterative-conditional-expectation (ICE) formula (also called the $g$-formula) says that under Assumptions A1--A3,
\begin{align}
\mu^{1,1} = E \{ m_0^*(L_0)\}
= E \left\{ E ( Y | A_0= 1, A_1=1, \overline{L}_1 ) |  A_0= 1, L_0 \right\},
\label{eq:iden-ICE}
\end{align}
where $m_1^* ( \overline{L}_1 ) = E ( Y | A_0= 1, A_1=1, \overline{L}_1 )$ and
$m_0^*(L_0) = E \{ m_1^* ( \overline{L}_1 ) | A_0= 1, L_0 \}$.
The inverse-probability-weighting (IPW) formula says that under Assumptions A1--A3,
\begin{align}
\mu^{1,1} = E \left\{ \frac{ A_0 A_1 Y } { \pi^*_0 (L_0)  \pi^*_1(\overline{L}_1 ) } \right\}, \label{eq:iden-IPW}
\end{align}
where 
$\pi^*_1(\overline{L}_1 )  = P(A_1 = 1 | A_0=1, \overline{L}_1 ) $, and
$\pi^*_0(L_0) =  P(A_0 = 1| L_0)$.
The identification formulas (\ref{eq:iden-ICE}) and (\ref{eq:iden-IPW}) are singly robust,
requiring that $\{m_1^* ( \overline{L}_1 ), m_0^*(L_0) \}$ or respectively $\{ \pi^*_1(\overline{L}_1 ), \pi^*_0(L_0) \}$
be consistently estimated from sample data.
Based on (\ref{eq:iden-ICE}) and (\ref{eq:iden-IPW}), doubly robust (or multiply robust) identification formulas can also be derived
via augmented IPW techniques,
such that $\mu^{1,1}$ is point identified if any of several subsets of $\{ m^*_1(\overline{L}_1), m^*_0 (L_0)\}$
and $\{ \pi^*_1(\overline{L}_1), \pi^*_0 (L_0)\}$ are consistently estimated.

We make two additional remarks about (\ref{eq:iden-ICE}) and (\ref{eq:iden-IPW}).
For clarity, we refer to the right-hand side of (\ref{eq:iden-ICE}) or (\ref{eq:iden-IPW})
as the ICE or IPW functional.
First, for any function $b_1(l_0,l_1, y)$,
the mean $E \{ b_1 (L_0, L^{1,1}_1, Y^{1,1}) \}$ is identified by
the ICE or IPW functional with $Y$ replaced by $b_1 (\overL_1, Y)$
under Assumptions A1$^\dag$, A2$^\dag$, and A3. As mentioned earlier,
a similar extension is not feasible under Assumptions A1--A3:
the ICE or IPW functional does not in general identify $E \{ b_1 (L_0, L_1, Y^{1,1}) \}$.
Second, more importantly, the ICE and IPW functionals
with $Y$ replaced by any function $b_1 (\overL_1, Y)$
remain identical to each other and may differ from $E \{ b_1 (L_0, L^{1,1}_1, Y^{1,1}) \}$,
even if Assumption A2 or  A2$^\dag$ is violated as in our development of sensitivity analysis.
See Supplement Lemma~\ref{lem:IPW-ICE} for a direct proof.

\section{Sensitivity analysis for 2 periods} \label{sec:2period}

To develop sensitivity analysis, we consider as a latent-variable model the counterfactual framework in Section \ref{sec:setup}, while relaxing Assumption A2 or A2$^\dag$.
For $K=2$, the full data consist of $(\overA_1,\overL_1, Y^{1,1}, Y^{1,0}, Y^{0,1}, Y^{0,0})$
or $( \overA_1, L_0,\{ (L^{a_0a_1}_1, Y^{a_0a_1}): a_0,a_1=0,1\})$.
The counterfactual outcomes and covariates are latent variables.
The observed data  $(\overA_1,\overL_1,Y )$ are then determined from the full data through the consistency assumption A1 or A1$^\dag$.
In other words, Assumption A1 or A1$^\dag$ is viewed to decide how the full data lead to the observed data.
The meaning of ``full data'' used here differs from the usual literature on missing-at-random problems, where
the full data refer to $(\overline{L}_1,Y^{1,1}, Y^{1,0}, Y^{0,1}, Y^{0,0})$, excluding $\overA_1$.
From a frequentist perspective, we denote by $P$ the true distribution of the full data, which is fixed but unknown.
The induced distribution of $P$ on $(\overA_1,\overL_1,Y )$
is the true distribution of the observed data.
The expectation with respect to $P$ is denoted as $E(\cdot)$.

\subsection{Primary sensitivity model and bounds} \label{sec:2period-primary}

Our starting point is to notice that sequential unmeasured confounding beyond \textit{primary} sequential unconfounding  Assumption A2
(i.e., possible departure from (\ref{eq:unconfounded-2period})) can be characterized by the following density ratios:
\begin{subequations} \label{eq:lam-def}
\begin{align}
 \lambda^*_1 (\overline{L}_1, y) & =  \frac{\dif P_{Y^{1,1}} (y | A_0=1, A_1=0, \overline{L}_1) }{\dif P_{Y^{1,1}} (y | A_0=1, A_1=1, \overline{L}_1) }, \label{eq:lam-def-1}\\
 \lambda^*_0 (L_0, y) & = \frac{\dif P_{Y^{1,1}} ( y | A_0=0, L_0) }{\dif P_{ Y^{1,1}} ( y | A_0=1, L_0 ) } , \label{eq:lam-def-0}
\end{align}
\end{subequations}
where $P_{Y^{1,1}} (y | A_0=1, A_1=a_1, \overline{L}_1)$ denotes the distribution of $Y^{1,1}$ given $(A_0=1, A_1=a_1, \overline{L}_1)$,
and $P_{ Y^{1,1}} ( y | A_0=a_0, L_0) $ denotes the distribution of $ Y^{1,1} $ given $(A_0=a_0, L_0) $.
The special case $\lambda^*_1 = \lambda^*_0 \equiv 1$ corresponds to sequential unconfounding (\ref{eq:unconfounded-2period}),
whereas any deviations of $\lambda^*_1$ and $\lambda^*_0$ from 1 indicate sequential unmeasured confounding, i.e.,
possible differences between $P_{Y^{1,1}} (y | A_0=1, A_1=0, \overline{L}_1)$ and $P_{Y^{1,1}} (y | A_0=1, A_1=1, \overline{L}_1)$
and between $P_{ Y^{1,1}} ( y | A_0=0, L_0) $ and $P_{ Y^{1,1}} ( y | A_0=1, L_0) $.
By Bayes' rule, $\lambda^*_1$ and $\lambda^*_0$ can be equivalently expressed as odds ratios:
\begin{subequations} \label{eq:lam-def-odds}
\begin{align}
 \lambda^*_1 (\overline{L}_1, Y^{1,1}) 
  & =  \frac{ \pi^*_1 ( \overline{L}_1 ) }{ 1-\pi^*_1 (\overline{L}_1 ) } \times
  \frac{ P(A_1=0 | A_0=1,  \overline{L}_1, Y^{1,1} ) }{ P(A_1=1 | A_0=1, \overline{L}_1, Y^{1,1} ) } , \label{eq:lam-def-odds-1} \\
 \lambda^*_0 (L_0, Y^{1,1}) 
 & = \frac{ \pi^*_0( L_0 ) }{ 1- \pi^*_0 ( L_0 ) } \times
 \frac{ P (A_0=0 | L_0, Y^{1,1} ) }{ P (A_0=1 | L_0,  Y^{1,1} ) } . \label{eq:lam-def-odds-0}
\end{align}
\end{subequations}
where $\pi^*_1(\overline{L}_1 )  = P(A_1 = 1 | A_0=1, \overline{L}_1 ) $, and
$\pi^*_0(L_0) =  P(A_0 = 1| L_0)$ as in Section \ref{sec:setup}.
The first expression above is  the odds of $A_1=0$ given the counterfactual outcome $Y^{1,1}$ and the covariate history $\overline{L}_1$
over that of $A_1=0$ given only $\overline{L}_1$ in period 1.
The second expression above is  the odds of $A_0=0$ given the counterfactual outcome $Y^{1,1}$ and the baseline covariate $L_0$
over that of $A_0=0$ given only $L_0$ in period 0.
Hereafter, we refer to $\lambda^*_1$ and $\lambda^*_0$ as the sensitivity ratios in periods 1 and 0 respectively.

We stress that the sensitivity ratios $\lambda^*_1$ and $\lambda^*_0$
satisfy inherent constraints, before any further restriction is introduced on unmeasured confounding.
In fact, the definition of $\lambda^*_1$ and $\lambda^*_0$ as density ratios in (\ref{eq:lam-def}) directly implies that
\begin{subequations} \label{eq:constr-2period-a}
\begin{align}
& E ( \lambda^*_1(\overline{L}_1,Y^{1,1}) | A_0=1, A_1=1, \overline{L}_1 ) \equiv 1, \label{eq:constr-2period-a1} \\
& E ( \lambda^*_0(L_0,Y^{1,1}) | A_0=1, L_0 ) \equiv 1 . \label{eq:constr-2period-a0}
\end{align}
\end{subequations}
Under consistency Assumption A1, the two constraints can be rewritten as
\begin{subequations}  \label{eq:constr-2period}
\begin{align}
& E \left( \lambda^*_1(\overline{L}_1,Y) | A_0=1, A_1=1, \overline{L}_1 \right) \equiv 1,  \label{eq:constr-2period-0}\\
& E \left\{ E ( \lambda^*_0(L_0,Y) \varrho_1 (\overline{L}_1, Y; \lambda^*_1 ) |A_0=1,A_1=1, \overline{L}_1)  | A_0=1, L_0 \right\}\equiv 1 , \label{eq:constr-2period-1}
\end{align}
\end{subequations}
depending only on the distribution of the observed data $(\overA_1,\overL_1,Y)$,
where $\varrho_1 (\overline{L}_1, y; \lambda_1 ) = \pi^*_1(\overline{L}_1) + (1-\pi^*_1(\overline{L}_1)) \lambda_1(\overline{L}_1, y) $.
While replacing $Y^{1,1}$ by $Y$ in (\ref{eq:constr-2period-a1}) easily yields (\ref{eq:constr-2period-0}),
more care is needed in converting (\ref{eq:constr-2period-a0}) to (\ref{eq:constr-2period-1}).
Notationally, we reserve $\lambda^*_1 (\overline{L}_1, y)$ and  $\lambda^*_0 (L_0, y)$ for indicating the true sensitivity ratios,
which are fixed but unknown, in the same way as $P$ for indicating the true distribution of the full data.
By Lemma~\ref{lem:density-ratio}, any nonnegative functions $\lambda_1 (\overline{L}_1, y)$ and  $\lambda_0 (L_0, y)$ satisfying (\ref{eq:constr-2period-b})
can be valid choices for $\lambda^*_1 (\overline{L}_1, y)$ and  $\lambda^*_0 (L_0, y)$,
while being fully compatible with the observed data $(\overline{A}_1,\overline{L}_1,Y)$.
Our proof of the sufficiency assertion in Lemma~\ref{lem:density-ratio} is constructive, as illustrated by Figure~\ref{fig:hist-combined-sen} later.

\begin{lem} \label{lem:density-ratio}
Suppose that Assumption A1 holds for $\overline{a}=(1,1)$.
The sensitivity ratios $\lambda^*_1$ and $\lambda^*_0$ satisfy (\ref{eq:constr-2period}).
Conversely, if nonnegative functions $\lambda_1 (\overline{L}_1, y)$ and  $\lambda_0 (L_0, y)$ satisfy the normalization constraints,
with $\varrho_1(\cdot)$ defined as in (\ref{eq:constr-2period}),
\begin{subequations} \label{eq:constr-2period-b}
\begin{align}
& E ( \lambda_1(\overline{L}_1,Y) | A_0=1, A_1=1, \overline{L}_1 ) \equiv 1,  \label{eq:constr-2period-b1} \\
& E\{ E ( \lambda_0(L_0,Y) \varrho_1 (\overline{L}_1, Y; \lambda_1 ) |A_0=1,A_1=1, \overline{L}_1)  | A_0=1, L_0 \}\equiv 1 ,  \label{eq:constr-2period-b0}
\end{align}
\end{subequations}
then there exists a probability distribution $Q$ for the full data $(\overline{A}_1,\overline{L}_1,Y^{1,1}, Y^{1,0}, Y^{0,1}, Y^{0,0})$ satisfying
the two properties: \vspace{-.05in}
\begin{itemize}\addtolength{\itemsep}{-.1in}
\item[(i)] the induced distribution of $Q$ on the observed data $(\overline{A}_1,\overline{L}_1,Y )$ under Assumption A1 coincides with
the true distribution of  $(\overline{A}_1,\overline{L}_1,Y )$,
\item[(ii)] $ \lambda^*_{1,Q} (\overline{L}_1, y) =\lambda_1 (\overline{L}_1, y) $ and
 $ \lambda^*_{0,Q} (L_0, y) = \lambda_0 (L_0, y)$,
\end{itemize} \vspace{-.05in}
where $\lambda^*_{1,Q}$ and $\lambda^*_{0,Q}$ are defined as $\lambda^*_1$ and $\lambda^*_0$ in (\ref{eq:lam-def})  with $P$ replaced by $Q$.
\end{lem}

For a sensitivity model, we postulate that
the sensitivity ratios $\lambda^*_1$ and $ \lambda^*_0$
may differ from 1 by at most a factor of $\Lambda_1$ and $\Lambda_0$ for any $y$ (almost surely in $\overline{L}_1$):
\begin{align} \label{eq:model-2period}
\Lambda_1^{-1} \le \lambda^*_1 (\overline{L}_1, y) \le \Lambda_1, \quad
\Lambda_0^{-1} \le \lambda^*_0 (L_0, y) \le \Lambda_0 ,
\end{align}
where $\Lambda_1 \ge 1$ and $\Lambda_0 \ge 1$ are sensitivity parameters, encoding the degrees of unmeasured confounding
in periods 1 and 0 respectively.
By design,
the extreme case of $\Lambda_0 = \Lambda_1 =1$ reduces to the primary sequential unconfounding assumption (\ref{eq:unconfounded-2period}),
where $\lambda^*_1(\overline{L}_1, y) = \lambda^*_0  (L_0, y) \equiv 1$ trivially satisfies the normalization constraints (\ref{eq:constr-2period-b}).

For sensitivity analysis, we are interested in deriving population bounds on $\mu^{1,1}$ under sensitivity model (\ref{eq:model-2period}).
First, we show that if the sensitivity ratios $\lambda^*_1$ and $ \lambda^*_0$ were known, then $\mu^{1,1}$ can be point identified
in two distinct ways, as a generalization of the ICE and IPW formulas (\ref{eq:iden-ICE}) and (\ref{eq:iden-IPW}) from sequential unconfounding.
For any nonnegative functions $\lambda_1 (\overline{L}_1, y)$ and  $\lambda_0 (L_0, y)$, we define
\begin{align} \label{eq:g-iden-2period}
\begin{split}
& \mu^{1,1}_{\mytext{ICE}} (\lambda_0,\lambda_1) =
E \left( \mathcal E_{L_0} \left[ \mathcal E_{\overL_1} \left\{ (\pi^*_0 + (1-\pi^*_0) \lambda_0) (\pi^*_1 + (1-\pi^*_1) \lambda_1) Y \right\} \right] \right),
 \\ 
& \mu^{1,1}_{\mytext{IPW}} (\lambda_0,\lambda_1) = E \left\{ A_0 A_1 \left( 1+ \frac{1-\pi^*_0 } { \pi^*_0 }\lambda_0 \right)
\left( 1+ \frac{1-\pi^*_1 } { \pi^*_1 }\lambda_1 \right) Y \right\} , 
\end{split}
\end{align}
where
$\pi^*_0 = \pi^*_0 (L_0) $,
$\pi^*_1 = \pi^*_1(\overline{L}_1 )$,
$\lambda_0 = \lambda_0 (L_0,Y)$, and $\lambda_1 = \lambda_1( \overL_1, Y)$ inside the expectations.
Throughout, $\mathcal E_{\overline{L}_1}(\cdot)$ denotes the conditional expectation $E(\cdot |A_0=1,A_1=1, \overline{L}_1)$
and $\mathcal E_{L_0}(\cdot)$ denotes the conditional expectation $E(\cdot |A_0=1, L_0)$.

\begin{lem} \label{lem:iden-mu-2period}
Suppose that Assumptions A1 and A3 hold for $\overline{a}=(1,1)$,
and $\lambda^*_1 (\overline{L}_1, Y^{1,1})$ and $\lambda^*_0 (L_0, Y^{1,1})$ are finite almost surely. Then
$\mu^{1,1} = \mu^{1,1}_{\mytext{ICE}} (\lambda^*_0,\lambda^*_1) = \mu^{1,1}_{\mytext{IPW}} (\lambda^*_0,\lambda^*_1)$
\end{lem}

The sharp upper bound on $\mu^{1,1}$ under model (\ref{eq:model-2period}) is defined as
\begin{align} \label{eq:upper-2period}
\mu^{1,1}_+ = \max_{Q} \; \mu^{1,1}_{\mytext{ICE}} (\lambda^*_{0,Q},\lambda^*_{1,Q})
= \max_{Q} \; \mu^{1,1}_{\mytext{IPW}} (\lambda^*_{0,Q},\lambda^*_{1,Q}),
\end{align}
over all possible distributions $Q$ for the full data such that (i) the induced distribution of $Q$ on the observed data $(\overA_1,\overL_1,Y )$
coincides with the true distribution, and (ii)
the sensitivity ratios $\lambda^*_{1,Q}$ and $\lambda^*_{0,Q}$ based on $Q$ as in Lemma~\ref{lem:density-ratio}
satisfy the range constraints as $\lambda^*_1$ and $\lambda^*_0$ in (\ref{eq:model-2period}).
In other words, $\mu^{1,1}_+$ is the maximum value of $\mu^{1,1}$ over all possible choices of the true distribution $P$ (denoted as $Q$)
under the sensitivity model (\ref{eq:model-2period}).
Lemma~\ref{lem:Q-lam-2period} gives an equivalent representation of $\mu^{1,1}_+$, based on
the characterization of sensitivity ratios in Lemma~\ref{lem:density-ratio}.
This result is important in
transforming the optimization from (\ref{eq:upper-2period}) over possible distributions $Q$ for the full data
to (\ref{eq:upper-2period-b}) over nonnegative functions $\lambda_1 (\overline{L}_1, y)$ and  $\lambda_0 (L_0, y)$ with constraints (\ref{eq:constr-2period-b})
and (\ref{eq:model-2period-b}) depending only on the observed data $(\overline{A}_1,\overline{L}_1,Y)$.

\begin{lem} \label{lem:Q-lam-2period}
The sharp upper bound $\mu^{1,1}_+$ in (\ref{eq:upper-2period}) can be
equivalently obtained as
\begin{align}  \label{eq:upper-2period-b}
\mu^{1,1}_+ = \max_{\lambda_0,\lambda_1} \; \mu^{1,1}_{\mytext{ICE}} (\lambda_0,\lambda_1)
= \max_{\lambda_0,\lambda_1} \; \mu^{1,1}_{\mytext{IPW}} (\lambda_0,\lambda_1),
\end{align}
over all possible nonnegative functions $ \lambda_1 = \lambda_1 (\overline{L}_1, y) $  and
$\lambda_0 = \lambda_0 (L_0, y)$
subject to the normalization constraints (\ref{eq:constr-2period-b})
and the range constraints that for any $y$,
\begin{align} \label{eq:model-2period-b}
\Lambda_1^{-1} \le \lambda_1 (\overline{L}_1, y) \le \Lambda_1, \quad
\Lambda_0^{-1} \le \lambda_0 (L_0, y) \le \Lambda_0 .
\end{align}
\end{lem}

As the final result of Section \ref{sec:2period-primary}, we show that the constrained optimization in (\ref{eq:upper-2period-b}) subject to (\ref{eq:constr-2period-b})
and (\ref{eq:model-2period-b})
can be solved through a dual relationship to an unconstrained optimization problem.
For any functions $q_0=q_0(L_0)$ and $q_1=q_1(\overline{L}_1)$, let
\begin{align*}
& \eta_{0+} (y,q_0) = y + (1-\pi^*_0) (\Lambda_0-\Lambda_0^{-1}) \rho_{\tau_0} ( y, q_0 ),\\
& \eta_{1+} (y, q_1) = y + (1-\pi^*_1) (\Lambda_1-\Lambda_1^{-1}) \rho_{\tau_1} (y, q_1 ),
\end{align*}
where $\tau_0 = \Lambda_0/(1+\Lambda_0)$, $\tau_1 = \Lambda_1/(1+\Lambda_1)$,
$\rho_\tau (y,q) = \tau (y-q)^+ + (1-\tau) (q-y)^+$, known as a ``check" function, and $c^+ = c$ for $c \ge 0$ or $0$ for $c <0$.

\begin{pro} \label{pro:upper-2period}
The sharp upper bound $\mu^{1,1}_+$ in (\ref{eq:upper-2period-b}) can be determined as
\begin{align}\label{eq:upper-2period-c}
\begin{split}
\mu^{1,1}_+ & = \min_{q_0,q_1}\; E \{ \mathcal E_{L_0} ( \mathcal E_{\overL_1} [ \eta_{1+} \{ \eta_{0+}(Y, q_0), q_1\} ] ) \} \\
& = \min_{q_0,q_1}\; E \left[ \frac{A_0A_1}{\pi^*_0 \pi^*_1} \eta_{1+} \{ \eta_{0+}(Y, q_0), q_1\}  \right],
\end{split}
\end{align}
over all possible functions $q_0 = q_0(L_0)$ and $q_1= q_1(\overline{L}_1)$.
Let $(\check q_0, \check q_1)$ be a solution to the optimization in (\ref{eq:upper-2period-c}). Then a solution, $(\check\lambda_0,\check\lambda_1)$,
to the optimization in (\ref{eq:upper-2period-b})
can be obtained such that $ \check\lambda_0 (L_0,y) = \Lambda_0$ if $ y > \check q_0(L_0)$ or $=\Lambda_0^{-1}$ if $y < \check q_0(L_0)$,
and $\check\lambda_1( \overline{L}_1, y) = \Lambda_1$ if $\eta_{0+}(y,\check q_0) > \check q_1(\overline{L}_1)$ or $= \Lambda_1^{-1}$ if $\eta_{0+}(y,\check q_0) < \check q_1(\overline{L}_1)$
\end{pro}

We provide some remarks on the implications of Proposition~\ref{pro:upper-2period}. First, as shown in Supplement Lemma~\ref{lem:eta-convex},
the function $\eta_{j+}(y,q_j)$ is increasing and convex in $y$ and convex in $q_j$ for $j=1,2$,
and the composite function $ \eta^1_+ (y, \overq_1) =\eta_{1+} ( \eta_{0+}(y,q_0), q_1)$ is increasing in $y$ and convex in $\overq_1=(q_0,q_1)$.
Then the sample analogs of the representations of $\mu^{1,1}_+$ in (\ref{eq:upper-2period-c})
can be computed via convex optimization, provided that $q_0=q_0(L_0)$ and $q_1=q_1(\overL_1)$ are
linearly parameterized such that $(\check q_0, \check q_1)$ are
contained in the parameterized function classes, in addition to the assumption that either the ICE or IPW functional can be consistently estimated
for the ``outcome'' $\eta^1_+ (Y, \overq_1)$.
If the parameterizations of $(q_0,q_1)$ are misspecified such that $(\check q_0, \check q_1)$ are not captured by
the parameterized function classes, then in the large-sample limit,
the computed upper bounds remain valid (or more precisely, conservative) in being no smaller than $\mu^{1,1}$.
Hence the representations of $\mu^{1,1}_+$ in (\ref{eq:upper-2period-c}) extend similar properties
associated with the sharp sensitivity bounds in single-period
cross-sectional studies (Dorn \& Guo 2022; Tan 2022). See Supplement Section~\ref{sec:tech-preparation} for a summary.
We leave further development of estimation methods with sample data to future work.

Second, the representations of $\mu^{1,1}_+$ in (\ref{eq:upper-2period-c}) also reveal how the sensitivity bound $\mu^{1,1}_+$ can be influenced by
certain aspects of the observed data, while the point identification of $\mu^{1,1}$ under sequential unconfounding is not affected.
See Section \ref{sec:numerical} for numerical illustration.
For example, if $\pi^*_0$ or $\pi^*_1$ is decreased in the observed-data distribution with all other aspects unchanged, then
the ICE or IPW formula (\ref{eq:iden-ICE}) or (\ref{eq:iden-IPW}) gives the same \textit{population value},
but the sharp upper bound $\mu^{1,1}_+$ is increased.
This phenomenon seems to be expected because smaller $\pi^*_0$ or $\pi^*_1$ indicates larger probabilities of
observing $(A_0,A_1)\not=(1,1)$ and hence missing $Y^{1,1}$.
For a nontrivial example, we point out that if the outcome $Y$ is distributed with more uncertainty (i.e., higher noise levels) while
its conditional mean $E(Y|A_0=A_1=1, \overL_1)$ is fixed, then
the same point value is identified by the ICE formula under sequential unconfounding,
but the sharp upper bound $\mu^{1,1}_+$ tends to be larger.
Formally, suppose $Y$ is distributed like $\alpha (\overL_1) + \sigma Z$
conditionally on $(A_0=A_1=1, \overL_1)$, where $\alpha(\cdot)$ is an arbitrary function, $\sigma \ge 0$, and $Z$ is a noise variable satisfying
$ E( Z | A_0=A_1=1, \overL_1) \equiv 0$. In Supplement Section \ref{sec:prf-mono-sig}, we show that for any $\sigma \ge \tilde \sigma \ge 0$,
\begin{align}
& \quad \min_{ q_0,q_1 }\;
 E \{ \mathcal E_{L_0} (  \mathcal E_{\overline{L}_1} [ \eta_{1+} \{ \eta_{0+} ( \alpha + \sigma Z,  q_0 ), q_1\} ] ) \} \nonumber \\
& \ge  \min_{ q_0,q_1 }\;
 E \{ \mathcal E_{L_0} (  \mathcal E_{\overline{L}_1} [ \eta_{1+} \{ \eta_{0+} ( \alpha + \tilde\sigma Z,  q_0 ), q_1\} ] ) \}  , \label{eq:mono-sig}
\end{align}
In the extreme case of $\tilde\sigma=0$, the right-hand side of (\ref{eq:mono-sig}) reduces to the ICE formula (\ref{eq:iden-ICE}).
The proof of (\ref{eq:mono-sig}) is lengthy, even though the interpretation may appear to be natural.
This observation sheds interesting light on the role of outcome predictive modeling in sensitivity analysis, echoing a similar discussion in Tan (2022).
For two observational studies, I and II, if the outcome can be more accurately predicted from the covariates in study I,
then the sensitivity upper (or lower) bounds are expected to be smaller (or larger), hence the sensitivity intervals
expected to be narrower,  from study I than from study II.

\subsection{Joint sensitivity model and bounds} \label{sec:2period-Jnt}

In Section \ref{sec:2period-primary}, sensitivity analysis is formulated by relaxing \textit{primary} sequential unconfounding Assumption A2.
In this section, we study sensitivity analysis by relaxing \textit{joint} sequential unconfounding Assumption A2$^\dag$
(i.e., possible departure from (\ref{eq:unconfounded-2period-Jnt})).
Such sequential unmeasured confounding can be characterized by the density ratios $\lambda^*_1 (\overline{L}_1, y)$ and
\begin{align}
 \lambda^*_{0,\mytext{Jnt}} (L_0, l_1, y) & = \frac{\dif P_{L_1^{1,1}, Y^{1,1}} (l_1, y | A_0=0, L_0) }{\dif P_{L_1^{1,1}, Y^{1,1}} (l_1, y | A_0=1, L_0 ) } ,
 \label{eq:lam-def-Jnt}
\end{align}
where $P_{L_1^{1,1}, Y^{1,1}} (l_1, y | A_0=a_0, L_0) $ denotes the distribution of $(L_1^{1,1}, Y^{1,1})$ given $(A_0=a_0, L_0) $.
By Bayes' rule, $\lambda^*_{0,\mytext{Jnt}}$ can be expressed as an odds ratio:
\begin{align*}
 \lambda^*_{0,\mytext{Jnt}} (L_0, L_1^{1,1}, Y^{1,1}) 
 & = \frac{ \pi^*_0( L_0 ) }{ 1- \pi^*_0 ( L_0 ) } \times
 \frac{ P (A_0=0 | L_0, L_1^{1,1}, Y^{1,1} ) }{ P (A_0=1 | L_0, L_1^{1,1}, Y^{1,1} ) } .
\end{align*}
Compared with $\lambda^*_0$ in (\ref{eq:lam-def-odds-0}), the preceding expression is the odds of $A_0=0$ given
$L_0$ and $(L_1^{1,1}, Y^{1,1})$, not just $Y^{1,1}$,
over that of $A_0=0$ given $L_0$ only in period 0.

Similarly to Lemma~\ref{lem:density-ratio}, the following result provides a necessary and sufficient condition for
nonnegative functions $\lambda_1 (\overline{L}_1,y)$ and $\lambda_{0,\mytext{Jnt}} (L_0, l_1, y)$
to serve as sensitivity ratios, while being fully compatible with the observed data $(\overline{A}_1, \overline{L}_1, Y)$.
The normalization condition (\ref{eq:constr-2period-Jnt-b1}) is the same as (\ref{eq:constr-2period-b1}) in Lemma~\ref{lem:density-ratio}
and is restated here for convenience.

\begin{lem} \label{lem:density-ratio-Jnt}
Suppose that Assumption A1$^\dag$ holds for $\overline{a}=(1,1)$.
The sensitivity ratios $\lambda^*_1$ and $\lambda^*_{0,\mytext{Jnt}}$ satisfy the normalization constraints
as $\lambda_1 $ and  $\lambda_{0,\mytext{Jnt}} $ in (\ref{eq:constr-2period-Jnt-b}).
Conversely, if nonnegative functions $\lambda_1 (\overline{L}_1, y)$ and  $\lambda_{0,\mytext{Jnt}} (L_0, l_1, y)$ satisfy
\begin{subequations} \label{eq:constr-2period-Jnt-b}
\begin{align}
& E ( \lambda_1(\overline{L}_1,Y) | A_0=1, A_1=1, \overline{L}_1 ) \equiv 1,  \label{eq:constr-2period-Jnt-b1} \\
& E\{ E ( \lambda_{0,\mytext{Jnt}}(\overline{L}_1,Y) \varrho_1 (\overline{L}_1, Y; \lambda_1 ) |A_0=1,A_1=1, \overline{L}_1)  | A_0=1, L_0 \}\equiv 1 ,
\label{eq:constr-2period-Jnt-b0}
\end{align}
\end{subequations}
then there exists a probability distribution $Q$ for the full data $(\overline{A}_1, L_0,\{(L_1^{a_0a_1},Y^{a_0a_1}): a_0,a_1=0,1\} )$ satisfying
the two properties: \vspace{-.05in}
\begin{itemize}\addtolength{\itemsep}{-.1in}
\item[(i)] the induced distribution of $Q$ on the observed data $(\overline{A}_1,\overline{L}_1,Y )$ under Assumption A1$^\dag$ coincides with
the true distribution of  $(\overline{A}_1,\overline{L}_1,Y )$,
\item[(ii)] $ \lambda^*_{1,Q} (\overline{L}_1, y) =\lambda_1 (\overline{L}_1, y) $ and
 $ \lambda^*_{0,\mytext{Jnt},Q} (L_0, l_1, y) = \lambda_{0,\mytext{Jnt}} (L_0, l_1, y)$,
\end{itemize} \vspace{-.05in}
where $\lambda^*_{1,Q}$ and $\lambda^*_{0,\mytext{Jnt},Q}$ are defined as $\lambda^*_1$ and $\lambda^*_{0,\mytext{Jnt}}$ in (\ref{eq:lam-def-1})
and (\ref{eq:lam-def-Jnt}) with $P$ replaced by $Q$.
\end{lem}

Similarly to Lemma~\ref{lem:iden-mu-2period}, the following result shows that if the sensitivity ratios $\lambda^*_1$ and $\lambda^*_{0,\mytext{Jnt}}$ were known,
then $\mu^{1,1}$ can be point identified by a generalization of the ICE and IPW formulas (\ref{eq:iden-ICE}) and (\ref{eq:iden-IPW}) from sequential unconfounding.

\begin{lem} \label{lem:iden-mu-2period-Jnt}
Suppose that Assumptions A1$^\dag$ and A3 hold for $\overline{a}=(1,1)$,
and $\lambda^*_1 (\overline{L}_1, Y^{1,1})$ and $\lambda^*_{0,\mytext{Jnt}} (L_0, L_1^{1,1}, Y^{1,1})$ are finite almost surely. Then
\begin{align*}
\mu^{1,1} = \mu^{1,1}_{\mytext{ICE}} (\lambda^*_{0,\mytext{Jnt}},\lambda^*_1) = \mu^{1,1}_{\mytext{IPW}} (\lambda^*_{0,\mytext{Jnt}},\lambda^*_1),
\end{align*}
where for any nonnegative functions $\lambda_1 (\overline{L}_1, y)$ and $\lambda_{0,\mytext{Jnt}} (L_0, l_1, y)$,
the functionals $\mu^{1,1}_{\mytext{ICE}}(\lambda_{0,\mytext{Jnt}},$  $\lambda_1)$ and $\mu^{1,1}_{\mytext{ICE}}(\lambda_{0,\mytext{Jnt}},\lambda_1)$ are defined as
 $\mu^{1,1}_{\mytext{ICE}}(\lambda_0,\lambda_1)$ and $\mu^{1,1}_{\mytext{ICE}}(\lambda_0,\lambda_1)$ in (\ref{eq:g-iden-2period})
 except with $\lambda_0=\lambda_0 (L_0, y)$ replaced by $\lambda_{0,\mytext{Jnt}}=\lambda_{0,\mytext{Jnt}} (L_0, l_1, y)$.
\end{lem}

For a sensitivity model similar to (\ref{eq:model-2period}), we postulate that
the sensitivity ratios $\lambda^*_1$ and $ \lambda^*_{0,\mytext{Jnt}}$
may differ from 1 by at most a factor of $\Lambda_1$ and $\Lambda_0$ for any $(l_1,y)$:
\begin{align} \label{eq:model-2period-Jnt}
\Lambda_1^{-1} \le \lambda^*_1 (\overline{L}_1, y) \le \Lambda_1, \quad
\Lambda_0^{-1} \le \lambda^*_{0,\mytext{Jnt}} (L_0, l_1, y) \le \Lambda_0 ,
\end{align}
where $\Lambda_1 \ge 1$ and $\Lambda_0 \ge 1$ are sensitivity parameters.
By design,
the extreme case of $\Lambda_0 = \Lambda_1 =1$ reduces to the joint sequential unconfounding assumption (\ref{eq:unconfounded-2period-Jnt}),
where $\lambda^*_1(\overline{L}_1, y) = \lambda^*_{0,\mytext{Jnt}}  (L_0, l_1, y) \equiv 1$
trivially satisfies the normalization constraints (\ref{eq:constr-2period-Jnt-b}).

The sharp upper bound on $\mu^{1,1}$ under model (\ref{eq:model-2period-Jnt}) is defined as
\begin{align}  \label{eq:upper-2period-Jnt}
\mu^{1,1}_{+,\mytext{Jnt}} = \max_{Q} \; \mu^{1,1}_{\mytext{ICE}} (\lambda^*_{0,\mytext{Jnt},Q},\lambda^*_{1,Q})
= \max_{Q} \; \mu^{1,1}_{\mytext{IPW}} (\lambda^*_{0,\mytext{Jnt},Q},\lambda^*_{1,Q}),
\end{align}
over all possible distributions $Q$ for the full data such that (i) the induced distribution of $Q$ on the observed data $(\overA_1,\overL_1,Y )$
coincides with the true distribution, and
(ii) the sensitivity ratios $\lambda^*_{1,Q}$ and $\lambda^*_{0,\mytext{Jnt},Q}$ based on $Q$ as in Lemma~\ref{lem:density-ratio-Jnt}
satisfy the range constraints as $\lambda^*_1$ and $\lambda^*_{0,\mytext{Jnt}}$ in (\ref{eq:model-2period-Jnt}).
In other words, $\mu^{1,1}_{+,\mytext{Jnt}}$ is the maximum value of $\mu^{1,1}$ over all possible choices of the true distribution $P$ (denoted as $Q$)
under the sensitivity model (\ref{eq:model-2period-Jnt}).
In parallel to Lemma~\ref{lem:Q-lam-2period}, the following result gives an equivalent representation of $\mu^{1,1}_{+,\mytext{Jnt}}$
in terms of only the observed data, by the characterization of sensitivity ratios in Lemma~\ref{lem:density-ratio-Jnt}.

\begin{lem} \label{lem:Q-lam-2period-Jnt}
The sharp upper bound $\mu^{1,1}_{+,\mytext{Jnt}}$ in (\ref{eq:upper-2period-Jnt}) can be equivalently obtained as
\begin{align}  \label{eq:upper-2period-Jnt-b}
\mu^{1,1}_{+,\mytext{Jnt}} = \max_{\lambda_{0,\mytext{Jnt}},\lambda_1} \; \mu^{1,1}_{\mytext{ICE}} (\lambda_{0,\mytext{Jnt}},\lambda_1)
= \max_{\lambda_{0,\mytext{Jnt}},\lambda_1} \; \mu^{1,1}_{\mytext{IPW}} (\lambda_{0,\mytext{Jnt}},\lambda_1),
\end{align}
over all possible nonnegative functions $ \lambda_1 = \lambda_1 (\overline{L}_1, y) $ and
$\lambda_{0,\mytext{Jnt}} = \lambda_{0,\mytext{Jnt}} (L_0, l_1, y)$
subject to the normalization constraints (\ref{eq:constr-2period-Jnt-b})
and the range constraints that for any $(l_1,y)$,
\begin{align} \label{eq:model-2period-Jnt-b}
\Lambda_1^{-1} \le \lambda_1 (\overline{L}_1, y ) \le \Lambda_1, \quad
\Lambda_0^{-1} \le \lambda_{0,\mytext{Jnt}} (L_0, l_1, y ) \le \Lambda_0 .
\end{align}
\end{lem}

To conclude Section \ref{sec:2period-Jnt}, we provide an interesting discovery that although joint sensitivity model (\ref{eq:model-2period-Jnt})
is more restrictive than (\ref{eq:model-2period}) [i.e.,
if model (\ref{eq:model-2period-Jnt}) holds, then model (\ref{eq:model-2period}) holds],
the sharp bounds on $\mu^{1,1}$ under the two models are identical to each other.

\begin{pro} \label{pro:upper-2period-Jnt}
The sharp upper bound $\mu^{1,1}_{+,\mytext{Jnt}}$ in  (\ref{eq:upper-2period-Jnt}) or (\ref{eq:upper-2period-Jnt-b})
coincides with $\mu^{1,1}_+$ in (\ref{eq:upper-2period}), (\ref{eq:upper-2period-b}) or (\ref{eq:upper-2period-c}),
i.e.,  $\mu^{1,1}_{+,\mytext{Jnt}} =\mu^{1,1}_+$.
The sharp bound $\mu^{1,1}_{+,\mytext{Jnt}}$ can be achieved by a probability distribution $Q$ for the full data such that
$ L_1^{1,1} \perp A_0 \,|\, L_0, Y^{1,1}$ under $Q$.
\end{pro}

\subsection{Product sensitivity model and bounds} \label{sec:2period-Prod}

As an alternative to Section \ref{sec:2period-Jnt}, we study sensitivity analysis by relaxing joint sequential unconfounding Assumption A2$^\dag$
through a decomposition of joint sensitivity ratios, which is reflected by how Assumption A2$^\dag$ can be equivalently stated as
(\ref{eq:unconfounded-2period-Jnt}) or (\ref{eq:unconfounded-2period-Jnt2}) in two periods.
The joint sensitivity ratio in (\ref{eq:lam-def-Jnt}) can be decomposed as
\begin{align*}
 \lambda^*_{0,\mytext{Jnt}} (L_0, l_1, y) & =\lambda^*_{0,L_1} (L_0, l_1) \, \lambda^*_{0,Y} (L_0, l_1, y),
\end{align*}
where
\begin{subequations} \label{eq:lam-def-0-Prod}
\begin{align}
& \lambda^*_{0,L_1} (L_0, l_1) = \frac{\dif P_{L_1^{1,1}} (l_1 | A_0=0, L_0) }{\dif P_{L_1^{1,1}} (l_1 | A_0=1, L_0 ) } ,  \label{eq:lam-def-0-L1} \\
& \lambda^*_{0,Y} (L_0, l_1, y) =   \frac{\dif P_{Y^{1,1}} (y| A_0=0, L_0, L_1^{1,1} = l_1) }{\dif P_{Y^{1,1}} (y| A_0=1, L_0, L_1^{1,1} = l_1) } . 
\end{align}
\end{subequations}
By Bayes' rule, $\lambda^*_{0,L_1}$ and $\lambda^*_{0,Y}$ can be expressed as odds ratios:
\begin{align*}
& \lambda^*_{0,L_1} (L_0, L_1^{1,1}) = \frac{ \pi^*_0( L_0 ) }{ 1- \pi^*_0 ( L_0 ) } \times
 \frac{ P (A_0=0 | L_0, L_1^{1,1} ) }{ P (A_0=1 | L_0,  L_1^{1,1} ) } , \\
& \lambda^*_{0,Y} (L_0, L_1^{1,1}, Y^{1,1} ) = \frac{ P (A_0=1 | L_0, L_1^{1,1} ) }{ P (A_0=0 | L_0,  L_1^{1,1} ) }  \times
 \frac{ P (A_0=0 | L_0, L_1^{1,1}, Y^{1,1} ) }{ P (A_0=1 | L_0,  L_1^{1,1}, Y^{1,1} ) } .
\end{align*}
The sensitivity ratio $\lambda^*_{0,L_1} (L_0, L_1^{1,1})$ is the odds of $A_0=0$ given $(L_0, L_1^{1,1})$ over that of $A_0=0$ given $L_0$ only in period 0, and
$\lambda^*_{0,Y} (L_0, L_1^{1,1}, Y^{1,1} )$ is the odds of $A_1=0$ given $(L_0, L_1^{1,1}, Y^{1,1})$ over that of $A_0=0$ given $(L_0, L_1^{1,1})$ in period 1.
It is worth noticing that, in contrast with $\lambda^*_0$, $\lambda^*_{0,\mytext{Jnt}}$, and $\lambda^*_{0,L_1}$,
the odds-ratio interpretation of $\lambda^*_{0,Y}$ above involves two odds of $A_1=0$ which both depend on counterfactual covariate $L_1^{1,1}$ or outcome $Y^{1,1}$.

In parallel to Lemmas~\ref{lem:density-ratio} and \ref{lem:density-ratio-Jnt}, the following result provides a necessary and sufficient condition for
nonnegative functions $\lambda_1 (\overline{L}_1,y)$, $\lambda_{0,Y} (L_0, l_1, y)$, and $\lambda_{0,L_1} (L_0, l_1)$
to serve as sensitivity ratios, while being fully compatible with the observed data $(\overline{A}_1, \overline{L}_1, Y)$.
The normalization condition (\ref{eq:constr-2period-Prod-b1}) for $\lambda_1$ is the same as (\ref{eq:constr-2period-b1}) in Lemma~\ref{lem:density-ratio}
or (\ref{eq:constr-2period-Prod-b1}) in Lemma~\ref{lem:density-ratio-Jnt}.
As a consistency check, it can be easily verified that if $\lambda_{0,Y} $ and $\lambda_{0,L_1}$
satisfy (\ref{eq:constr-2period-Prod-b0Y}) and (\ref{eq:constr-2period-Prod-b0L}), then
$\lambda_{0,\mytext{Jnt}} = \lambda_{0,L_1} \lambda_{0,Y} $ satisfies (\ref{eq:constr-2period-Jnt-b0}) for a joint sensitivity ratio in Lemma~\ref{lem:density-ratio-Jnt}.

\begin{lem} \label{lem:density-ratio-Prod}
Suppose that Assumption A1$^\dag$ holds for $\overline{a}=(1,1)$.
The sensitivity ratios $\lambda^*_1$, $\lambda^*_{0,Y}$, and $\lambda^*_{0,L_1}$ satisfy the normalization constraints
as $\lambda_1 $, $\lambda_{0,Y}$, and $\lambda_{0,L_1}$ in (\ref{eq:constr-2period-Prod-b}).
Conversely, if nonnegative functions $\lambda_1 (\overline{L}_1, y)$, $\lambda_{0,Y} (L_0, l_1, y)$, and $\lambda_{0,L_1} (L_0, l_1)$ satisfy
\begin{subequations} \label{eq:constr-2period-Prod-b}
\begin{align}
& E ( \lambda_1(\overline{L}_1,Y) | A_0=1, A_1=1, \overline{L}_1 ) \equiv 1,  \label{eq:constr-2period-Prod-b1} \\
& E ( \lambda_{0,Y}(\overline{L}_1,Y) \varrho_1 (\overline{L}_1, Y; \lambda_1 ) |A_0=1,A_1=1, \overline{L}_1) \equiv 1 ,
\label{eq:constr-2period-Prod-b0Y} \\
& E ( \lambda_{0,L_1}(\overline{L}_1) | A_0=1, L_0 ) \equiv 1 ,
\label{eq:constr-2period-Prod-b0L}
\end{align}
\end{subequations}
then there exists a probability distribution $Q$ for the full data $(\overline{A}_1, L_0,\{(L_1^{a_0a_1},Y^{a_0a_1}): a_0,a_1=0,1\} )$ satisfying
the two properties: \vspace{-.05in}
\begin{itemize}\addtolength{\itemsep}{-.1in}
\item[(i)] the induced distribution of $Q$ on the observed data $(\overline{A}_1,\overline{L}_1,Y )$ under Assumption A1$^\dag$ coincides with
the true distribution of  $(\overline{A}_1,\overline{L}_1,Y )$,
\item[(ii)] $ \lambda^*_{1,Q} (\overline{L}_1, y) =\lambda_1 (\overline{L}_1, y) $,
 $ \lambda^*_{0,Y,Q} (L_0, l_1, y) = \lambda_{0,Y} (L_0, l_1, y)$, and
 $ \lambda^*_{0,L_1,Q} (L_0, l_1) = \lambda_{0,L_1} $ $ (L_0, l_1)$,
\end{itemize} \vspace{-.05in}
where $\lambda^*_{1,Q}$, $\lambda^*_{0,Y,Q}$, and $\lambda^*_{0,L_1,Q}$ are defined as $\lambda^*_1$, $\lambda^*_{0,Y}$, and $\lambda^*_{0,L_1}$
in (\ref{eq:lam-def-1}) and (\ref{eq:lam-def-0-Prod}) with $P$ replaced by $Q$.
\end{lem}

For a sensitivity model, we postulate that
the sensitivity ratios $\lambda^*_1$, $\lambda^*_{0,Y}$, and $ \lambda^*_{0,L_1}$
may differ from 1 by at most a factor of $\Lambda_1$, $\Lambda_{0,Y}$, and $\Lambda_{0,L_1}$ for any $(l_1,y)$:
\begin{align}\label{eq:model-2period-Prod}
\begin{split}
& \Lambda_1^{-1} \le \lambda^*_1 (\overline{L}_1, y) \le \Lambda_1,  \\
& \Lambda_{0,L_1}^{-1} \le \lambda^*_{0,L_1} (L_0, l_1) \le \Lambda_{0,L_1} ,\quad
\Lambda_{0,Y}^{-1} \le \lambda^*_{0,Y} (L_0, l_1, y) \le \Lambda_{0,Y} ,
\end{split}
\end{align}
where $\Lambda_1 \ge 1$, $\Lambda_{0,Y} \ge 1$, and $\Lambda_{0,L_1} \ge 1$ are sensitivity parameters.
In the case of $\Lambda_1 =\Lambda_{0,Y} = \Lambda_{0,L_1} = 1$, model (\ref{eq:model-2period-Prod}) reduces to
the joint sequential unconfounding assumption (\ref{eq:unconfounded-2period-Jnt2}),
which is equivalent to assumption (\ref{eq:unconfounded-2period-Jnt}) recovered by
model (\ref{eq:model-2period-Jnt}) with $\Lambda_1 =\Lambda_{0,\mytext{Jnt}} = 1$.
However, if $\Lambda_{0,Y}>1$ or $ \Lambda_{0,L_1} > 1$, then model (\ref{eq:model-2period-Prod})
is more restrictive (i.e., placing more restrictions on the full data) than model (\ref{eq:model-2period-Jnt})
with $\Lambda_{0,\mytext{Jnt}} = \Lambda_{0,L_1} \Lambda_{0,Y}$.

For point identification of $\mu^{1,1}$ with known $(\lambda^*_1, \lambda^*_{0,Y}, \lambda^*_{0,L_1})$,
the generalized ICE and IPW formulas  in Lemma~\ref{lem:iden-mu-2period-Jnt}
remain valid by setting $\lambda^*_{0,\mytext{Jnt}} =\lambda^*_{0,L_1} \lambda^*_{0,Y}$.
Under sensitivity model (\ref{eq:model-2period-Prod}), the sharp upper bound on $\mu^{1,1}$ is defined as
\begin{align}  \label{eq:upper-2period-Prod}
\mu^{1,1}_{+,\mytext{Prod}} = \max_{Q} \; \mu^{1,1}_{\mytext{ICE}} (\lambda^*_{0,L_1,Q}\lambda^*_{0,Y,Q},\lambda^*_{1,Q})
= \max_{Q} \; \mu^{1,1}_{\mytext{IPW}} (\lambda^*_{0,L_1,Q}\lambda^*_{0,Y,Q},\lambda^*_{1,Q}),
\end{align}
over all possible distributions $Q$ for the full data such that the induced distribution of $Q$ on the observed data $(\overA_1,\overL_1,Y )$
coincides with the true distribution, and
the sensitivity ratios $\lambda^*_{1,Q}$, $\lambda^*_{0,Y,Q}$, and $\lambda^*_{0,L_1,Q}$ based on $Q$ as in Lemma~\ref{lem:density-ratio-Prod}
satisfy the range constraints as $\lambda^*_1$, $\lambda^*_{0,L_1}$, and $\lambda^*_{0,Y}$ in (\ref{eq:model-2period-Prod}).
By the characterization of sensitivity ratios in Lemma~\ref{lem:density-ratio-Prod},
the following result gives an equivalent representation of $\mu^{1,1}_{+,\mytext{Prod}}$
through optimization over nonnegative functions $(\lambda_1, \lambda_{0,Y},\lambda_{0,L_1})$
depending only on the observed data $(\overline{A}_1, \overline{L}_1, Y)$.

\begin{lem} \label{lem:Q-lam-2period-Prod}
The sharp upper bound $\mu^{1,1}_{+,\mytext{Prod}}$ in (\ref{eq:upper-2period-Prod}) can be equivalently obtained as
\begin{align}  \label{eq:upper-2period-Prod-b}
\mu^{1,1}_{+,\mytext{Prod}} = \max_{\lambda_{0,L_1},\lambda_{0,Y},\lambda_1} \; \mu^{1,1}_{\mytext{ICE}} (\lambda_{0,L_1} \lambda_{0,Y},\lambda_1)
= \max_{\lambda_{0,L_1},\lambda_{0,Y},\lambda_1} \; \mu^{1,1}_{\mytext{IPW}} (\lambda_{0,L_1} \lambda_{0,Y},\lambda_1),
\end{align}
over all nonnegative functions $ \lambda_1 = \lambda_1 (\overline{L}_1, y) $,
$\lambda_{0,Y} = \lambda_{0,Y} (L_0, l_1, y)$, and
$\lambda_{0,L_1} = \lambda_{0,L_1} (L_0, l_1)$
subject to the normalization constraints (\ref{eq:constr-2period-Prod-b})
and the range constraints that for any $(l_1,y)$,
\begin{align}\label{eq:model-2period-Prod-b}
\begin{split}
& \Lambda_1^{-1} \le \lambda_1 (\overline{L}_1, y) \le \Lambda_1 ,  \\  
& \Lambda_{0,L_1}^{-1} \le \lambda_{0,L_1} (L_0, l_1) \le \Lambda_{0,L_1},  \quad
 \Lambda_{0,Y}^{-1} \le \lambda_{0,Y} (L_0, l_1, y) \le \Lambda_{0,Y} .  
\end{split}
\end{align}
\end{lem}

Compared with (\ref{eq:upper-2period-b}) and (\ref{eq:upper-2period-Jnt-b}), the constrained optimization in (\ref{eq:upper-2period-Prod-b})
is complicated by the coupling of $\lambda_{0,L_1}$ and $\lambda_{0,Y}$ in the ICE and IPW formulas.
In general, there may not exist an exact representation of $\mu^{1,1}_{+,\mytext{Prod}}$ through unconstrained optimization. 
Nevertheless, we show that a conservative bound on $\mu^{1,1}_{+,\mytext{Prod}}$ can be determined through unconstrained optimization.
For any functions $q_{0,Y} = q_{0,Y}(\overline{L}_1)$ and $q_{0,L_1} = q_{0,L_1}( L_0)$, let
\begin{align*}
 \eta_{0+,Y} (y, q_{0,Y} )
& = y + (\Lambda_{0,Y}-\Lambda_{0,Y}^{-1}) \rho_{\tau_{0,Y}} ( y, q_{0,Y}  )  , \\
 \eta_{0+,L_1} (y, \tilde y, q_{0,L_1} )
& = \pi^*_0 y + (1-\pi^*_0) \left\{ \tilde y + (\Lambda_{0,L_1}-\Lambda_{0,L_1}^{-1}) \rho_{\tau_{0,L_1}} (\tilde y, q_{0,L_1} ) \right\} ,
\end{align*}
where $\tau_{0,Y} = \Lambda_{0,Y}/(1+\Lambda_{0,Y})$
and $\tau_{0,L_1} = \Lambda_{0,L_1}/(1+\Lambda_{0,L_1})$.
These definitions are distinct from $\eta_{0+} (\cdot)$, but the following relationship holds:
$ \eta_{0+,L_1} (y, y, q_{0,L_1} )= \eta_{0+} (y, q_{0,L_1} )$ and, if $\Lambda_{0,L_1}=1$, then
$ \eta_{0+,L_1} (y, \eta_{0+,Y} (y, q_{0,Y} ), q_{0,L_1} ) = \eta_{0+} (y, q_{0,Y} )$, independently of $q_{0,L_1}$.

\begin{pro} \label{pro:upper-2period-Prod}
The sharp upper bound $\mu^{1,1}_{+,\mytext{Prod}}$ in (\ref{eq:upper-2period-Prod-b}) satisfies
\begin{align*}
\mu^{1,1}_{+,\mytext{Prod}} \le \min\left( \mu^{1,1}_{+,\mytext{Prod,v1}} ,\; \mu^{1,1}_{+,\mytext{Prod,v2}} \right),
\end{align*}
where $\mu^{1,1}_{+,\mytext{Prod,v1}}=\mu^{1,1}_{+,\mytext{Prod,v1}} (\Lambda_{0,L_1},\Lambda_{0,Y},\Lambda_1)$
and  $\mu^{1,1}_{+,\mytext{Prod,v2}} = \mu^{1,1}_{+,\mytext{Prod,v2}} (\Lambda_{0,L_1},\Lambda_{0,Y},\Lambda_1)$ are defined as follows:
\begin{align*}
& \mu^{1,1}_{+,\mytext{Prod,v1}} = \min_{q_{0,L_1},q_{0,Y},q_1}
 E \{ \mathcal E_{L_0} (  \mathcal E_{\overline{L}_1} [ \eta_{1+} \{ \eta_{0+,L_1} ( Y, \eta_{0+,Y} (Y, q_{0,Y} ), q_{0,L_1} ), q_1\} ] ) \} ,  \\
& \mu^{1,1}_{+,\mytext{Prod,v2}} = \min_{q_{0,L_1},q_{0,Y},q_1}
 E \{ \mathcal E_{L_0} ( \eta_{0+,L_1} [ \mathcal E_{\overline{L}_1} \{\eta_{1+} (Y, q_1)\}, \mathcal E_{\overline{L}_1} \{\eta_{1+} ( \eta_{0+,Y} (Y, q_{0,Y}), q_1 )\},
q_{0,L_1} ] ) \} ,
\end{align*}
each over all possible functions $q_{0,L_1} = q_{0,L_1}(L_0)$, $q_{0,Y}=q_{0,Y}(\overline{L}_1)$,  and $q_1= q_1(\overline{L}_1)$.
\end{pro}

The conservative bound on $\mu^{1,1}_{+,\mytext{Prod}} $ in Lemma \ref{pro:upper-2period-Prod} is
the conjunction of two distinct conservative bounds, $\mu^{1,1}_{+,\mytext{Prod,v1}}$ and $\mu^{1,1}_{+,\mytext{Prod,v2}}$,
each of which is in general a conservative bound on $\mu^{1,1}_{+,\mytext{Prod}} $, based on unconstrained optimization.
However, the individual bound $\mu^{1,1}_{+,\mytext{Prod,v1}}$ or $\mu^{1,1}_{+,\mytext{Prod,v2}}$
becomes exact in the special case of $\Lambda_{0,L_1}=1$ or $\Lambda_{0,Y}=1$ respectively.
The two cases can be viewed as two restricted versions of the product sensitivity model.

\begin{cor} \label{cor:upper-2period-Prod}
(i) If $\Lambda_{0,L_1}=1$ (i.e., $L_1^{1,1} \perp A_0 \,|\, L_0 $) and $\Lambda_{0,Y}=\Lambda_0$, then $\mu^{1,1}_{+,\mytext{Prod}}$ can be determined as
\begin{align*}
& \mu^{1,1}_{+,\mytext{Prod}} = \mu^{1,1}_{+,\mytext{Prod,v1}} (1,\Lambda_0,\Lambda_1) = \min_{ q_{0,Y},q_1 }
 E \{ \mathcal E_{L_0} (  \mathcal E_{\overline{L}_1} [ \eta_{1+} \{ \eta_{0+} ( Y,  q_{0,Y} ), q_1\} ] ) \} ,
\end{align*}
over all possible functions $q_{0,Y}=q_{0,Y}(\overline{L}_1)$ and $q_1= q_1(\overline{L}_1)$,
where  $ \eta_{0+} (y, q_{0,Y} )$ is defined as $\eta_{0+} (y,q_0)$
in Proposition \ref{pro:upper-2period} with $q_0$ replaced by $q_{0,Y}$.\\
(ii) If $\Lambda_{0,Y} =1$ (i.e., $Y^{1,1} \perp A_0 \,|\, L_0, L_1^{1,1}$) and $\Lambda_{0,L_1}=\Lambda_0$, then $\mu^{1,1}_{+,\mytext{Prod}}$ can be determined as
\begin{align*}
& \mu^{1,1}_{+,\mytext{Prod}} =  \mu^{1,1}_{+,\mytext{Prod,v2}} (\Lambda_0,1,\Lambda_1) = \min_{q_{0,L_1},q_1}
 E \{ \mathcal E_{L_0} ( \eta_{0+} [ \mathcal E_{\overline{L}_1} \{\eta_{1+} (Y, q_1)\}, q_{0,L_1} ] ) \} ,
\end{align*}
over all possible functions $q_{0,L_1} = q_{0,L_1}(L_0)$ and $q_1= q_1(\overline{L}_1)$,
where $ \eta_{0+} (y, q_{0,L_1} )$ is defined as $\eta_{0+} (y,q_0)$
in Proposition \ref{pro:upper-2period} with $q_0$ replaced by $q_{0,L_1}$.
\end{cor}

We make some remarks to help understand our results. First, model (\ref{eq:model-2period-Prod}) with $\Lambda_{0,L_1}=1$ assumes that
$L_1^{1,1}$ is unconfounded with $A_0$ given $L_0$, but allows for unmeasured confounding of $Y^{1,1}$ with $A_0$ given $(L_0,L_1^{1,1})$
as well as that of $Y^{1,1}$ with $A_1$ given $(A_0=1, \overline{L}_1)$.
The sharp upper bound on $\mu^{1,1}$ in this case is no greater than that under primary sensitivity model (\ref{eq:model-2period})
or joint sensitivity model (\ref{eq:model-2period-Jnt}), after $\Lambda_{0,Y} = \Lambda_0$ being matched:
\begin{align*}
\mu^{1,1}_+ (\Lambda_0,\Lambda_1) = \mu^{1,1}_{+,\mytext{Jnt}} (\Lambda_0,\Lambda_1) \ge \mu^{1,1}_{+,\mytext{Prod,v1}} (1,\Lambda_0,\Lambda_1) ,
\end{align*}
where the dependency of $\mu^{1,1}_+$ and $\mu^{1,1}_{+,\mytext{Jnt}}$ on $(\Lambda_0,\Lambda_1) $ is made explicit in the notation.
This relationship can be deduced from the fact that model  (\ref{eq:model-2period-Prod}) with $(\Lambda_{0,L_1},\Lambda_{0,Y},\Lambda_1) = (1, \Lambda_0,\Lambda_1)$
is more restrictive than model (\ref{eq:model-2period}) or model (\ref{eq:model-2period-Jnt}) with parameters $(\Lambda_0,\Lambda_1)$.
The preceding inequality can also be established directly:
\begin{align}
\mu^{1,1}_+ (\Lambda_0,\Lambda_1) & = \min_{ q_0,q_1 }
 E \{ \mathcal E_{L_0} (  \mathcal E_{\overline{L}_1} [ \eta_{1+} \{ \eta_{0+} ( Y,  q_0 ), q_1\} ] ) \} \nonumber \\
& \ge \min_{ q_{0,Y},q_1 }
 E \{ \mathcal E_{L_0} (  \mathcal E_{\overline{L}_1} [ \eta_{1+} \{ \eta_{0+} ( Y,  q_{0,Y} ), q_1\} ] ) \} = \mu^{1,1}_{+,\mytext{Prod,v1}} (1,\Lambda_0,\Lambda_1).
\label{eq:comparison-Jnt-Prod1}
\end{align}
Any function $q_0 = q_0 (L_0)$ in the first line is allowed to be $q_{0,Y} = q_{0,Y} (\overline{L}_1)$ in the second line
and hence the minimum in the first line is no greater than in the second line.
This inequality reveals an interesting consequence of restricting $\Lambda_{0,L_1}=1$ for our sensitivity bounds.

Second, model (\ref{eq:model-2period-Prod}) with $\Lambda_{0,Y}=1$ assumes that
$Y^{1,1}$ is unconfounded with $A_0$ given $(L_0,L_1^{1,1})$, but allows for unmeasured confounding of $L_1^{1,1}$ with $A_0$ given $L_0$
as well as that of $Y^{1,1}$ with $A_1$ given $(A_0=1, \overline{L}_1)$.
The sharp upper bound on $\mu^{1,1}$ in this case is no greater than that under primary sensitivity model (\ref{eq:model-2period})
or joint sensitivity model (\ref{eq:model-2period-Jnt}), after $\Lambda_{0,L_1} = \Lambda_0$ being matched:
\begin{align}
\mu^{1,1}_+ (\Lambda_0,\Lambda_1) & = \min_{ q_0,q_1 }
 E \{ \mathcal E_{L_0} (  \mathcal E_{\overline{L}_1} [ \eta_{1+} \{ \eta_{0+} ( Y,  q_0 ), q_1\} ] ) \} \nonumber \\
& \ge \min_{q_{0,L_1},q_1}
 E \{ \mathcal E_{L_0} ( \eta_{0+} [ \mathcal E_{\overline{L}_1} \{\eta_{1+} (Y, q_1)\}, q_{0,L_1} ] ) \}
 = \mu^{1,1}_{+,\mytext{Prod,v2}} (\Lambda_0,1,\Lambda_1),
\label{eq:comparison-Jnt-Prod2}
\end{align}
because model  (\ref{eq:model-2period-Prod}) with $(\Lambda_{0,L_1},\Lambda_{0,Y},\Lambda_1) = (\Lambda_0,1,\Lambda_1)$
is more restrictive than model (\ref{eq:model-2period}) or model (\ref{eq:model-2period-Jnt}) with parameters $(\Lambda_0,\Lambda_1)$.
Compared with (\ref{eq:comparison-Jnt-Prod1}), inequality (\ref{eq:comparison-Jnt-Prod2}) is more intriguing.
For an independent verification, we provide in Supplement Section \ref{sec:prf-comparison-Jnt-Prod} a direct but lengthy proof of (\ref{eq:comparison-Jnt-Prod2})
without viewing the two sides of (\ref{eq:comparison-Jnt-Prod2}) as sensitivity bounds.
Nevertheless, a useful simplification from restricting $\Lambda_{0,Y}=1$ is that
the sharp upper bound $\mu^{1,1}_{+,\mytext{Prod,v2}} (\Lambda_0,1,\Lambda_1)$ can be sequentially computed as follows:
\begin{align*}
& \check\eta_{1+} ( \overline{L}_1) = \min_{q_1} \; E \{ \eta_{1+} (Y, q_1) | A_0=1,A_1=1,\overline{L}_1\} , \\
& \check\eta_{0+} (L_0) = \min_{q_0} \; E \{ \eta_{0+} ( \check\eta_{1+} ( \overline{L}_1), q_0)  | A_0=1, L_0 \} ,\\
& \mu^{1,1}_{+,\mytext{Prod,v2}} (\Lambda_0,1,\Lambda_1) = E ( \check\eta_{0+} (L_0) ) .
\end{align*}
By the definition of $\eta_{0+}$ and $\eta_{1+}$, 
each optimization step above amounts to finding conditional quantiles through  weighted quantile regression. 
In contrast, neither the sharp upper bound $\mu^{1,1}_+ (\Lambda_0,\Lambda_1)$ nor 
$\mu^{1,1}_{+,\mytext{Prod,v1}} (1,$ $\Lambda_0,\Lambda_1)$ can in general be computed
in such a sequential manner.

\section{Sensitivity analysis for $K$ periods} \label{sec:Kperiod}

We present sensitivity models and bounds for general $K \ge 2$ periods.
While the basic ideas are similar as in Section \ref{sec:2period},
the development requires additional considerations, including the discovery of monotonicity and convexity of the recursive functions in the sensitivity bounds
and careful use of induction arguments in the proofs.

\subsection{Primary sensitivity model and bounds} \label{sec:Kperiod-primary}

To relax the primary sequential unconfounding Assumption A2, we consider the following density ratios as sensitivity ratios for $k=0,1,\ldots,K-1$:
\begin{align*}
 \lambda^*_k (\overL_k, y) & =  \frac{\dif P_{Y^{\over1}} (y | A_k=0, \overA_{k-1}=\over1_{k-1},  \overL_k) }
 {\dif P_{Y^{\over1}} (y | A_k=1, \overA_{k-1}=\over1_{k-1}, \overL_k) },
\end{align*}
where $\overA_{-1} = \over1_{-1}$ can be dropped from the conditioning.
By Bayes' rule, each $\lambda^*_k$ can be equivalently expressed as an odds ratio:
\begin{align*}
 \lambda^*_k (\overL_k, Y^{\over1})
  & =  \frac{ \pi^*_k ( \overline{L}_k ) }{ 1-\pi^*_k (\overline{L}_k ) } \times
  \frac{ P(A_k=0 | \overA_{k-1}=\over1_{k-1}, \overL_k, Y^{\over1} ) }{ P(A_k=1 | \overA_{k-1}=\over1_{k-1}, \overL_k, Y^{\over1}  ) } .
\end{align*}
In the case of $K=2$, these expressions reduce to (\ref{eq:lam-def}) and (\ref{eq:lam-def-odds}) in Section \ref{sec:2period-primary}.

Similarly as in Lemma~\ref{lem:density-ratio}, a necessary and sufficient condition for
any nonnegative functions $\lambda_k (\overL_k, y)$, $k=0,1,\ldots,K-1$,
to serve as sensitivity ratios, while being fully compatible with the observed data $(\overA_{K-1},\overL_{K-1},Y )$,
is that for  $k=0,1,\ldots,K-1$,
\begin{align} \label{eq:constr-Kperiod-b}
\Big( \prod_{j=k}^ {K-1} \calE_{\overL_j} \Big) \left\{ \lambda_k (\overL_k, Y )  \varrho^{k+1} ( \underL_{k+1}, Y; \underlam_{k+1} )
 \right\} \equiv 1 ,
\end{align}
where
$\varrho^k ( \underL_k, y; \underlam_k ) = \prod_{j=k}^ {K-1}  \varrho_j (\overL_j, y; \lambda_j )$,
$\varrho_k (\overL_k, y; \lambda_k)  = \pi^*_k (\overL_k) + (1 - \pi^*_k (\overL_k) ) \lambda_k ( \overL_k, y)$,
$\calE_{\overL_k} (\cdot)$ denotes the conditional expectation $E( \cdot | \overA_k=\over1_k, \overL_k)$ as an operator,
and $( \prod_{j=k}^{K-1}  \calE_{\overL_j} ) (\cdot) $ $ = \calE_{\overL_k} \cdots \calE_{\overL_{K-1}} (\cdot)$ by composition.
In addition, we set $ \varrho^K (\underL_K, y; \underlam_K ) \equiv 1 $.

If the sensitivity ratios $\lambda^*_k$, $k=0,1,\ldots,K-1$, were known, then $\mu^{\over1}$ can be point identified by
generalized ICE and IPW formulas as follows. See Lemma~\ref{lem:iden-mu-2period} for the case of $K=2$.
For any nonnegative functions $\lambda_k (\overL_k, y)$, $k=0,1,\ldots,K-1$, we define
\begin{align} \label{eq:g-iden-Kperiod}
\begin{split}
& \mu^{\over1}_{\mytext{ICE}} (\overlam_{K-1}) = E \Big[ \Big( \prod_{j=0}^ {K-1} \calE_{\overL_j} \Big) \left\{ \varrho^0 ( \underL_0, Y; \underlam_0 ) Y \right\} \Big],
 \\ 
& \mu^{\over1}_{\mytext{IPW}} (\overlam_{K-1}) = E \Big[
 \Big\{ \prod_{j=0}^ {K-1} \frac{A_j}{\pi^*_j(\overL_j)} \Big\} \varrho^0 ( \underL_0, Y; \underlam_0 ) Y \Big], 
\end{split}
\end{align}

\begin{lem} \label{lem:iden-mu-Kperiod}
Suppose that Assumptions A1 and A3 hold,
and $\lambda^*_k (\overL_k, Y^{\over1} )$ is finite almost surely for $k=0,1,\ldots,K-1$. Then
$\mu^{\over1} = \mu^{\over1}_{\mytext{ICE}} (\overlam^*_{K-1}) =\mu^{\over1}_{\mytext{IPW}} (\overlam^*_{K-1})$.
\end{lem}

For a sensitivity model, we postulate that
the sensitivity ratio $\lambda^*_k$
may differ from 1 by at most a factor of $\Lambda_k$ for any $k=0,1,\ldots,K-1$ and $y$:
\begin{align} \label{eq:model-Kperiod}
\Lambda_k^{-1} \le \lambda^*_k (\overL_k, y) \le \Lambda_k,
\end{align}
where $\Lambda_k \ge 1$ is a sensitivity parameter, encoding the degrees of unmeasured confounding in periods $k$.
By design,
(\ref{eq:model-Kperiod}) in the extreme case of $\Lambda_0 = \Lambda_1 =\cdots=\Lambda_{K-1}$ reduces to the primary sequential unconfounding Assumption A2.

For sensitivity analysis,
the sharp upper bound on $\mu^{\over1}$ under model (\ref{eq:model-Kperiod}) is defined as
\begin{align} \label{eq:upper-Kperiod}
\mu^{\over1}_+ = \max_{Q} \; \mu^{\over1}_{\mytext{ICE}} (\overlam^*_{K-1, Q})
= \max_{Q} \; \mu^{\over1}_{\mytext{IPW}} (\overlam^*_{K-1, Q}),
\end{align}
over all possible distributions $Q$ for the full data such that (i) the induced distribution of $Q$ on the observed data $(\overA_{K-1},\overL_{K-1},Y )$
coincides with the true distribution, and (ii)
the sensitivity ratios $\lambda^*_{k,Q}$, $k=0,1,\ldots,K-1$, defined as $\lambda^*_k$ with $P$ replaced by $Q$,
satisfy the range constraints as $\lambda^*_k$ in (\ref{eq:model-Kperiod}).
In other words, $\mu^{\over1}_+$ is the maximum value of $\mu^{\over1}$ over all possible choices of the true distribution $P$ (denoted as $Q$)
under the sensitivity model (\ref{eq:model-Kperiod}).
By the characterization of sensitivity ratios through the constraints (\ref{eq:constr-Kperiod-b}),
the optimization in (\ref{eq:model-Kperiod}) over probability distributions $Q$ can be equivalently
converted to optimization over nonnegative functions $\lambda_k (\overline{L}_k, y) $ subject to constraints (\ref{eq:constr-Kperiod-b})
and (\ref{eq:model-Kperiod-b}).

\vspace{-.05in}
\begin{lem} \label{lem:Q-lam-Kperiod}
The sharp upper bound $\mu^{\over1}_+$ in (\ref{eq:upper-Kperiod}) can be
equivalently obtained as  \vspace{-.05in}
\begin{align}  \label{eq:upper-Kperiod-b}
\mu^{\over1}_+ = \max_{\overlam_{K-1} = (\lambda_0,\ldots,\lambda_{K-1}) } \; \mu^{\over1}_{\mytext{ICE}} (\overlam_{K-1})
= \max_{\overlam_{K-1}= (\lambda_0,\ldots,\lambda_{K-1}) } \; \mu^{\over1}_{\mytext{IPW}} (\overlam_{K-1}),
\end{align}
over all nonnegative functions $ \lambda_k = \lambda_k (\overline{L}_k, y) $, $k=0,1,\ldots,K-1$,
subject to the normalization constraints (\ref{eq:constr-Kperiod-b})
and the range constraints that for any $k=0,1,\ldots,K-1$ and $y$,
\begin{align} \label{eq:model-Kperiod-b}
 \Lambda_k^{-1} \le \lambda_k ( \overL_k, y) \le \Lambda_k .
\end{align}
\end{lem}

Finally, as an extension of Proposition \ref{pro:upper-2period}, the following result shows that the constrained optimization in (\ref{eq:upper-Kperiod-b})
can be solved through a dual relationship to an unconstrained optimization problem.
For $k=0, 1, \ldots, K-1$, we define recursively
\begin{align} \label{eq:eta-def12}
\begin{split}
& \eta_{k+} (y, q_k ) = y + (1-\pi^*_k) (\Lambda_k-\Lambda_k^{-1}) \rho_{\tau_k} (y, q_k  ), \\  
& \eta^k_+ (y, \overq_k ) =
\left\{ \begin{array}{lc}
\eta_{0+} ( y, q_0 ) , & k=0,\\
\eta_{k+} ( \eta^{k-1}_+ (y, \overq_{k-1} ), q_k ), & k \ge 1,
\end{array} \right. 
\end{split}
\end{align}
where $\tau_k = \Lambda_k / (1+\Lambda_k)$ and $\rho_\tau(\cdot)$ is the ``check'' function associated with $\tau$-quantile regression
as in Section \ref{sec:2period}. The second line can be stated as
$ \eta^k_+ (y, \overq_k ) = (\eta_{k+}\circ \cdots \circ \eta_{1+}\circ \eta_{0+}) (y, \overq_k ) $ by successive composition,
where $\eta_{k+}$ is considered a function of $y$ for fixed $q_k$.
Fortuitously, as shown in Lemma \ref{lem:eta-convex}, the function $\eta^k_+ (y, \overq_k)$ defined in (\ref{eq:eta-def12}) is convex in $\overq_k$
for any fixed $y$ and increasing and convex in $y$ for any fixed $\overq_k$.
The convexity in $\overq_k$ is instrumental to the use of Sion's minimax theorem in our proof.
The increasingness in $y$ implies that $\mu^{\over1}_+$ is nondecreasing if the observed outcome increases almost surely.

\begin{pro} \label{pro:upper-Kperiod}
The sharp upper bound $\mu^{\over1}_+$ in (\ref{eq:upper-Kperiod-b}) can be determined as
\begin{align}  \label{eq:upper-Kperiod-c}
\begin{split}
\mu^{\over1}_+ &= \min_{\overq_{K-1} = (q_0,\ldots,q_{K-1}) }\; E \Big[ \Big( \prod_{j=0}^ {K-1} \calE_{\overL_j} \Big)
\left\{ \eta^{K-1}_+ (Y, \overq_{K-1} ) \right\} \Big] \\
& = \min_{\overq_{K-1} = (q_0,\ldots,q_{K-1}) }\; E \Big[
\Big\{ \prod_{j=0}^ {K-1} \frac{A_j}{\pi^*_j(\overL_j)} \Big\} \eta^{K-1}_+ (Y, \overq_{K-1} )  \Big],
\end{split}
\end{align}
over all possible functions $q_k= q_k(\overline{L}_k)$, $k=0,1,\ldots,K-1$.
Let $(\check q_0, \ldots, \check q_{K-1})$ be a solution to the optimization in (\ref{eq:upper-Kperiod-c}).
Then a solution, $(\check\lambda_0, \cdots, \check\lambda_{K-1})$,
to the optimization in (\ref{eq:upper-Kperiod-b})
can be obtained such that $ \check\lambda_k (\overL_k,y) = \Lambda_k$ if $  \eta^{k-1}_+ (y, \check{\overq}_{k-1} ) > \check q_k(L_k)$
 or $=\Lambda_k^{-1}$ if $ \eta^{k-1}_+ (y, \check{\overq}_{k-1} ) < \check q_k(L_k)$ for $k=0,1,\ldots, K-1$,
 where $\eta^{(-1)}_+ (y, \check{\overq}_{-1}) =y $ and $\check{\overq}_{k-1} = (\check q_0, \ldots, \check q_{k-1})$ for $k\ge 1$.
\end{pro}

\vspace{-.15in}
\subsection{Joint sensitivity model and bounds} \label{sec:Kperiod-Jnt}

To relax the joint sequential unconfounding Assumption A2$^\dag$, we consider the following density ratios as sensitivity ratios for $k=0,1,\ldots,K-1$:
\begin{align*}
 \lambda^*_{k,\mytext{Jnt}} (\overL_k, \underl_{k+1}, y) & =  \frac{\dif P_{ \underL^{\over1}_{k+1}, Y^{\over1}} (\underl_{k+1}, y | A_k=0, \overA_{k-1}=\over1_{k-1},  \overL_k) }
 {\dif P_{\underL^{\over1}_{k+1},  Y^{\over1}} (\underl_{k+1}, y | A_k=1, \overA_{k-1}=\over1_{k-1}, \overL_k) },
\end{align*}
where $\underl_K$ is set to the null.
The sensitivity ratio $\lambda^*_{k,\mytext{Jnt}}$ differs from $\lambda^*_k$ in Section \ref{sec:Kperiod-primary}
except for $k=K-1$, where $\lambda^*_{K-1,\mytext{Jnt}} = \lambda^*_{K-1}$.
By Bayes' rule, each $\lambda^*_{k,\mytext{Jnt}}$ can be equivalently expressed as an odds ratio:
\begin{align*}
 \lambda^*_{k,\mytext{Jnt}} (\overL_k, \underL^{\over1}_{k+1}, Y^{\over1})
  & =  \frac{ \pi^*_k ( \overline{L}_k ) }{ 1-\pi^*_k (\overline{L}_k ) } \times
  \frac{ P(A_k=0 | \overA_{k-1}=\over1_{k-1}, \overL_k,  \underL^{\over1}_{k+1}, Y^{\over1} ) }
  { P(A_k=1 | \overA_{k-1}=\over1_{k-1},\overL_k, \underL^{\over1}_{k+1}, Y^{\over1}  ) } .
\end{align*}
In the case of $K=2$, these expressions reduce to those for $\lambda^*_{0,\mytext{Jnt}}$ and $\lambda^*_1$ in Section \ref{sec:2period-Jnt}.

Similarly as in Lemma~\ref{lem:density-ratio-Jnt}, a necessary and sufficient condition for
any nonnegative functions $\lambda_{k,\mytext{Jnt}} (\overL_k, \underl_{k+1}, y)$, $k=0,1,\ldots,K-1$,
to serve as sensitivity ratios, while being fully compatible with the observed data $(\overA_{K-1},\overL_{K-1},Y )$,
is that for  $k=0,1,\ldots,K-1$,
\begin{align} \label{eq:constr-Kperiod-Jnt-b}
\Big( \prod_{j=k}^ {K-1} \calE_{\overL_j} \Big) \left\{ \lambda_{k,\mytext{Jnt}} (\overL_{K-1}, Y )
\varrho^{k+1}_{\mytext{Jnt}} ( \overL_{K-1}, Y; \underlam_{k+1,\mytext{Jnt}} )
 \right\} \equiv 1 ,
\end{align}
where
$\varrho^k_{\mytext{Jnt}} ( \overL_{K-1}, y; \underlam_{k,\mytext{Jnt}} )
= \prod_{j=k}^ {K-1}  \varrho_{j,\mytext{Jnt}} (\overL_{K-1}, y; \lambda_{j,\mytext{Jnt}} )$ and
$\varrho_{k,\mytext{Jnt}} (\overL_{K-1}, y; \lambda_{k,\mytext{Jnt}})
 = \pi^*_k (\overL_k) + (1 - \pi^*_k (\overL_k) ) \lambda_{k,\mytext{Jnt}} ( \overL_{K-1}, y)$.
In addition, we set $ \varrho^K_{\mytext{Jnt}} (\overL_{K-1}, y; \underlam_{K,\mytext{Jnt}} ) \equiv 1 $.

If the sensitivity ratios $\lambda^*_{k,\mytext{Jnt}}$, $k=0,1,\ldots,K-1$, were known, then $\mu^{\over1}$ can be point identified by
generalized ICE and IPW formulas which are of the same form as those in Lemma~\ref{lem:iden-mu-Kperiod}.
See Lemma~\ref{lem:iden-mu-2period-Jnt} for the case of $K=2$.

\begin{lem} \label{lem:iden-mu-Kperiod-Jnt}
Suppose that Assumptions A1$^\dag$ and A3 hold,
and $\lambda^*_{k,\mytext{Jnt}} (\overL_k, \underL^{\over1}_{k+1}, Y^{\over1} )$ is finite almost surely for $k=0,1,\ldots,K-1$. Then
$\mu^{\over1} = \mu^{\over1}_{\mytext{ICE}} (\overlam^*_{K-1,\mytext{Jnt}}) =\mu^{\over1}_{\mytext{IPW}} (\overlam^*_{K-1,\mytext{Jnt}})$,
where for nonnegative functions $\lambda_{k,\mytext{Jnt}} (\overL_k, \underl_{k+1}, y)$, $k=0,1,\ldots,K-1$, the functionals
$\mu^{\over1}_{\mytext{ICE}} (\overlam_{K-1,\mytext{Jnt}})$ and $\mu^{\over1}_{\mytext{IPW}} (\overlam_{K-1,\mytext{Jnt}})$
are defined as
$\mu^{\over1}_{\mytext{ICE}} (\overlam_{K-1})$ and $\mu^{\over1}_{\mytext{IPW}} (\overlam_{K-1})$ in (\ref{eq:g-iden-Kperiod})
except with  $\varrho^0 ( \underL_0, Y; $ $\underlam_0 ) $ replaced by
$ \varrho^0_{\mytext{Jnt}} ( \overL_{K-1}, Y; \underlam_{0,\mytext{Jnt}} ) $.
\end{lem}

For a sensitivity model, we postulate that
the sensitivity ratio $\lambda^*_{k,\mytext{Jnt}}$
may differ from 1 by at most a factor of $\Lambda_k$ for any $k=0,1,\ldots,K-1$ and $ (\underl_{k+1},y)$:
\begin{align} \label{eq:model-Kperiod-Jnt}
\Lambda_k^{-1} \le \lambda^*_{k,\mytext{Jnt}} ( \overL_k, \underl_{k+1}, y) \le \Lambda_k,
\end{align}
where $\Lambda_k \ge 1$ is a sensitivity parameter, encoding the degrees of unmeasured confounding in periods $k$.
By design,
(\ref{eq:model-Kperiod-Jnt}) in the extreme case of $\Lambda_0 = \Lambda_1 =\cdots=\Lambda_{K-1}$ reduces to the joint sequential unconfounding Assumption A2$^\dag$.

For sensitivity analysis,
the sharp upper bound on $\mu^{\over1}$ under model (\ref{eq:model-Kperiod-Jnt}) is defined as
\begin{align} \label{eq:upper-Kperiod-Jnt}
\mu^{\over1}_{+,\mytext{Jnt}} = \max_{Q} \; \mu^{\over1}_{\mytext{ICE}} (\overlam^*_{K-1,\mytext{Jnt}, Q})
= \max_{Q} \; \mu^{\over1}_{\mytext{IPW}} (\overlam^*_{K-1,\mytext{Jnt}, Q}),
\end{align}
over all possible distributions $Q$ for the full data such that (i) the induced distribution of $Q$ on the observed data $(\overline{L}_{K-1},\overline{A}_{K-1},Y )$
coincides with the true distribution, and (ii)
the sensitivity ratios $\lambda^*_{k,\mytext{Jnt},Q}$, $k=0,1,\ldots,K-1$, defined as $\lambda^*_{k,\mytext{Jnt}}$ with $P$ replaced by $Q$,
satisfy the range constraints as $\lambda^*_{k,\mytext{Jnt}}$ in (\ref{eq:model-Kperiod-Jnt}).
In other words, $\mu^{\over1}_+$ is the maximum value of $\mu^{\over1}$ over all possible choices of the true distribution $P$ (denoted as $Q$)
under the sensitivity model (\ref{eq:model-Kperiod-Jnt}).
Similarly as in Lemma~\ref{lem:Q-lam-2period-Jnt},
the optimization in (\ref{eq:upper-Kperiod-Jnt}) over probability distributions $Q$ can be converted
to optimization over nonnegative functions $\lambda_{k,\mytext{Jnt}} (\overline{L}_k, \underl_{k+1}, y) $.

\begin{lem} \label{lem:Q-lam-Kperiod-Jnt}
The sharp upper bound $\mu^{\over1}_{+,\mytext{Jnt}}$ in (\ref{eq:upper-Kperiod-Jnt}) can be
equivalently obtained as
\begin{align}  \label{eq:upper-Kperiod-Jnt-b}
\mu^{\over1}_{+,\mytext{Jnt}} = \max_{\overlam_{K-1,\mytext{Jnt}}  } \;
\mu^{\over1}_{\mytext{ICE}} (\overlam_{K-1,\mytext{Jnt}})
= \max_{\overlam_{K-1,\mytext{Jnt}} } \; \mu^{\over1}_{\mytext{IPW}} (\overlam_{K-1,\mytext{Jnt}}),
\end{align}
over all possible nonnegative functions $ \lambda_{k,\mytext{Jnt}} = \lambda_{k,\mytext{Jnt}} (\overline{L}_k, \underl_{k+1}, y) $, $k=0,1,\ldots,K-1$,
subject to the normalization constraints (\ref{eq:constr-Kperiod-Jnt-b})
and the range constraints that for any $k=0,1,\ldots,K-1$ and $(\underl_{k+1},y)$, \vspace{-.2in}
\begin{align} \label{eq:model-Kperiod-Jnt-b}
 \Lambda_k^{-1} \le \lambda_{k,\mytext{Jnt}} ( \overL_k, \underl_{k+1}, y) \le \Lambda_k .
\end{align}
\end{lem}

Finally, as an extension of Proposition \ref{pro:upper-2period-Jnt}, the following result shows that
the sharp bound on $\mu^{\over1}$ under joint sensitivity model (\ref{eq:model-Kperiod-Jnt}) is identical to
that under primary sensitivity model (\ref{eq:model-Kperiod}), even though model (\ref{eq:model-Kperiod-Jnt})
is more restrictive than model (\ref{eq:model-Kperiod}).

\begin{pro} \label{pro:upper-Kperiod-Jnt}
The sharp upper bound $\mu^{\over1}_{+,\mytext{Jnt}}$ in  (\ref{eq:upper-Kperiod-Jnt}) or (\ref{eq:upper-Kperiod-Jnt-b})
coincides with $\mu^{\over1}_+$ in (\ref{eq:upper-Kperiod}), (\ref{eq:upper-Kperiod-b}), or (\ref{eq:upper-Kperiod-c}),
i.e.,  $\mu^{\over1}_{+,\mytext{Jnt}} =\mu^{\over1}_+$.
The sharp bound $\mu^{\over1}_{+,\mytext{Jnt}}$ can be achieved by a probability distribution $Q$ for the full data such that
$ \underL_{k+1}^{\over1} \perp A_k \,|\, \overL_k, Y^{\over1}$ under $Q$.
\end{pro}

\subsection{Product sensitivity model and bounds} \label{sec:Kperiod-Prod}

For an alternative approach to relax joint sequential unconfounding in Section \ref{sec:Kperiod-Jnt},
we consider a decomposition of the joint sensitivity ratios for $k =0,\ldots,K-2$:
\begin{align*}
 \lambda^*_{k,\mytext{Jnt}} (\overL_k, \underl_{k+1}, y) & =\lambda^*_{k,L_{k+1}} (\overL_k, l_{k+1} ) \, \lambda^*_{k,Y} (\overL_k, \underl_{k+1}, y) ,
\end{align*}
where
\begin{subequations}
\begin{align*}
 \lambda^*_{k,L_{k+1}} (\overL_k, l_{k+1} ) & =  \frac{\dif P_{ L^{\over1}_{k+1} } ( l_{k+1} | A_k=0, \overA_{k-1}=\over1_{k-1},  \overL_k) }
 {\dif P_{ L^{\over1}_{k+1} } ( l_{k+1} | A_k=1, \overA_{k-1}=\over1_{k-1}, \overL_k) },  \\ 
 \lambda^*_{k,Y} (\overL_k, \underl_{k+1}, y) & =  \frac{\dif P_{ \underL^{\over1}_{k+2}, Y^{\over1}} (\underl_{k+2}, y | A_k=0, \overA_{k-1}=\over1_{k-1}, \overL_k, L^{\over1}_{k+1}=l_{k+1}  ) }
 {\dif P_{\underL^{\over1}_{k+2},  Y^{\over1}} (\underl_{k+2}, y | A_k=1, \overA_{k-1}=\over1_{k-1}, \overL_k, L^{\over1}_{k+1}=l_{k+1}) }. 
\end{align*}
\end{subequations}
For convenience, we also set
$\lambda^*_{K-1,L_K} (\overL_{K-1}, l_K )  \equiv 1$
and  $\lambda^*_{K-1,Y} (\overL_{K-1}, \underl_K, y) = \lambda^*_{K-1} (\overL_{K-1}, y)$.
By Bayes' rule, $\lambda^*_{k,L_{k+1}}$ and $\lambda^*_{k,Y}$ can be expressed as odds ratios:
\begin{align*}
 \lambda^*_{k,L_{k+1}} (\overL_k, L_{k+1}^{\over1}) & = \frac{ \pi^*_k( \overL_k ) }{ 1- \pi^*_k ( \overL_k ) } \times
 \frac{ P (A_k=0 | \overL_k, L_{k+1}^{\over1} ) }{ P (A_k=1 | \overL_k,  L_{k+1}^{\over1} ) } , \\
\lambda^*_{k,Y} (\overL_k, L_{k+1}^{\over1}, Y^{\over1} ) & = \frac{ P (A_k=1 |  \overA_{k-1}
=\over1_{k-1}, \overL_k, L_{k+1}^{\over1} ) }{ P (A_k=0 |  \overA_{k-1}=\over1_{k-1}, \overL_k,  L_{k+1}^{\over1} ) } \\
& \quad \times
 \frac{ P (A_k=0 |  \overA_{k-1}=\over1_{k-1}, \overL_k, \underL_{k+1}^{\over1}, Y^{\over1} ) }
 { P (A_k=1 |  \overA_{k-1}=\over1_{k-1}, \overL_k,  \underL_{k+1}^{\over1}, Y^{\over1} ) } .
\end{align*}
In the case of $K=2$, these expressions reduce to those for $\lambda^*_{0,L_1}$ and $\lambda^*_{0,Y}$ in Section \ref{sec:2period-Prod}.

Similarly as in Lemma~\ref{lem:density-ratio-Prod}, a necessary and sufficient condition for
any nonnegative functions $\{ (\lambda_{k,L_{k+1}} (\overL_k, l_{k+1} ),  \lambda_{k,Y} (\overL_k, \underl_{k+1}, y)): k=0,1,\ldots,K-1\}$
to serve as sensitivity ratios, while being fully compatible with the observed data $(\overA_{K-1},\overL_{K-1},Y )$,
is that for $k=0,1,\ldots,$ $K-1$,
\begin{align}\label{eq:constr-Kperiod-Prod-b}
\begin{split}
& \qquad \qquad \calE_{\overL_k}  \left\{ \lambda_{k,L} (\overL_{k+1} ) \right\}  \equiv 1,  \\  
& \Big( \prod_{j=k+1}^ {K-1} \calE_{\overL_j} \Big) \left\{ \lambda_{k,Y} (\overL_{K-1}, Y )  \varrho^{k+1}_{\mytext{Jnt}}  (\overL_{K-1}, Y; \underlam_{k+1,\mytext{Jnt}} )
 \right\} \equiv 1, 
\end{split}
\end{align}
where
$\varrho^k_{\mytext{Jnt}} ( \overL_{K-1}, y; \underlam_{k,\mytext{Jnt}} )$
are defined as in Section \ref{sec:Kperiod-Jnt}
with $\lambda_{k,\mytext{Jnt}} ( \overL_k, \underl_{k+1}, y) =\lambda_{k,L_{k+1}}  (\overL_k, $ $l_{k+1} ) \lambda_{k,Y}(\overL_k, \underl_{k+1}, y)$.
For convenience, we set $\lambda_{K-1,L_K} (\overL_{K-1}, l_K )  \equiv 1$.
As a consistency check, it can be easily verified that if $\{( \lambda_{k,L_{k+1}}, \lambda_{k,Y}) : k=0,1,\ldots,K-1\} $
satisfy (\ref{eq:constr-Kperiod-Prod-b}), then
$\{ \lambda_{k,\mytext{Jnt}} = \lambda_{k,L_{k+1}} \lambda_{k,Y}: k=0,1,\ldots,K-1\} $ satisfy (\ref{eq:constr-Kperiod-Jnt-b})
for joint sensitivity ratios.

For a sensitivity model, we postulate that
the sensitivity ratios $\lambda^*_{k,L_{k+1}}$ and $\lambda^*_{k,Y}$
may differ from 1 by at most a factor of $\Lambda_{k,L_{k+1}}$ and $\Lambda_{k,Y}$ for any $k=0,1,\ldots,K-1$ and $ (\underl_{k+1},y)$:
\begin{align} \label{eq:model-Kperiod-Prod}
\Lambda_{k,L_{k+1}}^{-1} \le \lambda^*_{k,L_{k+1}}  (\overL_k, l_{k+1} ) \le \Lambda_{k,L_{k+1}}, \quad
\Lambda_{k,Y}^{-1} \le \lambda^*_{k,Y} ( \overL_k, \underl_{k+1}, y) \le \Lambda_{k,Y},
\end{align}
where $\Lambda_{k,L_{k+1}} \ge 1$ and $\Lambda_{k,Y} \ge 1$ are sensitivity parameters,
with $\Lambda_{K-1, L_K} =1$.
In the case of $(\Lambda_{k,L_{k+1}},\Lambda_{k,Y})= (1,1)$, $k=0,1,\ldots,K-1$, model (\ref{eq:model-Kperiod-Prod}) reduces to
 joint sequential unconfounding assumption A2$^\dag$,
which is also recovered by
model (\ref{eq:model-Kperiod-Jnt}) with $\Lambda_0=\cdots=\Lambda_{K-1}=1$.
However, if $\Lambda_{j,L_{j+1}}>1$ or $ \Lambda_{j,Y} > 1$ for some $j$, then model (\ref{eq:model-Kperiod-Prod})
is more restrictive than model (\ref{eq:model-Kperiod-Jnt})
with $\Lambda_{k,\mytext{Jnt}} = \Lambda_{k,L_{k+1}} \Lambda_{k,Y}$ with $k=0,1,\ldots,K-1$.

For point identification of $\mu^{\over1}$ with known $(\lambda^*_{k,L_{k+1}}, \lambda^*_{k,Y})$, $k=0,1,\ldots,K-1$,
the ICE and IPW formulas  in Lemma~\ref{lem:iden-mu-Kperiod-Jnt}
remain valid by setting $\lambda^*_{k,\mytext{Jnt}} =\lambda^*_{k,L_{k+1}} \lambda^*_{k,Y}$.
Under sensitivity model (\ref{eq:model-Kperiod-Prod}), the sharp upper bound on $\mu^{\over1}$ is defined as
\begin{align}  \label{eq:upper-Kperiod-Prod}
\begin{split}
\mu^{\over1}_{+,\mytext{Prod}} & = \max_{Q} \; \mu^{\over1}_{\mytext{ICE}} (\lambda^*_{0,L_1,Q}\lambda^*_{0,Y,Q},\ldots,\lambda^*_{K-1,Y,Q}) \\
& = \max_{Q} \; \mu^{\over1}_{\mytext{IPW}} (\lambda^*_{0,L_1,Q}\lambda^*_{0,Y,Q},\ldots,\lambda^*_{K-1,Y,Q}),
\end{split}
\end{align}
over all possible distributions $Q$ for the full data such that the induced distribution of $Q$ on the observed data $(\overA_{K-1},\overL_{K-1},Y )$
coincides with the true distribution, and
the sensitivity ratios $\lambda^*_{k,L_{k+1},Q}$ and $\lambda^*_{k,Y,Q}$, defined as $\lambda^*_{k,L_{k+1}}$ and $\lambda^*_{k,Y}$
with $P$ replaced by $Q$,
satisfy the range constraints as $\lambda^*_{k,L_{k+1}}$ and $\lambda^*_{k,Y}$ in (\ref{eq:model-Kperiod-Prod}).
Similarly as in Lemma~\ref{lem:Q-lam-2period-Prod},
the optimization in (\ref{eq:upper-Kperiod-Prod}) over probability distributions $Q$ can be converted
to optimization over nonnegative functions  $\lambda_{k,L_{k+1}} (\overL_k, l_{k+1} )$ and
$\lambda_{k,Y}( \overL_k, \underl_{k+1}, y)$, $k=0,1,\ldots,K-1$.

\begin{lem} \label{lem:Q-lam-Kperiod-Prod}
The sharp upper bound $\mu^{\over1}_{+,\mytext{Prod}}$ in (\ref{eq:upper-Kperiod-Prod}) can be
equivalently obtained as
\begin{align}  \label{eq:upper-Kperiod-Prod-b}
\begin{split}
\mu^{\over1}_{+,\mytext{Prod}} & = \max_{\overlam_{K-1,L_K},\overlam_{K-1,Y}} \; \mu^{\over1}_{\mytext{ICE}}
(\lambda_{0,L_1}\lambda_{0,Y},\ldots,\lambda_{K-1,Y}) \\
& = \max_{\overlam_{K-1,L_K},\overlam_{K-1,Y}} \; \mu^{\over1}_{\mytext{IPW}} (\lambda_{0,L_1}\lambda_{0,Y},\ldots,\lambda_{K-1,Y}),
\end{split}
\end{align}
over all possible nonnegative functions $\lambda_{k,L_{k+1}} =\lambda_{k,L_{k+1}} (\overL_k, l_{k+1} )$ and
$\lambda_{k,Y}=\lambda_{k,Y}( \overL_k, \underl_{k+1}, y)$, $k=0,1,\ldots,K-1$,
subject to the normalization constraints (\ref{eq:constr-Kperiod-Prod-b})
and the range constraints that for any $k=0,1,\ldots,K-1$ and $(\underl_{k+1},y)$,
\begin{align} \label{eq:model-Kperiod-Prod-b}
\Lambda_{k,L_{k+1}}^{-1} \le \lambda_{k,L_{k+1}}  (\overL_k, l_{k+1} ) \le \Lambda_{k,L_{k+1}}, \quad
\Lambda_{k,Y}^{-1} \le \lambda_{k,Y} ( \overL_k, \underl_{k+1}, y) \le \Lambda_{k,Y} .
\end{align}
\end{lem}

Finally, as an extension of Proposition~\ref{pro:upper-2period-Prod}, the following result establishes two conservative bounds on
the sharp upper bound $\mu^{\over1}_{+,\mytext{Prod}}$, which are determined through unconstrained optimization.
For $k=0, 1, \ldots, K-1$, let
\begin{align*}
& \eta_{k+,Y} (y, q_{k,Y}) = y + (\Lambda_{k,Y} - \Lambda_{k,Y}^{-1}) \rho_{\tau_{k,Y}} (y, q_{k,Y} ),\\
& \eta_{k+,L_{k+1}} (y, \tilde y, q_{k,L_{k+1}} )
 = \pi^*_k y + (1-\pi^*_k)
 \left\{\tilde y  + (\Lambda_{k,L_{k+1}}-\Lambda_{k,L_{k+1}}^{-1}) \rho_{\tau_{k,L_{k+1}}} (\tilde y , q_{k,L_{k+1}} ) \right\} ,
\end{align*}
where $ \eta_{(K-1)+,L_K} (y, \tilde y, q_{K-1,L_K} )  =\pi^*_k y + (1-\pi^*_k) \tilde y$ with $\Lambda_{K-1,L_K}=1$ and $q_{K-1,L_K}$ set to the null.
%
%
For the first conservative bound, we define recursively
\begin{align*}
& \quad \eta^k_{+,\mytext{Prod,v1}} (y, \overq_{k,L_{k+1}}, \overq_{k,Y} )\\
& = \left\{ \begin{array}{ll}
\eta_{0+,L_1} (y, \eta_{0+,Y} (y, q_{0,Y}), q_{0,L_1} )  , & k=0,\\
\eta_{k+,L_{k+1}} (\eta^{k-1}_{+,\mytext{Prod,v1}} , \eta^{k-1}_{+,\mytext{Prod,v1}}\circ \eta_{k+,Y} , q_{k,L_{k+1}} ) , & 1\le k\le K-1 ,
\end{array} \right.
\end{align*}
where, for brevity,
$\eta^{k-1}_{+,\mytext{Prod,v1}} = \eta^{k-1}_{+,\mytext{Prod,v1}} (y, \overq_{k-1,L_k}, \overq_{k-1,Y} )$ and
$ \eta^{k-1}_{+,\mytext{Prod,v1}} \circ  \eta_{k+,Y} = \eta^{k-1}_{+,\mytext{Prod,v1}}$ $ (\eta_{k+,Y} (y,q_{k,Y}), \overq_{k-1,L_k}, \overq_{k-1,Y} )$.
For the second conservative bound, we define the following functional for any function $b_{K-1} = b_{K-1} (\overL_{K-1},Y)$:
\begin{align*}
& \quad \eta^k_{+,\mytext{Prod,v2}} (b_{K-1}, \underq_{k,L_{k+1}}, \underq_{k,Y} ) \\
& = \left\{ \begin{array}{ll}
 \eta_{(K-1)+}  (b_{K-1}, q_{K-1,Y} ) , &  k = K-1, \\
 \eta_{k+,L_{k+1}} ( \mathcal E_{\overL_{k+1}} \eta^{k+1}_{+,\mytext{Prod,v2}},
 \mathcal E_{\overL_{k+1}} (\eta^{k+1}_{+,\mytext{Prod,v2}} \circ \eta_{k+,Y} ), q_{k,L_{k+1}} ), & 0 \le k \le K-2 ,
\end{array} \right.
\end{align*}
where, for brevity,
$ \eta^{k+1}_{+,\mytext{Prod,v2}} = \eta^{k+1}_{+,\mytext{Prod,v2}}(b_{K-1}, \underq_{k+1,L_{k+2}}, \underq_{k+1,Y} ) $ and
$\eta^{k+1}_{+,\mytext{Prod,v2}}\circ \eta_{k+,Y} = \eta^{k+1}_{+,\mytext{Prod,v2}} ( \eta_{k+,Y} ( b_{K-1}, q_{k,Y} ), \underq_{k+1,L_{k+2}}, \underq_{k+1,Y} )$.
For comparison of these definitions, $\eta^k_{+,\mytext{Prod,v1}}$ is defined forward for $k=0,1,\ldots,K-1$ as a function of $(y, \overq_{k,L_{k+1}}, \overq_{k,Y} )$
which can then be applied to the random variables $(\overL_{K-1},Y)$,
whereas $\eta^k_{+,\mytext{Prod,v2}}$ is defined backward for $k=K-1,\ldots,1,0$,
directly as a functional of $(\overL_{K-1},Y)$ with conditional expectations taken.
Similarly as $\eta^k_+ (y,\overq_k)$ in Lemma~\ref{lem:eta-convex}, $\eta^k_{+,\mytext{Prod,v1}} (y, \overq_{k,L_{k+1}}, \overq_{k,Y} )$ can be shown to be
convex in $(\overq_{k,L_{k+1}}, \overq_{k,Y} )$ and increasing and convex in $y$.
As shown in Lemma~\ref{lem:eta2-convex},  $\eta^k_{+,\mytext{Prod,v2}} (b_{K-1}, \underq_{k,L_{k+1}}, \underq_{k,Y} ) $
is convex in $(\underq_{k,L_{k+1}}, \underq_{k,Y} )$ and nondecreasing and convex in $b_{K-1}$.
Due to the expectations involved, $\eta^k_{+,\mytext{Prod,v2}}$ may not be strictly increasing in $b_{K-1}$.

\begin{pro} \label{pro:upper-Kperiod-Prod}
The sharp upper bound $\mu^{\over1}_{+,\mytext{Prod}}$ in (\ref{eq:upper-Kperiod-Prod-b}) satisfies
\begin{align*}
\mu^{\over1}_{+,\mytext{Prod}} \le \min\left( \mu^{\over1}_{+,\mytext{Prod,v1}} ,\; \mu^{\over1}_{+,\mytext{Prod,v2}} \right),
\end{align*}
where $\mu^{\over1}_{+,\mytext{Prod,v1}}=\mu^{\over1}_{+,\mytext{Prod,v1}} (\overLam_{K-1,L_K},\overLam_{K-1,Y})$
and  $\mu^{\over1}_{+,\mytext{Prod,v2}} = \mu^{\over1}_{+,\mytext{Prod,v2}} (\overLam_{K-1,L_K},\overLam_{K-1,Y})$ are defined as follows:
\begin{align*}
& \mu^{\over1}_{+,\mytext{Prod,v1}} = \min_{\overq_{K-1,L_K},\overq_{K-1,Y}}
 E \Big\{  \Big( \prod_{j=0}^{K-2} \calE_{\overL_j} \Big)  \eta^{K-1}_{+,\mytext{Prod,v1}} (Y, \overq_{K-1,L_K}, \overq_{K-1,Y} ) \Big\} \\
& \mu^{\over1}_{+,\mytext{Prod,v2}} = \min_{\underq_{0,L_1}, \underq_{0,Y}}
 E \left[ \mathcal E_{L_0} \left\{ \eta^0_{+,\mytext{Prod,v2}} (Y, \underq_{0,L_1}, \underq_{0,Y} ) \right\} \right],
\end{align*}
each over all possible functions $q _{k,L_{k+1}} =q _{k,L_{k+1}} (\overL_k )$ and
$q _{k,Y}=q _{k,Y}( \overL_{K-1} )$, $k=0,1,\ldots,K-1$, with $q_{K-1,L_K}$ set to the null.
\end{pro}

In the case of $K=2$, the two bounds
$\mu^{\over1}_{+,\mytext{Prod,v1}} ( (\Lambda_{0,L_1},1), (\Lambda_{0,Y},\Lambda_1) )$
and $\mu^{\over1}_{+,\mytext{Prod,v2}} $ $( (\Lambda_{0,L_1},1), (\Lambda_{0,Y},\Lambda_1) )$, with $\Lambda_{1,L_2}=1$  and $\Lambda_{1,Y}=\Lambda_1$,
reduce to $\mu^{\over1}_{+,\mytext{Prod,v1}} ( \Lambda_{0,L_1},\Lambda_{0,Y},\Lambda_1 )$
and $\mu^{\over1}_{+,\mytext{Prod,v2}} ( \Lambda_{0,L_1}, \Lambda_{0,Y},\Lambda_1 )$ in Proposition \ref{pro:upper-2period-Prod}
after suppressing $\Lambda_{1,L_2}=1$ in the notation. Similarly as in Corollary~\ref{cor:upper-Kperiod-Prod},
the two bounds in Proposition \ref{pro:upper-Kperiod-Prod}, although conservative in general, become exact in two special cases respectively.

\begin{cor} \label{cor:upper-Kperiod-Prod}
(i) If $\Lambda_{k,L_{k+1}}=1$ (i.e., $L_{k+1}^{\over1} \perp A_k \,|\, \overA_{k-1}=\over1_{k-1}, \overL_k $)
and $\Lambda_{k,Y} = \Lambda_k$ for $k=0,1,\ldots,K-1$, then $\mu^{\over1}_{+,\mytext{Prod}}$ can be determined as
\begin{align*}
& \mu^{\over1}_{+,\mytext{Prod}} = \mu^{\over1}_{+,\mytext{Prod,v1}} (\over1_{K-1},\overLam_{K-1})
= \min_{\overq_{K-1,Y} }  E \Big[ \Big( \prod_{j=0}^ {K-1} \calE_{\overL_j} \Big)
\left\{ \eta^{K-1}_+ (Y, \overq_{K-1,Y} ) \right\} \Big],
\end{align*}
over all possible functions $q _{k,Y}=q _{k,Y}( \overL_{K-1} )$, $k=0,1,\ldots,K-1$,
where  $ \eta^{K-1}_+ (y, \overq_{K-1,Y} )$ is defined as $\eta^{K-1}_+ (y, \overq_{K-1})$
in Proposition \ref{pro:upper-Kperiod} with $\overq_{K-1}$ replaced by $\overq_{K-1,Y}$.\\
(ii) If $\Lambda_{k,Y} =1$ (i.e., $Y^{\over1} \perp A_k \,|\, \overA_{k-1}=\over1_{k-1}, \overL_k, \overL_{k+1}^{\over1}$)
and $\Lambda_{k,L_{k+1}} = \Lambda_k$ for $k=0,1,\ldots,K-2$ and $\Lambda_{K-1,Y}=\Lambda_{K-1}$,
then $\mu^{\over1}_{+,\mytext{Prod}}$ can be determined as
\begin{align*}
& \mu^{\over1}_{+,\mytext{Prod}} =  \mu^{\over1}_{+,\mytext{Prod,v2}} ( (\overLam_{K-2},1), (\over1_{K-2},\Lambda_{K-1}) )
= \min_{ \underq_{0,L_1}, q_{K-1,Y} }
 E \left\{ \eta^0_{+,\mytext{v2}} (Y, \underq_{0,L_1}, q_{K-1,Y} ) \right\},
\end{align*}
over all possible functions $q _{k,L_{k+1}} =q _{k,L_{k+1}} (\overL_k )$, $k=0,1,\ldots,K-2$, and $q_{K-1,Y} = q_{K-1,Y} (\overL_{K-1})$,
where $\eta^0_{+,\mytext{v2}} (Y, \underq_{0,L_1}, q_{K-1,Y} )$ is defined recursively as follows:
\begin{align*}
& \quad \eta^k_{+,\mytext{v2}} (Y, \underq_{k,L_{k+1}}, q_{K-1,Y} ) \\
& = \left\{ \begin{array}{ll}
 \eta_{(K-1)+}  (Y, q_{K-1,Y} ) , &  k = K-1, \\
 \eta_{k+} ( \mathcal E_{\overL_{k+1}} \eta^{k+1}_{+,\mytext{v2}}(Y, \underq_{k+1,L_{k+2}}, q_{K-1,Y} ) , q_{k,L_{k+1}} ), & 0 \le k \le K-2 ,
\end{array} \right.
\end{align*}
and $\eta_{k+} (y, q_{k,L_{k+1}})$ or $\eta_{(K-1)+} (y, q_{K-1,Y})$  is defined as $\eta_{k+} (y, q_k)$ or $\eta_{(K-1)+} (y, q_{K-1})$
in Proposition \ref{pro:upper-Kperiod}, with $q_k$ replaced by $q_{k,L_{k+1}}$ or $q_{K-1}$ by $q_{K-1,Y}$ respectively.
\end{cor}

\section{Numerical study} \label{sec:numerical}

We present a numerical study with 2 periods to demonstrate various properties of the sharp upper bounds
$\mu^{1,1}_+$, $\mu^{1,1}_{+,\mytext{Prod,v1}} (1,\Lambda_0,\Lambda_1)$, and $\mu^{1,1}_{+,\mytext{Prod,v2}} (\Lambda_0,1,\Lambda_1)$
in Proposition~\ref{pro:upper-2period} and Corollary~\ref{cor:upper-2period-Prod},
including the dependency on the observed-data configuration, the comparison of the three sensitivity bounds,
and the construction of the implied distribution of $Y^{1,1}$.
Computational details are provided in Supplement Section \ref{sec:comp-details}.

\begin{figure}
\begin{minipage}[c]{.45\textwidth}
\includegraphics[width=3in, height=2.5in]{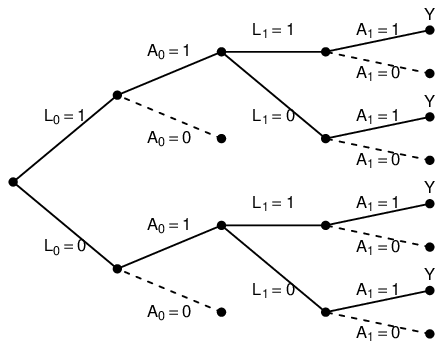} \vspace{-.3in}
\end{minipage}%
  \begin{minipage}[c]{0.6\textwidth} \scriptsize \vspace{-.2in}
  \begin{align*}
& P(L_0=1) = .5, \\
& P(A_0=1 |L_0=0) = P( A_0=1| L_0=1) = .5,\\
& P( L_1=1 | A_0=1,L_0=l_0) \\
& = \left\{ \begin{array}{ll}
.2, & l_0=0 ,\\
.8, & l_0=1 ,
\end{array} \right.\\
& P( A_1=1 | A_0=1, L_0, L_1=l_1 ) \\
& = \left\{ \begin{array}{ll}
 .4, &  l_1=0, \\
 .8, &  l_1=1,
\end{array} \right. \\
& Y \,|\, A_0= A_1=1, L_0=l_0, L_1= l_1 \sim \N ( \alpha_{l_0,l_1}, \sigma^2_{l_0}) .
 \end{align*}  \normalsize
 \end{minipage} \vspace{-.3in}
\caption{\small 
Data configuration in the numerical examples with 2 periods.
} \label{fig:data-diagram} 
\end{figure}

For our study, Figure~\ref{fig:data-diagram} gives the basic configuration of the observed data, which is modified from Table 21.2 in Hernan \& Robins (2022)
to include a baseline covariate $L_0$.
Our study is focused on sensitivity bounds for $\mu^{1,1}$, and hence the configuration is only shown for the observed data related to
the treatment sequence $(A_0,A_1) = (1,1)$.
Note that the configuration of the full data including counterfactual covariates and outcomes is unspecified, except for its implication on the observed data.
To accommodate various forms of heterogeneity, consider the following choices of the configuration. \vspace{-.1in}
\begin{itemize}\addtolength{\itemsep}{-.1in}
\item[(C1)] $\alpha_{l_0l_1} = 62, 58, 62, 58$ for $(l_0,l_1) =(0,0), (0,1), (1,0), (1,1)$, and $\sigma_0=\sigma_1 =1$.

\item[(C2)] $\alpha_{l_0l_1}$'s and $\sigma_{l_0}$'s are the same as in C1, but $P(A_0=1 |L_0=0) = P( A_0=1| L_0=1) = .25$
and $ P( A_1=1 | A_0=1, L_0, L_1=l_1 ) = .2, .4$ for $l_0=0,1$, i.e., the conditional probabilities of $A_0$ and $A_1$ are halved.

\item[(C3)] $\sigma_{l_0}$'s are the same as in C1, but $\alpha_{l_0l_1} = 62, 58, 66, 57$ for $(l_0,l_1) =(0,0), (0,1), (1,0), (1,1)$.

\item[(C4)] $\alpha_{l_0l_1}$'s are the same as in C1, but $\sigma_0=1$ and $\sigma_1=2$.
\end{itemize} \vspace{-.1in}
For each of the preceding configuration, it can be verified that under the assumption of sequential unconfounding, $\mu^{1,1}$ is point identified as
\begin{align*}
(.5)\left\{ (.8) \alpha_{00} + (.2) \alpha_{01} \right\}+(.5)\left\{ (.2) \alpha_{10} + (.8) \alpha_{11} \right\} = 60 .
\end{align*}
Hence $\mu^{1,1}$ remains the same under sequential unconfounding,
regardless of the specified differences in  $\alpha_{l_0l_1}$'s, $\sigma_{l_0}$'s, or the conditional probabilities of $A_0$ and $A_1$.
This is distinct from the fact that
the aforementioned choices may affect the sampling variation of an estimator for $\mu^{1,1}$ for any sample size $n$
and its asymptotic variance as $n\to\infty$ under sequential unconfounding.
However, by definition, the probability limit of any consistent estimator of $\mu^{1,1}$ remains $60$ for the preceding configurations
under sequential unconfounding.

\begin{figure} [t!]
\centering
\includegraphics[width=6in, height=2in]{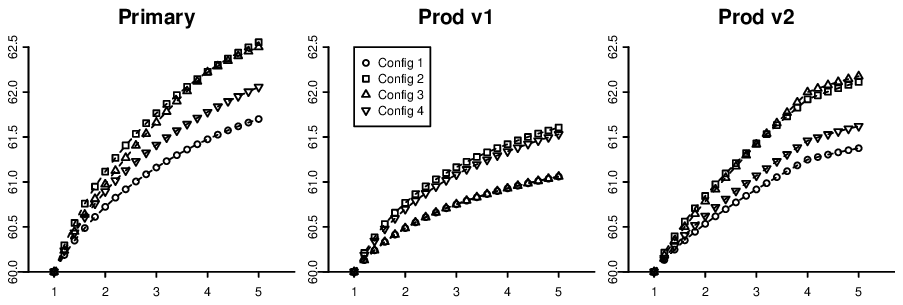} \vspace{-.1in}
\caption{\small 
Sharp upper bounds $\mu^{1,1}_+$ (Primary), $\mu^{1,1}_{+,\mytext{Prod,v1}} (1,\Lambda_0,\Lambda_1)$,
and $\mu^{1,1}_{+,\mytext{Prod,v2}} (\Lambda_0,1,\Lambda_1)$ over a range of $\Lambda_0 = \Lambda_1$ in four configurations, estimated using a sample size $10^8$.
} \label{fig:upper-bounds-method} 
\end{figure}

Figure~\ref{fig:upper-bounds-method} shows the sensitivity bounds
$\mu^{1,1}_+$, $\mu^{1,1}_{+,\mytext{Prod,v1}} (1,\Lambda_0,\Lambda_1)$, and $\mu^{1,1}_{+,\mytext{Prod,v2}} (\Lambda_0,1,\Lambda_1)$,
referred to as three methods, over a range of sensitivity parameter $\Lambda_0=\Lambda_1$ with the data configurations C1--C4.
The main findings can be summarized as follows. \vspace{-.1in}
\begin{itemize}\addtolength{\itemsep}{-.1in}
\item For each method, the sensitivity bounds exhibit notable differences between configurations C1--C4, even though all
sensitivity bounds reduce to the same value $60$ at $\Lambda_0=\Lambda_1=1$.
The sensitivity upper bounds in configurations C2--C4 are consistently larger than those in C1 at the same sensitivity parameter,
which corroborates our discussion after Proposition~\ref{pro:upper-2period}.
Compared with C1, the conditional probabilities $\pi^*_0$ and $\pi^*_1$ are reduced in C2,
indicating that there is a smaller probability of observing $A_0=A_1=1$ and hence observing $Y^{1,1}$.
In addition, the conditional means, $(\alpha_{10},\alpha_{11})$, of $Y$ given $A_0=A_1=1$, $L_0=1$, and $L_1=l_1$ for $l_1=0,1$ show greater heterogeneity in C2
(the larger of $ \alpha_{10}$ and $\alpha_{11}$ is increased to $66$ from $62$),
and the conditional standard deviations, $\sigma_1$, of $Y$ given $A_0=A_1=1$, $L_0=1$, and $L_1=l_1$ for $l_1=0,1$  are doubled in C3.
Both differences increase the probabilities of $Y^{1,1}$ taking larger values than in C1, hence leading to larger sensitivity upper bounds.

\item The increases of the upper bounds in configurations C2--C4 over C1 appear to vary among the three methods.
In particular, the upper bounds from $\mu^{1,1}_{+,\mytext{Prod,v1}} (1,\Lambda_0,\Lambda_1)$ under the restricted product sensitivity model
such that $L^{1,1}_1 \perp A_0 |L_0$ are virtually the same in C1 and C3, whereas
those from  $\mu^{1,1}_+$ under the primary sensitivity model and
from $\mu^{1,1}_{+,\mytext{Prod,v2}} (1,\Lambda_0,\Lambda_1)$ under the restricted product sensitivity model
such that $Y^{1,1} \perp A_0 |L_0, L^{1,1}_1$ show substantial differences between C1 and C3.
Further study is desired to understand this phenomenon and its implications.

\item The upper bounds from $\mu^{1,1}_+$ are consistently larger than those from
$\mu^{1,1}_{+,\mytext{Prod,v1}} (1,\Lambda_0,\Lambda_1)$, and $\mu^{1,1}_{+,\mytext{Prod,v2}} (\Lambda_0,1,\Lambda_1)$
at the same sensitivity parameter, as theoretically shown in (\ref{eq:comparison-Jnt-Prod1}) and (\ref{eq:comparison-Jnt-Prod2}).
See Supplement Figure~\ref{fig:upper-bounds-config} for the same results as Figure~\ref{fig:upper-bounds-method}, but re-grouped
to facilitate comparison between different methods.
\end{itemize} \vspace{-.1in}

\begin{sidewaysfigure} 
\includegraphics[angle=270, totalheight=6.5in]{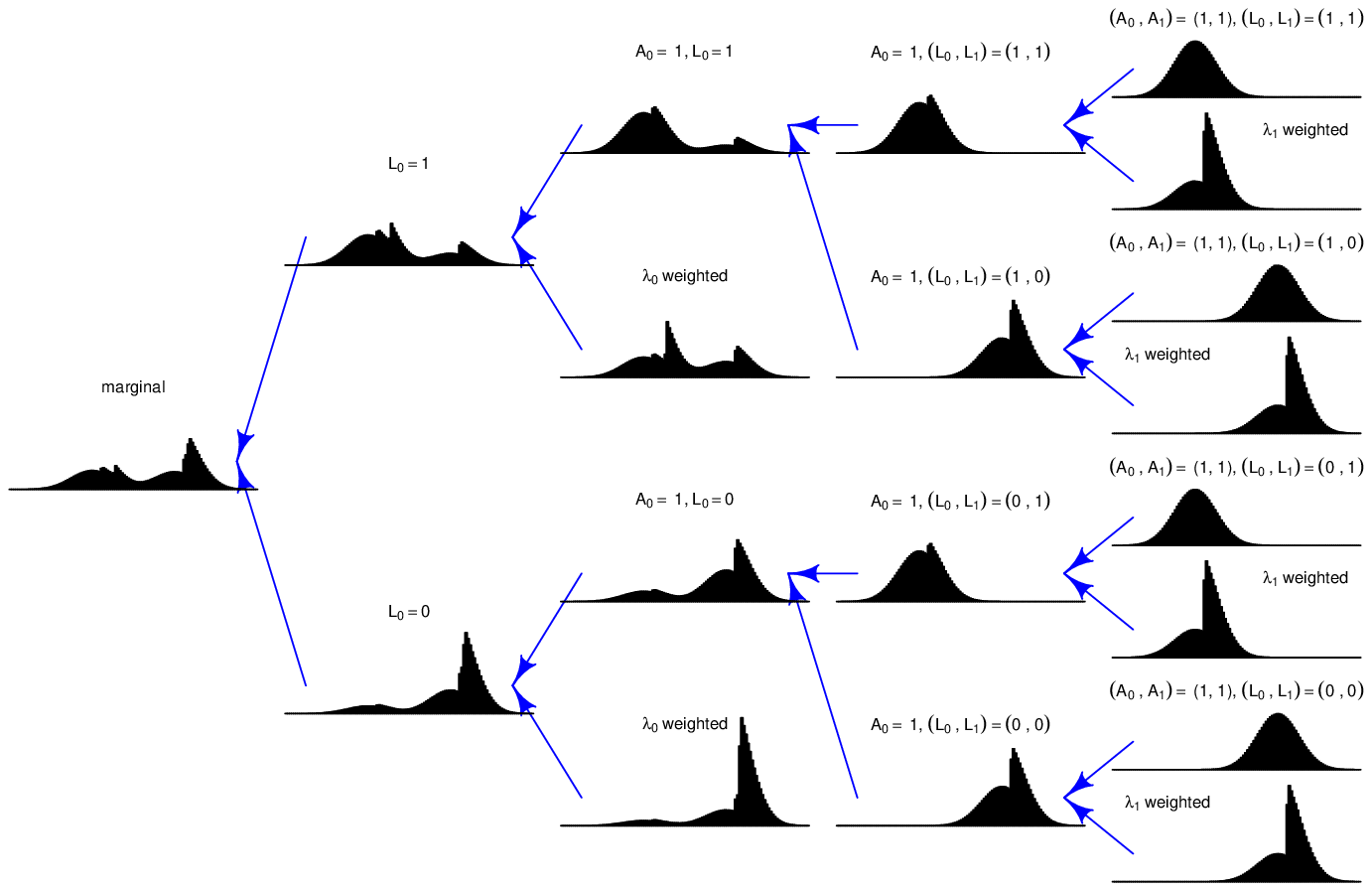} \vspace{-.1in}
\caption{\scriptsize
Distributions of $Y^{1,1}$ under various conditions (labeled on the top) based on the
sharp upper bound $\mu^{1,1}_+$ at $\Lambda_0=\Lambda_1=2$ in Proposition \ref{pro:upper-2period}, using a sample size $10^8$ in configuration C1.
On the rightmost are the observed distributions of $Y^{1,1}$ given $(A_0,A_1)=(1,1)$ and $(L_0,L_1)$ in the odd rows,
and the $\lambda_1$ weighted distributions given $(A_0,A_1)=(1,0)$ and $(L_0,L_1)$ in the even rows.
On the leftmost is the marginal distribution of $Y^{1,1}$. Each pair of arrows indicates a mixture of two distributions.
} \label{fig:hist-combined-sen} 
\end{sidewaysfigure}

Finally, to highlight the constructiveness of our methods, Figure~\ref{fig:hist-combined-sen} shows the conditional and marginal distributions
of $Y^{1,1}$ achieving the sharp upper bound $\mu^{1,1}_+$ at $\Lambda_0=\Lambda_1=2$ under the primary sensitivity model in configuration C1.
The distributions, approximated using a sample size $10^8$, are constructed step-by-step as in the proof of Lemma~\ref{lem:density-ratio},
based on the sensitivity ratios $(\check\lambda_0, \check\lambda_1)$ associated with the solution $(\check q_0,\check q_1)$
in Proposition~\ref{pro:upper-2period}.
The marginal distribution of $Y^{1,1}$ in Figure~\ref{fig:hist-combined-sen} yields the largest possible value of
$\mu^{1,1}= E (Y^{1,1})$ under the stated sensitivity model, while being fully compatible with the observed-data distribution (i.e., configuration C1).
For comparison, the corresponding distributions of $Y^{1,1}$ under sequential unconfounding are shown in Supplement Figure~\ref{fig:hist-combined-nosen}.

\section{Further discussion}  \label{sec:further-discussion}

\textbf{Sensitivity lower bounds.}\;
Our development is focused on sensitivity upper bounds, but can be readily extended to sensitivity lower bounds.
For example, for two periods, the sharp lower bound, $\mu^{1,1}_-$, on $\mu^{1,1}$ under model (\ref{eq:model-2period}) is defined as
(\ref{eq:upper-2period}) and Lemma~\ref{lem:Q-lam-2period} is applicable, provided
that $\mu^{1,1}_+$ is replaced by $\mu^{1,1}_-$ and the maximization is replaced by the minimization.
Moreover, Proposition~\ref{pro:upper-2period} can be extended as follows. For $k=0,1$, let
$\eta_{k-} (y,q_0) = y - (1-\pi^*_k) (\Lambda_k-\Lambda_k^{-1}) \rho_{1-\tau_k} ( y, q_k )$,
where $\tau_k = \Lambda_k/(1+\Lambda_k)$ as before.
Then the sharp lower bound $\mu^{1,1}_-$ can be determined as
\begin{align*}
\mu^{1,1}_- & = \min_{q_0,q_1}\; E \{ \mathcal E_{L_0} ( \mathcal E_{\overL_1} [ \eta_{1-} \{ \eta_{0-}(Y, q_0), q_1\} ] ) \} \\
& = \min_{q_0,q_1}\; E \left[ \frac{A_0A_1}{\pi^*_0 \pi^*_1} \eta_{1-} \{ \eta_{0-}(Y, q_0), q_1\}  \right],
\end{align*}
over all possible functions $q_0 = q_0(L_0)$ and $q_1= q_1(\overline{L}_1)$.

\vspace{.05in}
\textbf{Multiple treatment strategies.}\; Our sensitivity models and bounds are so far presented for the mean of
the counterfactual outcome under a single treatment strategy, for example, $\overline{a}=(1,1)$.
We discuss how our results can be extended to the means of counterfactual outcomes under multiple treatment strategies,
from which sensitivity bounds can be derived about ATEs.
For simplicity, our discussion deals with two periods, but the main ideas are applicable to general $K$ periods.

First, the primary sensitivity ratios for any single treatment strategy $\overa=(a_0,a_1)$ with $a_0,a_1 \in \{0,1\}$ are defined as
\begin{align*}
 \lambda^{\overa *}_1 (\overline{L}_1, y) & =  \frac{\dif P_{Y^{\overa}} (y | A_0=a_0, A_1=1-a_1, \overline{L}_1) }{\dif P_{Y^{\overa}} (y | A_0=a_0, A_1=a_1, \overline{L}_1) }, \\
 \lambda^{\overa *}_0 (L_0, y) & = \frac{\dif P_{Y^{\overa}} ( y | A_0=1-a_0, L_0) }{\dif P_{ Y^{\overa}} ( y | A_0=a_0, L_0 ) } .
\end{align*}
The previous definition (\ref{eq:lam-def}) corresponds to $\lambda^*_1 = \lambda^{1,1 *}_1$ and $\lambda^*_0 = \lambda^{1,1 *}_0$.
For a primary sensitivity model, we postulate that for each $\overa \in \{0,1\}\times \{0,1\}$,
\begin{align} \label{eq:model-2period-multi}
\Lambda_1^{-1} \le \lambda^{\overa *}_1 (\overline{L}_1, y) \le \Lambda_1, \quad
\Lambda_0^{-1} \le \lambda^{\overa *}_0 (L_0, y) \le \Lambda_0 ,
\end{align}
where $\Lambda_1 \ge 1$ and $\Lambda_0 \ge 1$ are sensitivity parameters, which may also vary over different $\overa$.
Model (\ref{eq:model-2period-multi}) differs from the previous model (\ref{eq:model-2period}) in that
the range constraints are placed simultaneously for all $\overa \in \{0,1\}\times \{0,1\}$, not just a single choice $\overa=(1,1)$.
The sharp upper bound on $\mu^{1,1}$ under model (\ref{eq:model-2period-multi}) is defined as
\begin{align} \label{eq:upper-2period-multi}
\mu^{1,1}_+ = \max_{Q} \; \mu^{1,1}_{\mytext{ICE}} (\lambda^{1,1 *}_{0,Q},\lambda^{1,1 *}_{1,Q})
= \max_{Q} \; \mu^{1,1}_{\mytext{IPW}} (\lambda^{1,1 *}_{0,Q},\lambda^{1,1 *}_{1,Q}),
\end{align}
over all possible distributions $Q$ for the full data as in the previous definition (\ref{eq:upper-2period}) except that
the sensitivity ratios $\lambda^{\overa *}_{1,Q}$ and $\lambda^{\overa *}_{0,Q}$ based on $Q$ are subject to
the range constraints (\ref{eq:model-2period-multi}) for all $\overa \in \{0,1\}\times \{0,1\}$, not just $\overa=(1,1)$.
The results in Section~\ref{sec:2period-primary} can be extended in a simultaneous manner, including
Lemmas~\ref{lem:density-ratio} and \ref{lem:Q-lam-2period} and Proposition~\ref{pro:upper-2period}.
For $\overline{t} \in \{0,1\}\times \{0,1\}$, the sharp upper bound, $\mu^{\overline{t}}_+$, on $\mu^{\overline{t}}$ under model (\ref{eq:model-2period})
or (\ref{eq:model-2period-multi})
is defined as $\mu^{1,1}_+$ in (\ref{eq:upper-2period}) or (\ref{eq:upper-2period-multi}) except with $(1,1)$ replaced by $\overline{t}$,
and the sharp lower bound $\mu^{\overline{t}}_+$ can be defined accordingly.
From Proposition~\ref{pro:upper-2period-multi}(ii),
$\mu^{1,1}_+ -\mu^{\overline{t}}_-$ is the sharp upper bound on the ATE between $(1,1)$ and $\overline{t}$,
and the sharp bounds can be simultaneously achieved for all $\overline{t} \not=(1,1)$.

\begin{pro} \label{pro:upper-2period-multi}
(i) The sharp upper bound $\mu^{1,1}_+$ in (\ref{eq:upper-2period-multi}) is identical to
$\mu^{1,1}_+$ in (\ref{eq:upper-2period}) and can be equivalently obtained as (\ref{eq:upper-2period-b}) or (\ref{eq:upper-2period-c}).
The same relationship holds for $\mu^{\overline{t}}_+$ under model (\ref{eq:model-2period}) and (\ref{eq:model-2period-multi}).
(ii) The sharp bounds $\{\mu^{\overline{t}}_+: \overline{t} \in\{0,1\}\times \{0,1\} \}$ can be simultaneously achieved by
a single choice of $Q$ allowed in (\ref{eq:upper-2period-multi}).
The result also holds with a different choice of $Q$ if $\mu^{\overline{t}}_+$ is replaced by the corresponding lower sharp bound $\mu^{\overline{t}}_-$
for some treatment strategies $\overline{t}$.
\end{pro}

Second, the joint sensitivity ratios for any single treatment strategy $\overa=(a_0,a_1)$ with $a_0,a_1 \in \{0,1\}$ are defined as
$\lambda^{\overa *}_1 (\overline{L}_1, y)$ and
\begin{align*}
\lambda^{\overa *}_{0,\mytext{Jnt}} (L_0, l_1, y) & = \frac{\dif P_{L_1^{\overa}, Y^{\overa}} (l_1, y | A_0=1-a_0, L_0) }
{\dif P_{L_1^{\overa}, Y^{\overa}} (l_1, y | A_0=a_0, L_0 ) } .
\end{align*}
For a joint sensitivity model, we postulate that for each $\overa \in \{0,1\}\times \{0,1\}$,
\begin{align} \label{eq:model-2period-Jnt-multi}
\Lambda_1^{-1} \le \lambda^{\overa *}_1 (\overline{L}_1, y) \le \Lambda_1, \quad
\Lambda_0^{-1} \le \lambda^{\overa *}_{0,\mytext{Jnt}} (L_0, l_1, y) \le \Lambda_0 ,
\end{align}
The sharp upper bound on $\mu^{1,1}$ under model (\ref{eq:model-2period-Jnt-multi}) is defined as
\begin{align}  \label{eq:upper-2period-Jnt-multi}
\mu^{1,1}_{+,\mytext{Jnt}} = \max_{Q} \; \mu^{1,1}_{\mytext{ICE}} (\lambda^{1,1 *}_{0,\mytext{Jnt},Q},\lambda^{1,1 *}_{1,Q})
= \max_{Q} \; \mu^{1,1}_{\mytext{IPW}} (\lambda^{1,1 *}_{0,\mytext{Jnt},Q},\lambda^{1,1 *}_{1,Q}),
\end{align}
over all possible distributions $Q$ for the full data as in the previous definition (\ref{eq:upper-2period-Jnt}) except that
$\lambda^{\overa *}_{1,Q}$ and $\lambda^{\overa *}_{0,\mytext{Jnt},Q}$ are subject to the range constraints (\ref{eq:model-2period-Jnt-multi})
for all $\overa \in \{0,1\}\times \{0,1\}$.
The results in Section~\ref{sec:2period-Jnt} can also be extended.
The sharp bounds on $\mu^{\overline{t}}$ under the primary and joint sensitivity models coincide for each treatment strategy $\overline{t}$,
and these bounds can be simultaneously attained for all $\overline{t}\in \{0,1\}\times\{0,1\}$.
A caveat is that Proposition~\ref{pro:upper-2period-Jnt-multi}(i) remains valid,
but the simultaneous attainment in Proposition~\ref{pro:upper-2period-Jnt-multi}(ii) may fail if it is required that
$L^{1,1}_1 = L^{1,0}_1$ regardless of $A_0=0$ or $1$, in addition to
$L^{1,1}_1 = L^{1,0}_1 =L_1$ if $A_0=1$ by Assumption A1$^\dag$. We defer further discussion to
Supplement Section \ref{sec:pro-upper-2period-Jnt-multi}.

\begin{pro} \label{pro:upper-2period-Jnt-multi}
(i) The sharp upper bound $\mu^{1,1}_{+,\mytext{Jnt}}$ in (\ref{eq:upper-2period-Jnt-multi})
is identical to $\mu^{1,1}_+$ in (\ref{eq:upper-2period-multi}) and hence to
$\mu^{1,1}_+$ in (\ref{eq:upper-2period}), (\ref{eq:upper-2period-b}), or (\ref{eq:upper-2period-c}),
and can be achieved by a probability distribution $Q$ for the full data such that
$ L_1^{1,1} \perp A_0 \,|\, L_0, Y^{1,1}$ under $Q$.
The result also hold if $(1,1)$ is replaced by any $\overline{t} \in\{0,1\}\times \{0,1\}$.
(ii) The sharp bounds $\{\mu^{\overline{t}}_{+,\mytext{Jnt}}: \overline{t} \in\{0,1\}\times \{0,1\} \}$ can be simultaneously achieved by
a single choice of $Q$ allowed in (\ref{eq:upper-2period-Jnt-multi}), where
$\mu^{\overline{t}}_{+,\mytext{Jnt}}$ may be replaced by the corresponding lower sharp bound
for some treatment strategies $\overline{t}$.
\end{pro}

Third, the product sensitivity ratios (\ref{eq:lam-def-0-Prod}) and sensitivity model (\ref{eq:model-2period-Prod}) can also be defined
with multiple treatment strategies.
The results in Section~\ref{sec:2period-Jnt} can be extended accordingly.
The bound in Proposition~\ref{pro:upper-2period-Prod} remains conservative for each treatment strategy, and
hence the question of simultaneous attainment is not relevant.
The bounds in Corollary~\ref{cor:upper-2period-Prod} remain to be sharp and can be simultaneously attained for all treatment strategies.
We defer to Supplement Section \ref{sec:discussion-multi} further discussion about the extension under the
unconditional consistency that $L^{1,1}_1 = L^{1,0}_1$ and $L^{0,1}_1=L^{0,0}_1$, instead of only Assumption A1$^\dag$.

\vspace{.05in}
\textbf{Sample sensitivity analysis.}\;
Our paper is currently focused on population sensitivity analysis, i.e., sensitivity models and bounds at the population level,
assuming access to an infinite amount of sample data.
There are various computational and statistical problems to be addressed.
For example, although the sharp bounds in Proposition~\ref{pro:upper-Kperiod} and Corollary~\ref{cor:upper-Kperiod-Prod} are
represented through convex optimization, it is desirable to develop specialized algorithms for efficient implementation.
Moreover, the ICE and IPW functionals used to represent the sensitivity bounds need to be estimated with sample data,
depending on modeling assumptions on the conditional expectations and probabilities involved.
As mentioned in Section~\ref{sec:setup}, doubly or multiply robust functionals can also be derived and exploited.

\vspace{.3in}
\centerline{\bf\Large References}

\begin{description}\addtolength{\itemsep}{-.15in}

\item Bonvini, M., Kennedy, E.H., Ventura, V., and Wasserman, L. (2022)
Sensitivity analysis for marginal structural models, arXiv:2210.04681.

\item Dorn, J. and Guo, K. (2022) Sharp sensitivity analysis for inverse propensity weighting via quantile balancing, {\em Journal of the American Statistical Association}, to appear. 

\item Dorn, J., Guo, K., and Kallus, N. (2021) Doubly-valid/doubly-sharp sensitivity analysis for causal inference with unmeasured confounding, arXiv:2112.11449.

\item Franks, A., D'Amour, A., and Feller, A. (2020) Flexible sensitivity analysis for observational
studies without observable implications, {\em Journal of the American Statistical Association}, 115, 1730-1746.

\item Hernan, M.A. and Robins, J.M. (2022)  {\em Causal Inference}, book draft.

\item Kallus, N. and Zhou, A. (2020) Confounding-robust policy evaluation in infinite-horizon reinforcement learning,
{\em Conference on Neural Information Processing Systems}.

\item Robins J.M. (1986) A new approach to causal inference in mortality studies
with sustained exposure periods --- Application to control of the healthy
worker survivor effect, {\em Mathematical Modelling}, 7, 1393-1512.

\item Robins, J.M., Rotnitzky, A., and Scharfstein, D.O. (2000) Sensitivity analysis for selection bias and unmeasured confounding in missing
data and causal inference models, in {\em Statistical Models in Epidemiology, the Environment, and Clinical Trials}, Berlin: Springer, 1-94.

\item Rosenbaum, P.R. (2002) {\em Observational Studies} (2nd edition). New York: Springer.

\item Scharfstein, D.,  Nabi, R., Kennedy, E.H., Huang, M.-Y., Bonvini, M., and Smid, M. (2021) Semiparametric sensitivity analysis:
Unmeasured confounding in observational studies, \\ arXiv:2104.08300.

\item Tan, Z. (2006) A distributional approach for causal inference using propensity scores,
{\em Journal of the American Statistical Association}, 101, 1619-1637.

\item Tan, Z. (2022) Model-assisted sensitivity analysis for treatment effects under unmeasured confounding via regularized calibrated estimation, arXiv:2209.11383.

\item Yadlowsky, S., Namkoong, H., Basu, S., Duchi, J., and Tian, L. (2022) Bounds on the conditional and average treatment effect
with unobserved confounding factors, {\em Annals of Statistics}, to appear.

\item Zhao, Q., Small, D.S., and Bhattacharya, B.B. (2019) Sensitivity analysis
for inverse probability weighting estimators via the percentile bootstrap,
{\em Journal of the Royal Statistical Society}, Ser.~B, 81, 735-761.

\end{description}

\clearpage

\setcounter{page}{1}

\setcounter{section}{0}
\setcounter{equation}{0}

\setcounter{table}{0}
\setcounter{figure}{0}

\setcounter{pro}{0}
\renewcommand{\thepro}{S\arabic{pro}}

\setcounter{lem}{0}
\renewcommand{\thelem}{S\arabic{lem}}

\setcounter{thm}{0}
\renewcommand{\thethm}{S\arabic{thm}}

\setcounter{ass}{0}
\renewcommand{\theass}{S\arabic{ass}}

\renewcommand{\thesection}{\Roman{section}}
\renewcommand{\theequation}{S\arabic{equation}}

\renewcommand\thetable{S\arabic{table}}
\renewcommand\thefigure{S\arabic{figure}}

\begin{center}
{\Large Supplementary Material for}\\
{\Large ``Sensitivity models and bounds under sequential unmeasured confounding in longitudinal studies''}

\vspace{.1in} Zhiqiang Tan
\end{center}

\section{Additional discussion of sensitivity analysis with multiple treatment strategies} \label{sec:discussion-multi}

We discuss the extension of the product sensitivity models and bounds with multiple treatment strategies in two periods.
For any single treatment strategy $\overa$, the joint sensitivity ratio $\lambda^{\overa *}_{0,\mytext{Jnt}} (L_0, l_1, y)$ can be decomposed as
\begin{align*}
 \lambda^{\overa *}_{0,\mytext{Jnt}} (L_0, l_1, y) & =\lambda^{\overa *}_{0,L_1} (L_0, l_1) \, \lambda^{\overa *}_{0,Y} (L_0, l_1, y),
\end{align*}
where
\begin{align*}
& \lambda^{\overa *}_{0,L_1} (L_0, l_1) = \frac{\dif P_{L_1^{\overa}} (l_1 | A_0=1-a_0, L_0) }{\dif P_{L_1^{\overa}} (l_1 | A_0=a_0, L_0 ) } , \\
& \lambda^{\overa *}_{0,Y} (L_0, l_1, y) =   \frac{\dif P_{Y^{\overa}} (y| A_0=1-a_0, L_0, L_1^{\overa} = l_1) }{\dif P_{Y^{\overa}} (y| A_0=a_0, L_0, L_1^{\overa} = l_1) } .
\end{align*}
For a product sensitivity model, we postulate that for each $\overa \in \{0,1\}\times \{0,1\}$,
\begin{align}\label{eq:model-2period-Prod-multi}
\begin{split}
& \Lambda_1^{-1} \le \lambda^{\overa *}_1 (\overline{L}_1, y) \le \Lambda_1,  \\
& \Lambda_{0,L_1}^{-1} \le \lambda^{\overa *}_{0,L_1} (L_0, l_1) \le \Lambda_{0,L_1} ,\quad
\Lambda_{0,Y}^{-1} \le \lambda^{\overa *}_{0,Y} (L_0, l_1, y) \le \Lambda_{0,Y} ,
\end{split}
\end{align}
where $\Lambda_1 \ge 1$, $\Lambda_{0,Y} \ge 1$, and $\Lambda_{0,L_1} \ge 1$ are sensitivity parameters.
Under model (\ref{eq:model-2period-Prod-multi}), the sharp upper bound on $\mu^{1,1}$ is defined as
\begin{align}  \label{eq:upper-2period-Prod-multi}
\mu^{1,1}_{+,\mytext{Prod}} = \max_{Q} \; \mu^{1,1}_{\mytext{ICE}} (\lambda^{1,1 *}_{0,L_1,Q}\lambda^{1,1 *}_{0,Y,Q},\lambda^{1,1 *}_{1,Q})
= \max_{Q} \; \mu^{1,1}_{\mytext{IPW}} (\lambda^{1,1 *}_{0,L_1,Q}\lambda^{1,1 *}_{0,Y,Q},\lambda^{1,1 *}_{1,Q}),
\end{align}
over all possible distributions $Q$ for the full data as in the previous definition (\ref{eq:upper-2period-Prod}) except that
$\lambda^*_{0,L_1,Q}$, $\lambda^*_{0,Y,Q}$, and $\lambda^*_{1,Q}$ are subject to the range constraints (\ref{eq:model-2period-Prod-multi})
for all $\overa \in \{0,1\}\times \{0,1\}$.
Similarly to the related results in Section~\ref{sec:2period-Jnt},
the results in Section~\ref{sec:2period-Prod} can be extended accordingly, including
the simultaneous attainment of the sharp bounds in Corollary~\ref{cor:upper-2period-Prod}
for all treatment strategies, which depends on the allowance that $L^{1,1}_1 \not= L^{1,0}_1$ given $A_0=0$ or
$L^{0,1}_1 \not= L^{0,0}_1$ given $A_0=1$ under Assumption A1$^\dag$.

There is an interesting difference between the sharp bounds in the two special cases of model (26) in Corollary 1,
under the unconditional consistency assumption that $L^{1,1}_1 = L^{1,0}_1$ and $L^{0,1}_1 = L^{0,0}_1$ regardless of $A_0=0$ or $1$.
\vspace{-.05in}
\begin{itemize}\addtolength{\itemsep}{-.1in}
\item For model (\ref{eq:model-2period-Prod-multi}) with $\Lambda_{0,L_1}=1$
(i.e., $L_1^{\overa} \perp A_0 \,|\, L_0 $ for each $\overa\in\{0,1\}\times \{0,1\}$) and $\Lambda_{0,Y}=\Lambda_0$,
the bound $\mu^{1,1}_{+,\mytext{Prod,v1}} (1,\Lambda_0,\Lambda_1) $ in Corollary~\ref{cor:upper-2period-Prod}(i) remains
a sharp upper bound on $\mu^{1,1}$.
The result also holds if $(1,1)$ is replaced by any $\overline{t} \in\{0,1\}\times \{0,1\}$.
The sharp bounds $\{\mu^{\overline{t}}_{+,\mytext{Jnt}}(1,\Lambda_0,\Lambda_1) : \overline{t} \in\{0,1\}\times \{0,1\} \}$ can be simultaneously achieved by
a single choice of $Q$ allowed in (\ref{eq:upper-2period-Prod-multi}), where
$\mu^{\overline{t}}_{+,\mytext{Jnt}}(1,\Lambda_0,\Lambda_1) $ may be replaced by the corresponding lower sharp bound
for some treatment strategies $\overline{t}$.

\item For model (\ref{eq:model-2period-Prod-multi}) with $\Lambda_{0,Y} =1$ (i.e., $Y^{\overa} \perp A_0 \,|\, L_0, L_1^{\overa}$ for each $\overa\in\{0,1\}\times \{0,1\}$)
and $\Lambda_{0,L_1}=\Lambda_0$,
the upper bound $\mu^{1,1}_{+,\mytext{Prod,v2}} (\Lambda_0,1,\Lambda_1) $ in Corollary~\ref{cor:upper-2period-Prod}(ii)
remains a sharp upper bound on $\mu^{1,1}$. The result also holds if $(1,1)$ is replaced by any $\overline{t} \in\{0,1\}\times \{0,1\}$.
However, the sharp bounds in general are not simultaneously achieved for all treatment strategies by a single choice of the joint distribution $Q$ on the full data.
\end{itemize} \vspace{-.05in}
The reason for the simple extension in the first case above is that the additional constraints
introduced by the unconditional consistency assumption are automatically satisfied.
For example, the constraint (\ref{eq:upper-2periond-Jnt-multi-prf6}) is satisfied,
$ \lambda^{1,1 *}_{0,L_1} (L_0, l_1) = \lambda^{1,0 *}_{0,L_1} (L_0, l_1) \equiv 1$,
because it is assumed that $L_1^{1,1} \perp A_0 \,|\, L_0 $ and  $L_1^{1,0} \perp A_0 \,|\, L_0 $ in the first case.

\section{Technical details: Preparation} \label{sec:tech-preparation}


In the case of $K=1$ (single period), let
\begin{align*}
 \lambda^*_0 (L_0, y) & = \frac{\dif P_{Y^1} ( y | A_0=0, L_0) }{\dif P_{ Y^1} ( y | A_0=1, L_0 ) } .
\end{align*}
Then $E(Y^1) = \mu^1_{\mytext{ICE}} (\lambda^*_0)
= \mu^1_{\mytext{IPW}} (\lambda^*_0 )$, where
\begin{align*}
& \mu^1_{\mytext{ICE}} (\lambda_0 ) =   E\left[ E \left\{ (\pi^*_0 + (1-\pi^*_0) \lambda_0) Y | A_0=1, L_0\right\} \right],  \\
& \mu^1_{\mytext{IPW}} (\lambda_0 ) = E \left\{ A_0 \left( 1+ \frac{1-\pi^*_0 } { \pi^*_0 }\lambda_0 \right) Y \right\} .
\end{align*}
The sharp upper bound on $\mu^1$ is obtained as
\begin{align}
\mu^1_+ = \max_{\lambda_0 } \, \mu^1_{\mytext{ICE}} (\lambda_0)
= \max_{\lambda_0} \, \mu^1_{\mytext{IPW}} (\lambda_0 ) \label{eq:upper-lam-1period}
\end{align}
over all possible nonnegative functions $\lambda_0 = \lambda_0(L_0,y)$ subject to
\begin{align} \label{eq:contr-1period}
\begin{split}
 & E ( \lambda_0(L_0,Y) | A_0=1, L_0 ) \equiv 1,\\
 & \Lambda_0^{-1} \le \lambda_0(L_0,y) \le \Lambda_0 \;\text{for any}\; y .
\end{split}
\end{align}
The following result can be derived from Proposition~1 in Tan (2022).
This result is only about optimization (independently of counterfactual variables): the maximum in (\ref{eq:upper-lam-1period})
equals the minimum in (\ref{eq:upper-1period}), depending only on the observed data $(A_0,L_0,Y)$.

\begin{lem} \label{lem:upper-1period}
The sharp upper bound $\mu^1_+$ can be determined as
\begin{align} \label{eq:upper-1period}
\begin{split}
\mu^1_+ & = \min_{q_0}\, E \left[ E \left\{ Y + (1-\pi^*_0) (\Lambda_0-\Lambda_0^{-1}) \rho_{\tau_0} ( Y, q_0) | A_0=1, L_0 \right\} \right] \\
& = \min_{q_0}\, E \left\{ \frac{A_0 Y}{\pi^*_0} + A_0 \frac{1-\pi^*_0}{\pi^*_0} (\Lambda_0-\Lambda_0^{-1}) \rho_{\tau_0} ( Y, q_0) \right\}
\end{split}
\end{align}
over all possible functions $q_0 = q_0(L_0)$, where $\tau_0 = \Lambda_0/(1+\Lambda_0)$.
Let $\check q_0 (L_0)$ be a solution to the optimization in (\ref{eq:upper-1period}). Then a solution, $\check \lambda_0(L_0,y)$,
to the optimization in (\ref{eq:upper-lam-1period})
can be obtained such that $ \check \lambda_0 (L_0,y) = \Lambda_0$ if $ y > \check q_0(L_0)$ or $=\Lambda_0^{-1}$ if $y < \check q_0(L_0)$.
\end{lem}

We make two remarks about extensions of Lemma~\ref{lem:upper-1period}.
First, a similar characterization can also be obtained about sharp bounds of $\mu_b^1 = E( b(L_0, Y^1) )$ for any function $b(l_0,y)$.
The sharp upper bound on $\mu_b^1$ is obtained as
\begin{align*}
\mu^1_{b+} = \max_{\lambda_0 } \, \mu^1_{b,\mytext{ICE}} (\lambda_0)
= \max_{\lambda_0} \, \mu^1_{b,\mytext{IPW}} (\lambda_0 )  
\end{align*}
over all possible nonnegative functions $\lambda_0 = \lambda_0(L_0,y)$ subject to (\ref{eq:contr-1period}),
where $\mu^1_{b,\mytext{ICE}} (\lambda_0)$
and $\mu^1_{b,\mytext{IPW}} (\lambda_0 ) $ are defined as $\mu^1_{\mytext{ICE}} (\lambda_0)$
and $\mu^1_{\mytext{IPW}} (\lambda_0 ) $ except with $Y$ replaced by $b(L_0,Y)$. Then (\ref{eq:upper-1period}) holds
with the ``outcome'' $Y$ replaced by $b(L_0,Y)$.
By removing the $\pi^*_0 b( L_0,Y)$ term, the ICE version of this extension of (\ref{eq:upper-1period}) can be stated as
\begin{subequations}  \label{eq:upper-1period-b}
\begin{align}
& \quad \max_{\lambda_0 } \, E\left[ E \left\{ (1-\pi^*_0) \lambda_0 b(L_0,Y) | A_0=1, L_0\right\} \right] \label{eq:upper-1period-b1} \\
& = \min_{q_0}\, E \left[ E \left\{ (1-\pi^*_0) \{ b(L_0,Y) + (\Lambda_0-\Lambda_0^{-1}) \rho_{\tau_0} ( b(L_0,Y), q_0) \} | A_0=1, L_0 \right\} \right] \label{eq:upper-1period-b2}
\end{align}
\end{subequations}
over all possible nonnegative functions $\lambda_0 = \lambda_0(L_0,y)$ subject to (\ref{eq:contr-1period}) and
over all possible functions $q_0 = q_0(L_0)$.
Let $\check q_0 (L_0)$ be a solution to the optimization in (\ref{eq:upper-1period-b2}).
Then a solution, $\check \lambda_0(L_0,y)$, to the optimization in (\ref{eq:upper-1period-b1})
can be obtained such that $ \check \lambda_0 (L_0,y) = \Lambda_0$ if $ b(L_0,y) > \check q_0(L_0)$ or $=\Lambda_0^{-1}$ if $ b(L_0,y) < \check q_0(L_0)$.

Second, 
by absorbing $1-\pi^*_0 (L_0)$ into $b(L_0,Y)$ and $q_0(L_0)$ and dropping the conditioning on $A_0=1$,
the identity (\ref{eq:upper-1period-b})
can be transformed into
\begin{align}  \label{eq:upper-1period-c}
\begin{split}
& \quad \max_{\lambda_0 } \, E\left[ E \left\{  \lambda_0  b(L_0,Y) | L_0\right\} \right]\\
& =\min_{q_0}\, E \left[ E \left\{ b(L_0,Y)  + (\Lambda_0-\Lambda_0^{-1}) \rho_{\tau_0} ( b(L_0,Y), q_0) |  L_0 \right\} \right]\\
\end{split}
\end{align}
over all possible nonnegative functions $\lambda_0$ subject to the following constraints modified from (\ref{eq:contr-1period}) and over all possible functions $q_0$,
\begin{align*}
 & E ( \lambda_0(L_0,Y) | L_0 ) \equiv 1,\\
 & \Lambda_0^{-1} \le \lambda_0(L_0,y) \le \Lambda_0 \;\text{for any}\; y .
\end{align*}
This form is convenient when applied in our proof of Proposition~\ref{pro:upper-2period-Prod}.

\vspace{.1in}
For $K=2$ (2 periods), the following result shows the equivalence of ICE and IPW functionals, regardless of any counterfactual outcomes or covariates,
as mentioned in Section \ref{sec:setup}.
Under Assumptions A1$^\dag$, A2$^\dag$, and A3, both the ICE and IPW functionals reduce to $E\{ b_1( L_0, L^{1,1}_1, Y^{1,1}) \}$.
Under Assumptions A1--A3, both the ICE and IPW functionals reduce to $\mu^{1,1} = E(Y^{1,1})$ for $b_1 (\overl_1,y)=y$.

\begin{lem} \label{lem:IPW-ICE}
Under Assumption A3 for $\overline{a}=(1,1)$, it holds that
\begin{align*}
& \quad E \left\{ \frac{ A_0 A_1}{  \pi^*_0(L_0) \pi^*_1 (\overline{L}_1) } b_1(\overline{L}_1, Y) \right\} \\
& = E \left( E\left[ E \left\{ b_1(\overline{L}_1, Y) |A_0=1, A_1=1, \overline{L}_1 \right\} | A_0=1, L_0\right] \right),
\end{align*}
for any function $b_1(\overl_1, y)$ such that all (conditional) expectations are well defined above.
\end{lem}

\begin{prf}
The result follows from repeated applications of the law of iterated expectations:
\begin{align*}
& \quad E \left\{ \frac{ A_0 A_1}{  \pi^*_0(L_0) \pi^*_1 (\overline{L}_1) } b_1(\overline{L}_1, Y) \right\}
 = E \left\{ \frac{ A_0 E( A_1 b_1(\overline{L}_1, Y) |A_0, \overline{L}_1)  }{  \pi^*_0(L_0) \pi^*_1 (\overline{L}_1) } \right\} \\
& = E \left\{ \frac{ A_0 E(b_1(\overline{L}_1, Y) |A_0=1,A_1=1, \overline{L}_1)  }{  \pi^*_0(L_0) } \right\} \\
& = E \left\{ \frac{ E ( A_0   E(b_1(\overline{L}_1, Y) |A_0=1,A_1=1, \overline{L}_1)  | L_0)  }{  \pi^*_0(L_0) } \right\} \\
& = E \left( E\left[ E \left\{ b_1(\overline{L}_1, Y) |A_0=1, A_1=1, \overline{L}_1 \right\} | A_0=1, L_0\right] \right).
\end{align*}
The first line is from the law of iterated expectations.
The second line follows because
\begin{align*}
& \quad A_0 E( A_1 b_1(\overline{L}_1, Y) |A_0, \overline{L}_1)
 = A_0  E( A_1 b_1(\overline{L}_1, Y) |A_0=1, \overline{L}_1) \\
& = A_0 P(A_1=1| A_0=1, \overline{L}_1) E\{ b_1(\overline{L}_1, Y) |A_0=1, A_1=1, \overline{L}_1 \} .
\end{align*}
The third line is also from the law of iterated expectations.
The fourth line follows because
\begin{align*}
& \quad  E [ A_0  E\{ b_1(\overline{L}_1, Y) |A_0=1,A_1=1, \overline{L}_1 \}  | L_0 ] \\
& =  P(A_0=1|L_0)  E [ E \{ b_1(\overline{L}_1, Y) |A_0=1,A_1=1, \overline{L}_1 \}  | A_0=1 , L_0 ]  .
\end{align*}
\end{prf}

\section{Technical details: Proofs for Section \ref{sec:2period-primary}}

\subsection{Proof of Lemma~\ref{lem:density-ratio}}

\textbf{Necessity.}\; Equation (\ref{eq:constr-2period-a1}) easily leads to (\ref{eq:constr-2period-0}) by  consistency Assumption A1
($Y= Y^{1,1}$ if $A_0=A_1=1$). To convert (\ref{eq:constr-2period-a0}) to (\ref{eq:constr-2period-1}),
it suffices to show that
\begin{align*}
& \quad E ( \lambda^*_0(L_0,Y^{1,1}) | A_0=1, L_0 ) \nonumber \\
& = E \left\{ E ( \lambda^*_0(L_0,Y) \varrho_1 (\overline{L}_1, Y; \lambda^*_1 ) |A_0=1,A_1=1, \overline{L}_1)  | A_0=1, L_0 \right\} .
\end{align*}
This holds by direct calculation using the law of iterated expectations:
\begin{align*}
& \quad E ( \lambda^*_0(L_0,Y^{1,1}) | A_0=1, L_0 )
= E \left\{  E ( \lambda^*_0(L_0,Y^{1,1}) | A_0=1, \overline{L}_1 ) | A_0=1, L_0 \right\} ,
\end{align*}
and
\begin{align*}
& \quad E ( \lambda^*_0(L_0,Y^{1,1}) | A_0=1, \overline{L}_1 )\\
& = \pi^*_1(\overline{L}_1) E \left\{\lambda^*_0(L_0,Y^{1,1}) |A_0=1, A_1=1, \overline{L}_1 \right\} \\
& \quad  + (1-\pi^*_1(\overline{L}_1)) E \left\{\lambda^*_0(L_0,Y^{1,1}) |A_0=1, A_1=0, \overline{L}_1 \right\} \\
& = \pi^*_1(\overline{L}_1) E \left\{\lambda^*_0(L_0,Y^{1,1}) |A_0=1, A_1=1, \overline{L}_1 \right\} \\
& \quad  + (1-\pi^*_1(\overline{L}_1)) E \left\{\lambda^*_0(L_0,Y^{1,1}) \lambda^*_1(\overline{L}_1,Y^{1,1})  |A_0=1, A_1=1, \overline{L}_1 \right\} \\
& = E \left\{ \lambda^*_0(L_0,Y) \varrho_1 (\overline{L}_1, Y; \lambda^*_1 ) |A_0=1, A_1=1, \overline{L}_1 \right\},
\end{align*}
where $\varrho_1 (\overline{L}_1, Y; \lambda^*_1 ) = \pi^*_1(\overline{L}_1) + (1-\pi^*_1(\overline{L}_1)) \lambda^*_1(\overline{L}_1, Y) $, and
the second step above follows from the definition of $\lambda^*_1$ in (\ref{eq:lam-def}).

\textbf{Sufficiency.}\; We construct a desired probability distribution $Q$ for the full data $(\overline{A}_1,\overline{L}_1,$ $ Y^{1,1}, Y^{1,0}, Y^{0,1}, Y^{0,0})$.
First, we let $\dif Q_{\overline{A}_1,\overline{L}_1} = \dif P_{\overline{A}_1,\overline{L}_1}$,
i.e., the marginal distribution of $(\overline{A}_1,\overline{L}_1)$ under $Q$ is the same as that under $P$.
Second, we define the conditional distribution of $Y^{1,1}$ given $A_0=1$, $A_1=1$ or $0$, and $\overline{L}_1$ as follows:
\begin{align}\label{eq:density-ratio-Q-YL1}
\begin{split}
& \dif Q_{Y^{1,1}} (y | A_0=1,A_1=1, \overline{L}_1) = \dif P_{Y } (y | A_0=1,A_1=1, \overline{L}_1) ,\\
& \dif Q_{Y^{1,1}} (y | A_0=1,A_1=0, \overline{L}_1) = \lambda_1(\overline{L}_1, y) \dif P_{Y } (y | A_0=1,A_1=1, \overline{L}_1),
\end{split}
\end{align}
where the second line yields a proper distribution (i.e., integrated to 1) due to (\ref{eq:constr-2period-b1}).
By the law of iterated expectations, the above definition determines the following conditional distributions:
\begin{align*}
& \quad \dif Q_{Y^{1,1}} (y | A_0=1, \overline{L}_1) \\
& = \pi^*_1(\overline{L}_1) \dif Q_{Y^{1,1}} (y | A_0=1,A_1=1, \overline{L}_1) +
(1-\pi^*_1(\overline{L}_1)) \dif Q_{Y^{1,1}} (y | A_0=1,A_1=0, \overline{L}_1) \\
& =  \varrho_1 (\overline{L}_1, y; \lambda_1 ) \dif P_{Y} (y | A_0=1,A_1=1, \overline{L}_1),
\end{align*}
and
\begin{align}
& \quad \dif Q_{Y^{1,1}} (y | A_0=1, L_0 ) \nonumber \\
& =\int \left\{ \varrho_1 (\overline{L}_1, y; \lambda_1 ) \dif P_{Y} (y | A_0=1,A_1=1, L_0, L_1=l_1) \right\}\, \dif P_{L_1} (l_1 | A_0=1, L_0),
\label{eq:density-ratio-Q-YL0}
\end{align}
where the integration is taken over all possible $l_1$.
Third,  we define the conditional distribution of $Y^{1,1}$ given $A_0=0$ and $L_0$ as follows:
\begin{align*}
& \dif Q_{Y^{1,1}} (y | A_0=0, L_0) = \lambda_0 (L_0, y) \dif Q_{Y^{1,1}} (y | A_0=1, L_0),
\end{align*}
which yields a proper distribution (i.e., integrated to 1) due to (\ref{eq:constr-2period-b0}).
Fourth, we define
\begin{align*}
& \quad \dif Q_{Y^{1,1}} (y | A_0=0,A_1= 1, \overline{L}_1) = \dif Q_{Y^{1,1}} (y | A_0=0,A_1=0, \overline{L}_1)  \\
& 
= \dif Q_{Y^{1,1}} (y | A_0=0, L_0),
\end{align*}
i.e., $Y^{1,1}$ is conditionally independent of $(A_1,L_1)$ given $A_0=0$ and $L_0$ under $Q$ [notationally, $Y^{1,1} \perp (A_1,L_1) | A_0 =0, L_0$ under $Q$].
These choices for $\dif Q_{Y^{1,1}} (y | A_0=0,A_1= 1, \overline{L}_1)$ and
$\dif Q_{Y^{1,1}} (y | A_0=0,A_1=0, \overline{L}_1)$ ensure that they are compatible with the marginalization $\dif Q_{Y^{1,1}} (y | A_0=0, L_0)$ defined earlier.
Combining the preceding four steps completes the definition of $Q$ for $(\overline{A}_1, \overline{L}_1, Y^{1,1})$.
The marginal distribution of $(\overA_1, \overL_1)$ under $Q$ is fixed the same as under $P$.
Finally, we define
\begin{align*}
& \quad \dif Q_{Y^{1,0}, Y^{0,1}, Y^{0,0}} (\cdot | \overline{A}_1, \overline{L}_1, Y^{1,1}) =
\dif Q_{Y^{1,0}, Y^{0,1}, Y^{0,0}} (\cdot | \overline{A}_1, \overline{L}_1) \\
& = \dif P_{Y^{1,0}}(\cdot | \overline{A}_1, \overline{L}_1)
\dif P_{Y^{0,1}}(\cdot | \overline{A}_1, \overline{L}_1) \dif P_{Y^{0,0}}(\cdot | \overline{A}_1, \overline{L}_1)  ,
\end{align*}
i.e., $(Y^{1,1},Y^{1,0}, Y^{0,1}, Y^{0,0})$ are conditionally independent of each other given $(\overline{A}_1, \overline{L}_1)$ under $Q$,
and the conditional distribution of $Y^{1,0}$, $Y^{0,1}$, or $Y^{0,0}$ respectively given $(\overline{A}_1, \overline{L}_1)$ under $Q$ is the same as that under $P$.
This completes our construction of $Q$ for the full data $(\overline{A}_1, \overline{L}_1, Y^{1,1}, Y^{1,0}, Y^{0,1}, Y^{0,0})$.
By design, $Q$ can be easily verified to satisfy the desired properties (i) and (ii) in Lemma~\ref{lem:density-ratio}.

\vspace{.1in}
\textbf{Remark.}\; From our construction of $Q$,
the primary sensitivity ratios for treatment strategy $(1,1)$ under $Q$ are $\lambda_1$ and $\lambda_0$ as pre-specified,
and, in a separable manner, the primary sensitivity ratios for each treatment strategy $\overa\not=(1,1)$ under $Q$
are the same as under $P$ (or any probability distribution for the full data, compatible with the observed data).
The separability of the primary density ratios between treatment strategies is made explicit in Lemma~\ref{lem:density-ratio-multi} later
as an extension of Lemma~\ref{lem:density-ratio}.

\subsection{Proof of Lemma~\ref{lem:iden-mu-2period}}

Denote $\varpi_0 ( L_0, y ) = P (A_0=1 | L_0,  Y^{1,1}=y )$
and $\varpi_1 (\overline{L}_1, y  ) = P (A_1=1 | A_0=1, \overline{L}_1, Y^{1,1}=y )$.
First, by definition of $\lambda^*_1$ and $\lambda^*_0$, direct calculation yields
\begin{align} \label{eq:iden-lam-prf1}
\begin{split}
& \frac{1}{ \varpi_0 ( L_0, Y^{1,1} ) } = 1 + \frac{1-\pi^*_0 (L_0)} { \pi^*_0(L_0) } \lambda^*_0 ( L_0, Y^{1,1} ) ,\\
& \frac{1}{ \varpi_1 (\overline{L}_1, Y^{1,1}  ) }= 1 + \frac{1-\pi^*_1(\overline{L}_1 ) } { \pi^*_1(\overline{L}_1 ) } \lambda^*_1 (\overline{L}_1, Y^{1,1}) .
\end{split}
\end{align}
By Assumption A3, $\pi^*_0(L_0) >0$ and $\pi^*_1(\overline{L}_1 )>0$ almost surely. In addition, it is assumed that
$\lambda^*_0 ( L_0, Y^{1,1} ) < \infty$ and $\lambda^*_1 (\overline{L}_1, Y^{1,1}) <\infty$ almost surely.
Applying these properties in (\ref{eq:iden-lam-prf1}) shows that
$\varpi_0 ( L_0, Y^{1,1} )>0$ and $\varpi_1 (\overline{L}_1, Y^{1,1}  )>0$ almost surely.

Second, we show that
\begin{align}
\mu^{1,1} = E \left\{ \frac{A_0 A_1 Y }{\varpi_0 ( L_0, Y^{1,1} ) \varpi_1 (\overline{L}_1, Y^{1,1} ) } \right\}.\label{eq:iden-lam-prf2}
\end{align}
The right-hand side of (\ref{eq:iden-lam-prf2}) can be calculated as
\begin{align*}
& \quad E \left\{ \frac{A_0 A_1 Y }{\varpi_0 ( L_0,  Y^{1,1} ) \varpi_1 (\overline{L}_1, Y^{1,1} ) } \right\}
 = E \left\{ \frac{E( A_0 A_1 |\overline{L}_1, Y^{1,1}) }{\varpi_0 ( L_0,  Y^{1,1} ) \varpi_1 (\overline{L}_1, Y^{1,1} ) }  Y^{1,1} \right\} \\
& = E \left\{ \frac{ P(A_0=1| \overline{L}_1, Y^{1,1} ) }{ \varpi_0 (L_0, Y^{1,1} ) }  Y^{1,1} \right\}  \\
& = E \left\{ \frac{ E( P(A_0=1| \overline{L}_1, Y^{1,1} ) | L_0,Y^{1,1}) }{ \varpi_0 (L_0, Y^{1,1} ) }  Y^{1,1} \right\} \\
& = E(Y^{1,1}).
\end{align*}
The first line follows from the fact that $Y = Y^{1,1}$ if $A_0=A_1=1$ by Assumption A1 and then the law of iterated expectations.
The second line follows because $E (A_0 A_1 |\overline{L}_1, Y^{1,1} ) = P(A_0=1| \overline{L}_1, Y^{1,1} ) P(A_1=1 | A_0=1, \overline{L}_1, Y^{1,1} )$
and $P (A_1=1 | A_0=1, L_0, \overline{L}_1, Y^{1,1} ) = \varpi_1 (\overline{L}_1, Y^{1,1} ) >0$ almost surely.
The third line follows from the law of iterated expectations.
The fourth line follows because
$E( P(A_0=1| \overline{L}_1, Y^{1,1} ) | L_0,Y^{1,1}) = P(A_0=1 | L_0, Y^{1,1})$ by the law of iterated expectations
and $P (A_0=1 | L_0, Y^{1,1} ) = \varpi_0 ( L_0, Y^{1,1} ) >0 $ almost surely.

\textbf{Cautionary note.}\; The right-hand side of (\ref{eq:iden-lam-prf2}) can also be calculated as
\begin{align*}
& E \left\{ \frac{A_0 A_1 Y }{\varpi_0 ( L_0,  Y^{1,1} ) \varpi_1 (\overline{L}_1, Y^{1,1} ) } \right\}
 = E \left\{ \frac{E( A_0 A_1 |\overline{L}_1, Y) Y }{\varpi_0 ( L_0,  Y ) \varpi_1 (\overline{L}_1, Y ) } \right\} ,
\end{align*}
using the fact that $Y = Y^{1,1}$ if $A_0=A_1=1$ by Assumption A1 and then the law of iterated expectations.
It is tempting to say that a similar calculation to the above proof of (\ref{eq:iden-lam-prf2}), with $Y^{1,1}$ replaced by $Y$, would show
that the right-hand side of (\ref{eq:iden-lam-prf2}) equals $E(Y)$, a contradiction.
This reasoning is incorrect. Although
$E (A_0 A_1 |\overline{L}_1, Y ) = P(A_0=1| \overline{L}_1, Y ) P(A_1=1 | A_0=1, \overline{L}_1, Y )$ holds,
the conditional probability $P(A_1=1 | A_0=1, \overline{L}_1, Y )$ cannot be canceled out with  $\varpi_1 (\overline{L}_1, Y ) $.
This is because  $\varpi_1 (\overline{L}_1, y  ) = P (A_1=1 | A_0=1, \overline{L}_1, Y^{1,1}=y )$,
different from $P(A_1=1 | A_0=1, \overline{L}_1, Y=y)$.

Combining (\ref{eq:iden-lam-prf1})--(\ref{eq:iden-lam-prf2}) yields the IPW identity $\mu^{1,1} = \mu^{1,1}_{\mytext{IPW}} (\lambda^*_0,\lambda^*_1)$:
\begin{align*}
& \quad  \mu^{1,1} = E \left\{ A_0 A_1 \left( 1+ \frac{1-\pi^*_0 (L_0) } { \pi^*_0 (L_0) }\lambda^*_0( L_0, Y^{1,1} ) \right)
\left( 1+ \frac{1-\pi^*_1 (\overline{L}_1 ) } { \pi^*_1(\overline{L}_1 )  }\lambda^*_1 (\overline{L}_1, Y^{1,1}) \right) Y \right\} \\
& = E \left\{ A_0 A_1 \left( 1+ \frac{1-\pi^*_0 (L_0) } { \pi^*_0 (L_0) }\lambda^*_0( L_0, Y ) \right)
\left( 1+ \frac{1-\pi^*_1 (\overline{L}_1 ) } { \pi^*_1(\overline{L}_1 )  }\lambda^*_1 (\overline{L}_1, Y) \right) Y \right\},
\end{align*}
where the second line follows by Assumption A1.

Finally, the ICE identity, $\mu^{1,1} = \mu^{1,1}_{\mytext{ICE}}(\lambda^*_0,\lambda^*_1)$, follows from the IPW identity by applying Lemma \ref{lem:IPW-ICE}
with $b (\overline L_1,Y) = (\pi^*_0 + (1-\pi^*_0) \lambda^*_0) (\pi^*_1 + (1-\pi^*_1) \lambda^*_1) Y$, regardless of sequential unmeasured confounding.
Alternatively, we give a direct proof of the ICE identity, which is useful for the proof of Proposition~\ref{pro:upper-2period}.
The first step is to show that
\begin{align}
& \quad E \left\{ \varrho_0(L_0,Y^{1,1}; \lambda^*_0) Y^{1,1} | A_0=1, \overline{L}_1 \right\} \nonumber \\
& = E \left\{ \varrho_0(L_0,Y ; \lambda^*_0) (\pi^*_1 + (1-\pi^*_1) \lambda^*_1) Y |A_0=1, A_1=1, \overline{L}_1 \right\} , \label{eq:iden-lam-prf3}
\end{align}
where $\varrho_0(L_0, y; \lambda_0) =\pi^*_0(L_0) + (1-\pi^*_0(L_0)) \lambda_0(L_0,y)$. 
The left-hand side of (\ref{eq:iden-lam-prf3}) can be calculated using the law of iterated expectations as
\begin{align*}
& \quad E \left\{\varrho_0(L_0,Y^{1,1}; \lambda^*_0) Y^{1,1} |A_0=1, \overline{L}_1 \right\} \\
& = \pi^*_1 E \left\{\varrho_0(L_0,Y^{1,1}; \lambda^*_0) Y^{1,1} |A_0=1, A_1=1, \overline{L}_1 \right\} \\
& \quad  + (1-\pi^*_1) E \left\{\varrho_0(L_0,Y^{1,1}; \lambda^*_0)  Y^{1,1} |A_0=1, A_1=0, \overline{L}_1 \right\} \\
& = \pi^*_1 E \left\{\varrho_0(L_0,Y^{1,1}; \lambda^*_0) Y^{1,1} |A_0=1, A_1=1, \overline{L}_1 \right\} \\
& \quad + (1-\pi^*_1) E \left\{\varrho_0(L_0,Y^{1,1}; \lambda^*_0) \lambda^*_1( \overline{L}_1,Y^{1,1})  Y^{1,1} |A_0=1, A_1=1, \overline{L}_1 \right\} ,
\end{align*}
which is the right-hand side of (\ref{eq:iden-lam-prf3}) by Assumption A1.
The second equality above follows from the definition of $\lambda^*_1$ in (\ref{eq:lam-def}).
Similarly, the second step proceeds as follows:
\begin{align}
& \quad E (Y^{1,1} | L_0)
= E \left[ (\pi^*_0 + (1-\pi^*_0) \lambda^*_0 (L_0,Y^{1,1})) Y^{1,1} | A_0=1, L_0\right] \nonumber \\
& = E \left[ E \left\{ (\pi^*_0 + (1-\pi^*_0) \lambda^*_0 (L_0,Y^{1,1})) Y^{1,1} | A_0=1, \overline{L}_1 \right\} | A_0=1, L_0\right]. \label{eq:iden-lam-prf4}
\end{align}
Combining (\ref{eq:iden-lam-prf3})--(\ref{eq:iden-lam-prf4}) shows that
\begin{align}
& \quad E (Y^{1,1} | L_0)
= E \left[ (\pi^*_0 + (1-\pi^*_0) \lambda^*_0 (L_0,Y^{1,1})) Y^{1,1} | A_0=1, L_0\right] \nonumber \\
& = E \left[  E \left\{(\pi^*_0 + (1-\pi^*_0) \lambda^*_0)(\pi^*_1 + (1-\pi^*_1) \lambda^*_1) Y |A_0=1, A_1=1, \overline{L}_1 \right\} | A_0=1, L_0\right].
\label{eq:iden-lam-prf5}
\end{align}
Using $\mu^{1,1} = E \{ E (Y^{1,1} | L_0)\}$ then leads to the ICE identity.

\subsection{Proof of Lemma~\ref{lem:Q-lam-2period}}

First, we show that $\mu^{1,1}_+ \le \max_{\lambda_0,\lambda_1}  \mu^{1,1}_{\mytext{ICE}} (\lambda_0,\lambda_1)$
over all possible nonnegative functions $ \lambda_1$  and $\lambda_0 $
subject to  (\ref{eq:constr-2period-b}) and (\ref{eq:model-2period-b}). In fact,
for any allowed $Q$ in (\ref{eq:upper-2period}),
$\lambda^*_{1,Q}$ and $\lambda^*_{0,Q}$ satisfy the normalization constraints (\ref{eq:constr-2period-b})
by the sufficiency assertion in Lemma~\ref{lem:density-ratio}, and satisfy the range constraints (\ref{eq:model-2period-b})
by the assumption on $Q$.
Then $(\lambda^*_{0,Q},\lambda^*_{1,Q})$ are allowed as $(\lambda_0,\lambda_1)$ in the optimization in (\ref{eq:upper-2period-b}),
and hence $ \mu^{1,1}_{\mytext{ICE}} (\lambda^*_{0,Q},\lambda^*_{1,Q}) \le \max_{\lambda_0,\lambda_1}  \mu^{1,1}_{\mytext{ICE}} (\lambda_0,\lambda_1)$.
Taking the maximum over allowed $Q$ gives the desired result.

Next, we show that $\mu^{1,1}_+ \ge \max_{\lambda_0,\lambda_1}  \mu^{1,1}_{\mytext{ICE}} (\lambda_0,\lambda_1)$
over all possible nonnegative functions $ \lambda_1$  and $\lambda_0 $
subject to  (\ref{eq:constr-2period-b}) and (\ref{eq:model-2period-b}).
For any such  $(\lambda_0,\lambda_1)$, there exists a probability distribution $Q$ for the full data satisfying properties (i) and (ii)
in Lemma~\ref{lem:density-ratio}.
By these properties and the range assumption about $(\lambda_0,\lambda_1)$,
the distribution $Q$ is allowed in the optimization in (\ref{eq:upper-2period}),
and hence $ \mu^{1,1}_{\mytext{ICE}} (\lambda_0,\lambda_1) =
\mu^{1,1}_{\mytext{ICE}} (\lambda^*_{0,Q},\lambda^*_{1,Q}) \le \max_{\tilde Q}  \mu^{1,1}_{\mytext{ICE}} (\lambda^*_{0,\tilde Q},\lambda^*_{1,\tilde Q})$.
Taking the maximum over allowed $(\lambda_0,\lambda_1)$ gives the desired result.

\subsection{Proof of Proposition~\ref{pro:upper-2period}}

By Lemma~\ref{lem:density-ratio}, for $\lambda_0(L_0,y)$ and $\lambda_1 (\overline{L}_1,y)$ satisfying (\ref{eq:constr-2period-b}),
there exists a probability distribution $Q$ for the full data satisfying properties (i) and (ii).
In other words, $Q$ can be a valid choice for $P$.
Applying (\ref{eq:iden-lam-prf5}) but with $P$ replaced by $Q$ yields
\begin{align*}
& \quad E \left[  E \left\{(\pi^*_0 + (1-\pi^*_0) \lambda_0)(\pi^*_1 + (1-\pi^*_1) \lambda_1) Y |A_0=1, A_1=1, \overline{L}_1 \right\} | A_0=1, L_0\right] \\
& = E_{\lambda_1} \left[ (\pi^*_0 + (1-\pi^*_0) \lambda_0 (L_0,Y^{1,1})) Y^{1,1} | A_0=1, L_0\right] ,
\end{align*}
where $E_{\lambda_1} (\cdot | A_0=1, L_0)$ denotes the conditional expectation under $Q$, depending on $\lambda_1$.
This can also be obtained from the definition of $\dif Q_{Y^{1,1}} (\cdot | A_0=1, L_0)$ in (\ref{eq:density-ratio-Q-YL0}), which indicates that
for any function $b_0(L_0,y)$,
\begin{align*}
& \quad E_{\lambda_1} \left[ b_0(L_0,Y^{1,1}) | A_0=1, L_0 \right] \\
& = E\left[ E \left\{
(\pi^*_1 + (1-\pi^*_1) \lambda_1 ) b_0(L_0,Y ) |A_0=1,A_1=1,\overline{L}_1 \right\} | A_0=1, L_0 \right].
\end{align*}
Hence $\mu^{1,1}_{\mytext{ICE}} (\lambda_0,\lambda_1)$ can be expressed as
\begin{align*}
\mu^{1,1}_{\mytext{ICE}} (\lambda_0,\lambda_1)
=E \left[ E_{\lambda_1} \left\{ (\pi^*_0 + (1-\pi^*_0) \lambda_0 (L_0,Y^{1,1})) Y^{1,1} | A_0=1, L_0\right\} \right].
\end{align*}

For any fixed $\lambda_1$, the normalization constraint (\ref{eq:constr-2period-b0}) on $\lambda_0$ is equivalent to
\begin{align*}
E_{\lambda_1} \left( \lambda_0 (L_0,Y^{1,1})  | A_0=1, L_0 \right) \equiv 1.
\end{align*}
Applying Lemma~\ref{lem:upper-1period} with $Y$ replaced by  $Y^{1,1}$ yields
\begin{align*}
& \max_{\lambda_0} \;  \mu^{1,1}_{\mytext{ICE}} (\lambda_0,\lambda_1)
= \min_{q_0}\;  E \left[ E_{\lambda_1} \left\{ \eta_{0+}(Y^{1,1}, q_0) | A_0=1, L_0 \right\} \right] ,
\end{align*}
where $\eta_{0+}(Y^{1,1}, q_0) = Y^{1,1} + (1-\pi^*_0) (\Lambda_0-\Lambda_0^{-1}) \rho_{\tau_0} ( Y^{1,1}, q_0) $.
For any fixed $q_0$, we apply  single-period  identity (\ref{eq:upper-1period-b}), with $L_0$, $A_0$, $\lambda_0$, and $b(L_0,Y)$ set to
$L_1$, $A_1$, $\lambda_1$, and $\eta_{0+}(Y , q_0)$ respectively,
conditionally on $(A_0=1,L_0)$, and obtain
\begin{align*}
& \quad \max_{\lambda_1} \;  E_{\lambda_1} \left\{ \eta_{0+}(Y^{1,1}, q_0)  | A_0=1, L_0 \right\}\\
& = \max_{\lambda_1} \; E\left[ E \left\{
(\pi^*_1 + (1-\pi^*_1) \lambda_1 (\overline{L}_1,Y )) \eta_{0+}(Y , q_0) |A_0=1,A_1=1,\overline{L}_1 \right\} | A_0=1, L_0 \right] \\
& = \min_{q_1} \; E\left[ E \left\{ \eta_{1+} ( \eta_{0+} (Y, q_0), q_1) | A_0=1, A_1=1, \overline{L}_1 \right\} | A_0=1, L_0 \right],
\end{align*}
where $ \eta_{1+} ( \eta_{0+} (Y, q_0), q_1)=
\eta_{0+}(Y, q_0) + (1-\pi^*_1) (\Lambda_1-\Lambda_1^{-1}) \rho_{\tau_1} ( \eta_{0+}(Y, q_0), q_1)$.
Combining the preceding two displays yields
\begin{align*}
& \quad \max_{\lambda_1} \max_{\lambda_0} \; \mu^{1,1}_{\mytext{ICE}} (\lambda_0,\lambda_1) \\
& = \max_{\lambda_1}  \min_{q_0} \; E \left( E\left[ E \left\{
(\pi^*_1 + (1-\pi^*_1) \lambda_1 ) \eta_{0+}(Y , q_0) |A_0=1,A_1=1,\overline{L}_1 \right\} | A_0=1, L_0 \right] \right) \\
& =  \min_{q_0} \max_{\lambda_1}  \; E \left( E\left[ E \left\{
(\pi^*_1 + (1-\pi^*_1) \lambda_1  ) \eta_{0+}(Y , q_0) |A_0=1,A_1=1,\overline{L}_1 \right\} | A_0=1, L_0 \right] \right) \\
& =  \min_{q_0} \min_{q_1} \; E \left( E\left[ E \left\{
\eta_{1+} ( \eta_{0+} (Y, q_0), q_1) | A_0=1, A_1=1, \overline{L}_1 \right\} | A_0=1, L_0 \right] \right) .
\end{align*}
The third line above follows from Sion's minimax theorem, because
the functional involved,
$$
E \left( E\left[ E \left\{
(\pi^*_1 + (1-\pi^*_1) \lambda_1  ) \eta_{0+}(Y , q_0) |A_0=1,A_1=1,\overline{L}_1 \right\} | A_0=1, L_0 \right] \right),
$$
is linear in $\lambda_1$ and convex in $q_0$.

\subsection{Proof of inequality (\ref{eq:mono-sig})}  \label{sec:prf-mono-sig}

We give a proof of inequality (\ref{eq:mono-sig}): if $ \mathcal E_{\overL_1} Z \equiv 0$, then for any $\sigma \ge \tilde \sigma \ge 0$,
\begin{align*}
& \quad \min_{ q_0,q_1 }\;
 E \{ \mathcal E_{L_0} (  \mathcal E_{\overline{L}_1} [ \eta_{1+} \{ \eta_{0+} ( \alpha + \sigma Z,  q_0 ), q_1\} ] ) \} \\
& \ge  \min_{ q_0,q_1 }\;
 E \{ \mathcal E_{L_0} (  \mathcal E_{\overline{L}_1} [ \eta_{1+} \{ \eta_{0+} ( \alpha + \tilde\sigma Z,  q_0 ), q_1\} ] ) \}  .
\end{align*}
In fact, we establish the following stronger result (Lemma~\ref{lem:mono-sig}), in a similar setting as Lemma~\ref{lem:bound-ineq} later.
The functionals $\eta_{0+}$ and $\eta_{1+}$ are defined as
\begin{align*}
& \eta_{0+} (y, q_0) = y + c_0 \varrho_{\tau_0} (y - q_0), \\
& \eta_{1+} (y, q_1) = y + c_1 \varrho_{\tau_1} (y-q_1),
\end{align*}
where $c_0, c_1 \ge 0$, $\tau_0,\tau_1 \in (0,1)$,
and $\varrho_\tau ( u ) = \tau u^+ + (1-\tau)(-u)^+$.
It is allowed that $\tau_0, \tau_1 \in (0, 1/2]$, although the definition $\tau_j = \Lambda_j/(1+\Lambda_j)$
satisfies $\tau_j \in (1/2,1)$ for $\Lambda_j >1$.
By the notation of the main paper, $\rho_\tau (y,q) = \varrho_\tau ( y-q)$.
[Note that $\varrho_\tau(\cdot)$ should be distinguished from $\varrho_1(\cdot)$ in (\ref{eq:constr-2period-b0}).]
Inequality (\ref{eq:mono-sig}) can be deduced by applying Lemma~\ref{lem:mono-sig} to
$(L_1,Z)$ under $\dif P_{L_1} (l_1 |A_0=1,L_0) \times \dif P_{Z}(z | A_0=1,A_1=1, L_0, L_1= l_1)$,
conditionally on $L_0$, and then taking expectations over $L_0$. The range assumptions $c_0(1-\tau_0) \in [0,1]$ and
$c_1 (1-\tau_1) \in [0,1]$ are satisfied due to (\ref{eq:eta-convex-coef}) in the proof of Lemma~\ref{lem:eta-convex}.
It is allowed that $c_0(1-\tau_0)= 1$ or $c_1(1-\tau_1) = 1$ in Lemma~\ref{lem:mono-sig},
although this is ruled out for any finite $\Lambda_0$ or $\Lambda_1$ by (\ref{eq:eta-convex-coef}).

\begin{lem} \label{lem:mono-sig}
Let $(X_1,Z)$ be random variables under a joint distribution with $\mathcal E_1 Z \equiv 0$.
Suppose that $c_0(1-\tau_0) \in [0,1]$ and
$c_1(X_1) (1-\tau_1) \in [0,1]$ almost surely. Then for any $\sigma \ge \tilde \sigma \ge 0$,
\begin{align}
 \min_{ q_0,q_1 }\; \mathcal E_0 ( \mathcal E_1 [ \eta_{1+} \{ \eta_{0+} (\alpha + \sigma Z,  q_0 ), q_1\} ] )
& \ge \min_{ q_0,q_1 }\;
  \mathcal E_0 ( \mathcal E_1 [ \eta_{1+} \{ \eta_{0+} (\alpha + \tilde\sigma Z,  q_0 ), q_1\} ] ) , \label{eq:mono-sig-b}
\end{align}
where the minimization is over any real number $q_0$ and a real-valued function $q_1 = q_1(X_1)$,
and $\eta_{1+}(\cdot)$ is applied with $c_1 = c_1(X_1)$.
Throughout, $\mathcal E_1(\cdot)$ denotes conditional expectation of $Z$ given $X_1$,
and $\mathcal E_0 (\cdot)$ denotes the (unconditional) expectation over $X_1$.
\end{lem}

\begin{prf}
If $\tilde\sigma=0$, then (\ref{eq:mono-sig-b}) follows directly from Jensen's inequality by the convexity of $\eta^1_+(\cdot,\overq_1)$:
\begin{align*}
\mathcal E_0 ( \mathcal E_1 [ \eta_{1+} \{ \eta_{0+} (\alpha + \sigma Z,  q_0 ), q_1\} ] )
& \ge
  \mathcal E_0 ( \mathcal E_1 [ \eta_{1+} \{ \eta_{0+} (\alpha ,  q_0 ), q_1\} ] ) .
\end{align*}
Hence it remains to show (\ref{eq:mono-sig-b}) for $\sigma \ge \tilde\sigma> 0$.
By resetting $\alpha(X_1)-q_0$ to $\alpha(X_1)$ and $q_1(X_1) - q_0$ to $q_1(X_1)$,
it suffices to show that for any $x_1$, the conditional version of (\ref{eq:mono-sig-b}) holds given $X_1=x_1$:
\begin{align*}
& \quad \min_{ q_1 }\; \mathcal E_1 \big\{ \alpha + \sigma Z + c_0 \varrho_{\tau_0} (\alpha+ \sigma Z) +
  c_1 \varrho_{\tau_1} ( \alpha + \sigma Z + c_0 \varrho_{\tau_0} (\alpha+ \sigma Z) - q_1 ) \big\} \\
& \ge \min_{ q_1 }\; \mathcal E_1 \big\{ \alpha + \tilde\sigma Z + c_0 \varrho_{\tau_0} (\alpha+ \tilde\sigma Z) +
  c_1 \varrho_{\tau_1} ( \alpha + \tilde\sigma Z  + c_0 \varrho_{\tau_0} (\alpha+ \tilde\sigma Z) - q_1 ) \big\}  ,
\end{align*}
or equivalently, with $\mathcal E_1 ( \alpha + \sigma Z ) =\mathcal E_1 ( \alpha + \tilde \sigma Z ) =\alpha $ removed from both sides,
\begin{align} \label{eq:mono-sig-prf1}
\begin{split}
& \quad \min_{ q_1 }\;  \mathcal E_1 \big\{ c_0 \varrho_{\tau_0} (\alpha+ \sigma Z) +
  c_1 \varrho_{\tau_1} ( \alpha + \sigma Z + c_0 \varrho_{\tau_0} (\alpha+ \sigma Z) - q_1 ) \big\} \\
& \ge \min_{ q_1 }\;  \mathcal E_1 \big\{ c_0 \varrho_{\tau_0} (\alpha+ \tilde\sigma Z) +
  c_1 \varrho_{\tau_1} ( \alpha + \tilde\sigma Z  + c_0 \varrho_{\tau_0} (\alpha+ \tilde\sigma Z) - q_1 ) \big\}  ,
\end{split}
\end{align}
where $\alpha=\alpha(x_1)$, $q_1= q_1 (x_1)$, $c_1 = c_1(x_1)$, and $\mathcal E_1(\cdot)$ denotes the conditional expectation given $X_1=x_1$.
In the following, the dependency of $\alpha$, $q_1$, $c_1$, and $\mathcal E_1(\cdot)$ on $x_1$ is suppressed.
Then (\ref{eq:mono-sig-prf1}) involves the distribution of only a single random variable $Z$.

First, we show that
\begin{align*}
& \quad \underbrace{ \min_{ q_1 }\;  \mathcal E_1 \big\{ c_0 \varrho_{\tau_0} (\alpha+ \sigma Z) +
  c_1 \varrho_{\tau_1} ( \alpha + \sigma Z + c_0 \varrho_{\tau_0} (\alpha+ \sigma Z) - q_1 ) \big\} }_{\one(\sigma)} \\
& =  \mathcal E_1 \big\{ c_0 \varrho_{\tau_0} (\alpha+ \sigma Z) +
  c_1 \varrho_{\tau_1} ( \alpha + \sigma Z  + c_0 \varrho_{\tau_0} (\alpha+ \sigma Z) - (\alpha + \sigma q_\myZ + c_0\varrho_{\tau_0}(\alpha+\sigma q_\myZ) ) ) \big\} \\
& = \underbrace{ \mathcal E_1 \big\{ c_0 \varrho_{\tau_0} (\alpha+ \sigma Z) \big\} }_{\one_1(\sigma)} +
 \underbrace{ \mathcal E_1 \big\{ c_1 \varrho_{\tau_1} ( \sigma Z - \sigma q_\myZ  + c_0 \varrho_{\tau_0} (\alpha+ \sigma Z) - c_0\varrho_{\tau_0}(\alpha+\sigma q_\myZ) ) \big\} }_{\one_2(\sigma)},
\end{align*}
where $q_\myZ$ is a $\tau_1$-quantile of $Z$. In fact, the minimum in the first line is achieved by
a $\tau_1$-quantile of $ \alpha + \sigma Z  + c_0 \varrho_{\tau_0} (\alpha+ \sigma Z) $.
By Lemma~\ref{lem:eta-convex}, $y + c_0\varrho_{\tau_0}(y)$ is increasing in $y$.
Applying this with $q_\myZ$ being a $\tau_1$-quantile of $Z$ implies that
$ \alpha + \sigma q_\myZ + c_0\varrho_{\tau_0}(\alpha+\sigma q_\myZ)$ is a $\tau_1$-quantile
of $ \alpha + \sigma Z  + c_0 \varrho_{\tau_0} (\alpha+ \sigma Z) $. Hence the preceding display holds.

Second, we show that  $\mathcal E_1 \varrho_{\tau_0} (\alpha+ \sigma Z) $ is nondecreasing in $\sigma >0$.
This can be seen as an monotonic extension of Jensen's inequality, which only asserts that
$\mathcal E_1 \varrho_{\tau_0} (\alpha+ \sigma Z) \ge \varrho_{\tau_0} (\alpha)$.
Note that $ \varrho_{\tau_0} (\alpha+ \sigma z) $ is convex in $\sigma>0$.
By Bertsekas (1973),  $\mathcal E_1 \varrho_{\tau_0} (\alpha+ \sigma Z) $ is convex in $\sigma>0$
and the sub-gradient of $\mathcal E_1 \varrho_{\tau_0} (\alpha+ \sigma Z)$ can be calculated as
\begin{align*}
\partial \mathcal E_1 \varrho_{\tau_0} (\alpha+ \sigma Z) = \mathcal E_1 \big\{ (\partial \varrho_{\tau_0} ) Z  \big\},
\end{align*}
where
\begin{align*}
\partial \varrho_{\tau_0}
\left\{ \begin{array}{ll}
= -(1-\tau_0),  & \alpha+ \sigma Z <0,\\
\in [ -(1-\tau_0), \tau_0 ], & \alpha+ \sigma Z =0, \\
= \tau_0 , & \alpha+ \sigma Z > 0.
\end{array} \right.
\end{align*}
It suffices to show that $ \partial \mathcal E_1 \varrho_{\tau_0} (\alpha+ \sigma Z) \ge 0$.
Consider the following two cases: $-\alpha \ge 0$ or $-\alpha <0$. In the case of $-\alpha \ge 0$,
\begin{align*}
\mathcal E_1 \big\{ (\partial \varrho_{\tau_0} ) Z \big\}  & \ge
\mathcal E_1 \big\{ ( -(1-\tau_0) 1_{B_+^c} + \tau_0 1_{B_+} )Z  \big\} \\
& = \mathcal E_1 \big\{ 1_{B_+} Z  \big\} \ge 0,
\end{align*}
where $B_+= \{ \alpha + \sigma Z >0\}$.
The first inequality follows because if $\alpha+ \sigma Z =0$, then $Z = -\alpha/\sigma \ge 0$ and hence
 $(\partial \varrho_{\tau_0} ) Z \ge -(1-\tau_0) Z $.
The second inequality follows because $0 = \mathcal E_1 Z = \mathcal E_1 ( 1_{B_+} Z  + 1_{B_+^c} Z) $.
The third inequality follows because if $\alpha+ \sigma Z >0$, then $Z > -\alpha/\sigma > 0$.
In the case of $-\alpha < 0$,
\begin{align*}
\mathcal E_1 \big\{ (\partial \varrho_{\tau_0} ) Z \big\}  & \ge
\mathcal E_1 \big\{ ( -(1-\tau_0) 1_{B_-} + \tau_0 1_{B_-^c} )Z  \big\} \\
& = \mathcal E_1 \big\{- 1_{B_- } Z  \big\} \ge 0,
\end{align*}
where $B_-= \{ \alpha + \sigma Z < 0\}$.
The first inequality follows because if $\alpha+ \sigma Z =0$, then $Z = -\alpha/\sigma < 0$ and hence
 $(\partial \varrho_{\tau_0} ) Z \ge \tau_0 Z $.
The second inequality follows because $0 = \mathcal E_1 Z = \mathcal E_1 ( 1_{B_-} Z  + 1_{B_-^c} Z) $.
The third inequality follows because if $\alpha+ \sigma Z < 0$, then $Z < -\alpha/\sigma < 0$.
Combining the two cases gives the desired result.

Third, we show that $\one(\sigma)$ is nondecreasing in $\{\sigma>0: \alpha + \sigma q_\myZ \le 0\}$, if nonempty.
Denote by $F(z)$ the cumulative distribution function of $Z$.
If $\alpha + \sigma q_\myZ \le 0$, then
by splitting the real line into $(-\alpha/\sigma, \infty)$, $(q_\myZ, -\alpha/\sigma]$, and $(-\infty, q_\myZ]$, $\one_2(\sigma)$ can be calculated as follows:
\begin{align*}
& \quad \one_2 (\sigma) \\
& = c_1 \int_{(-\alpha/\sigma,\infty)} \tau_1 \big\{ \sigma z - \sigma q_\myZ + c_0 \tau_0 (\alpha + \sigma z) + c_0 (1-\tau_0) (\alpha + \sigma q_\myZ) \big\} \,\dif F (z) \\
& \quad + c_1 \int_{(q_\myZ, -\alpha/\sigma]} \tau_1 \big\{ \sigma z - \sigma q_\myZ + c_0 (-1+\tau_0) (\alpha + \sigma z) + c_0 (1-\tau_0) (\alpha + \sigma q_\myZ) \big\} \,\dif F (z) \\
& \quad + c_1 \int_{(-\infty, q_\myZ]} (-1+\tau_1) \big\{ \sigma z - \sigma q_\myZ + c_0 (-1+\tau_0) (\alpha + \sigma z) + c_0 (1-\tau_0) (\alpha + \sigma q_\myZ) \big\} \,\dif F (z) \\
& = \underbrace{
c_1 \int_{(-\infty,\infty)} \tau_1 \big\{ \sigma z - \sigma q_\myZ + c_0 (-1+\tau_0) (\alpha + \sigma z) + c_0 (1-\tau_0) (\alpha + \sigma q_\myZ) \big\} \,\dif F (z)
}_{\two_1(\sigma)}\\
& \quad + \underbrace{
c_1 \int_{(-\infty,q_\myZ]} (-1) \big\{ \sigma z - \sigma q_\myZ + c_0 (-1+\tau_0) (\alpha + \sigma z) + c_0 (1-\tau_0) (\alpha + \sigma q_\myZ) \big\} \,\dif F (z)
 }_{\two_2 (\sigma)} \\
& \quad + \underbrace{
c_1 \int_{(-\alpha/\sigma,\infty) } \tau_1 c_0 (\alpha + \sigma z)  \,\dif F (z)
 }_{\two_3(\sigma)} .
\end{align*}
Note that $ \mathcal E_1 \varrho_{\tau_0} (\alpha+ \sigma Z) $ can be calculated as
\begin{align}
& \quad  \mathcal E_1 \varrho_{\tau_0} (\alpha+ \sigma Z) \nonumber \\
& =  \int_{(-\alpha/z,\infty)} \tau_0 (\alpha + \sigma z)\,\dif F(z) +
    \int_{(-\infty, -\alpha/z]} (-1+\tau_0) (\alpha + \sigma z)\,\dif F(z) \nonumber \\
& =   \int_{(-\alpha/z,\infty)} (\alpha + \sigma z)\,\dif F(z) +
    \int_{(-\infty, \infty)} (-1+\tau_0) (\alpha + \sigma z)\,\dif F(z) \nonumber \\
& =  \int_{(-\alpha/z,\infty)} (\alpha + \sigma z)\,\dif F(z) + (-1+\tau_0)\alpha , \label{eq:rho0-expression}
\end{align}
where the last equality follows because $\mathcal E_1 Z=0$. Hence
$\two_3 (\sigma) = c_0 c_1 \tau_1 \mathcal E_1 \varrho_{\tau_0} (\alpha+ \sigma Z) + c_0 c_1 \tau_1 (-1+\tau_0) \alpha$.
The gradients of $ \two_1(\sigma ) + \two_2(\sigma)$ can be calculated as
\begin{align*}
& \quad \frac{\dif }{\dif \sigma} ( \two_1(\sigma ) + \two_2(\sigma) ) \\
& = c_1 \int_{(-\infty,\infty)} \tau_1 \big\{ z - q_\myZ + c_0 (-1+\tau_0) z + c_0 (1-\tau_0) q_\myZ \big\} \,\dif F (z)\\
& \quad +
c_1 \int_{(-\infty,q_\myZ]} (-1) \big\{ z - q_\myZ + c_0 (-1+\tau_0) z + c_0 (1-\tau_0) q_\myZ \big\} \,\dif F (z) \\
& = c_1 \int_{(-\infty,q_\myZ]} (-1) ( 1 - c_0(1-\tau_0)) z \,\dif F (z)\\
& \ge 0.
\end{align*}
The second equality uses $\mathcal E_1 Z= 0$ and $q_\myZ$ being a $\tau_1$-quantile of $Z$, i.e., $\int_{ (-\infty, q_\myZ]} \dif F(z) = \tau_1$.
The last inequality uses $1-c_0 (1-\tau_0) \ge 0$
and $\int_{(-\infty,q_\myZ]} z \,\dif F (z) \le 0$,
which can be shown by combining the following two cases.
If $q_\myZ \le 0$, then
$\int_{(-\infty,q_\myZ]} z \,\dif F (z) \le 0$.
If $q_\myZ >0$, then $\int_{(-\infty,q_\myZ]} z \,\dif F (z) = \int_{(q_\myZ, \infty)} (-z) \,\dif F (z) \le 0$ using $\mathcal E_1 Z=0$.
Hence $ \two_1(\sigma ) + \two_2(\sigma) $ is nondecreasing in $\sigma>0 $,
which together with $\mathcal E_1 \varrho_{\tau_0} (\alpha+ \sigma Z)$ being nondecreasing implies that
$\one(\sigma)= \one_1(\sigma) + \one_2 (\sigma)
= \two_1(\sigma) + \two_2(\sigma) + c_0(1 + c_1\tau_1) \mathcal E_1 \varrho_{\tau_0} (\alpha+ \sigma Z) + c_0 c_1 \tau_1 (-1+\tau_0) \alpha$
is nondecreasing in $\sigma >0$ provided $\alpha + \sigma q_\myZ \le 0$.

Fourth, we show that $\one(\sigma)$ is nondecreasing in $\{\sigma>0: \alpha + \sigma q_\myZ > 0\}$, if nonempty.
If $\alpha + \sigma q_\myZ > 0$, then
by splitting the real line into $(q_\myZ, \infty)$, $(-\alpha/\sigma, q_\myZ]$, and $(-\infty, -\alpha/\sigma]$, $\one_2(\sigma)$ can be calculated
similarly as in the third step:
\begin{align*}
& \quad \one_2 (\sigma) \\
& = c_1 \int_{(q_\myZ,\infty)} \tau_1 \big\{ \sigma z - \sigma q_\myZ + c_0 \tau_0 (\alpha + \sigma z) - c_0 \tau_0 (\alpha + \sigma q_\myZ) \big\} \,\dif F (z) \\
& \quad + c_1 \int_{(-\alpha/\sigma, q_\myZ]} (-1+\tau_1) \big\{ \sigma z - \sigma q_\myZ + c_0 \tau_0 (\alpha + \sigma z) - c_0 \tau_0 (\alpha + \sigma q_\myZ) \big\} \,\dif F (z) \\
& \quad + c_1 \int_{(-\infty, -\alpha/\sigma]} (-1+\tau_1) \big\{ \sigma z - \sigma q_\myZ + c_0 (-1+\tau_0) (\alpha + \sigma z) - c_0 \tau_0 (\alpha + \sigma q_\myZ) \big\} \,\dif F (z) \\
& = \underbrace{
c_1 \int_{(-\infty,\infty)} \tau_1 \big\{ \sigma z - \sigma q_\myZ + c_0 \tau_0 (\alpha + \sigma z) - c_0 \tau_0 (\alpha + \sigma q_\myZ) \big\} \,\dif F (z)
}_{\three_1(\sigma)}\\
& \quad + \underbrace{
c_1 \int_{(-\infty, q_\myZ]} (-1) \big\{ \sigma z - \sigma q_\myZ + c_0 \tau_0 (\alpha + \sigma z) - c_0 \tau_0 (\alpha + \sigma q_\myZ)  \big\} \,\dif F (z)
 }_{\three_2 (\sigma)} \\
& \quad + \underbrace{
c_1 \int_{(-\infty,-\alpha/\sigma] } (1-\tau_1) c_0 (\alpha + \sigma z)  \,\dif F (z)
 }_{\three_3(\sigma)} .
\end{align*}
From  $\mathcal E_1 Z=0$ and (\ref{eq:rho0-expression}), $ \mathcal E_1 \varrho_{\tau_0} (\alpha+ \sigma Z) $ can also be expressed as
\begin{align*}
&  \mathcal E_1 \varrho_{\tau_0} (\alpha+ \sigma Z)
 =  \int_{(-\infty, -\alpha/z]} (-1) (\alpha + \sigma z)\,\dif F(z)  +  (-1+\tau_0) \alpha ,
\end{align*}
Hence
$\three_3 (\sigma) = - c_0 c_1 (1-\tau_1) \mathcal E_1 \varrho_{\tau_0} (\alpha+ \sigma Z)$.
The gradients of $ \three_1(\sigma ) + \three_2(\sigma)$ can be calculated as
\begin{align*}
& \quad \frac{\dif }{\dif \sigma} ( \three_1(\sigma ) + \three_2(\sigma) ) \\
& = c_1 \int_{(-\infty,\infty)} \tau_1 \big\{ z - q_\myZ + c_0 \tau_0 z - c_0 \tau_0 q_\myZ \big\} \,\dif F (z)\\
& \quad +
c_1 \int_{(-\infty, q_\myZ]} (-1) \big\{z - q_\myZ + c_0 \tau_0 z - c_0 \tau_0 q_\myZ \big\} \,\dif F (z) \\
& = c_1 \int_{(-\infty,q_\myZ]} (-1) ( 1 + c_0\tau_0) z \,\dif F (z)\\
& \ge 0.
\end{align*}
The second equality uses $\mathcal E_1 Z= 0$ and $q_\myZ$ being a $\tau_1$-quantile of $Z$, and
the last inequality uses $\int_{(-\infty,q_\myZ]} z \,\dif F (z) \le 0$,  as shown in the third step above.
Hence $ \three_1(\sigma ) + \three_2(\sigma) $ is nondecreasing in $\sigma>0 $,
which together with $1 - c_1(1-\tau_1)\ge0$ and $\mathcal E_1 \varrho_{\tau_0} (\alpha+ \sigma Z)$ being nondecreasing implies that
$\one(\sigma)= \one_1(\sigma) + \one_2 (\sigma) =
\three_1(\sigma) + \three_2(\sigma) + c_0(1 - c_1(1-\tau_1)) \mathcal E_1 \varrho_{\tau_0} (\alpha+ \sigma Z) + c_0 c_1 (1-\tau_1) (-1+\tau_0) \alpha$
is nondecreasing in $\sigma >0$ provided $\alpha + \sigma q_\myZ > 0$.

Combining the third and fourth steps shows that $\one(\sigma)$ is nondecreasing in $\sigma>0$.
\end{prf}

\vspace{-.1in}
\section{Technical details: Proofs for Section \ref{sec:2period-Jnt}}

\vspace{-.1in}
\subsection{Proof of Lemma~\ref{lem:density-ratio-Jnt}}  \label{sec:prf-lem-density-ratio-Jnt}

The necessity assertion can be proved similarly as that in Lemma~\ref{lem:density-ratio}.
In the following, we prove the sufficiency assertion, which involves subtle differences from that in Lemma~\ref{lem:density-ratio}.
We construct a desired probability distribution $Q$ for the full data $(\overline{A}_1, L_0,\{(L_1^{a_0a_1},Y^{a_0a_1}): a_0,a_1=0,1\} )$.
First, we let $\dif Q_{A_0,L_0} = \dif P_{A_0,L_0}$,
i.e., the marginal distribution of $(A_0,L_0)$ under $Q$ is the same as that under $P$.
Second, we let
\begin{align}
\dif Q_{L_1^{1,1}} (l_1 | A_0=1, L_0) = \dif P_{L_1} (l_1 | A_0=1, L_0), \label{eq:density-ratio-Q-Jnt-L1}
\end{align}
i.e., the conditional distribution of $L_1^{1,1}$ given $A_0=1$ and $L_0$ under $Q$
(equivalently, that of $L_1$ given $A_0=1$ and $L_0$ under $Q$) is the same as that of $L_1$ given $A_0=1$ under $P$.
We also define for $a_1=0,1$ and any $l_1$,
\begin{align}
& Q ( A_1= a_1 | A_0=1, L_0, L_1^{1,1} = l_1) = P ( A_1= a_1 | A_0=1, L_0, L_1 = l_1),\label{eq:density-ratio-Q-Jnt-A1}
\end{align}
Third, we define the conditional distribution of $Y^{1,1}$ given $A_0=1$, $A_1=1$ or $0$, and $(L_0, L_1^{1,1})$ as follows:
\begin{align} \label{eq:density-ratio-Q-Jnt-YL1}
\begin{split}
& \dif Q_{Y^{1,1}} (y | A_0=1,A_1=1, L_0, L_1^{1,1}=l_1) = \dif P_{Y } (y | A_0=1,A_1=1, L_0, L_1=l_1) ,\\
& \dif Q_{Y^{1,1}} (y | A_0=1,A_1=0, L_0, L_1^{1,1}=l_1) = \lambda_1(L_0, l_1, y) \dif P_{Y } (y | A_0=1,A_1=1, L_0, L_1=l_1),
\end{split}
\end{align}
where the second line yields a proper distribution (i.e., integrated to 1) due to (\ref{eq:constr-2period-b1}).
This agrees with the definition (\ref{eq:density-ratio-Q-YL1}) for $\dif Q_{Y^{1,1}} (y | A_0=1,A_1, \overline{L}_1)$ in the proof of Lemma~\ref{lem:density-ratio},
under Assumption A1$^\dag$, $L_1=L_1^{1,1}$ if $A_0=1$ under $Q$ as well as $P$.
Combining the proceeding three steps yields
the conditional distribution of $ (A_1, L_0, L_1^{1,1},Y^{1,1} )$ given $A_0=1$
(as well as the marginal probabilities of $A_0$) under $Q$.
In particular, applying the law of iterated expectations with (\ref{eq:density-ratio-Q-Jnt-A1}) and (\ref{eq:density-ratio-Q-Jnt-YL1})
determines the conditional distribution
\begin{align*}
& \quad \dif Q_{Y^{1,1}} (y | A_0=1, L_0, L_1^{1,1}=l_1) \\
& =  \varrho_1 (L_0, l_1, y; \lambda_1 ) \dif P_{Y} (y | A_0=1,A_1=1, L_0, L_1 =l_1),
\end{align*}
which, in conjunction with $\dif Q_{L_1^{1,1}} (l_1 | A_0=1,L_0)$ in (\ref{eq:density-ratio-Q-Jnt-L1}), yields
\begin{align}
& \quad \dif Q_{L_1^{1,1},Y^{1,1}} (l_1, y | A_0=1, L_0 ) \nonumber \\
& = \dif Q_{Y^{1,1}} (y | A_0=1,L_0, L_1^{1,1}=l_1) \times \dif Q_{L_1^{1,1}} (l_1 | A_0=1,L_0) \nonumber \\
& =  \varrho_1 (L_0, l_1, y; \lambda_1 ) \dif P_{Y} (y | A_0=1,A_1=1, L_0, L_1 =l_1)  \dif P_{L_1} (l_1 | A_0=1,L_0) . \label{eq:density-ratio-Q-Jnt-YL1-b}
\end{align}
Fourth,  we define the conditional distribution of $(L_1^{1,1},Y^{1,1})$ given $A_0=0$ and $L_0$ as follows:
\begin{align}
& \dif Q_{L_1^{1,1},Y^{1,1}} (l_1, y | A_0=0, L_0) = \lambda_{0,\mytext{Jnt}} (L_0, l_1, y) \dif Q_{L_1^{1,1}, Y^{1,1}} (l_1, y | A_0=1, L_0),
 \label{eq:density-ratio-Q-Jnt-A0}
\end{align}
which yields a proper distribution (i.e., integrated to 1) due to (\ref{eq:constr-2period-Jnt-b0}).
Fifth, we define the conditional probabilities of $A_1$ given $A_0=0$ and $L_0$ under $Q$ the same as under $P$:
\begin{align*}
& Q ( A_1=a_1 | A_0=0, L_0) = P( A_1=a_1 | A_0=0, L_0),
\end{align*}
for $a_1=0,1$. We also define
\begin{align*}
& \dif Q_{L_1^{1,1},Y^{1,1}} (l_1, y | A_0=0, A_1 , L_0) = \dif Q_{L_1^{1,1},Y^{1,1}} (l_1, y | A_0=0, L_0),
\end{align*}
i.e., $(L_1^{1,1},Y^{1,1})$ is conditionally independent of $A_1$ given $A_0=0$ and $L_0$ under $Q$
[notationally, $ (L_1^{1,1},Y^{1,1}) \perp A_1 | A_0=0, L_0$ under $Q$].
Combining the proceeding five steps completes the definition of $Q$ for
$ (\overline{A}_1, L_0, L_1^{1,1},Y^{1,1} )$.
The marginal distribution of $(\overA_1, L_0)$ under $Q$ is fixed the same as under $P$.
Finally, we define
\begin{align} \label{eq:density-ratio-Q-Jnt-final1}
\begin{split}
& \quad \dif Q_{\{(L_1^{a_0a_1},Y^{a_0a_1}): a_0,a_1=0,1, (a_0,a_1)\not=(1,1)\}} (\cdot | \overline{A}_1, L_0, L_1^{1,1},Y^{1,1})\\
& = \dif Q_{\{(L_1^{a_0a_1},Y^{a_0a_1}): a_0,a_1=0,1, (a_0,a_1)\not=(1,1)\}} (\cdot | \overline{A}_1, L_0, L_1^{1,1} ) \\
& = \dif Q_{ L_1^{1,0},Y^{1,0}} (\cdot | \overline{A}_1, L_0, L_1^{1,1} )
\times \dif P_{\{(L_1^{a_0a_1},Y^{a_0a_1}): a_0=0, a_1=0,1 \}} (\cdot | \overline{A}_1, L_0) ,
\end{split}
\end{align}
i.e., $ (L_1^{1,0},Y^{1,0})$ and $\{(L_1^{a_0a_1},Y^{a_0a_1}): a_0=0, a_1=0,1 \}$ are conditionally independent of each other
and of $Y^{1,1}$ given $(\overline{A}_1, L_0, L_1^{1,1})$ under $Q$, and
$\{(L_1^{a_0a_1},Y^{a_0a_1}): a_0=0, a_1=0,1 \}$ given $(\overline{A}_1, L_0)$ is conditionally independent of $L_1^{1,1}$ under $Q$
and is identically distributed under $Q$ and $P$.
Moreover, we define
\begin{subequations} \label{eq:density-ratio-Q-Jnt-final2}
\begin{align}
& \quad \dif Q_{ L_1^{1,0},Y^{1,0}} (l_1^\prime, y |A_0=1,A_1, L_0, L_1^{1,1}=l_1 ) \nonumber \\
& = \left\{ \begin{array}{ll}
\dif P_{Y^{1,0}} (y |A_0=1,A_1, L_0, L_1^{1,0}=l_1), & l_1=l_1^\prime, \label{eq:density-ratio-Q-Jnt-final2-1} \\
0, & l_1\not=l_1^\prime,
\end{array} \right. \\
& \quad \dif Q_{ L_1^{1,0},Y^{1,0}} (\cdot |A_0=0, A_1, L_0, L_1^{1,1} ) \nonumber \\
& = \dif Q_{ L_1^{1,0},Y^{1,0}} (\cdot |A_0=0, A_1, L_0) = \dif P_{ L_1^{1,0},Y^{1,0}} (\cdot |A_0=0, L_0 ), \label{eq:density-ratio-Q-Jnt-final2-0}
\end{align}
\end{subequations}
i.e., $L_1^{1,1}=L_1^{1,0}$ almost surely given $A_0=1$ under $Q$,
the conditional distribution of $Y^{1,0}$ given $(A_0=1,A_1, L_0, L_1^{1,0})$ under $Q$ is the same as under $P$,
and $(L_1^{1,0},Y^{1,0})$ given $( A_0=0, L_0 )$ is conditionally independent of $(A_1, L_1^{1,1})$ under $Q$
and is identically distributed under $Q$ and $P$.
This completes our construction of $Q$ for the full data $(\overline{A}_1, L_0,\{(L_1^{a_0a_1},Y^{a_0a_1}): a_0,a_1=0,1\} )$.
By design, $Q$ can be easily verified to satisfy the desired properties (i) and (ii) in Lemma~\ref{lem:density-ratio-Jnt}.

\vspace{.1in}
\textbf{Remark.}\;
(i) From our construction of $Q$, the joint sensitivity ratios for
treatment strategy $(1,1)$ are $\lambda_1$ and $\lambda_{0,\mytext{Jnt}}$ as pre-specified, and, in a separable manner,
the joint sensitivity ratios for treatment strategy $(1,0)$, $(0,1)$, or $(0,0)$ are the same as under $P$
(or any probability distribution for the full data, compatible with the observed data).
The separability of the joint density ratios between treatment strategies is made explicit in Lemma~\ref{lem:density-ratio-Jnt-multi} later
as an extension of Lemma~\ref{lem:density-ratio-Jnt}.

(ii) By the definition (\ref{eq:density-ratio-Q-Jnt-final1}),
the consistency that $L_1^{0,1} = L_1^{0,0}$ given $A_0=0$ is preserved under $Q$, the same as under $P$.
By the definition (\ref{eq:density-ratio-Q-Jnt-final2}),
the consistency that $L_1^{1,1} = L_1^{1,0}$ given $A_0=1$ is also satisfied under $Q$.
Such consistency is required by Assumption A1$^\dag$.
However, Assumption A1$^\dag$ in principle allows that under $Q$ as well as $P$,
$L_1^{0,1} \not= L_1^{0,0}$ given $A_0=1$ or $L_1^{1,1} \not= L_1^{1,0}$ given $A_0=0$.
In fact, enforcing the consistency that $L_1^{1,1} = L_1^{1,0}$ given $A_0=0$ under $Q$ would imply that
%
%
the joint sensitivity ratio $\lambda^{1,0 *}_{0,\mytext{Jnt},Q}$ for treatment strategy $(1,0)$
is tied with (no longer separable from) the pre-specified sensitivity ratio $\lambda^{1,1 *}_{0,\mytext{Jnt},Q}=\lambda_{0,\mytext{Jnt}}$ for treatment strategy $(1,1)$,
as shown in (\ref{eq:upper-2periond-Jnt-multi-prf7}).
See Remark (iv) in Section \ref{sec:pro-upper-2period-Jnt-multi}.

(iii) To satisfy $L_1^{1,1} = L_1^{1,0}$ given $A_0=0$, the definition (\ref{eq:density-ratio-Q-Jnt-final2-0}) can be modified as follows:
\begin{align*}
& \quad \dif Q_{ L_1^{1,0},Y^{1,0}} (l_1^\prime, y  |A_0=0, A_1, L_0, L_1^{1,1}=l_1 )  \\
& = \dif Q_{ L_1^{1,0},Y^{1,0}} (l_1^\prime, y  |A_0=0, L_0, L_1^{1,1} = l_1)  \\
& = \left\{ \begin{array}{ll}
\dif Q_{ Y^{1,0}} ( y  |A_0=1, L_0, L_1^{1,0} = l_1) , & l_1 = l_1^\prime,\\
0, & l_1 \not= l_1^\prime,
\end{array} \right.
\end{align*}
i.e., $L_1^{1,1}=L_1^{1,0}$ almost surely given $A_0=0$ under $Q$, and
$ Y^{1,0}$ is conditionally independent of $(A_0,A_1)$ given $(L_0, L_1^{1,0})$ under $Q$
[notatinally, $Y^{1,0} \perp (A_0,A_1) | L_0,L_1^{1,0}$ under $Q$].
In this case, the joint sensitivity ratio $\lambda^{1,0 *}_{0,\mytext{Jnt},Q}$ for treatment strategy $(1,0)$
is tied with  $\lambda^{1,1 *}_{0,\mytext{Jnt},Q}=\lambda_{0,\mytext{Jnt}}$ for treatment strategy $(1,1)$ as mentioned above:
\begin{align}\label{eq:rem-lem-density-ratio-Jnt}
\begin{split}
& \quad \lambda^{1,0 *}_{0,\mytext{Jnt},Q} (L_0,l_1,y)
= \frac{\dif Q_{L_1^{1,0},Y^{1,0}} ( l_1,y | A_0=0, L_0) }{ \dif Q_{L_1^{1,0},Y^{1,0}} ( l_1,y | A_0=1, L_0) }\\
& = \frac{\dif Q_{L_1^{1,0} } ( l_1 | A_0=0, L_0) }{ \dif Q_{L_1^{1,0} } ( l_1 | A_0=1, L_0) }
 = \frac{\dif Q_{L_1^{1,1} } ( l_1 | A_0=0, L_0) }{ \dif Q_{L_1^{1,1} } ( l_1 | A_0=1, L_0) } \\
& = E_Q \left\{ \lambda^{1,1 *}_{0,\mytext{Jnt},Q} (L_0,L_1^{1,1},Y^{1,1}) | A_0=1, L_0, L_1^{1,1} = l_1 \right\},
\end{split}
\end{align}
where $E_Q (\cdot)$ denotes the expectation taken under $Q$.
The second equality follows from $Y^{1,0} \perp A_0 | L_0,L_1^{1,0}$ under $Q$,
the third equality follows from $L_1^{1,1} = L_1^{1,0}$ given $A_0=0$,
and the last equality follows from the identity (\ref{eq:importance-iden-b}).
In addition, by the definition (\ref{eq:density-ratio-Q-Jnt-final1}),
the consistency that $L_0^{0,1} = L_1^{0,0}$ given $A_0=1$ is preserved under $Q$, if the same is assumed under $P$.
Therefore, Lemma~\ref{lem:density-ratio-Jnt} can be extended such that the unconditional consistency assumption
 that $L_1^{1,1} = L_1^{1,0}$  and  $L_1^{0,1} = L_1^{0,0}$ is satisfied under $Q$,
 at the cost of allowing $\lambda^{1,0 *}_{0,\mytext{Jnt},Q}$ to depend on $\lambda^{1,1 *}_{0,\mytext{Jnt},Q}=\lambda_{0,\mytext{Jnt}}$,
 as indicated by (\ref{eq:rem-lem-density-ratio-Jnt}).

\subsection{Proof of Lemma~\ref{lem:iden-mu-2period-Jnt}}

Denote $\varpi_{0,\mytext{Jnt}} ( L_0, l_1, y) = P (A_0=1 | L_0, L_1^{1,1}=l_1, Y^{1,1}=y )$
and $\varpi_1 (\overline{L}_1, y ) = P (A_1=1 | A_0=1, \overline{L}_1, Y^{1,1}=y )$.
For the IPW identity, by similar reasoning as in the proof of Lemma~\ref{lem:iden-mu-2period}, it suffices to show that
\begin{align}
\mu^{1,1} = E \left\{ \frac{A_0 A_1 Y }{\varpi_{0,\mytext{Jnt}} ( L_0, L_1^{1,1}, Y^{1,1} ) \varpi_1 (\overline{L}_1, Y^{1,1} ) } \right\}.
\label{eq:iden-lam-Jnt-prf1}
\end{align}
However, the proof of (\ref{eq:iden-lam-Jnt-prf1}) differs in a subtle way from that of (\ref{eq:iden-lam-prf2}).
The right-hand side of (\ref{eq:iden-lam-Jnt-prf1}) can be calculated as
\begin{align*}
& \quad E \left\{ \frac{A_0 A_1 Y }{\varpi_{0,\mytext{Jnt}} ( L_0, L_1^{1,1}, Y^{1,1} ) \varpi_1 (\overline{L}_1, Y^{1,1} ) } \right\}
 = E \left\{ \frac{E(A_0 A_1| L_0, L_1^{1,1},Y^{1,1}) }{\varpi_{0,\mytext{Jnt}} ( L_0, L_1^{1,1}, Y^{1,1} ) \varpi_1 (\overline{L}_1, Y^{1,1} ) } Y^{1,1} \right\} \\
& = E \left\{ \frac{ P(A_0=1| L_0, L_1^{1,1}, Y^{1,1} ) }{ \varpi_{0,\mytext{Jnt}} ( L_0, L_1^{1,1}, Y^{1,1} ) }  Y^{1,1} \right\}  \\
& = E(Y^{1,1}).
\end{align*}
The first line follows because $Y = Y^{1,1}$ if $A_0=A_1=1$ by Assumption A1 and then the law of iterated expectations.
The second line follows because $E (A_0 A_1 |L_0, L_1^{1,1}, Y^{1,1} ) = P(A_0=1| L_0, L_1^{1,1}, Y^{1,1} ) P(A_1=1 | A_0=1, L_0, L_1^{1,1}, Y^{1,1} )$ and,
 by Assumption A1$^\dag$ ($L_1 = L_1^{1,1}$ if $A_0=1$),
$P (A_1=1 | A_0=1, L_0, L_1^{1,1}, Y^{1,1} ) = P (A_1=1 | A_0=1, L_0, L_1, Y^{1,1} ) = \varpi_1 (\overline{L}_1, Y^{1,1} ) >0$ almost surely.
The third line follows because
$P (A_0=1 | L_0, L_1^{1,1}, Y^{1,1} ) = \varpi_{0,\mytext{Jnt}} ( L_0, L_1^{1,1}, Y^{1,1} ) >0 $ almost surely.

The ICE identity follows from the IPW identity and Lemma \ref{lem:IPW-ICE}.
Alternatively, a direct proof of the ICE identity can be obtained similarly as in the proof of Lemma~\ref{lem:iden-mu-2period}.
The corresponding two steps proceed as follows:
\begin{align*}
& \quad E \left\{ (\pi^*_0 + (1-\pi^*_0) \lambda^*_{0,\mytext{Jnt}} (\overline{L}_1, Y^{1,1})) Y^{1,1} | A_0=1, \overline{L}_1 \right\} \nonumber \\
& = E \left\{ (\pi^*_0 + (1-\pi^*_0) \lambda^*_{0,\mytext{Jnt}} (\overline{L}_1, Y )) (\pi^*_1 + (1-\pi^*_1) \lambda^*_1) Y |A_0=1, A_1=1, \overline{L}_1 \right\} ,
\end{align*}
and
\begin{align*}
& \quad E (Y^{1,1} | L_0)
= E \left[ (\pi^*_0 + (1-\pi^*_0) \lambda^*_{0,\mytext{Jnt}} (L_0,L_1^{1,1}, Y^{1,1}))Y^{1,1}  | A_0=1, L_0\right] \nonumber \\
& = E \left[ E \left\{ (\pi^*_0 + (1-\pi^*_0) \lambda^*_{0,\mytext{Jnt}} (\overline{L}_1, Y^{1,1}))Y^{1,1}  | A_0=1, \overline{L}_1 \right\} | A_0=1, L_0\right] .
\end{align*}
Combining these two displays gives
\begin{align}
& \quad E (Y^{1,1} | L_0)
= E \left[ (\pi^*_0 + (1-\pi^*_0) \lambda^*_{0,\mytext{Jnt}} (L_0,L_1^{1,1}, Y^{1,1}))Y^{1,1}  | A_0=1, L_0\right] \nonumber \\
& = E \left[ E \left\{ (\pi^*_0 + (1-\pi^*_0) \lambda^*_{0,\mytext{Jnt}} (\overline{L}_1, Y^{1,1}))Y^{1,1}  | A_0=1, \overline{L}_1 \right\} | A_0=1, L_0\right] \nonumber \\
& = E \left[
E \left\{ (\pi^*_0 + (1-\pi^*_0) \lambda^*_{0,\mytext{Jnt}} ) (\pi^*_1 + (1-\pi^*_1) \lambda^*_1) Y |A_0=1, A_1=1, \overline{L}_1 \right\}
| A_0=1, L_0 \right] .  \label{eq:iden-lam-Jnt-prf2}
\end{align}
Using $\mu^{1,1} = E \{ E (Y^{1,1} | L_0)\}$ then leads to the ICE identity.

\subsection{Proof of Proposition~\ref{pro:upper-2period-Jnt}}  \label{sec:prf-pro-upper-Kperiod-Jnt}

First, we show that $\mu^{1,1}_{+,\mytext{Jnt}} \le \mu^{1,1}_+$. It suffices to show that
model (\ref{eq:model-2period-Jnt}) is more restrictive than model (\ref{eq:model-2period}),
i.e., if model (\ref{eq:model-2period-Jnt}) holds, then model (\ref{eq:model-2period}) holds.
This then indicates that for any distribution $Q$ on $(\overline{A}_1, L_0,\{(L_1^{a_0a_1},Y^{a_0a_1}): a_0,a_1=0,1\} )$
allowed in the optimization in (\ref{eq:upper-2period-Jnt}),
the induced distribution of $Q$ on $(\overline{A}_1, \overline{L}_1, Y^{1,1},Y^{1,0},Y^{0,1},Y^{0,0})$, denoted as $\tilde Q$,
is allowed in the optimization in (\ref{eq:upper-2period-c}).
Moreover, the ICE (or equivalently IPW) functionals in (\ref{eq:upper-2period-Jnt}) and (\ref{eq:upper-2period-c}), evaluated
at $Q$ and $\tilde Q$ respectively, are the same by Lemmas~\ref{lem:iden-mu-2period-Jnt} and \ref{lem:iden-mu-2period}:
\begin{align*}
\mu^{1,1}_Q = \mu^{1,1}_{\mytext{ICE}} (\lambda^*_{0,\mytext{Jnt},Q},\lambda^*_{1,Q})
= \mu^{1,1}_{\mytext{ICE}} (\lambda^*_{0,\tilde Q},\lambda^*_{1,\tilde Q})
\end{align*}
where $\mu^{1,1}_Q$ is the mean of $Y^{1,1}$ under $Q$ or under $\tilde Q$.
Taking the maximum of $\mu^{1,1}_{\mytext{ICE}} (\lambda^*_{0,\mytext{Jnt},Q},\lambda^*_{1,Q})$ in
the preceding display over allowed $Q$
shows that $\mu^{1,1}_{+,\mytext{Jnt}} \le \mu^{1,1}_+$.

To show that model (\ref{eq:model-2period-Jnt}) implies model (\ref{eq:model-2period}),
we use the following result. If $p_0(u,v)$ and $p_1(u,v)$ are two joint densities with $p_0$ absolutely continuous with respect to $p_1$, then
\begin{align*}
E_{p_1} \left\{ \frac{p_0(U,V) }{p_1(U,V)} \Big| V=v \right\} = \frac{p_0(v)}{p_1(v)},
\end{align*}
where $p_0(v)$ and $p_1(v)$ are the associated marginal densities,
and $E_{p_1}(\cdot)$ denotes the expectation under $p_1(u,v)$.
Then $\lambda^*_{0,\mytext{Jnt}}$ and $\lambda^*_0$ are related as follows:
\begin{align}
E \left\{ \lambda^*_{0,\mytext{Jnt}} (L_0, L_1^{1,1}, Y^{1,1}) | A_0=1, L_0, Y^{1,1}=y \right\}
=  \lambda^*_0 (L_0, y). \label{eq:importance-iden}
\end{align}
By the monotonicity of expectations, it follows that if  $\lambda^*_{0,\mytext{Jnt}}$ satisfies the range constraints (\ref{eq:model-2period-Jnt}), then
$\lambda^*_0$ satisfies the range constraints (\ref{eq:model-2period}), i.e., model (\ref{eq:model-2period-Jnt}) implies model (\ref{eq:model-2period}).

Second, we show that $\mu^{1,1}_{+,\mytext{Jnt}}  \ge \mu^{1,1}_+$.
For pedagogical purposes, we give two distinct proofs.
The first proof uses the representation of the sharp bounds through optimization over allowed full-data distributions.
The second proof uses the representation of the sharp bounds through optimization over allowed sensitivity ratios.

\textbf{First proof of $\mu^{1,1}_{+,\mytext{Jnt}}  \ge \mu^{1,1}_+$.}\;
We show that for any distribution $Q$ on $(\overline{A}_1, \overline{L}_1, Y^{1,1},Y^{1,0},$ $Y^{0,1},Y^{0,0})$
allowed in the optimization in (\ref{eq:upper-2period}),
there exists a probability distribution $\tilde Q$ on $(\overline{A}_1, L_0,\{(L_1^{a_0a_1},$ $Y^{a_0a_1}): a_0,a_1=0,1\} )$
such that the induced distribution of $\tilde Q$ on  $(\overline{A}_1, \overline{L}_1, Y )$ is the same as that of $P$,
and $\lambda^*_{1,\tilde Q} (\overline{L}_1, y) = \lambda^*_{1,Q} (\overline{L}_1,y) $ and
$\lambda^*_{0,\mytext{Jnt},\tilde Q} (L_0, l_1, y) = \lambda^*_{0,Q} (L_0,y)$ for any $l_1$ and $y$.
(Note that the meanings of $Q$ and $\tilde Q$ are interchanged, compared with those in the proof of $\mu^{1,1}_{+,\mytext{Jnt}}  \le \mu^{1,1}_+$.)
This then indicates that $\tilde Q$ is allowed in the optimization in (\ref{eq:upper-2period-Jnt}), with $\lambda^*_{1,\tilde Q}$
and $\lambda^*_{0,\mytext{Jnt},\tilde Q}$
satisfying the range constraints, and
\begin{align*}
\mu^{1,1}_Q
= \mu^{1,1}_{\mytext{ICE}} (\lambda^*_{0,Q},\lambda^*_{1,Q})
= \mu^{1,1}_{\mytext{ICE}} (\lambda^*_{0,\mytext{Jnt},\tilde Q},\lambda^*_{1,\tilde Q}),
\end{align*}
where $\mu^{1,1}_Q$ is the mean of $Y^{1,1}$ under $Q$ or under $\tilde Q$.
Taking the maximum of $\mu^{1,1}_{\mytext{ICE}} (\lambda^*_{0,Q},\lambda^*_{1,Q}) $ in
the preceding display over  allowed $Q$
shows that $\mu^{1,1}_+ \le \mu^{1,1}_{+,\mytext{Jnt}}$.

We construct the aforementioned $\tilde Q$ as follows. In the first step, we let
\begin{align}
\dif \tilde Q_{\overline{A}_1, L_0, Y^{1,1} } = \dif Q_{\overline{A}_1, L_0, Y^{1,1} } .  \label{eq:upper-2periond-Jnt-prf}
\end{align}
In the second step, we define
\begin{align*}
& \quad \dif \tilde Q_{L_1^{1,1}} (l_1 | A_0=1, A_1, L_0, Y^{1,1} )
= \dif Q_{L_1} (l_1 | A_0=1, A_1, L_0, Y^{1,1} ),
\end{align*}
i.e., the conditional distribution of $L_1^{1,1}$ given $(A_0=1,A_1,L_0, Y^{1,1})$ under $\tilde Q$
is the same as that of $L_1$ given $(A_0=1,A_1,L_0, Y^{1,1})$ under $Q$.
Combining these two steps yields the conditional distribution of $(A_1, L_0, L_1^{1,1}, Y^{1,1})$ given $A_0=1$
under $\tilde Q$.
The conditional distribution of $(A_1, L_0, L_1, Y^{1,1})$ given $A_0=1$
under $\tilde Q$, with $L_1 = L_1^{1,1}$ if $A_0=1$ by Assumption A1$^\dag$,
is the same as that of $(A_1, L_0, L_1, Y^{1,1})$ given $A_0=1$ under $Q$. Hence
$\lambda^*_{1,\tilde Q} (\overline{L}_1, y) = \lambda^*_{1,Q} (\overline{L}_1,y) $ for any $y$.
In the third step, we define
\begin{align*}
& \quad \dif \tilde Q_{L_1^{1,1}} (l_1 | A_0=0, A_1, L_0, Y^{1,1} ) \\
& = \dif \tilde Q_{L_1^{1,1}} (l_1 | A_0=0 , L_0, Y^{1,1} ) = \dif \tilde Q_{L_1^{1,1}} (l_1 | A_0=1 , L_0, Y^{1,1} ),
\end{align*}
where $\dif \tilde Q_{L_1^{1,1}} (l_1 | A_0=1 , L_0, Y^{1,1} )$ is determined from the first two steps.
That is, $L_1^{1,1}$ is conditionally independent of $A_1$ given $(A_0=0, L_0, Y^{1,1})$ under $\tilde Q$
[notationally, $L^{1,1}_1 \perp A_1 | A_0=0, L_0, Y^{1,1} $ under $\tilde Q$], and
$L_1^{1,1}$ is conditionally independent of $A_0$ given $(L_0, Y^{1,1})$ under $\tilde Q$
[notationally, $L^{1,1}_1 \perp A_0 | L_0, Y^{1,1}$ under $\tilde Q$].
The latter conditional independence together with the fact that by (\ref{eq:upper-2periond-Jnt-prf})
the conditional distribution of $Y^{1,1}$ given $(A_0,L_0)$ under $\tilde Q$ is the same as under $Q$ implies that for any $l_1$ and $y$,
\begin{align*}
\lambda^*_{0,\mytext{Jnt},\tilde Q} (L_0, l_1, y) = \lambda^*_{0,\tilde Q} (L_0,y) = \lambda^*_{0,Q} (L_0,y) .
\end{align*}
Combining the proceeding three steps  completes the definition of $\tilde Q$ for
$ (\overline{A}_1, L_0, L_1^{1,1},Y^{1,1} )$.
In the final step, we define
$\dif \tilde Q_{\{(L_1^{a_0a_1},Y^{a_0a_1}): a_0,a_1=0,1, (a_0,a_1)\not=(1,1)\}} (\cdot | \overline{A}_1, L_0, L_1^{1,1},Y^{1,1})$
in the same manner as (\ref{eq:density-ratio-Q-Jnt-final1})--(\ref{eq:density-ratio-Q-Jnt-final2}) with $Q$ replaced by $\tilde Q$.
By design, the induced distribution of $\tilde Q$ on  $(\overline{A}_1, \overline{L}_1, Y )$
can be easily verified to be the same as that of $P$.

From this proof, it also follows that if the sharp upper bound $\mu^{1,1}_+$ is achieved by a distribution $Q$ on
$(\overline{A}_1, \overline{L}_1, Y^{1,1},Y^{1,0},Y^{0,1},Y^{0,0})$,
then $\mu^{1,1}_{+,\mytext{Jnt}}$ is achieved by the distribution $\tilde Q$ constructed above, for which  $L_1^{1,1} \perp A_0 | L_0, Y^{1,1}$.

\textbf{Second proof of $\mu^{1,1}_{+,\mytext{Jnt}}  \ge \mu^{1,1}_+$.}\;
We use the representation (\ref{eq:upper-2period-b}) for $\mu^{1,1}_+$ and the representation (\ref{eq:upper-2period-Jnt-b})
for $\mu^{1,1}_{+,\mytext{Jnt}}$. For any nonnegative functions $\lambda_1 (\overL_1,y)$ and $\lambda_0 (L_0, y)$ allowed in
(\ref{eq:upper-2period-b}), i.e., satisfying the constraints (\ref{eq:constr-2period-b}) and (\ref{eq:model-2period-b}),
the functions $\lambda _1 (\overL_1, y)$ and $\lambda_{0,\mytext{Jnt}} (L_0, l_1, y) =\lambda_0 (L_0, y)$, which is independent of $l_1$,
can be easily shown to be allowed in (\ref{eq:upper-2period-Jnt-b}), i.e., satisfying the constraints (\ref{eq:constr-2period-Jnt-b}) and (\ref{eq:model-2period-Jnt-b}).
For these choices, $\mu^{1,1}_{\mytext{ICE}} (\lambda_0,\lambda_1) = \mu^{1,1}_{\mytext{ICE}} (\lambda_{0,\mytext{Jnt}},\lambda_1)$.
Taking the maximum over allowed $(\lambda_0,\lambda_1)$ implies that $\mu^{1,1}_+ \le \mu^{1,1}_{+,\mytext{Jnt}}$.

From the second proof, it also follows that if the sharp upper bound $\mu^{1,1}_+$ is achieved by $\lambda_1 (\overL_1,y)$ and $\lambda_0 (L_0, y)$,
then  $\mu^{1,1}_{+,\mytext{Jnt}}$ is achieved by $\lambda _1 (\overL_1, y)$ and $\lambda_{0,\mytext{Jnt}} (L_0, l_1, y)=\lambda_0 (L_0, y)$ discussed above.
For the distribution $Q$ associated with these $\lambda_1$ and $\lambda_{0,\mytext{Jnt}}$ from Lemma~\ref{lem:density-ratio-Jnt},
applying the identity (\ref{eq:importance-iden}) shows that $ \lambda^*_{0,\mytext{Jnt},Q} (L_0, l_1, y) = \lambda^*_{0,Q} (L_0, y) $
and hence
$\dif Q_{L^{1,1}_1} (l_1 | A_0=0, L_0, Y^{1,1}) / \dif Q_{L^{1,1}_1} (l_1 | A_0=1, L_0, Y^{1,1}) \equiv 1$, i.e., $L^{1,1}_1 \perp A_0 | L_0, Y^{1,1}$.

\section{Technical details: Proofs for Section \ref{sec:2period-Prod}}

\subsection{Proof of Lemma~\ref{lem:density-ratio-Prod}}

The necessity assertion can be proved similarly as that in Lemma~\ref{lem:density-ratio}.
In the following, we prove the sufficiency assertion.
The proof differs from that of the sufficiency assertion in Lemma \ref{lem:density-ratio-Jnt}
in only one of the steps involved.

We construct a desired probability distribution $Q$ for the full data $(\overline{A}_1, L_0,\{(L_1^{a_0a_1},Y^{a_0a_1}): a_0,a_1=0,1\} )$.
The first three steps and the fifth and final (sixth) steps
are the same as in the proof of the sufficiency in Lemma \ref{lem:density-ratio-Jnt} (Section \ref{sec:prf-lem-density-ratio-Jnt}).
For the fourth step,  we define the conditional distribution of $Y^{1,1}$ given $(A_0=0, L_0, L_1^{1,1})$ as follows:
\begin{align}
& \dif Q_{Y^{1,1}} ( y | A_0=0, L_0, L_1^{1,1}=l_1 ) = \lambda_{0,Y} (L_0, l_1, y) \dif Q_{Y^{1,1}} ( y | A_0=1, L_0, L_1^{1,1}=l_1 ) \label{eq:density-ratio-Q-Prod-Y} \\
& = \lambda_{0,Y} (L_0, l_1, y)  \varrho_1 (L_0, l_1, y; \lambda_1 ) \dif P_{Y} (y | A_0=1,A_1=1, L_0, L_1 =l_1) , \nonumber
\end{align}
which yields a proper distribution (i.e., integrated to 1) due to (\ref{eq:constr-2period-Prod-b0Y}),
and define the conditional distribution of $L_1^{1,1}$ given $(A_0=0, L_0)$ as follows:
\begin{align}
& \dif Q_{L_1^{1,1}} (l_1 | A_0=0, L_0) = \lambda_{0,L_1} (L_0, l_1 ) \dif Q_{L_1^{1,1}} (l_1 | A_0=1, L_0),\label{eq:density-ratio-Q-Prod-L1}
\end{align}
which yields a proper distribution (i.e., integrated to 1) due to (\ref{eq:constr-2period-Prod-b0L}).
Combining the two definitions determines the conditional distribution of $(L_1^{1,1}, Y^{1,1})$ given $(A_0=0, L_0)$ as
\begin{align}
& \quad \dif Q_{L_1^{1,1},Y^{1,1}} (l_1, y | A_0=0, L_0 ) \nonumber \\
& = \dif Q_{Y^{1,1}} (y | A_0=0,L_0, L_1^{1,1}=l_1) \times \dif Q_{L_1^{1,1}} (l_1 | A_0=0,L_0) \nonumber \\
& =  \lambda_{0,L_1} (L_0, l_1 )\lambda_{0,Y} (L_0, l_1, y) \dif Q_{L_1^{1,1},Y^{1,1}} (l_1, y | A_0=1, L_0). \label{eq:density-ratio-Q-Prod-A0}
\end{align}
As indicated by (\ref{eq:density-ratio-Q-Jnt-A0}) and (\ref{eq:density-ratio-Q-Prod-A0}),
the distribution $Q$ is the same as that in the proof of the sufficiency assertion in Lemma \ref{lem:density-ratio-Jnt},
provided that $\lambda_{0,\mytext{Jnt}} (L_0, l_1,y) = \lambda_{0,L_1} (L_0, l_1 )\lambda_{0,Y} (L_0, l_1, y) $.
This choice of $\lambda_{0,\mytext{Jnt}} (L_0, l_1,y)$ can be easily verified to satisfy constraint (\ref{eq:constr-2period-Jnt-b0}) by using
the constraints (\ref{eq:constr-2period-Prod-b0Y})--(\ref{eq:constr-2period-Prod-b0L}) satisfied by
$\lambda_{0,Y}$ and $\lambda_{0,L_1}$.
However, by the refined
definitions (\ref{eq:density-ratio-Q-Prod-Y}) and (\ref{eq:density-ratio-Q-Prod-L1}) behind (\ref{eq:density-ratio-Q-Prod-A0}) based on
 $\lambda_{0,Y}$ and $\lambda_{0,L_1}$,
the distribution $Q$ defined here further satisfies $ \lambda^*_{0,Y,Q} (L_0, l_1, y) = \lambda_{0,Y} (L_0, l_1, y)$ and
 $ \lambda^*_{0,L_1,Q} (L_0, l_1) = \lambda_{0,L_1} (L_0, l_1)$,

\subsection{Proof of Proposition~\ref{pro:upper-2period-Prod}}  \label{sec:prf-pro-upper-2period-Prod}

We show that
$ \mu^{1,1}_{+,\mytext{Prod}} \le  \mu^{1,1}_{+,\mytext{Prod,v1}} $
and
$ \mu^{1,1}_{+,\mytext{Prod}} \le  \mu^{1,1}_{+,\mytext{Prod,v2}} $.
By Lemma~\ref{lem:density-ratio-Prod}, for $\lambda_{0,L_1}(L_0,l_1)$, $\lambda_{0,Y}(L_0,l_1,y)$, and $\lambda_1 (\overline{L}_1,y)$
satisfying (\ref{eq:constr-2period-Prod-b}),
there exists a probability distribution $Q$ for the full data satisfying properties (i) and (ii).
In other words, $Q$ can be a valid choice for $P$.
For $\lambda_{0,\mytext{Jnt}}=\lambda_{0,L_1}\lambda_{0,Y}$, applying (\ref{eq:iden-lam-Jnt-prf2}) with $P$ replaced by $Q$ yields
\begin{align*}
& \quad E \left[  E \left\{(\pi^*_0 + (1-\pi^*_0) \lambda_{0,L_1}\lambda_{0,Y})(\pi^*_1 + (1-\pi^*_1) \lambda_1) Y |A_0=1, A_1=1, \overline{L}_1 \right\} | A_0=1, L_0\right] \\
& = E_{\lambda_1} \Big[ E_{\lambda_1} \Big\{ \underbrace{  (\pi^*_0 + (1-\pi^*_0)  \lambda_{0,L_1}(L_0,L_1^{1,1}) \lambda_{0,Y} (L_0,L_1^{1,1},Y^{1,1})) Y^{1,1}  }_{\one}  |
 A_0=1, \overline{L}_1 \Big\} | A_0=1, L_0  \Big] \\
& = E_{\lambda_1} \left[ \;\one\; | A_0=1, L_0\right] ,
\end{align*}
where $E_{\lambda_1} (\cdot | A_0=1, \overline{L}_1)$ and $E_{\lambda_1} (\cdot | A_0=1, L_0)$
denote the conditional expectations under $Q$, depending on $\lambda_1$.
Hence $\mu^{1,1}_{\mytext{ICE}} (\lambda_{0,L_1}\lambda_{0,Y}, \lambda_1)$ can be expressed as
\begin{align}
& \quad \mu^{1,1}_{\mytext{ICE}} (\lambda_{0,L_1}\lambda_{0,Y},\lambda_1) = E \left( E_{\lambda_1} \left[ E_{\lambda_1} \left\{ \;\one\; |
 A_0=1, \overline{L}_1 \right\} | A_0=1, L_0\right]  \right) \nonumber \\
& =E \left( E_{\lambda_1} \left[ \;\one\; | A_0=1, L_0\right] \right). \label{eq:pro-upper-2period-Prod-prf1}
\end{align}

\textbf{First upper bound.}\; We show that
$ \mu^{1,1}_{+,\mytext{Prod}} \le  \mu^{1,1}_{+,\mytext{Prod,v1}} $.
For any fixed $\lambda_1$, the normalization constraint (\ref{eq:constr-2period-Prod-b0Y}) on $\lambda_{0,Y}$ is equivalent to
\begin{align*}
E_{\lambda_1} \left( \lambda_{0,Y} (\overline{L}_1,Y^{1,1})  | A_0=1, \overline{L}_1 \right) \equiv 1.
\end{align*}
For any fixed $\lambda_1$ and $\lambda_{0,L_1}$,
we apply identity (\ref{eq:upper-1period-c}), with $L_0$, $\lambda_0$, and $b(L_0,Y)$ set to
$L_1$, $\lambda_{0,Y}$, and $ (1-\pi^*_0) \lambda_{0,L_1}(L_0,L_1) Y^{1,1}$ respectively,
conditionally on $(A_0=1,L_0)$, and obtain
\begin{align}
& \quad \max_{\lambda_{0,Y}} \; E_{\lambda_1} \left[ E_{\lambda_1} \left\{ \; \one \; |
 A_0=1, \overline{L}_1 \right\}  | A_0=1, L_0\right] \nonumber  \\
& = \min_{q_{0,Y}} \; E_{\lambda_1} \Big[  E_{\lambda_1} \Big\{ \underbrace{ \pi^*_0 Y^{1,1} + (1-\pi^*_0)  \lambda_{0,L_1}
 ( Y^{1,1} + \tilde \Lambda_{0,Y} \rho_{\tau_{0,Y}} (Y^{1,1}, q_{0,Y}) ) }_{\two}  |
 A_0=1, \overline{L}_1 \Big\}  | A_0=1, L_0 \Big] \label{eq:pro-upper-2period-Prod-prf2} \\
& = \min_{q_{0,Y}} \; E_{\lambda_1} \Big[ \underbrace{ \pi^*_0 \mathcal E_{0,\overline{L}_1} (Y^{1,1})
+ (1-\pi^*_0)  \lambda_{0,L_1} \mathcal E_{0,\overline{L}_1} ( Y^{1,1} + \tilde \Lambda_{0,Y}
\rho_{\tau_{0,Y}} (Y^{1,1}, q_{0,Y}) ) }_{\three}  | A_0=1, L_0 \Big]  ,\nonumber
\end{align}
where $\tilde \Lambda_{0,Y} = \Lambda_{0,Y} - \Lambda_{0,Y}^{-1}$,
and $\mathcal E_{0,\overline{L}_1}(\cdot)$ denotes the conditional expectation $E_{\lambda_1} (\cdot |A_0=1,\overline{L}_1)$.
For any fixed $\lambda_1$ and $q_{0,Y}$, we apply identity (\ref{eq:upper-1period-b}) with
$L_0$, $A_0$, $\lambda_0$, and $ b(L_0,Y)$ set to $L_0$, $A_0$, $\lambda_{0,L_1}$ and $\mathcal E_{0,\overline{L}_1} ( Y^{1,1} + \tilde \Lambda_{0,Y}
\rho_{\tau_{0,Y}} (Y^{1,1}, q_{0,Y}) )$, and obtain
\begin{align*}
& \quad \max_{\lambda_{0,L_1}} \; E \left(  E_{\lambda_1} [ \;\three\; | A_0=1, L_0 ] \right) \\
& = \min_{q_{0,L_1}} \; E \Big( E_{\lambda_1} \Big[ \underbrace{ \pi^*_0 \mathcal E_{0,\overline{L}_1} (Y^{1,1})
+ (1-\pi^*_0)  \mathcal E_{0,\overline{L}_1} ( Y^{1,1} + \tilde \Lambda_{0,Y}
\rho_{\tau_{0,Y}} (Y^{1,1}, q_{0,Y}) ) }_{\four} \\
& \qquad  + \underbrace{ (1-\pi^*_0)  \tilde \Lambda_{0,L_1} \rho_{\tau_{0,L_1}} \left\{ \mathcal E_{0,\overline{L}_1} ( Y^{1,1} + \tilde \Lambda_{0,Y}
\rho_{\tau_{0,Y}} (Y^{1,1}, q_{0,Y}) ) , q_{0,L_1} \right\} }_{\five} | A_0=1, L_0 \Big] \Big),
\end{align*}
where $\tilde \Lambda_{0,L_1} = \Lambda_{0,L_1} - \Lambda_{0,L_1}^{-1}$.
For any fixed $\lambda_1$, combining the preceding three displays yields
\begin{align}
& \quad \max_{\lambda_{0,L_1}} \max_{\lambda_{0,Y}} \; \mu^{1,1}_{\mytext{ICE}} (\lambda_{0,L_1}\lambda_{0,Y},\lambda_1) \nonumber \\
& =  \max_{\lambda_{0,L_1}}  \min_{q_{0,Y}} \; E \left( E_{\lambda_1} \left[  \;\three\;  | A_0=1, L_0\right] \right) \nonumber \\
& =  \min_{q_{0,Y}}  \max_{\lambda_{0,L_1}} \; E \left( E_{\lambda_1} \left[  \;\three\;  | A_0=1, L_0\right] \right) \nonumber \\
& =  \min_{q_{0,Y}}  \min_{q_{0,L_1}} \; E \left( E_{\lambda_1} \left[  \;\four + \five \;   | A_0=1, L_0\right] \right) . \label{eq:pro-upper-2period-Prod-prf3}
\end{align}
The third line follows from Sion's minimax theorem, because the functional $\three$ is linear in $\lambda_{0,L_1}$ and convex in $q_{0,Y}$.
By Jensen's inequality with the convexity of $\rho_{\tau_{0,L_1}} (y, q)$ in $y$,
\begin{align*}
& \quad \rho_{\tau_{0,L_1}} \left\{ \mathcal E_{0,\overline{L}_1} ( Y^{1,1} + \tilde \Lambda_{0,Y}
\rho_{\tau_{0,Y}} (Y^{1,1}, q_{0,Y}) ) , q_{0,L_1} \right\} \\
& \le \mathcal E_{0,\overline{L}_1} \left\{ \rho_{\tau_{0,L_1}} ( Y^{1,1} + \tilde \Lambda_{0,Y}
\rho_{\tau_{0,Y}} (Y^{1,1}, q_{0,Y}) , q_{0,L_1} )  \right\},
\end{align*}
and hence
\begin{align*}
& \five \le (1-\pi^*_0)  \tilde \Lambda_{0,L_1}
\mathcal E_{0,\overline{L}_1} \left\{ \rho_{\tau_{0,L_1}} ( Y^{1,1} + \tilde \Lambda_{0,Y}
\rho_{\tau_{0,Y}} (Y^{1,1}, q_{0,Y}) , q_{0,L_1} )  \right\} .
\end{align*}
Substituting this into (\ref{eq:pro-upper-2period-Prod-prf3}) yields
\begin{align}
& \quad \max_{\lambda_{0,L_1}} \max_{\lambda_{0,Y}} \; \mu^{1,1}_{\mytext{ICE}} (\lambda_{0,L_1}\lambda_{0,Y},\lambda_1) \nonumber \\
& \le \min_{q_{0,Y}}  \min_{q_{0,L_1}} \; E \Big( E_{\lambda_1} \Big[ \pi^*_0 \mathcal E_{0,\overline{L}_1} (Y^{1,1})
 + (1-\pi^*_0)  \mathcal E_{0,\overline{L}_1} ( Y^{1,1} + \tilde \Lambda_{0,Y} \rho_{\tau_{0,Y}} (Y^{1,1}, q_{0,Y}) ) \nonumber \\
& \qquad + (1-\pi^*_0)  \tilde \Lambda_{0,L_1}
\mathcal E_{0,\overline{L}_1} \left\{ \rho_{\tau_{0,L_1}} ( Y^{1,1} + \tilde \Lambda_{0,Y}
\rho_{\tau_{0,Y}} (Y^{1,1}, q_{0,Y}) , q_{0,L_1} )  \right\}  | A_0=1, L_0\Big] \Big) \nonumber \\
& = \min_{q_{0,Y}}  \min_{q_{0,L_1}} \; E \Big( E_{\lambda_1} \Big[ \pi^*_0  Y^{1,1}
 + (1-\pi^*_0)  ( Y^{1,1} + \tilde \Lambda_{0,Y} \rho_{\tau_{0,Y}} (Y^{1,1}, q_{0,Y}) ) \nonumber  \\
& \qquad + (1-\pi^*_0)  \tilde \Lambda_{0,L_1}
 \rho_{\tau_{0,L_1}} ( Y^{1,1} + \tilde \Lambda_{0,Y}
\rho_{\tau_{0,Y}} (Y^{1,1}, q_{0,Y}) , q_{0,L_1} )   | A_0=1, L_0\Big] \Big) \nonumber  \\
& = \min_{q_{0,Y}}  \min_{q_{0,L_1}} \; E \Big( E_{\lambda_1} \Big[
\underbrace{ \eta_{0+,L_1} ( Y^{1,1}, \eta_{0+,Y} (Y^{1,1}, q_{0,Y} ), q_{0,L_1} ) }_{\six}
 | A_0=1, L_0\Big] \Big)  , \label{eq:pro-upper-2period-Prod-prf4}
\end{align}
where the third step follows from the law of iterated expectations,
and the last step follows from the definitions of $\eta_{0+,Y}$ and $\eta_{0+,L_1}$.
For any fixed $q_{0,Y}$ and $q_{0,L_1}$, we apply
identity (\ref{eq:upper-1period-b}), with $L_0$, $A_0$, $\lambda_0$, and $b(L_0, Y)$ set to
$L_1$, $A_1$, $\lambda_1$, and $\six$ respectively,
conditionally on $(A_0=1,L_0)$, and obtain
\begin{align*}
& \quad \max_{\lambda_1} \;  E_{\lambda_1} \left\{ \;\six\; | A_0=1, L_0 \right\}\\
& = \max_{\lambda_1} \; E\left[ E \left\{
(\pi^*_1 + (1-\pi^*_1) \lambda_1 (\overline{L}_1,Y )) \six |A_0=1,A_1=1,\overline{L}_1 \right\} | A_0=1, L_0 \right] \\
& = \min_{q_1} \; E\left[ E \left\{ \eta_{1+} (\six, q_1) | A_0=1, A_1=1, \overline{L}_1 \right\} | A_0=1, L_0 \right] .
\end{align*}
Combining the preceding two displays yields
\begin{align*}
& \quad \max_{\lambda_1} \max_{\lambda_{0,L_1}} \max_{\lambda_{0,Y}} \; \mu^{1,1}_{\mytext{ICE}} (\lambda_{0,L_1}\lambda_{0,Y},\lambda_1)  \\
& \le \max_{\lambda_1} \min_{q_{0,Y}}  \min_{q_{0,L_1}} \; E\left[ E \left\{
(\pi^*_1 + (1-\pi^*_1) \lambda_1 (\overline{L}_1,Y )) \six |A_0=1,A_1=1,\overline{L}_1 \right\} | A_0=1, L_0 \right] \\
& \le \min_{q_{0,Y}}  \min_{q_{0,L_1}} \max_{\lambda_1}   \; E\left[ E \left\{
(\pi^*_1 + (1-\pi^*_1) \lambda_1 (\overline{L}_1,Y )) \six |A_0=1,A_1=1,\overline{L}_1 \right\} | A_0=1, L_0 \right] \\
& = \min_{q_{0,Y}}  \min_{q_{0,L_1}} \min_{q_1} \;  E\left[ E \left\{ \eta_{1+} (\six, q_1) | A_0=1, A_1=1, \overline{L}_1 \right\} | A_0=1, L_0 \right] .
\end{align*}
The third line follows from two applications of
the inequality that $ \max_u \min_v f(u,v) \le \min_v \max_u f(u,v)$ for any function $f$.
This concludes that $ \mu^{1,1}_{+,\mytext{Prod}} \le  \mu^{1,1}_{+,\mytext{Prod,v1}} $.

\textbf{Second upper bound.}\; We show that
$ \mu^{1,1}_{+,\mytext{Prod}} \le  \mu^{1,1}_{+,\mytext{Prod,v2}} $, through a different path departing from (\ref{eq:pro-upper-2period-Prod-prf2}).
For any fixed $\lambda_{0,L_1}$ and $q_{0,Y}$, we apply
identity (\ref{eq:upper-1period-b}), with $L_0$, $A_0$, $\lambda_0$, and $Y$ set to
$L_1$, $A_1$, $\lambda_1$, and $\two$ respectively,
conditionally on $(A_0=1,L_0)$, and obtain
\begin{align*}
& \quad \max_{\lambda_1} \;  E_{\lambda_1} \left\{ \;\two\; | A_0=1, L_0 \right\}\\
& = \max_{\lambda_1} \; E\left[ E \left\{
(\pi^*_1 + (1-\pi^*_1) \lambda_1 (\overline{L}_1,Y )) \two |A_0=1,A_1=1,\overline{L}_1 \right\} | A_0=1, L_0 \right] \\
& = \min_{q_1} \; E\left[ E \left\{ \eta_{1+} (\two, q_1) | A_0=1, A_1=1, \overline{L}_1 \right\} | A_0=1, L_0 \right] ,
\end{align*}
where $\two = \pi^*_0 Y^{1,1} + (1-\pi^*_0)  \lambda_{0,L_1} \eta_{0+,Y} (Y^{1,1}, q_{0,Y} )$ as defined earlier and
$\eta_{1+} (\two, q_1) = \two + (1-\pi^*_1) \tilde \Lambda_1 \rho_{\tau_1} (\two, q_1)$ with $\tilde \Lambda_1 = \Lambda_1 - \Lambda_1^{-1}$.
For any fixed $\lambda_{0,L_1}$, combining the preceding display with (\ref{eq:pro-upper-2period-Prod-prf1}) and (\ref{eq:pro-upper-2period-Prod-prf2}) yields
\begin{align}
& \quad \max_{\lambda_1} \max_{\lambda_{0,Y}} \; \mu^{1,1}_{\mytext{ICE}} (\lambda_{0,L_1}\lambda_{0,Y},\lambda_1) \nonumber \\
& =  \max_{\lambda_1}  \min_{q_{0,Y}} \; E \left( E_{\lambda_1} \left[  \;\two\;  | A_0=1, L_0\right] \right) \nonumber \\
& =  \min_{q_{0,Y}}  \max_{\lambda_1} \; E \left( E_{\lambda_1} \left[  \;\two\;  | A_0=1, L_0\right] \right) \nonumber \\
& =  \min_{q_{0,Y}} \min_{q_1} \; E \left( E\left[ E \left\{ \eta_{1+} (\two, q_1) | A_0=1, A_1=1, \overline{L}_1 \right\} | A_0=1, L_0 \right] \right).
\label{eq:pro-upper-2period-Prod-prf5}
\end{align}
The third line follows from Sion's minimax theorem, because the functional $\two$ is linear in $\lambda_1$ and convex in $q_{0,Y}$.
By the convexity of $\rho_{\tau_1}(y,q)$ in $y$,
\begin{align*}
& \quad  \rho_{\tau_1} (\two, q_1)
 = (\pi^*_0 + (1-\pi^*_0) \lambda_{0,L_1}) \rho_{\tau_1} (\two /(\pi^*_0 + (1-\pi^*_0) \lambda_{0,L_1}), q_{1, \lambda_{0,L_1}}) \\
& \le  \pi^*_0 \rho_{\tau_1} (Y^{1,1}, q_{1, \lambda_{0,L_1}}) + (1-\pi^*_0) \lambda_{0,L_1} \rho_{\tau_1} ( \eta_{0+,Y} (Y^{1,1}, q_{0,Y} ), q_{1, \lambda_{0,L_1}}) ,
\end{align*}
and hence
\begin{align*}
\eta_{1+} (\two, q_1) & \le  \pi^*_0  \left\{ Y^{1,1}  +  (1-\pi^*_1) \tilde \Lambda_1 \rho_{\tau_1} (Y^{1,1}, q_{1, \lambda_{0,L_1}})\right\} \\
& \quad + (1 - \pi^*_0)  \lambda_{0,L_1}  \left\{ \eta_{0+,Y} (Y^{1,1}, q_{0,Y} )
  +  (1-\pi^*_1) \tilde \Lambda_1 \rho_{\tau_1} ( \eta_{0+,Y} (Y^{1,1}, q_{0,Y} ), q_{1, \lambda_{0,L_1}}) \right\} \\
& =   \pi^*_0  \eta_{1+} ( Y^{1,1}, q_{1, \lambda_{0,L_1}} ) + (1-\pi^*_0) \lambda_{0,L_1}\eta_{1+} ( \eta_{0+,Y} (Y^{1,1}, q_{0,Y} ), q_{1, \lambda_{0,L_1}} ) ,
\end{align*}
where $q_{1, \lambda_{0,L_1}} = q_1 / (\pi^*_0 + (1-\pi^*_0) \lambda_{0,L_1})$.
Substituting this into (\ref{eq:pro-upper-2period-Prod-prf5}) yields
\begin{align}
& \quad \max_{\lambda_1} \max_{\lambda_{0,Y}} \; \mu^{1,1}_{\mytext{ICE}} (\lambda_{0,L_1}\lambda_{0,Y},\lambda_1) \nonumber \\
& \le \min_{q_{0,Y}} \min_{q_1} \; E \Big( E\Big[  \pi^*_0 \mathcal E_{\overline{L}_1} \left\{ \eta_{1+} ( Y, q_{1, \lambda_{0,L_1}} ) \right\} \nonumber \\
& \qquad + (1-\pi^*_0)\lambda_{0,L_1}
\mathcal E_{\overline{L}_1} \left\{ \eta_{1+} ( \eta_{0+,Y} (Y, q_{0,Y} ), q_{1, \lambda_{0,L_1}} ) \right\} | A_0=1, L_0 \Big] \Big) \nonumber \\
& = \min_{q_{0,Y}} \min_{\tilde q_1} \; \nonumber \\
& \qquad E \Big( E\Big[ \underbrace{ \pi^*_0 \mathcal E_{\overline{L}_1} \left\{ \eta_{1+} ( Y, \tilde q_1 ) \right\}
 + (1-\pi^*_0)\lambda_{0,L_1}
\mathcal E_{\overline{L}_1} \left\{ \eta_{1+} ( \eta_{0+,Y} (Y, q_{0,Y} ), \tilde q_1 ) \right\} }_{\seven} | A_0=1, L_0 \Big] \Big) ,
\label{eq:pro-upper-2period-Prod-prf6}
\end{align}
where the last step follows because the minimization over all possible $q_1$ is equivalent to that over all possible
$\tilde q_1$, with the dependency on $\lambda_{0,L_1}$ absorbed.
For any fixed $\lambda_1$ and $q_{0,Y}$, we apply identity (\ref{eq:upper-1period-b}) with
$L_0$, $A_0$, $\lambda_0$, and $b(L_0,Y)$ set to $L_0$, $A_0$, $\lambda_{0,L_1}$, and
$\mathcal E_{\overline{L}_1} \{ \eta_{1+} ( \eta_{0+,Y} (Y, q_{0,Y} ), \tilde q_1 ) \}$, and obtain
\begin{align}
& \quad \max_{\lambda_{0,L_1}} \; E \Big( E\Big[  \;\seven\; | A_0=1, L_0 \Big] \Big) \nonumber \\
& = \min_{q_{0,L_1}} \; E \Big( E\Big[  \pi^*_0 \mathcal E_{\overline{L}_1} \left\{ \eta_{1+} ( Y, \tilde q_1 ) \right\}
+ (1-\pi^*_0) \mathcal E_{\overline{L}_1} \left\{ \eta_{1+} ( \eta_{0+,Y} (Y, q_{0,Y} ), \tilde q_1 ) \right\}  \nonumber \\
& \qquad + (1-\pi^*_0) \tilde\Lambda_{0,L_1} \rho_{\tau_{0,L_1}} \left(
\mathcal E_{\overline{L}_1} \left\{ \eta_{1+} ( \eta_{0+,Y} (Y, q_{0,Y} ), \tilde q_1 ) \right\}, q_{0,L_1} \right) | A_0=1, L_0 \Big] \Big)  \nonumber \\
& = \min_{q_{0,L_1}} \; E \Big( E\Big[  \eta_{0+,\mytext{Prod}} ( \mathcal E_{\overline{L}_1} \left\{ \eta_{1+} ( Y, \tilde q_1 ) \right\},
\mathcal E_{\overline{L}_1} \left\{ \eta_{1+} ( \eta_{0+,Y} (Y, q_{0,Y} ), \tilde q_1 ) \right\}, q_{0,L_1} )  | A_0=1, L_0 \Big] \Big) .
\label{eq:pro-upper-2period-Prod-prf7}
\end{align}
Combining the preceding two displays yields
\begin{align*}
& \quad  \max_{\lambda_{0,L_1}} \max_{\lambda_1}\max_{\lambda_{0,Y}} \; \mu^{1,1}_{\mytext{ICE}} (\lambda_{0,L_1}\lambda_{0,Y},\lambda_1)  \\
& \le  \max_{\lambda_{0,L_1}} \min_{q_{0,Y}} \min_{\tilde q_1} \; E \Big( E\Big[  \;\seven\; | A_0=1, L_0 \Big] \Big) \\
& \le \min_{q_{0,Y}} \min_{\tilde q_1}   \max_{\lambda_{0,L_1}} \; E \Big( E\Big[  \;\seven\; | A_0=1, L_0 \Big] \Big) \\
& = \min_{q_{0,Y}} \min_{\tilde q_1} \min_{q_{0,L_1}} \; \\
& \qquad E \Big( E\Big[  \eta_{0+,\mytext{Prod}} ( \mathcal E_{\overline{L}_1} \left\{ \eta_{1+} ( Y, \tilde q_1 ) \right\},
\mathcal E_{\overline{L}_1} \left\{ \eta_{1+} ( \eta_{0+,Y} (Y, q_{0,Y} ), \tilde q_1 ) \right\}, q_{0,L_1} ) | A_0=1, L_0  \Big] \Big) .
\end{align*}
The third line follows from two applications of
the inequality that $ \max_u \min_v f(u,v) \le \min_v \max_u f(u,v)$ for any function $f$.
This concludes that $ \mu^{1,1}_{+,\mytext{Prod}} \le  \mu^{1,1}_{+,\mytext{Prod,v2}} $.

\subsection{Proof of Corollary \ref{cor:upper-2period-Prod}}

For results (i) and (ii), we show that in the special case considered, the corresponding conservative bound in Proposition \ref{pro:upper-2period-Prod}
becomes exact and simplified as stated.

(i) Suppose that $\Lambda_{0,L_1}=1$. Then by definition,
$ \eta_{0+,L_1} (y, \tilde y, q_{0,L_1} ) = \pi^*_0 y + (1-\pi^*_0) \tilde y  $
and
\begin{align*}
\eta_{0+,L_1} ( y, \eta_{0+,Y} (y, q_{0,Y} ), q_{0,L_1} ) = \eta_{0+} (y, q_{0,Y}),
\end{align*}
independent of $q_{0,L_1}$. Hence $\mu^{1,1}_{+,\mytext{Prod,v1}}$ is simplified as stated in Corollary \ref{cor:upper-2period-Prod}.
In the proof of Proposition \ref{pro:upper-2period-Prod} in Section \ref{sec:prf-pro-upper-2period-Prod},
it can be easily verified that $\five=0$ and
\begin{align*}
& \quad \six  = \pi^*_0 Y^{1,1} + (1-\pi^*_0) \eta_{0+,Y} (Y^{1,1}, q_{0,Y} ) \\
& = Y^{1,1} + (1-\pi^*_0) \tilde \Lambda_{0,Y} \rho_{\tau_{0,Y}} (Y^{1,1}, q_{0,Y} (\overline{L}_1) )
= \eta_{0,\Lambda_{0,Y}} (Y^{1,1}, q_{0,Y} ).
\end{align*}
The final step in the proof of First upper bound can be modified as follows:
\begin{align*}
& \quad \max_{\lambda_1}  \max_{\lambda_{0,Y}} \; \mu^{1,1}_{\mytext{ICE}} ( \lambda_{0,Y},\lambda_1)  \\
& = \max_{\lambda_1} \min_{q_{0,Y}} \; E\left[ E \left\{
(\pi^*_1 + (1-\pi^*_1) \lambda_1 (\overline{L}_1,Y )) \six |A_0=1,A_1=1,\overline{L}_1 \right\} | A_0=1, L_0 \right] \\
& = \min_{q_{0,Y}}  \max_{\lambda_1}   \; E\left[ E \left\{
(\pi^*_1 + (1-\pi^*_1) \lambda_1 (\overline{L}_1,Y )) \six |A_0=1,A_1=1,\overline{L}_1 \right\} | A_0=1, L_0 \right] \\
& = \min_{q_{0,Y}} \min_{q_1} \;  E\left[ E \left\{ \eta_{1+} (\six, q_1) | A_0=1, A_1=1, \overline{L}_1 \right\} | A_0=1, L_0 \right] .
\end{align*}
The second line follows because the use of Jensen's inequality to upper bound $\five$ is no longer needed.
The third line follows from Sion's minimax theorem, because $\six$ is linear in $\lambda_1$ and convex in $q_{0,Y}$.
This completes the proof of result (i).

(ii) Suppose that $\Lambda_{0,Y}=1$. Then by definition, $\eta_{0+,Y} (y, q_{0,Y}) = y$ and
\begin{align*}
\eta_{0+,L_1} ( y, \eta_{0+,Y} (y, q_{0,Y} ), q_{0,L_1} ) = \eta_{0+} (y, q_{0,L_1}),
\end{align*}
independent of $q_{0,Y}$.  Hence $\mu^{1,1}_{+,\mytext{Prod,v2}}$ is simplified as stated in Corollary \ref{cor:upper-2period-Prod}.
In the proof of Proposition \ref{pro:upper-2period-Prod} in Section \ref{sec:prf-pro-upper-2period-Prod},
it can be easily verified that $ \two = (\pi^*_0 + (1-\pi^*_0) \lambda_{0,L_1} ) Y^{1,1}$ and
\begin{align*}
& \quad \seven = \pi^*_0 \mathcal E_{\overline{L}_1} \left\{ \eta_{1+} ( Y, \tilde q_1 ) \right\}
 + (1-\pi^*_0)\lambda_{0,L_1}
\mathcal E_{\overline{L}_1} \left\{ \eta_{1+} (Y , \tilde q_1 ) \right\}\\
& = (\pi^*_0 + (1-\pi^*_0) \lambda_{0,L_1} ) \mathcal E_{\overline{L}_1} \left\{  \eta_{1+} (Y , \tilde q_1 ) \right\} .
\end{align*}
The final step in the proof of Second upper bound can be modified as follows:
\begin{align*}
& \quad  \max_{\lambda_{0,L_1}} \max_{\lambda_1}\; \mu^{1,1}_{\mytext{ICE}} (\lambda_{0,L_1} ,\lambda_1)  \\
& =  \max_{\lambda_{0,L_1}}  \min_{\tilde q_1} \; E \Big( E\Big[  \;\seven\; | A_0=1, L_0 \Big] \Big) \\
& = \min_{\tilde q_1}  \max_{\lambda_{0,L_1}}  \; E \Big( E\Big[  \;\seven\; | A_0=1, L_0 \Big] \Big) \\
& = \min_{\tilde q_1} \min_{q_{0,L_1}} \; E \Big( E\Big[  \eta_{0+,\Lambda_{0,L_1}} ( \mathcal E_{\overline{L}_1} \left\{ \eta_{1+} ( Y, \tilde q_1 ) \right\},
 q_{0,L_1} )  | A_0=1, L_0 \Big] \Big) .
\end{align*}
The second line follows because the use of the convexity of $\rho_{\tau_1}$ to upper bound $\rho_{\tau_1}(\two,q_1)$ is no longer needed.
The third line follows from Sion's minimax theorem, because $\seven$ is linear in $\lambda_{0,L_1}$ and convex in $\tilde q_1$.
This completes the proof of result (ii).

\subsection{Direct proof of inequality (\ref{eq:comparison-Jnt-Prod2})}  \label{sec:prf-comparison-Jnt-Prod}

We give a direct proof of inequality (\ref{eq:comparison-Jnt-Prod2}),
\begin{align*}
\min_{ q_0,q_1 }\;
 E \{ \mathcal E_{L_0} (  \mathcal E_{\overline{L}_1} [ \eta_{1+} \{ \eta_{0+} ( Y,  q_0 ), q_1\} ] ) \}
& \ge \min_{q_{0,L_1},q_1}\;
 E \{ \mathcal E_{L_0} ( \eta_{0+} [ \mathcal E_{\overline{L}_1} \{\eta_{1+} (Y, q_1)\}, q_{0,L_1} ] ) \},
\end{align*}
without viewing the two sides of the inequality as sensitivity bounds.
In fact, we establish the following stronger result (Lemma~\ref{lem:bound-ineq}). The functionals $\eta_{0+}$ and $\eta_{1+}$ are defined as
\begin{align*}
& \eta_{0+} (y, q_0) = y + c_0 \varrho_{\tau_0} (y - q_0), \\
& \eta_{1+} (y, q_1) = y + c_1 \varrho_{\tau_1} (y-q_1),
\end{align*}
where $c_0, c_1 \ge 0$, $\tau_0,\tau_1 \in (0,1)$,
and $\varrho_\tau ( u ) = \tau u^+ + (1-\tau)(-u)^+$.
It is allowed that $\tau_0, \tau_1 \in (0, 1/2]$, although the definition $\tau_j = \Lambda_j/(1+\Lambda_j)$
satisfies $\tau_j \in (1/2,1)$ for $\Lambda_j >1$.
By the notation of the main paper, $\rho_\tau (y,q) = \varrho_\tau ( y-q)$.
[Note that $\varrho_\tau(\cdot)$ should be distinguished from $\varrho_1(\cdot)$ in (\ref{eq:constr-2period-b0}).]
For $j=0,1$, $\varrho_{\tau_j} (\cdot)$ can be easily shown to satisfy the triangle inequality: for any $u_1,u_2 \in \bbR$,
\begin{align*}
 \varrho_{\tau_j} ( u_1 + u_2 ) \le \varrho_{\tau_j} ( u_1 ) + \varrho_{\tau_j} (u_2 ).
\end{align*}
Inequality (\ref{eq:comparison-Jnt-Prod2}) can be deduced by applying Lemma~\ref{lem:bound-ineq} to
$(L_1,Y)$ under $\dif P_{L_1} (l_1 |A_0=1,L_0) \times \dif P_{Y}(y | A_0=1,A_1=1, L_0, L_1= l_1)$,
conditionally on $L_0$, and then taking expectations over $L_0$. The range assumptions $c_0(1-\tau_0) \in [0,1]$ and
$c_1 (1-\tau_1) \in [0,1]$ are satisfied due to (\ref{eq:eta-convex-coef}) in the proof of Lemma~\ref{lem:eta-convex}.
It is allowed that $c_0(1-\tau_0)= 1$ or $c_1(1-\tau_1) = 1$ in Lemma~\ref{lem:bound-ineq},
although this is ruled out for any finite $\Lambda_0$ or $\Lambda_1$ by (\ref{eq:eta-convex-coef}).

\begin{lem} \label{lem:bound-ineq}
Let $(X_1,Y)$ be random variables under a joint distribution. Suppose that $c_0(1-\tau_0) \in [0,1]$ and
$c_1(X_1) (1-\tau_1) \in [0,1]$ almost surely. Then
\begin{align}
\min_{ q_0,q_1 } \;
 \mathcal E_0 ( \mathcal E_1 [ \eta_{1+} \{ \eta_{0+} ( Y,  q_0 ), q_1\} ] )
& \ge \min_{q_0,q_1}\;
 \mathcal E_0 ( \eta_{0+} [ \mathcal E_1 \{\eta_{1+} (Y, q_1)\}, q_0 ]  ) , \label{eq:bound-ineq}
\end{align}
where the minimization is over any real number $q_0$ and a real-valued function $q_1 = q_1(X_1)$,
and $\eta_{1+}(\cdot)$ is applied with $c_1 = c_1(X_1)$. Throughout,
$\mathcal E_1(\cdot)$ denotes conditional expectation of $Y$ given $X_1$,
and $\mathcal E_0 (\cdot)$ denotes the (unconditional) expectation over $X_1$.
\end{lem}

\begin{prf}
By resetting $Y-q_0$ to $Y$ and $q_1(X_1) - q_0$ to $q_1(X_1)$, it suffices to show that for any $x_1$, the conditional version of (\ref{eq:bound-ineq}) holds given $X_1=x_1$:
\begin{align*}
& \quad \min_{ q_1 } \;
 \mathcal E_1 \big\{ Y + c_0 \varrho_{\tau_0} (Y) + c_1 \varrho_{\tau_1} ( Y - q_1 + c_0 \varrho_{\tau_0} ( Y) ) \big\}  \\
& \ge \min_{ q_1}\;
 \big\{ \mathcal E_1 Y + c_1 \mathcal E_1 \varrho_{\tau_1} ( Y -q_1 ) + c_0 \varrho_{\tau_0} ( \mathcal E_1 Y + c_1\mathcal E_1 \varrho_{\tau_1} (Y- q_1 )) \big\} ,
\end{align*}
or equivalently, with $\mathcal E_1 Y$ removed from both sides,
\begin{align}
\begin{split}
& \quad \min_{ q_1 } \;
 \mathcal E_1 \big\{ c_0 \varrho_{\tau_0} (Y) + c_1 \varrho_{\tau_1} ( Y - q_1 + c_0 \varrho_{\tau_0} ( Y) ) \big\}  \\
& \ge \min_{ q_1}\;
 \big\{ c_1 \mathcal E_1 \varrho_{\tau_1} ( Y -q_1 ) + c_0 \varrho_{\tau_0} ( \mathcal E_1 Y + c_1\mathcal E_1 \varrho_{\tau_1} (Y- q_1 ) ) \big\} ,
\end{split}  \label{eq:bound-ineq-prf1}
\end{align}
where $q_1= q_1 (x_1)$, $c_1 = c_1(x_1)$, and $\mathcal E_1(\cdot)$ denotes the conditional expectation given $X_1=x_1$.
We proceed in five steps. In the following, the dependency of $q_1$, $c_1$, and $\mathcal E_1(\cdot)$ on $x_1$ is suppressed.
Then (\ref{eq:bound-ineq-prf1}) involves the distribution of only a single random variable $Y$.

First, we show that for any $y \in \bbR$ and $q_1 \in \bbR$,
\begin{align}
\begin{split}
&  \quad \underbrace{ c_0 \varrho_{\tau_0} (y) + c_1 \varrho_{\tau_1} ( y - q_1 - c_0\tau_0 q_1 + c_0 \varrho_{\tau_0} ( y) ) }_{\one}\\
& \ge \underbrace{ c_1 \varrho_{\tau_1} ( y -q_1 ) + c_0 \varrho_{\tau_0} (  y + c_1 \varrho_{\tau_1} (y- q_1 ) ) }_{\two},
\end{split} \label{eq:bound-ineq-prf2}
\end{align}
in any of the following cases: (i) $y \ge 0$, (ii) $y<0, q_1\ge 0$,
(iii) $q_1 \le y <0$, and (iv) $ y \le q_1 <0$ and $c_0(1-\tau_0)=1$.
Note that $\one$ can be obtained from the left-hand side of (\ref{eq:bound-ineq-prf1}) with $q_1$ reset to $(1+c_0\tau_0) q_1$.

\begin{itemize}
\item[(i)] If $y \ge 0$, then $\varrho_{\tau_0} ( y)  = \tau_0 y$
and $\varrho_{\tau_0} (  y + c_1 \varrho_{\tau_1} (y- q_1 ) )=
\tau_0 (  y + c_1 \varrho_{\tau_1} (y- q_1 ) )$, and hence by direct calculation
\begin{align*}
& \one = \two= c_0\tau_0 y + c_1 ( 1 + c_0 \tau_0 ) \varrho_{\tau_1} ( y - q_1).
\end{align*}

\item[(ii)] If $y<0, q_1\ge 0$, then
$  y - q_1 - c_0\tau_0 q_1 + c_0 \varrho_{\tau_0} ( y) = (1- c_0 (1-\tau_0)) y -(1+c_0\tau_0)q_1 \le 0$ with $1- c_0 (1-\tau_0) \ge 0$, and hence
\begin{align*}
 &  \quad \one - \two = - c_0 (1-\tau_0) y +
 c_1(1-\tau_1) ( - (1- c_0 (1-\tau_0)) y  + (1+c_0\tau_0)q_1 )  \\
 & \qquad \qquad - \big\{ c_1(1-\tau_1) (-y + q_1) + c_0 \varrho_{\tau_0} (  y + c_1 (1-\tau_1) (-y+ q_1 ) ) \big\} \\
 & = - c_0 (1-\tau_0) y +   c_1(1-\tau_1) ( c_0(1-\tau_0) y + c_0 \tau_0 q_1 )-  c_0 \varrho_{\tau_0} (  y + c_1 (1-\tau_1) (-y+ q_1 ) ) \\
 & \ge 0 .
\end{align*}
The last step uses the triangle inequality
$\varrho_{\tau_0} (  y + c_1 (1-\tau_1) (-y+ q_1 ) )
\le \varrho_{\tau_0} (  (1-c_1 (1-\tau_1)) y ) + \varrho_{\tau_0} ( c_1(1-\tau_1) q_1 )
= (1-\tau_0) (1-c_1 (1-\tau_1)) (-y) + \tau_0 c_1(1-\tau_1) q_1 $,
with $1- c_1 (1-\tau_1) \ge 0$.

\item[(iii)] If $q_1 \le y <0$, then
$  y - q_1 - c_0\tau_0 q_1 + c_0 \varrho_{\tau_0} ( y) =  y - q_1 - c_0\tau_0 q_1 + c_0 (1-\tau_0) (- y) \ge 0$, and hence
\begin{align*}
 &  \quad \one - \two = - c_0 (1-\tau_0) y +
 c_1 \tau_1 ( y - q_1 - c_0\tau_0 q_1 - c_0 (1-\tau_0) y )  \\
 & \qquad \qquad - \big\{ c_1 \tau_1 (y - q_1) + c_0 \varrho_{\tau_0} (  y + c_1 \tau_1 (y- q_1 ) ) \big\} \\
 & = - c_0 (1-\tau_0) y +   c_1 \tau_1 ( - c_0(1-\tau_0) y - c_0 \tau_0 q_1 )-  c_0 \varrho_{\tau_0} (  y + c_1 \tau_1 (y- q_1 ) ) \\
 & \ge 0.
\end{align*}
The last step uses the triangle inequality
$ \varrho_{\tau_0} (  y + c_1 \tau_1 (y- q_1 ) )
\le \varrho_{\tau_0} (  (1+c_1 (1-\tau_1)) y ) + \varrho_{\tau_1} ( -c_1\tau_1 q_1 )
= (1-\tau_0) (1+c_1 \tau_1 ) (-y) + \tau_0 c_1 \tau_1  (-q_1) $.

\item[(iv)] If $ y \le q_1 <0$, then
 $ y + c_1 \varrho_{\tau_1} (y- q_1 ) = y + c_1(1-\tau_1) (-y + q_1) = (1- c_1 (1-\tau_1) ) y + c_1(1-\tau_1)q_1 \le 0$
 with $1- c_1 (1-\tau_1) \ge 0$, and hence
\begin{align*}
 &  \quad \one - \two = - c_0 (1-\tau_0) y +
 c_1 \varrho_{\tau_1} ( y - q_1 - c_0\tau_0 q_1 - c_0 (1-\tau_0) y )  \\
 & \qquad \qquad - \big\{ c_1(1-\tau_1) (-y +q_1) + c_0  (1-\tau_0) (- y + c_1(1-\tau_1) (y-q_1) ) \big\} \\
 & = c_1 \varrho_{\tau_1} ( y - q_1 - c_0\tau_0 q_1 - c_0 (1-\tau_0) y ) - c_1 (1-\tau_1) (1- c_0(1-\tau_0)) (-y + q_1) \\
 & = c_1 \varrho_{\tau_1} ( (1- c_0(1-\tau_0)) (y - q_1) - c_0 q_1 ) - c_1 (1-\tau_1) (1- c_0(1-\tau_0)) (-y + q_1) ,
\end{align*}
which may be negative.
If further $c_0(1-\tau_0)=1$, then $\one-\two = c_1  \varrho_{\tau_1} ( - c_0 q_1 ) \ge 0$.
\end{itemize}

Second, we show that if $q_1 \ge 0$ or if $q_1 <0$ and $c_0(1-\tau_0)=1$, then
\begin{align}
\begin{split}
&  \quad \mathcal E_1 \big\{ c_0 \varrho_{\tau_0} (Y) + c_1 \varrho_{\tau_1} ( Y - q_1 - c_0\tau_0 q_1 + c_0 \varrho_{\tau_0} ( Y) ) \big\} \\
& \ge c_1 \mathcal E_1 \varrho_{\tau_1} (   Y -q_1 ) + c_0 \varrho_{\tau_0} ( \mathcal E_1 Y + c_1 \mathcal E_1 \varrho_{\tau_1} (Y- q_1 ) ) .
\end{split} \label{eq:bound-ineq-prf3}
\end{align}
In fact, combining cases (i) and (ii) above indicates that  (\ref{eq:bound-ineq-prf2}) holds for any $q_1\ge 0$ and $y\in\bbR$ and hence
(\ref{eq:bound-ineq-prf3}) holds for $q_1 \ge 0$:
\begin{align}
&  \quad \mathcal E_1 \big\{ c_0 \varrho_{\tau_0} (Y) + c_1 \varrho_{\tau_1} ( Y - q_1 - c_0\tau_0 q_1 + c_0 \varrho_{\tau_0} ( Y) ) \big\} \nonumber \\
& \ge \mathcal E_1 \big \{ c_1 \varrho_{\tau_1} ( Y -q_1 ) + c_0 \varrho_{\tau_0} (  Y + c_1 \varrho_{\tau_1} (Y- q_1 ) ) \big\} \label{eq:bound-ineq-prf4} \\
& \ge c_1  \mathcal E_1 \varrho_{\tau_1} ( Y -q_1 ) + c_0 \varrho_{\tau_0} ( \mathcal E_1 Y + c_1 \mathcal E_1 \varrho_{\tau_1} (Y- q_1 ) ) , \nonumber
\end{align}
where the last step follows from Jensen's inequality by the convexity of $\varrho_{\tau_0}(\cdot)$.
Similarly, combining cases (i), (iii), and (iv) above indicates that
(\ref{eq:bound-ineq-prf3}) holds if $q_1 < 0$ and $c_0(1-\tau_0)=1$.

The following remark is helpful to motivating the subsequent steps of our proof.
If (\ref{eq:bound-ineq-prf2}) were also valid for any $q_1 <0$ and $y\in\bbR$ with $c_0(1-\tau_0)<1$,
then (\ref{eq:bound-ineq-prf2}) would hold for any $q_1\in\bbR$ and $y\in\bbR$, the preceding display would hold for any $q_1\in\bbR$,
and hence minimizing over $q_1 \in \bbR$ would imply the desired result (\ref{eq:bound-ineq-prf1}), because
\begin{align*}
&  \quad \min_{ q_1 } \;
 \mathcal E_1 \big\{ c_0 \varrho_{\tau_0} (Y) + c_1 \varrho_{\tau_1} ( Y - q_1- c_0\tau_0 q_1 + c_0 \varrho_{\tau_0} ( Y) ) \big\}  \\
& = \min_{ q_1 } \;
 \mathcal E_1 \big\{ c_0 \varrho_{\tau_0} (Y) + c_1 \varrho_{\tau_1} ( Y - q_1 + c_0 \varrho_{\tau_0} ( Y) ) \big\} .
\end{align*}
Unfortunately, this direct reasoning cannot be realized. There are examples where
 inequality (\ref{eq:bound-ineq-prf2}) is violated for some $q_1 <0$ and $y\in \bbR$, and
  (\ref{eq:bound-ineq-prf3}) is violated for some $q_1 <0$
(although satisfied for $q_1 \ge 0$). A more careful reasoning is needed than requiring (\ref{eq:bound-ineq-prf3}) for $q_1\in\bbR$.

Third, we show that (\ref{eq:bound-ineq-prf3}) holds if $q_1 < 0$ and $\mathcal E_1 Y \ge 0$, through a different reasoning than in (\ref{eq:bound-ineq-prf4}).
In fact, if $\mathcal E_1 Y \ge 0$, then
the right-hand side of (\ref{eq:bound-ineq-prf3}) can be rewritten as
\begin{align*}
& \quad c_1  \mathcal E_1  \varrho_{\tau_1} (Y -q_1 ) + c_0 \varrho_{\tau_0} ( \mathcal E_1 Y + c_1 \mathcal E_1 \varrho_{\tau_1} (Y- q_1 ) ) \\
& = c_1 \mathcal E_1  \varrho_{\tau_1} (Y -q_1 ) + c_0 \tau_0 ( \mathcal E_1 Y + c_1 \mathcal E_1 \varrho_{\tau_1} (Y- q_1 ) ) \\
& = \mathcal E_1 \big\{ c_1 \varrho_{\tau_1} (Y -q_1 ) + c_0 \tau_0 ( Y + c_1 \varrho_{\tau_1} (Y- q_1 ) ) \big\} ,
\end{align*}
and  (\ref{eq:bound-ineq-prf3})  is equivalent to
\begin{align}
\begin{split}
&  \quad \mathcal E_1 \big\{ c_0 \varrho_{\tau_0} (Y) + c_1 \varrho_{\tau_1} ( Y - q_1 - c_0\tau_0 q_1 + c_0 \varrho_{\tau_0} ( Y) ) \big\} \\
& \ge \mathcal E_1 \big\{ c_1 \varrho_{\tau_1} (Y -q_1 ) + c_0 \tau_0 ( Y + c_1 \varrho_{\tau_1} (Y- q_1 ) ) \big\} .
\end{split} \label{eq:bound-ineq-prf5}
\end{align}
Hence it suffices to show that for any $y \in \bbR$ and $q_1 <0 $,
\begin{align}
\begin{split}
&  \quad \underbrace{ c_0 \varrho_{\tau_0} (y) + c_1 \varrho_{\tau_1} ( y - q_1 - c_0\tau_0 q_1 + c_0 \varrho_{\tau_0} ( y) ) }_{\one} \\
& \ge \underbrace{ c_1 \varrho_{\tau_1} (y -q_1 ) + c_0 \tau_0 ( y + c_1 \varrho_{\tau_1} (y- q_1 ) ) }_{\three} .
\end{split}  \label{eq:bound-ineq-prf6}
\end{align}
If $y \ge 0$, then (\ref{eq:bound-ineq-prf6}) is the same as (\ref{eq:bound-ineq-prf2}) in case (i).
If $q_1 \le y <0$, then (\ref{eq:bound-ineq-prf6}) is implied by (\ref{eq:bound-ineq-prf2}) in case (iii),
because $ \varrho_{\tau_0} ( y + c_1 \varrho_{\tau_1} (y- q_1 ) ) \ge \tau_0 ( y + c_1 \varrho_{\tau_1} (y- q_1 ) )$.
If $y \le q_1 <0$, then by the calculation of $\one-\two$ in case (iv) above,
\begin{align*}
 &  \quad \one - \three = \one-\two - c_0 ( y + c_1(1-\tau_1) (-y+q_1) )  \\
 & = c_1 \varrho_{\tau_1} ( (1- c_0(1-\tau_0)) (y - q_1) - c_0 q_1 ) - c_1 (1-\tau_1) (1- c_0(1-\tau_0)) (-y + q_1) \\
 & \qquad \qquad - c_0 ( y + c_1(1-\tau_1) (-y+q_1) ) \\
 & = c_1 \varrho_{\tau_1} ( (1- c_0(1-\tau_0)) (y - q_1) - c_0 q_1 ) - c_1 (1-\tau_1)\big\{ (1- c_0(1-\tau_0)) (-y + q_1) +c_0q_1 \big\} \\
 & \qquad \qquad - c_0 (1- c_1(1-\tau_1)) y \\
 &  \ge - c_0 (1- c_1(1-\tau_1)) y  \ge 0,
\end{align*}
where the second last step uses $\varrho_{\tau_1} (u) \ge -(1-\tau_1) u$ for $u = (1- c_0(1-\tau_0)) (y - q_1) - c_0 q_1$.
Combining the three cases shows that (\ref{eq:bound-ineq-prf6}) holds for $y \in \bbR$ and $q_1 <0 $.

Fourth, we show that (\ref{eq:bound-ineq-prf3}) holds if $q_1 < 0$,
$\mathcal E_1 Y < 0$, but $ | \mathcal E_1 Y| \le c_1\mathcal E_1 \varrho_{\tau_1} (Y-q_1) $.
In this case, (\ref{eq:bound-ineq-prf3}) remains the same as (\ref{eq:bound-ineq-prf5}),
and (\ref{eq:bound-ineq-prf5}) follows from (\ref{eq:bound-ineq-prf6}) for any $y \in \bbR$ and $q_1 <0 $.

The implications of the second, third, and four steps can be seen as follows. Let
\begin{align*}
\check q_1 = \argmin_{q_1}\;
\mathcal E_1 \big\{ c_0 \varrho_{\tau_0} (Y) + c_1 \varrho_{\tau_1} ( Y - q_1 - c_0\tau_0 q_1 + c_0 \varrho_{\tau_0} ( Y) ) \big\} .
\end{align*}
Applying (\ref{eq:bound-ineq-prf3}) with $q_1$ set to $\check q_1$ shows that
if $\check q_1\ge 0$, or if $\check q_1 <0$ and $c_0(1-\tau_0)=1$, or if $\check q_1 < 0$ and $\mathcal E_1 Y \ge 0$, or if
 $\check q_1 < 0$, $\mathcal E_1 Y < 0$, but $ | \mathcal E_1 Y| \le c_1\mathcal E_1 \varrho_{\tau_1} (Y-\check q_1) $, then
\begin{align}
\begin{split}
&  \quad \mathcal E_1 \big\{ c_0 \varrho_{\tau_0} (Y) + c_1 \varrho_{\tau_1} ( Y - \check q_1 - c_0\tau_0 \check q_1 + c_0 \varrho_{\tau_0} ( Y) ) \big\} \\
& \ge c_1 \mathcal E_1 \varrho_{\tau_1} (   Y - \check q_1 ) + c_0 \varrho_{\tau_0} ( \mathcal E_1 Y + c_1 \mathcal E_1 \varrho_{\tau_1} (Y- \check q_1 ) ) ,
\end{split} \label{eq:bound-ineq-prf7}
\end{align}
which implies the desired result (\ref{eq:bound-ineq-prf1}).

Finally, to complete the proof of (\ref{eq:bound-ineq-prf1}),  we show that given $c_0(1-\tau_0) \in [0, 1)$,
if $\check q_1 < 0$,
$\mathcal E_1 Y < 0$, but $ | \mathcal E_1 Y| > c_1\mathcal E_1 \varrho_{\tau_1} (Y-\check q_1) $, then
\begin{align}
\begin{split}
&  \quad \mathcal E_1 \big\{ c_0 \varrho_{\tau_0} (Y) + c_1 \varrho_{\tau_1} ( Y - \check q_1 - c_0\tau_0 \check q_1 + c_0 \varrho_{\tau_0} ( Y) ) \big\} \\
& \ge c_1 \mathcal E_1 \varrho_{\tau_1} (   Y - \delta \check q_1 ) + c_0 \varrho_{\tau_0} ( \mathcal E_1 Y + c_1 \mathcal E_1 \varrho_{\tau_1} (Y- \delta \check q_1 ) ) ,
\end{split} \label{eq:bound-ineq-prf8}
\end{align}
where $ \delta= \frac{1+c_0\tau_0}{ 1- c_0(1-\tau_0)} >0$.
In fact, by definition, $(1+c_0 \tau_0) \check q_1$ is a $\tau$-quantile of $Y + c_0 \varrho_{\tau_0}(Y)$,
which is an increasing function of $Y$ by Lemma~\ref{lem:eta-convex}.
If $\check q_1<0$, then $\delta \check q_1$ is a $\tau$-quantile of $Y$, satisfying
$ \delta \check q_1 + c_0 \varrho_{\tau_0} ( \delta \check q_1)  = (1- c_0(1-\tau_0)) (\delta \check q_1) = (1+c_0 \tau_0) \check q_1$.
Moreover, if $ | \mathcal E_1 Y| > c_1\mathcal E_1 \varrho_{\tau_1} (Y-\check q_1) $, then $\delta \check q_1$ being a $\tau$-quantile of $Y$ implies that
\begin{align*}
 | \mathcal E_1 Y| > c_1\mathcal E_1 \varrho_{\tau_1} (Y-\check q_1)
 \ge c_1\mathcal E_1 \varrho_{\tau_1} (Y- \delta \check q_1) .
\end{align*}
Hence if   $\check q_1 < 0$,
$\mathcal E_1 Y < 0$, but $ | \mathcal E_1 Y| > c_1\mathcal E_1 \varrho_{\tau_1} (Y-\check q_1) $, then
the right-hand side of
(\ref{eq:bound-ineq-prf8}) can be calculated as
\begin{align*}
&  \quad  c_1 \mathcal E_1 \varrho_{\tau_1} (   Y - \delta \check q_1 ) + c_0 \varrho_{\tau_0} ( \mathcal E_1 Y + c_1 \mathcal E_1 \varrho_{\tau_1} (Y- \delta \check q_1 ) ) \\
& = c_1 \mathcal E_1 \varrho_{\tau_1} (   Y - \delta \check q_1 ) + c_0 (1-\tau_0) ( -\mathcal E_1 Y - c_1 \mathcal E_1 \varrho_{\tau_1} (Y- \delta \check q_1 ) ) \\
& = \mathcal E_1 \big\{ c_1 \varrho_{\tau_1} (   Y - \delta \check q_1 ) + c_0 (1-\tau_0) ( -  Y - c_1  \varrho_{\tau_1} (Y- \delta \check q_1 ) ) \big\}.
\end{align*}
Hence it suffices to show that for any $y \in \bbR$ and $q_1 <0 $,
\begin{align}
\begin{split}
&  \quad \underbrace{ c_0 \varrho_{\tau_0} (y) + c_1 \varrho_{\tau_1} ( y - q_1 - c_0\tau_0 q_1 + c_0 \varrho_{\tau_0} ( y) ) }_{\one} \\
& \ge \underbrace{ c_1 \varrho_{\tau_1} (   y - \delta q_1 ) + c_0 (1-\tau_0) ( -  y - c_1  \varrho_{\tau_1} (y- \delta q_1 ) ) }_{\four} .
\end{split}  \label{eq:bound-ineq-prf9}
\end{align}
If $y \ge 0$ and $q_1 <0$, then (\ref{eq:bound-ineq-prf9}) is implied by (\ref{eq:bound-ineq-prf2}) in case (i),
because $\varrho_{\tau_1} (   y - q_1 ) \ge \varrho_{\tau_1} (   y - \delta q_1 )$ and
$ \varrho_{\tau_0} ( y + c_1 \varrho_{\tau_1} (y- q_1 ) ) \ge \varrho_{\tau_0} ( y + c_1 \varrho_{\tau_1} (y- \delta q_1 ) )  \ge
(1-\tau_0) ( -y - c_1 \varrho_{\tau_1} (y- q_1 ) )$.
If $y <0$, then  $\varrho_{\tau_0} ( y) = -(1-\tau_0)y$  and direct calculation using $y - (1+c_0 \tau_0) q_1 = y - (1-c_0(1-\tau_0)) \delta q_1 $ yields
\begin{align*}
 & \one = \four =  - c_0 (1-\tau_0) y + c_1 (1-c_0(1-\tau_0)) \varrho_{\tau_1} ( y - \delta q_1  ).
\end{align*}
Combining the two cases shows that (\ref{eq:bound-ineq-prf9}) holds for $y \in \bbR$ and $q_1 <0 $.
\end{prf}

\section{Technical details: Proofs for Section \ref{sec:Kperiod}}

\subsection{Proof of Proposition \ref{pro:upper-Kperiod}}

We apply induction to show a strengthened version of the desired result as stated in (\ref{eq:pro-upper-Kperiod-prf3}).
Applying induction in $K$ to the original result does not seem to work.
The proof involves conceptually similar steps as in the proof of Lemma~\ref{pro:upper-2period} for 2 periods,
but a crucial property is the convexity of $\eta^k_+ (y, \overq_k)$ in $\overq_k$ as shown in Lemma~\ref{lem:eta-convex}.

By an extension of Lemma~\ref{lem:density-ratio}, for nonnegative functions $(\lambda_0(L_0,y), \ldots, \lambda_{K-1}(\overL_{K-1},y))$
satisfying (\ref{eq:constr-Kperiod-b}),
there exists a probability distribution $Q$ for the full data satisfying the following two properties in parallel to those in Lemma~\ref{lem:density-ratio}: \vspace{-.05in}
\begin{itemize}\addtolength{\itemsep}{-.05in}
\item[(i)] the induced distribution of $Q$ on the observed data $(\overline{A}_{K-1},\overline{L}_{K-1},Y )$ under Assumption~A1  coincides with
the true distribution of  $(\overline{A}_{K-1},\overline{L}_{K-1},Y )$,
\item[(ii)] $ \lambda^*_{k,Q} (\overline{L}_k, y) =\lambda_k (\overline{L}_k, y) $ for $k=0,1,\ldots,K-1$.
\end{itemize} \vspace{-.05in}
where $\lambda^*_{k,Q}$ is defined as $\lambda^*_k$ with $P$ replaced by $Q$.
In other words, $Q$ can be a valid choice for $P$.
By the law of iterated expectations under $Q$ similarly as in (\ref{eq:iden-lam-prf5}),
it can be shown that for $k=1,\ldots, K-1$ and any function $b_{k-1} (\overL_{k-1}, y)$,
\begin{align}
& \quad E_{\underlam_k} \left\{ b_{k-1} (\overL_{k-1}, Y^{\over1} ) | \overA_{k-1} = \over1_{k-1}, \overL_{k-1} \right\} \nonumber \\
& =  \calE_{\overL_{k-1}} E_{\underlam_{k+1}} \left\{  (\pi^*_k + (1-\pi^*_k) \lambda_k (\overL_k, Y^{\over1}))
b_{k-1} (\overL_{k-1}, Y^{\over1} ) | \overA_k = \over1_k, \overL_k \right\} , \label{eq:pro-upper-Kperiod-prf2}
\end{align}
where $E_{\underlam_k } (\cdot | \overA_{k-1} = \over1_{k-1}, \overL_{k-1})$ denotes the conditional expectation  under $Q$, depending on
$\underlam_k $, and $\underlam_K$ is set to the null and $E_{\underlam_K} (\cdot |  \overA_{K-1}= \over1_{K-1}, \overL_{K-1} ) $ reduces to
$E(\cdot |  \overA_{K-1}= \over1_{K-1}, \overL_{K-1})$.
Then repeated application of (\ref{eq:pro-upper-Kperiod-prf2}) yields
\begin{align}
& \quad E_{\underlam_k} \left\{ b_{k-1} (\overL_{k-1}, Y^{\over1} ) | \overA_{k-1} = \over1_{k-1}, \overL_{k-1} \right\} \nonumber \\
& = \Big( \prod_{j=k-1}^ {K-1} \calE_{\overL_j} \Big) \left\{ \varrho^k ( \underL_k, Y; \underlam_k ) b_{k-1} (\overL_{k-1}, Y^{\over1}) \right\} .
\label{eq:pro-upper-Kperiod-prf1}
\end{align}
For $k=0,1,\ldots,K-1$, we show by induction that for any fixed $\underlam_{k+1}$,
\begin{align}
& \max_{\overlam_k} \;  \mu^{\over1}_{\mytext{ICE}} (\overlam_{K-1} )
= \min_{\overq_k}\;  E \Big[  \Big( \prod_{j=0}^{k-1} \calE_{\overL_j} \Big)
E_{\underlam_{k+1} } \left\{ \eta^k_+ (Y^{\over1}, \overq_k ) | \overA_k= \over1_k, \overL_k \right\} \Big] .\label{eq:pro-upper-Kperiod-prf3}
\end{align}
The case of $k=K-1$ yields the desired result.

First, we show that (\ref{eq:pro-upper-Kperiod-prf3}) holds for $k=0$.
By (\ref{eq:pro-upper-Kperiod-prf1}) with $k=1$ and $b_0(L_0,y)=(\pi^*_0 + (1-\pi^*_0) \lambda_0 (L_0,y) )y$,
the generalized ICE formula in Lemma~\ref{lem:iden-mu-Kperiod}
can be written as
\begin{align*}
\mu^{\over1}_{\mytext{ICE}} (\overlam_{K-1} )
= E \left[ E_{\underlam_1} \left\{ (\pi^*_0 + (1-\pi^*_0) \lambda_0 (L_0,Y^{\over1})) Y^{\over1} | A_0=1, L_0\right\} \right].
\end{align*}
For any fixed $\underlam_1$, the normalization constraint (\ref{eq:constr-Kperiod-b}) on $\lambda_0$ is equivalent to
\begin{align*}
E_{\underlam_1} \left( \lambda_0 (L_0,Y^{\over1})  | A_0=1, L_0 \right) \equiv 1.
\end{align*}
Applying Lemma~\ref{lem:upper-1period} with $Y$ replaced by  $Y^{\over1}$ yields
\begin{align*}
& \max_{\lambda_0} \;  \mu^{\over1}_{\mytext{ICE}} (\overlam_{K-1} )
= \min_{q_0}\;  E \left[ E_{\underlam_1} \left\{ \eta_{0+}(Y^{\over1}, q_0) | A_0=1, L_0 \right\} \right] ,
\end{align*}
where $\eta_{0+}(y, q_0) = y + (1-\pi^*_0) (\Lambda_0-\Lambda_0^{-1}) \rho_{\tau_0} (y, q_0) $.
Hence (\ref{eq:pro-upper-Kperiod-prf3}) holds for $k=0$.

Suppose that for some $1\le k \le K-1$, equation (\ref{eq:pro-upper-Kperiod-prf3}) holds with $k$ replaced by $k-1$ for any fixed $\underlam_k$:
\begin{align}
& \max_{\overlam_{k-1} } \;  \mu^{\over1}_{\mytext{ICE}} (\overlam_{K-1} )
= \min_{\overq_{k-1} }\;  E \Big[  \Big( \prod_{j=0}^{k-2} \calE_{\overL_j} \Big)
E_{\underlam_k } \left\{ \eta^{k-1}_+ (Y^{\over1}, \overq_{k-1} ) | \overA_{k-1} = \over1_{k-1}, \overL_{k-1} \right\} \Big] .\label{eq:pro-upper-Kperiod-prf4}
\end{align}
It is desired to show that
(\ref{eq:pro-upper-Kperiod-prf3}) holds. For any fixed $\underlam_{k+1}$,
the normalization constraint (\ref{eq:constr-Kperiod-b}) on $\lambda_k$ is equivalent to
\begin{align*}
E_{\underlam_{k+1} } \left( \lambda_k (L_k ,Y^{\over1})  | \overA_k= \over1_k, \overL_k \right) \equiv 1.
\end{align*}
We apply single-period identity (\ref{eq:upper-1period-b}), with $L_0$, $A_0$, $\lambda_0$, and $b(L_0,Y)$ replaced by
$L_k$, $A_k$, $\lambda_k$, and $\eta^{k-1}_+ (Y^{\over1},$ $ \overq_{k-1} ) $ respectively,
conditionally on $( \overA_{k-1}= \over1_{k-1}, \overL_{k-1} )$, and obtain
\begin{align}
& \quad \max_{\lambda_k} \;
E_{\underlam_k } \left\{ \eta^{k-1}_+ (Y^{\over1}, \overq_{k-1} ) | \overA_{k-1} = \over1_{k-1}, \overL_{k-1} \right\} \nonumber \\
& = \max_{\lambda_k} \; \calE_{\overL_{k-1}}
E_{\underlam_{k+1} } \left\{ (\pi^*_k + (1-\pi^*_k) \lambda_k (\overL_k, Y^{\over1}))
\eta^{k-1}_+ (Y^{\over1}, \overq_{k-1} ) | \overA_k = \over1_k, \overL_k \right\} \nonumber \\
& = \max_{q_k} \; \calE_{\overL_{k-1}}
E_{\underlam_{k+1} } \left\{ \eta_{k+}(
\eta^{k-1}_+ (Y^{\over1}, \overq_{k-1} ), q_k) | \overA_k = \over1_k, \overL_k  \right\} \nonumber \\
& = \max_{q_k} \; \calE_{\overL_{k-1}}
E_{\underlam_{k+1} } \left\{
\eta^k_+ (Y^{\over1}, \overq_k ) | \overA_k = \over1_k, \overL_k  \right\} ,  \label{eq:pro-upper-Kperiod-prf5}
\end{align}
where the first line follows from (\ref{eq:pro-upper-Kperiod-prf2}). Combining (\ref{eq:pro-upper-Kperiod-prf4}) and (\ref{eq:pro-upper-Kperiod-prf5}) yields
\begin{align*}
& \quad \max_{\lambda_k} \max_{\overlam_{k-1} } \;  \mu^{\over1}_{\mytext{ICE}} (\overlam_{K-1} )  \\
& = \max_{\lambda_k} \min_{\overq_{k-1} }\;  E \Big[  \Big( \prod_{j=0}^{k-2} \calE_{\overL_j} \Big)
E_{\underlam_k } \left\{ \eta^{k-1}_+ (Y^{\over1}, \overq_{k-1} ) | \overA_{k-1} = \over1_{k-1}, \overL_{k-1} \right\} \Big]  \\
& = \min_{\overq_{k-1} } \max_{\lambda_k} \;  E \Big[  \Big( \prod_{j=0}^{k-2} \calE_{\overL_j} \Big)
E_{\underlam_k } \left\{ \eta^{k-1}_+ (Y^{\over1}, \overq_{k-1} ) | \overA_{k-1} = \over1_{k-1}, \overL_{k-1} \right\} \Big]  \\
& = \min_{\overq_{k-1} }  \min_{ q_k} \;  E \Big[  \Big( \prod_{j=0}^{k-1} \calE_{\overL_j} \Big)
E_{\underlam_{k+1} } \left\{ \eta^k_+ (Y^{\over1}, \overq_k ) | \overA_k= \over1_k, \overL_k \right\} \Big] .
\end{align*}
The third line above follows from Sion's minimax theorem, because
the functional involved,
\begin{align*}
& \quad E \Big[  \Big( \prod_{j=0}^{k-2} \calE_{\overL_j} \Big)
E_{\underlam_k } \left\{ \eta^{k-1}_+ (Y^{\over1}, \overq_{k-1} ) | \overA_{k-1} = \over1_{k-1}, \overL_{k-1} \right\} \Big]  \\
& = E \Big[  \Big( \prod_{j=0}^{k-1} \calE_{\overL_j} \Big)
E_{\underlam_{k+1} } \left\{ (\pi^*_k + (1-\pi^*_k) \lambda_k (\overL_k, Y^{\over1}))
\eta^{k-1}_+ (Y^{\over1}, \overq_{k-1} ) | \overA_k = \over1_k, \overL_k \right\}   \Big]
\end{align*}
is linear in $\lambda_k$ and, by Lemma~\ref{lem:eta-convex}, convex in $\overq_{k-1}$.
This completes the proof of (\ref{eq:pro-upper-Kperiod-prf3}).

\begin{lem} \label{lem:eta-convex}
For any $k=0,1,\ldots, K-1$, the function $\eta^k_+ (y, \overq_k)$ defined in (\ref{eq:eta-def12}) is convex in $\overq_k$
for any fixed $y$ and increasing and convex in $y$ for any fixed $\overq_k$.
\end{lem}

\begin{prf}
We first demonstrate an auxiliary result. By definition,
\begin{align*}
& \quad \eta^k_+ (y, \overq_k) = \eta^{k-1}_+ + c_k \rho_{\tau_k} ( \eta^{k-1}_+ , q_k) \\
& =  \eta^{k-1}_+ + c_k \left\{ \tau_k ( \eta^{k-1}_+ - q_k)^+ + (1-\tau_k) ( q_k - \eta^{k-1}_+)^+ \right\},
\end{align*}
where $\eta^{k-1}_+ = y$ if $k=0$ or $\eta^{k-1}_+ (y, \overq_{k-1})$ if $k\ge 1$ and $c_k = \pi^*_k ( \Lambda_k - \Lambda_k^{-1})$. Hence $\eta^k_+ (y, \overq_k)$
can be expressed as $\eta^k_+ (y, \overq_k) = q_k + f_k \circ g_k $, where $g_k = \eta^{k-1}_+ - q_k$ and
\begin{align*}
& \quad f_k (x ) = x + c_k \left\{ \tau_k x^+ + (1-\tau_k) (-x)^+ \right\} \\
& = \left\{ \begin{array}{lc}
 (1+ c_k \tau_k) x , & x \ge 0,\\
 (1 - c_k (1-\tau_k)) x, & x < 0.
\end{array} \right.
\end{align*}
Then $f_k(x)$ is convex in $x$. Moreover, $f_k(x)$ is increasing in $x$, because $1+ c_k \tau_k \ge 1$ and
$0 < 1 - c_k (1-\tau_k) \le 1$ by the following calculation:
\begin{align}
0 \le c_k (1- \tau_k) \le (\Lambda_k - \Lambda_k^{-1}) (1-\tau_k) = (\Lambda_k-1)/\Lambda_k  < 1 ,  \label{eq:eta-convex-coef}
\end{align}
where $ 0 \le c_k  \le \Lambda_k - \Lambda_k^{-1}$ and $\tau_k = \Lambda_k / (1+\Lambda_k)$.

The proof of Lemma~\ref{lem:eta-convex} is by induction. For $k=0$, $\eta^0_+ (y, q_0) = \eta_{0+} (y, q_0)$ is trivially convex in $q_0$ for fixed $y$
and convex in $y$ for fixed $q_0$. Moreover, $\eta^0_+ (y, q_0) = q_0 + f_0( y-q_0)$ is increasing in $y$ by
the increasingness of $f_0$ demonstrated above.

Suppose that for some $1\le k \le K-1$,  $\eta^{k-1}_+ (y, \overq_{k-1})$ is convex in $\overq_{k-1}$ for fixed $y$
and increasing and convex in $y$ for fixed $\overq_{k-1}$.
It is desired to show that $\eta^k_+ (y, \overq_k)= q_k + f_k \circ g_k$ is convex in $\overq_k$ for fixed $y$ and increasing and convex in $y$ for fixed $\overq_k$.
By assumption, $g_k = \eta^{k-1}_+ - q_k$ is convex in $\overq_k = (\overq_{k-1}, q_k)$ for fixed $y$ and
increasing and convex in $y$ for fixed $\overq_k$.
As shown above, $f_k$ is increasing and convex.
Then $\eta^k_+ (y, \overq_k)$ is increasing in $y$ by
the composition rule that $f_k \circ g_k$ is increasing if $f_k$  and $g_k$ are increasing.
The convexity of $\eta^k_+ (y, \overq_k)$ in $q_k$ or in $y$ follows
from the composition rule that $f_k \circ g_k$ is convex if $f_k$ is convex and nondecreasing and $g_k$ is convex.
\end{prf}

%

\subsection{Proof of Proposition \ref{pro:upper-Kperiod-Jnt}}

The two directions of the result, $\mu^{\over1}_{+,\mytext{Jnt}} \le \mu^{\over1}_+$ and $\mu^{\over1}_{+,\mytext{Jnt}} \ge \mu^{\over1}_+$,
can be shown by directly extending the proof of Proposition \ref{pro:upper-2period-Jnt} in Section \ref{sec:prf-pro-upper-Kperiod-Jnt}.

\subsection{Proof of Proposition \ref{pro:upper-Kperiod-Prod}}

By an extension of Lemma~\ref{lem:density-ratio-Prod}, for nonnegative functions
$\{(\lambda_{k,L_{k+1}} (\overL_k, l_{k+1} ),  \lambda_{k,Y} (\overL_k, \underl_{k+1}, y)): k=0,1,\ldots,K-1\}$
satisfying (\ref{eq:constr-Kperiod-Prod-b}),
there exists a probability distribution $Q$ for the full data satisfying the following two properties in parallel to those in Lemma~\ref{lem:density-ratio-Prod}: \vspace{-.05in}
\begin{itemize}\addtolength{\itemsep}{-.05in}
\item[(i)] the induced distribution of $Q$ on the observed data $(\overline{A}_{K-1},\overline{L}_{K-1},Y )$ under Assumption~A1$^\dag$ coincides with
the true distribution of  $(\overline{A}_{K-1},\overline{L}_{K-1},Y )$,
\item[(ii)] $ \lambda^*_{k,L_{k+1},Q} (\overL_k, l_{k+1} ) =\lambda_{k,L_{k+1}} (\overL_k, l_{k+1} )$ and
$\lambda_{k,Y,Q} (\overL_k, \underl_{k+1}, y) = \lambda_{k,Y} (\overL_k, \underl_{k+1}, y)$ for $k=0,1,\ldots,K-1$.
\end{itemize} \vspace{-.05in}
where $\lambda^*_{k,L_{k+1},Q}$ and $\lambda^*_{k,Y,Q}$ are defined as $\lambda^*_{k,L_{k+1}}$ and $\lambda^*_{k,Y}$ with $P$ replaced by $Q$.
In other words, $Q$ can be a valid choice for $P$.
The joint sensitivity ratios under $Q$ are
$\lambda_{k,\mytext{Jnt}} ( \overL_k, \underl_{k+1}, y) =\lambda_{k,L_{k+1}}  (\overL_k, $ $l_{k+1} ) \lambda_{k,Y}(\overL_k, \underl_{k+1}, y)$
for $k=0,1,\ldots,K-1$.
By the law of iterated expectations under $Q$,
it can be shown that for $k=1,\ldots, K-1$ and any function $b_{K-1} (\overL_{k-1}, \underl_k, y)$,
\begin{align}
& \quad E_{\underlam_{k,\mytext{Jnt}}} \left\{ b_{K-1} (\overL_{k-1}, \underL^{\over1}_k, Y^{\over1} ) | \overA_{k-1} = \over1_{k-1}, \overL_{k-1} \right\} \nonumber \\
& =  \calE_{\overL_{k-1}} E_{\underlam_{k+1,\mytext{Jnt}}} \left\{  (\pi^*_k + (1-\pi^*_k) \lambda_{k,\mytext{Jnt}} (\overL_k, \underL^{\over1}_{k+1}, Y^{\over1}))
b_{K-1} (\overL_k, \underL^{\over1}_{k+1} ,  Y^{\over1} ) | \overA_k = \over1_k, \overL_k \right\} , \label{eq:pro-upper-Kperiod-Prod-prf2}
\end{align}
where $E_{\underlam_{k,\mytext{Jnt}} } (\cdot | \overA_{k-1} = \over1_{k-1}, \overL_{k-1})$ denotes the conditional expectation  under $Q$, depending on
$\underlam_{k,\mytext{Jnt}} $, and $\underlam_{K,\mytext{Jnt}}$ is set to the null and $E_{\underlam_{K,\mytext{Jnt}}} (\cdot |  \overA_{K-1}= \over1_{K-1}, \overL_{K-1} ) $ reduces to
$E(\cdot |  \overA_{K-1}= \over1_{K-1}, \overL_{K-1})$.
Repeated application of (\ref{eq:pro-upper-Kperiod-Prod-prf2}) then yields
\begin{align}
& \quad E_{\underlam_{k,\mytext{Jnt}} } \left\{ b_{K-1} (\overL_{k-1},\underL^{\over1}_k, Y^{\over1} ) | \overA_{k-1} = \over1_{k-1}, \overL_{k-1} \right\} \nonumber \\
& = \Big( \prod_{j=k-1}^ {K-1} \calE_{\overL_j} \Big) \left\{ \varrho^k_{\mytext{Jnt}} ( \overL_{K-1}, Y; \underlam_{k,\mytext{Jnt}} ) b_{K-1} (\overL_{K-1}, Y ) \right\} .
\label{eq:pro-upper-Kperiod-Prod-prf1}
\end{align}
In the following, we show that $\mu^{\over1}_{+,\mytext{Prod}} \le \mu^{\over1}_{+,\mytext{Prod,v1}}$ (First upper bound)
and $\mu^{\over1}_{+,\mytext{Prod}} \le \mu^{\over1}_{+,\mytext{Prod,v2}}$ (Second upper bound).

\textbf{First upper bound.}\;
For $k=0,1,\ldots,K-1$, we show by induction that for any fixed $\underlam_{k+1,L_{k+2}}$ and $\underlam_{k+1,Y}$ (hence fixed $\underlam_{k+1,\mytext{Jnt}}$),
\begin{align}
& \quad \max_{\overlam_{k,L_{k+1}}, \overlam_{k,Y}} \;  \mu^{\over1}_{\mytext{ICE}} (\lambda_{0,L_1}\lambda_{0,Y},\ldots,\lambda_{K-1,Y}) \nonumber \\
& \le \min_{\overq_{k,L_{k+1}}, \overq_{k,Y} }\;  E \Big[  \Big( \prod_{j=0}^{k-1} \calE_{\overL_j} \Big)
E_{\underlam_{k+1,\mytext{Jnt}} } \left\{ \eta^k_{+,\mytext{Prod,v1}} (Y^{\over1}, \overq_{k,L_{k+1}}, \overq_{k,Y} ) | \overA_k= \over1_k, \overL_k \right\} \Big] ,
\label{eq:pro-upper-Kperiod-Prod-prf3}
\end{align}
where $\overq_{k,L_{k+1}} = \overq_{k,L_{k+1}} (\overL_k, L^{\over1}_{k+1} )$ and $\overq_{k,Y}  = \overq_{k,Y} (\overL_k, \underL^{\over1}_{k+1})$
inside $E_{\underlam_{k+1,\mytext{Jnt}} } (\cdot | \overA_k= \over1_k, \overL_k)$, and $\prod_{j=0}^{-1} \calE_{\overL_j}$ is set to the null.
The case of $k=K-1$ yields the desired result.

First, we show that (\ref{eq:pro-upper-Kperiod-Prod-prf3}) holds for $k=0$.
By (\ref{eq:pro-upper-Kperiod-Prod-prf1}) with $k=1$ and $b_{K-1}(L_0,\underl_1, y)=(\pi^*_0 + (1-\pi^*_0) \lambda_{0,\mytext{Jnt}} (L_0,\underl_1,y) )y$,
the generalized ICE formula in Lemma~\ref{lem:iden-mu-Kperiod-Jnt}
can be written as
\begin{align*}
& \quad \mu^{\over1}_{\mytext{ICE}} (\lambda_{0,L_1}\lambda_{0,Y},\ldots,\lambda_{K-1,Y}) \\
& = E \left[ E_{\underlam_{1,\mytext{Jnt}}} \left\{ (\pi^*_0 + (1-\pi^*_0) \lambda_{0,\mytext{Jnt}} (L_0, L^{\over1}_1, Y^{\over1})) Y^{\over1} | A_0=1, L_0\right\} \right],
\end{align*}
which matches (\ref{eq:pro-upper-2period-Prod-prf1}) in the case of 2 periods with $\lambda_1$ and $(L^{1,1}_1, Y^{1,1})$
replaced by $\underlam_{1,\mytext{Jnt}}$ and $(L^{\over1}_1, Y^{\over1})$. Then applying inequality (\ref{eq:pro-upper-2period-Prod-prf4}) yields
\begin{align}
& \quad \max_{\lambda_{0,L_1}, \lambda_{0,Y}} \;  \mu^{\over1}_{\mytext{ICE}} (\lambda_{0,L_1}\lambda_{0,Y},\ldots,\lambda_{K-1,Y}) \nonumber \\
& \le \min_{ q_{0,L_1}, q_{0,Y} }\;  E \Big[
E_{\underlam_{1,\mytext{Jnt}} } \left\{ \eta_{0+,L_1} (Y^{\over1}, \eta_{0+,Y} ( Y^{\over1}, q_{0,Y}), q_{0,L_1} ) | A_0= 1, L_0 \right\} \Big]
\label{eq:pro-upper-Kperiod-Prod-prf3-init} \\
& = \min_{ q_{0,L_1}, q_{0,Y} }\;  E \Big[
E_{\underlam_{1,\mytext{Jnt}} } \left\{ \eta_{0+,\mytext{Prod,v1}} (Y^{\over1}, q_{0,L_1}, q_{0,Y} ) | A_0= 1, L_0 \right\} \Big] , \nonumber
\end{align}
i.e., (\ref{eq:pro-upper-Kperiod-Prod-prf3}) holds for $k=0$.

Suppose that for some $1\le k \le K-1$, equation (\ref{eq:pro-upper-Kperiod-Prod-prf3}) holds with $k$ replaced by $k-1$ for any fixed $\underlam_k$:
\begin{align}
& \quad \max_{\overlam_{k-1,L_{k}}, \overlam_{k-1,Y}} \;  \mu^{\over1}_{\mytext{ICE}} (\lambda_{0,L_1}\lambda_{0,Y},\ldots,\lambda_{K-1,Y}) \nonumber \\
& \le \min_{\overq_{k-1,L_{k}}, \overq_{k-1,Y} }\;  E \Big[  \Big( \prod_{j=0}^{k-2} \calE_{\overL_j} \Big)
E_{\underlam_{k,\mytext{Jnt}} } \left\{ \eta^{k-1}_{+,\mytext{Prod,v1}} (Y^{\over1}, \overq_{k-1,L_{k}}, \overq_{k-1,Y} )
| \overA_{k-1}= \over1_{k-1}, \overL_{k-1} \right\} \Big] .
\label{eq:pro-upper-Kperiod-Prod-prf4}
\end{align}
It is desired to show that
(\ref{eq:pro-upper-Kperiod-Prod-prf3}) holds. For any fixed $\underlam_{k+1,L_{k+1}}$ and $\underlam_{k+1,Y}$ (hence fixed $\underlam_{k+1,\mytext{Jnt}}$),
we apply an extension of (\ref{eq:pro-upper-2period-Prod-prf4}) [similar to the extension (\ref{eq:upper-1period-b}) of (\ref{eq:upper-1period})],
with $L_0$, $A_0$, $(\lambda_{0,L_1},\lambda_{0,Y},$ $ \lambda_1)$, and $(L^{1,1}_1, Y^{1,1})$ replaced by
$L_k$, $A_k$, $(\lambda_{k,L_{k+1}},\lambda_{k,Y}, \underlam_{k+1,\mytext{Jnt}})$, and $(L^{\over1}_{k+1}, (\underL^{\over1}_{k+2}, Y^{\over1}) )$ and
with the ``outcome'' $Y^{1,1}$ replaced by
$\eta^{k-1}_{+,\mytext{Prod,v1}} = \eta^{k-1}_{+,\mytext{Prod,v1}} (Y^{\over1}, \overq_{k-1,L_{k}}, \overq_{k-1,Y} ) $,
conditionally on $( \overA_{k-1}= \over1_{k-1}, \overL_{k-1} )$, and obtain
\begin{align}
& \quad \max_{\lambda_{k,L_{k+1}},\lambda_{k,Y}} \;
E_{\underlam_{k,\mytext{Jnt}}} \left\{ \eta^{k-1}_{+,\mytext{Prod,v1}} (Y^{\over1}, \overq_{k-1,L_{k}}, \overq_{k-1,Y} )
 | \overA_{k-1} = \over1_{k-1}, \overL_{k-1} \right\} \nonumber \\
& = \max_{\lambda_{k,L_{k+1}},\lambda_{k,Y}} \; \calE_{\overL_{k-1}}
E_{\underlam_{k+1,\mytext{Jnt}} } \left\{ (\pi^*_k + (1-\pi^*_k) \lambda_{k,\mytext{Jnt}} (\overL_k, \underL^{\over1}_{k+1}, Y^{\over1}))
\eta^{k-1}_{+,\mytext{Prod,v1}} | \overA_k = \over1_k, \overL_k \right\} \nonumber \\
& \le \max_{q_{k,L_{k+1}},q_{k,Y}} \; \calE_{\overL_{k-1}}
E_{\underlam_{k+1,\mytext{Jnt}} } \left\{ \eta_{k+}(
\eta^{k-1}_{+,\mytext{Prod,v1}}, q_{k,L_{k+1}},q_{k,Y}) | \overA_k = \over1_k, \overL_k  \right\} \label{eq:pro-upper-Kperiod-Prod-prf5}  \\
& = \max_{q_{k,L_{k+1}},q_{k,Y}} \; \calE_{\overL_{k-1}}
E_{\underlam_{k+1,\mytext{Jnt}} } \left\{
 \eta^k_{+,\mytext{Prod,v1}} (Y^{\over1}, \overq_{k,L_{k+1}}, \overq_{k,Y} ) | \overA_k = \over1_k, \overL_k  \right\} ,   \nonumber
\end{align}
where the first line follows from (\ref{eq:pro-upper-Kperiod-Prod-prf2}). Combining (\ref{eq:pro-upper-Kperiod-Prod-prf4}) and (\ref{eq:pro-upper-Kperiod-Prod-prf5}) yields
\begin{align*}
& \quad \max_{\lambda_{k,L_{k+1}},\lambda_{k,Y}} \max_{\overlam_{k-1,L_{k}}, \overlam_{k-1,Y}} \;  \mu^{\over1}_{\mytext{ICE}} (\lambda_{0,L_1}\lambda_{0,Y},\ldots,\lambda_{K-1,Y}) \\
& \le \max_{\lambda_{k,L_{k+1}},\lambda_{k,Y}} \min_{\overq_{k-1,L_{k}}, \overq_{k-1,Y} } \;  E \Big[  \Big( \prod_{j=0}^{k-2} \calE_{\overL_j} \Big)
E_{\underlam_{k,\mytext{Jnt}} } \left\{ \eta^{k-1}_{+,\mytext{Prod,v1}}
| \overA_{k-1}= \over1_{k-1}, \overL_{k-1} \right\} \Big] \\
& \le \min_{\overq_{k-1,L_{k}}, \overq_{k-1,Y} } \max_{\lambda_{k,L_{k+1}},\lambda_{k,Y}} \;  E \Big[  \Big( \prod_{j=0}^{k-2} \calE_{\overL_j} \Big)
E_{\underlam_{k,\mytext{Jnt}} } \left\{ \eta^{k-1}_{+,\mytext{Prod,v1}}
| \overA_{k-1}= \over1_{k-1}, \overL_{k-1} \right\} \Big]  \\
& \le \min_{\overq_{k-1,L_{k}}, \overq_{k-1,Y} }  \min_{q_{k,L_{k+1}},q_{k,Y}} \;  E \Big[  \Big( \prod_{j=0}^{k-1} \calE_{\overL_j} \Big)
E_{\underlam_{k+1,\mytext{Jnt}} } \left\{
 \eta^k_{+,\mytext{Prod,v1}} (Y^{\over1}, \overq_{k,L_{k+1}}, \overq_{k,Y} ) | \overA_k = \over1_k, \overL_k  \right\} \Big] .
\end{align*}
The third line above follows from the max-min inequality.
This completes the proof of (\ref{eq:pro-upper-Kperiod-Prod-prf3}).

\textbf{Second upper bound.}\;
We show by induction that for $k=K-1,\ldots, 1,0$ (downward) and any function $b_{K-1} (\overL_k, \underl_{k+1}, y)$,
\begin{align}
& \quad \max_{\underlam_{k,L_{k+1}}, \underlam_{k,Y}} \; E \Big[  \Big( \prod_{j=0}^{k-1} \calE_{\overL_j} \Big)
E_{\underlam_{k,\mytext{Jnt}}} \left\{ b_{K-1} (\overL_k, \underL^{\over1}_{k+1}, Y^{\over1} ) | \overA_{k-1} = \over1_{k-1}, \overL_k  \right\} \Big]  \nonumber \\
& \le \min_{\underq_{k,L_{k+1}}, \underq_{k,Y} }\; E \Big[  \Big( \prod_{j=0}^k \calE_{\overL_j} \Big)
 \left\{ \eta^k_{+,\mytext{Prod,v2}} ( b_{K-1} (\overL_{K-1}, Y), \underq_{k,L_{k+1}}, \underq_{k,Y} ) \right\} \Big] ,
\label{eq:pro-upper-Kperiod-Prod-prf6}
\end{align}
where  $E_{\underlam_{k,\mytext{Jnt}} } (\cdot | \overA_{k-1} = \over1_{k-1}, \overL_{k-1})$ denotes the conditional expectation  under $Q$
as in  (\ref{eq:pro-upper-Kperiod-Prod-prf1}),
and $\overA_{-1} = \over1_{-1}$ and $\overL_{-1}$ are set to the null and $E_{\underlam_{0,\mytext{Jnt}}} (\cdot | \overA_{-1} = \over1_{-1}, \overL_{-1} )$
reduces to $E_{\underlam_{0,\mytext{Jnt}}} (\cdot )$, the unconditional expectation under $Q$.
The case of (\ref{eq:pro-upper-Kperiod-Prod-prf6}) with $k=0$ and $b_{K-1} (\overL_0, \underl_1, y)=y$ yields the desired result.

First, we show that (\ref{eq:pro-upper-Kperiod-Prod-prf6}) holds for $k=K-1$.
Note that $\lambda_{K-1,L_K} \equiv 1$. By the law of iterated expectations,
\begin{align*}
& \quad E_{\lambda_{K-1,Y}} \left\{ b_{K-1} (\overL_{K-1}, Y^{\over1} ) | \overA_{K-2} = \over1_{K-2}, \overL_{K-1} \right\}   \nonumber \\
& = E \left\{ (\pi^*_{K-1} + (1-\pi^*_{K-1}) \lambda_{K-1,Y})
b_{K-1} (\overL_{K-1}, Y ) | \overA_{K-1} = \over1_{K-1}, \overL_{K-1} \right\} .
\end{align*}
We apply single-period identity (\ref{eq:upper-1period-b}),
with $L_0$, $A_0$, and $\lambda_0$ replaced by
$L_{K-1}$, $A_{K-1}$, and $\lambda_{K-1,Y}$ and with the ``outcome'' $b(L_0,Y)$ set to $b_{K-1} (\overL_{K-1}, Y)$,
conditionally on $(\overA_{K-2} = \over1_{K-2}, \overL_{K-2})$, and obtain
\begin{align*}
& \quad \max_{ \lambda_{K-1,Y}} \;
\mathcal E_{\overL_{K-2}} E_{\lambda_{K-1,Y}} \left\{ b_{K-1} (\overL_{K-1}, Y^{\over1} ) | \overA_{K-2} = \over1_{K-2}, \overL_{K-1} \right\}   \nonumber \\
& = \max_{ \lambda_{K-1,Y}} \;
 \mathcal E_{\overL_{K-2}}  E \left\{ (\pi^*_{K-1} + (1-\pi^*_{K-1}) \lambda_{K-1,Y})
b_{K-1} (\overL_{K-1}, Y ) | \overA_{K-1} = \over1_{K-1}, \overL_{K-1} \right\} \\
& = \min_{ q_{K-1,Y} }\;
 \mathcal E_{\overL_{K-2}} \mathcal E_{\overL_{K-1}} \left\{ \eta_{(K-1)+} ( b_{K-1} (\overL_{K-1}, Y), q_{K-1,Y} ) \right\} \\
& = \min_{ q_{K-1,Y} }\;
  \mathcal E_{\overL_{K-2}} \mathcal E_{\overL_{K-1}} \left\{  \eta_{(K-1)+} ( b_{K-1} (\overL_{K-1}, Y), q_{K-1,Y} ) \right\}.
\end{align*}
Hence (\ref{eq:pro-upper-Kperiod-Prod-prf6}) holds (with equality) for $k=K-1$.

Suppose that for some $1 \le k \le K-2$, equation (\ref{eq:pro-upper-Kperiod-Prod-prf6}) holds with $k$ replaced by $k+1$:
\begin{align}
& \quad \max_{\underlam_{k+1,L_{k+2}}, \underlam_{k+1,Y}} \;E \Big[  \Big( \prod_{j=0}^k \calE_{\overL_j} \Big)
E_{\underlam_{k+1,\mytext{Jnt}}} \left\{ b_{K-1} (\overL_{k+1}, \underL^{\over1}_{k+2}, Y^{\over1} ) | \overA_k = \over1_k, \overL_{k+1}  \right\}  \Big]  \nonumber \\
& \le \min_{\underq_{k+1,L_{k+2}}, \underq_{k+1,Y} }\; E \Big[  \Big( \prod_{j=0}^{k+1} \calE_{\overL_j} \Big)
 \left\{ \eta^{k+1}_{+,\mytext{Prod,v2}} ( b_{K-1} (\overL_{K-1}, Y), \underq_{k+1,L_{k+2}}, \underq_{k+1,Y} ) \right\} \Big] .
\label{eq:pro-upper-Kperiod-Prod-prf7}
\end{align}
It is desired to show that
(\ref{eq:pro-upper-Kperiod-Prod-prf6}) holds.
By the law of iterated expectations under $Q$,
\begin{align*}
& \quad
 E_{\underlam_{k,\mytext{Jnt}}} \left\{ b_{K-1} (\overL_{k-1}, \underL^{\over1}_k, Y^{\over1} ) | \overA_{k-1} = \over1_{k-1}, \overL_k \right\}   \\
& =  E_{\underlam_{k+1,\mytext{Jnt}}} \Big\{ \underbrace{ (\pi^*_k + (1-\pi^*_k) \lambda_{k,\mytext{Jnt}} (\overL_k, \underL^{\over1}_{k+1}, Y^{\over1}))
b_{K-1} (\overL_{k-1}, \underL^{\over1}_k, Y^{\over1} )  }_{\one}  | \overA_k = \over1_k, \overL_k \Big\} \\
& = E \left\{ E_{\underlam_{k+1,\mytext{Jnt}}} ( \;\one\;| \overA_k=\over1_k, \overL_{k+1} )  | \overA_k = \over1_k, \overL_k \right\} .
\end{align*}
Note that taking the conditional expectation $\mathcal E_{\overL_{k-1}}$ on both sides of the first equality above
leads to (\ref{eq:pro-upper-Kperiod-Prod-prf2}).
Then we apply identity (\ref{eq:upper-1period-c}) conditionally on $(\overA_k = \over1_k, \overL_k )$,
similarly as in (\ref{eq:pro-upper-2period-Prod-prf2}) for the proof of Second upper bound in Proposition~\ref{pro:upper-2period-Prod}:
\begin{align*}
& \quad \max_{\lambda_{k,Y}} \;
 E_{\underlam_{k,\mytext{Jnt}}} \left\{ b_{K-1} (\overL_{k-1}, \underL^{\over1}_k, Y^{\over1} ) | \overA_{k-1} = \over1_{k-1}, \overL_k \right\}   \\
& = \max_{\lambda_{k,Y}}\;
 E \left\{ E_{\underlam_{k+1,\mytext{Jnt}}} ( \;\one\;| \overA_k=\over1_k, \overL_{k+1} )
 | \overA_k = \over1_k, \overL_k \right\} \\
& = \min_{q_{k,Y}} \;
\mathcal E_{\overL_k}  \Big[ E_{\underlam_{k+1,\mytext{Jnt}}} \Big\{  \underbrace{ \pi^*_k  b_{K-1}  + (1-\pi^*_k) \lambda_{k,L_{k+1}}
 \eta_{k+,Y} ( b_{K-1}, q_{k,Y} )  }_{\two} | \overA_k=\over1_k, \overL_{k+1}  \Big\} \Big] .
\end{align*}
Moreover, we apply (\ref{eq:pro-upper-Kperiod-Prod-prf7}) with $b_{K-1}$ replaced by $\two$:
\begin{align*}
& \quad \max_{\underlam_{k+1,L_{k+2}}, \underlam_{k+1,Y}} \; E \Big[  \Big( \prod_{j=0}^{k-1} \calE_{\overL_j} \Big)
 \mathcal E_{\overL_k}  E_{\underlam_{k+1,\mytext{Jnt}}} ( \;\two\; | \overA_k = \over1_k, \overL_{k+1} ) \Big] \\
& \le \min_{\underq_{k+1,L_{k+2}}, \underq_{k+1,Y} }\;  E \Big[  \Big( \prod_{j=0}^k \calE_{\overL_j} \Big)
 \mathcal E_{\overL_{k+1}}  \left\{ \eta^{k+1}_{+,\mytext{Prod,v2}} ( \;\two , \underq_{k+1,L_{k+2}}, \underq_{k+1,Y} ) \right\} \Big]
\end{align*}
Combining the preceding two displays yields
\begin{align}
& \quad  \max_{\underlam_{k+1,L_{k+2}}, \underlam_{k+1,Y}} \max_{\lambda_{k,Y}} \; E \Big[  \Big( \prod_{j=0}^{k-1} \calE_{\overL_j} \Big)
E_{\underlam_{k,\mytext{Jnt}}} \left\{ b_{K-1} (\overL_k, \underL^{\over1}_{k+1}, Y^{\over1} ) | \overA_{k-1} = \over1_{k-1}, \overL_k  \right\} \Big] \nonumber \\
& =  \max_{\underlam_{k+1,L_{k+2}}, \underlam_{k+1,Y}} \min_{q_{k,Y}} \; E \Big[  \Big( \prod_{j=0}^{k-1} \calE_{\overL_j} \Big)
 \mathcal E_{\overL_k}  E_{\underlam_{k+1,\mytext{Jnt}}} ( \;\two\; | \overA_k = \over1_k, \overL_{k+1} ) \Big]  \nonumber \\
& \le \min_{q_{k,Y}}   \max_{\underlam_{k+1,L_{k+2}}, \underlam_{k+1,Y}}  \; E \Big[  \Big( \prod_{j=0}^{k-1} \calE_{\overL_j} \Big)
 \mathcal E_{\overL_k}  E_{\underlam_{k+1,\mytext{Jnt}}} ( \;\two\; | \overA_k = \over1_k, \overL_{k+1} ) \Big]
 \label{eq:pro-upper-Kperiod-Prod-prf8} \\
& \le \min_{q_{k,Y}}  \min_{\underq_{k+1,L_{k+2}}, \underq_{k+1,Y} }\;  E \Big[  \Big( \prod_{j=0}^k  \calE_{\overL_j} \Big)
 \mathcal E_{\overL_{k+1}}  \left\{ \eta^{k+1}_{+,\mytext{Prod,v2}} ( \;\two , \underq_{k+1,L_{k+2}}, \underq_{k+1,Y} ) \right\} \Big] , \nonumber
\end{align}
where the third line follows from the max-min inequality.
By the convexity and homogeneity of $\eta^{k+1}_{+,\mytext{Prod,v2}} (\cdot)$
in Lemma~\ref{lem:eta2-convex} (i) and (iii), it can be shown similarly as in (\ref{eq:pro-upper-2period-Prod-prf6}) that
\begin{align}
& \quad  \min_{\underq_{k+1,L_{k+2}}, \underq_{k+1,Y} }\;  E \Big[  \Big( \prod_{j=0}^k \calE_{\overL_j} \Big)
 \mathcal E_{\overL_{k+1}}  \left\{ \eta^{k+1}_{+,\mytext{Prod,v2}} ( \;\two , \underq_{k+1,L_{k+2}}, \underq_{k+1,Y} ) \right\} \Big] \nonumber \\
& \le  \min_{\underq_{k+1,L_{k+2}}, \underq_{k+1,Y} }\;  E \Big[  \Big( \prod_{j=0}^k \calE_{\overL_j} \Big)
  \underbrace{ \mathcal E_{\overL_{k+1}}  \Big\{  \pi^*_k  \eta^{k+1}_{+,\mytext{Prod,v2}}
    + (1-\pi^*_k) \lambda_{k,L_{k+1}}  (\eta^{k+1}_{+,\mytext{Prod,v2}}\circ \eta_{k+,Y}) \Big\} }_{\three} \Big], \label{eq:pro-upper-Kperiod-Prod-prf9}
\end{align}
where, as abbreviations, $\eta^{k+1}_{+,\mytext{Prod,v2}}=\eta^{k+1}_{+,\mytext{Prod,v2}} (b_{K-1} , \underq_{k+1,L_{k+2}}, \underq_{k+1,Y} )$ and
$\eta^{k+1}_{+,\mytext{Prod,v2}}\circ \eta_{k+,Y} = \eta^{k+1}_{+,\mytext{Prod,v2}} ( \eta_{k+,Y} ( b_{K-1}, q_{k,Y} ), \underq_{k+1,L_{k+2}}, \underq_{k+1,Y} )$.
We apply identity (\ref{eq:upper-1period-b}) conditionally on $(\overA_{k-1} = \over1_{k-1}, \overL_{k-1} )$, similarly as in (\ref{eq:pro-upper-2period-Prod-prf7}):
\begin{align*}
& \quad \max_{\lambda_{k,L_{k+1}}} \;
\mathcal E_{\overL_k}   (\;\three\;)  \\
& = \min_{q_{k,L_{k+1}}}  \;
\mathcal E_{\overL_k} \left\{ \eta_{k+,L_{k+1}} \left(
\mathcal E_{\overL_{k+1}}  \eta^{k+1}_{+,\mytext{Prod,v2}}  ,
\mathcal E_{\overL_{k+1}} (\eta^{k+1}_{+,\mytext{Prod,v2}} \circ \eta_{k+,Y}) , q_{k,L_{k+1}} \right) \right\} \\
& = \min_{q_{k,L_{k+1}}}  \;
\mathcal E_{\overL_k} \left\{  \eta^k_{+,\mytext{Prod,v2}} (b_{K-1} , \underq_{k,L_{k+1}}, \underq_{k,Y} ) \right\}.
\end{align*}
Combining the preceding three displays yields
\begin{align*}
& \quad  \max_{\lambda_{k,L_{k+1}}} \max_{\underlam_{k+1,L_{k+2}}, \underlam_{k+1,Y}} \max_{\lambda_{k,Y}} \; E \Big[  \Big( \prod_{j=0}^{k-1} \calE_{\overL_j} \Big)
E_{\underlam_{k,\mytext{Jnt}}} \left\{ b_{K-1} (\overL_k, \underL^{\over1}_{k+1}, Y^{\over1} ) | \overA_{k-1} = \over1_{k-1}, \overL_k  \right\} \Big] \\
& \le  \max_{\lambda_{k,L_{k+1}}} \min_{q_{k,Y}}  \min_{\underq_{k+1,L_{k+2}}, \underq_{k+1,Y} }\;  E \Big[  \Big( \prod_{j=0}^k \calE_{\overL_j} \Big)
  (\;\three\;) \Big] \\
&  \le  \min_{q_{k,Y}}  \min_{\underq_{k+1,L_{k+2}}, \underq_{k+1,Y} } \max_{\lambda_{k,L_{k+1}}} \;  E \Big[  \Big( \prod_{j=0}^k \calE_{\overL_j} \Big)
  (\;\three\;) \Big] \\
& = \min_{q_{k,Y}}  \min_{\underq_{k+1,L_{k+2}}, \underq_{k+1,Y} } \min_{q_{k,L_{k+1}}}  \; E \Big[  \Big( \prod_{j=0}^k \calE_{\overL_j} \Big)
 \left\{ \eta^k_{+,\mytext{Prod,v2}} (b_{K-1} , \underq_{k,L_{k+1}}, \underq_{k,Y} ) \right\} \Big] ,
\end{align*}
where the third line follows from the max-min inequality.
This completes the proof of (\ref{eq:pro-upper-Kperiod-Prod-prf6}).

\begin{lem} \label{lem:eta2-convex}
(i) For $k=K-1,\ldots,1,0$,  the functional
$\eta^k_{+,\mytext{Prod,v2}} (b_{K-1} , \underq_{k,L_{k+1}}, \underq_{k,Y} )$ is nondecreasing and convex in $b_{K-1}$: for any functions $b_{K-1}$ and $\tilde b_{K-1}$
and $\alpha = \alpha (\overL_k ) \in [0,1] $,
\begin{align*}
& \quad \eta^k_{+,\mytext{Prod,v2}} ( \alpha b_{K-1}+(1-\alpha)\tilde b_{K-1} , \underq_{k,L_{k+1}}, \underq_{k,Y} ) \\
& \le \alpha \eta^k_{+,\mytext{Prod,v2}} (  b_{K-1} , \underq_{k,L_{k+1}}, \underq_{k,Y} )  +
 (1- \alpha) \eta^k_{+,\mytext{Prod,v2}} ( \tilde b_{K-1} , \underq_{k,L_{k+1}}, \underq_{k,Y} ) ,
\end{align*}
and if $b_K \le \tilde b_K$, then $\eta^k_{+,\mytext{Prod,v2}} (  b_{K-1} , \underq_{k,L_{k+1}}, \underq_{k,Y} )
\le \eta^k_{+,\mytext{Prod,v2}} ( \tilde b_{K-1} , \underq_{k,L_{k+1}}, \underq_{k,Y} )$.\\
(ii) $\eta^k_{+,\mytext{Prod,v2}} (b_{K-1} , \underq_{k,L_{k+1}}, \underq_{k,Y} )$ is convex in $(\underq_{k,L_{k+1}}, \underq_{k,Y} )$:
for any functions $(\underq_{k,L_{k+1}},$ $ \underq_{k,Y} )$ and $(\tilde\underq_{k,L_{k+1}}, \tilde\underq_{k,Y} )$
and scalar $\alpha\in [0,1] $,
\begin{align*}
& \quad \eta^k_{+,\mytext{Prod,v2}} (b_{K-1} , \alpha \underq_{k,L_{k+1}}+(1-\alpha)\tilde \underq_{k,L_{k+1}}, \alpha\underq_{k,Y} + (1-\alpha) \tilde\underq_{k,Y} ) \\
& \le \alpha \eta^k_{+,\mytext{Prod,v2}} (  b_{K-1} , \underq_{k,L_{k+1}}, \underq_{k,Y} )  +
 (1- \alpha) \eta^k_{+,\mytext{Prod,v2}} ( b_{K-1} , \tilde \underq_{k,L_{k+1}}, \tilde \underq_{k,Y} )  .
\end{align*}
(iii) $\eta^k_{+,\mytext{Prod,v2}} (b_{K-1} , \underq_{k,L_{k+1}}, \underq_{k,Y} )$ is homogeneous:
for any function $\gamma = \gamma (\overL_k ) >0$,
\begin{align*}
\eta^k_{+,\mytext{Prod,v2}} (\gamma b_{K-1}  , \gamma \underq_{k,L_{k+1}}, \gamma \underq_{k,Y} )
 = \gamma \eta^k_{+,\mytext{Prod,v2}} ( b_{K-1}  ,  \underq_{k,L_{k+1}},  \underq_{k,Y} )  .
\end{align*}
\end{lem}

\begin{prf}
The proof is by induction. For $k=K-1$, the results (i)--(iii) hold directly; the nondecreasingness of $\eta_{(K-1)+} (b_{K-1}; q_{K-1,Y})$ in $b_{K-1}$
follows from that of $\eta_{k+}(y,q_k)$ in $y$ by the proof of Lemma~\ref{lem:eta-convex}.
Suppose that for some $0 \le k \le K-2$, the results hold for $\eta^{k+1}_{+,\mytext{Prod,v2}}$.
It is desired to show that the results hold for $\eta^k_{+,\mytext{Prod,v2}}$. By definition,
\begin{align*}
& \quad \eta^k_{+,\mytext{Prod,v2}} (b_{K-1} , \underq_{k,L_{k+1}}, \underq_{k,Y} ) \\
& = \eta_{k+,L_{k+1}} \big\{ \mathcal E_{\overL_{k+1}} \eta^{k+1}_{+,\mytext{Prod,v2}} (b_{K-1} , \underq_{k+1,L_{k+2}}, \underq_{k+1,Y} ), \\
& \qquad \mathcal E_{\overL_{k+1}} \eta^{k+1}_{+,\mytext{Prod,v2}} ( \eta_{k+,Y} ( b_{K-1}, q_{k,Y} ), \underq_{k+1,L_{k+2}}, \underq_{k+1,Y} ), q_{k,L_{k+1}} \big\}.
\end{align*}

The convexity of $\eta^k_{+,\mytext{Prod,v2}}(b_{K-1} , \underq_{k,L_{k+1}}, \underq_{k,Y} )$ in $b_{K-1}$ can be shown as follows:
\begin{align*}
& \quad \eta^k_{+,\mytext{Prod,v2}} ( \alpha b_{K-1}+(1-\alpha)\tilde b_{K-1} , \underq_{k,L_{k+1}}, \underq_{k,Y} ) \\
& \le \eta_{k+,L_{k+1}} \big\{ \mathcal E_{\overL_{k+1}} \eta^{k+1}_{+,\mytext{Prod,v2}} (\alpha b_{K-1}+(1-\alpha)\tilde b_{K-1} , \underq_{k+1,L_{k+2}}, \underq_{k+1,Y} ), \\
& \quad \mathcal E_{\overL_{k+1}} \eta^{k+1}_{+,\mytext{Prod,v2}} ( \alpha \eta_{k+,Y}(b_{K-1}, q_{k,Y} )
+(1-\alpha) \eta_{k+,Y}(\tilde b_{K-1}, q_{k,Y} ), \underq_{k+1,L_{k+2}}, \underq_{k+1,Y} ), q_{k,L_{k+1}} \big\} \\
& \le \eta_{k+,L_{k+1}} \big\{ \alpha  \mathcal E_{\overL_{k+1}} \eta^{k+1}_{+,\mytext{Prod,v2}} (b_{K-1}, \underq_{k+1,L_{k+2}}, \underq_{k+1,Y} ) \\
& \quad + (1-\alpha) \mathcal E_{\overL_{k+1}} \eta^{k+1}_{+,\mytext{Prod,v2}} (\tilde b_{K-1}, \underq_{k+1,L_{k+2}}, \underq_{k+1,Y} ) , \\
& \quad \alpha\mathcal E_{\overL_{k+1}} \eta^{k+1}_{+,\mytext{Prod,v2}} ( \eta_{k+,Y}(b_{K-1}, q_{k,Y} ) , \underq_{k+1,L_{k+2}}, \underq_{k+1,Y} ) \\
& \quad + (1-\alpha) \mathcal E_{\overL_{k+1}} \eta^{k+1}_{+,\mytext{Prod,v2}} ( \eta_{k+,Y}(\tilde b_{K-1}, q_{k,Y} ),
\underq_{k+1,L_{k+2}}, \underq_{k+1,Y} ), q_{k,L_{k+1}} \big\}  \\
& \le \alpha \eta_{k+,L_{k+1}} \big\{ \mathcal E_{\overL_{k+1}} \eta^{k+1}_{+,\mytext{Prod,v2}} (b_{K-1} , \underq_{k+1,L_{k+2}}, \underq_{k+1,Y} ), \\
& \quad \mathcal E_{\overL_{k+1}} \eta^{k+1}_{+,\mytext{Prod,v2}} ( \eta_{k+,Y} ( b_{K-1}, q_{k,Y} ), \underq_{k+1,L_{k+2}}, \underq_{k+1,Y} ), q_{k,L_{k+1}} \big\} \\
& \quad + (1-\alpha) \eta_{k+,L_{k+1}} \big\{ \mathcal E_{\overL_{k+1}} \eta^{k+1}_{+,\mytext{Prod,v2}} (\tilde b_{K-1} , \underq_{k+1,L_{k+2}}, \underq_{k+1,Y} ), \\
& \quad \mathcal E_{\overL_{k+1}} \eta^{k+1}_{+,\mytext{Prod,v2}} ( \eta_{k+,Y} (\tilde b_{K-1}, q_{k,Y} ), \underq_{k+1,L_{k+2}}, \underq_{k+1,Y} ), q_{k,L_{k+1}} \big\}.
\end{align*}
The first step uses
the convexity of $\eta_{k+,Y} ( y, q_{k,Y} )$ in $y$,
the nondecreasingness of $\eta^k_{+,\mytext{Prod,v2}} (b_{K-1} ,$ $ \underq_{k,L_{k+1}}, \underq_{k,Y} )$ in $b_{K-1}$, and
the nondecreasingness of $\eta_{k+,L_{k+1}} ( y, \tilde y,  q_{k,L_{k+1}} )$ in $\tilde y$ by the proof of Lemma~\ref{lem:eta-convex}.
The second step uses
the convexity of $\eta^k_{+,\mytext{Prod,v2}} (b_{K-1} , \underq_{k,L_{k+1}}, \underq_{k,Y} )$ in $b_{K-1}$, and
the nondecreasingness of $\eta_{k+,L_{k+1}} ( y, \tilde y,  q_{k,L_{k+1}} )$ in $y$ and in $\tilde y$.
The third step uses
the convexity of $\eta_{k+,L_{k+1}} ( y, \tilde y,  q_{k,L_{k+1}} )$ in $(y,\tilde y)$ jointly.
Similarly, the convexity of $\eta^k_{+,\mytext{Prod,v2}} (b_{K-1} , \underq_{k,L_{k+1}}, \underq_{k,Y} )$ in $(\underq_{k,L_{k+1}}, \underq_{k,Y} )$
can be shown.

The nondecreasingness of $\eta^k_{+,\mytext{Prod,v2}}(b_{K-1} , \underq_{k,L_{k+1}}, \underq_{k,Y} )$ in $b_{K-1}$
can be shown by using the nondecreasingness of $\eta_{k+,Y} ( y, q_{k,Y} )$ in $y$,
the nondecreasingness of $\eta^k_{+,\mytext{Prod,v2}} (b_{K-1} ,$ $ \underq_{k,L_{k+1}}, \underq_{k,Y} )$ in $b_{K-1}$, and
the nondecreasingness of $\eta_{k+,L_{k+1}} ( y, \tilde y,  q_{k,L_{k+1}} )$ in $y$ and in $\tilde y$.

The homogeneity of  $\eta^k_{+,\mytext{Prod,v2}} (b_{K-1} , \underq_{k,L_{k+1}}, \underq_{k,Y} )$  in its arguments
can be shown by using the homogeneity of $\eta_{k+,Y} ( y, q_{k,Y} )$,
the homogeneity of $ \eta^{k+1}_{+,\mytext{Prod,v2}} (b_{K-1} , \underq_{k+1,L_{k+2}}, \underq_{k+1,Y} )$,
and the homogeneity of $\eta_{k+,L_{k+1}}(y, \tilde y, q_{k,L_{k+1}} )$ in their corresponding arguments.
\end{prf}

\subsection{Proof of Corollary \ref{cor:upper-Kperiod-Prod}}

For results (i) and (ii), we show that in the special case considered, the corresponding conservative bound in Proposition \ref{pro:upper-Kperiod-Prod}
becomes exact and takes the simplified form.

(i) Suppose that $\Lambda_{k,L_{k+1}}=1$ (i.e., $L_{k+1}^{\over1} \perp A_k \,|\, \overA_{k-1}=\over1_{k-1}, \overL_k $)
and $\Lambda_{k,Y} = \Lambda_k$ for $k=0,1,\ldots,K-1$. Then by definition,
$ \eta_{k+,L_{k+1}} (y, \tilde y, q_{k,L_{k+1}} )
 = \pi^*_k y + (1-\pi^*_k) \tilde y$ and
\begin{align*}
 & \eta_{k+, L_{k+1}} (y,  \eta_{k+,Y} (y, q_{k,Y}) , q_{k,L_{k+1}} ) = \eta_{k+} (y, q_{k,Y}) ,
\end{align*}
independently of $q_{k,L_{k+1}}$. Hence $ \eta^k_{+,\mytext{Prod,v1}} (y, \overq_{k,L_{k+1}}, \overq_{k,Y} ) = \eta^k_+ ( y, \overq_{k,Y}) $ by induction,
and $\mu^{\over1}_{+,\mytext{Prod,v1}} (\over1_{K-1},\overLam_{K-1})$ can be simplified as stated in Corollary~\ref{cor:upper-Kperiod-Prod}.
In the induction proof of (\ref{eq:pro-upper-Kperiod-Prod-prf3}), inequalities (\ref{eq:pro-upper-Kperiod-Prod-prf3-init})
and (\ref{eq:pro-upper-Kperiod-Prod-prf5}) become equalities, because  (\ref{eq:pro-upper-2period-Prod-prf4}) used becomes
equality, as explained in the proof of Corollary \ref{cor:upper-2period-Prod}.
The max-min inequality in the final step also becomes equality by Sion's minimax theorem because
the functional involved,
\begin{align*}
& \quad E \Big[  \Big( \prod_{j=0}^{k-2} \calE_{\overL_j} \Big)
E_{\underlam_{k,\mytext{Jnt}} } \left\{ \eta^{k-1}_{+,\mytext{Prod,v1}}
| \overA_{k-1}= \over1_{k-1}, \overL_{k-1} \right\} \Big] \\
& = E \Big[ \Big( \prod_{j=0}^{k-1} \calE_{\overL_j} \Big)
E_{\underlam_{k+1,\mytext{Jnt}} } \left\{ (\pi^*_k + (1-\pi^*_k) \lambda_{k,Y} )
\eta^{k-1}_+ ( Y^{\over1}, \overq_{k-1,Y}) | \overA_k = \over1_k, \overL_k \right\} \Big]
\end{align*}
is linear in $\lambda_{k,Y}$ and convex in $\overq_{k-1,Y}$, in the absence of $\lambda_{k,L_{k+1}}$ and $\overq_{k-1,L_k}$.
Hence (\ref{eq:pro-upper-Kperiod-Prod-prf3}) holds with equality, which implies that the bound $\mu^{\over1}_{+,\mytext{Prod,v1}} (\over1_{K-1},\overLam_{K-1})$ becomes exact.

(ii) Suppose that $\Lambda_{k,Y} =1$ (i.e., $Y^{\over1} \perp A_k \,|\, \overA_{k-1}=\over1_{k-1}, \overL_k, \overL_{k+1}^{\over1}$)
and $\Lambda_{k,L_{k+1}} = \Lambda_k$ for $k=0,1,\ldots,K-2$ and $\Lambda_{K-1,Y}=\Lambda_{K-1}$.
Then by definition,
$ \eta_{k+,Y} (y,  q_{k,Y} )= y $  and
\begin{align*}
 & \eta_{k+, L_{k+1}} (y,  \eta_{k+,Y} (y, q_{k,Y}) , q_{k,L_{k+1}} ) = \eta_{k+} (y, q_{k,L_{k+1}}) ,
\end{align*}
independently of $q_{k,Y}$.
Hence $\eta^k_{+,\mytext{Prod,v2}} (Y, \underq_{k,L_{k+1}}, \underq_{k,Y} )$ and
$\mu^{\over1}_{+,\mytext{Prod,v2}} ( (\overLam_{K-2},1), (\over1_{K-2},\Lambda_{K-1}) )$ can be simplified as stated in Corollary~\ref{cor:upper-Kperiod-Prod}.
In the induction proof of (\ref{eq:pro-upper-Kperiod-Prod-prf7}), optimization over $\lambda_{k,Y}$ for $k \le K-2$ is no longer needed and hence
the max-min inequality (\ref{eq:pro-upper-Kperiod-Prod-prf8}) and the convexity-based upper bound (\ref{eq:pro-upper-Kperiod-Prod-prf9}) can be dropped.
The max-min inequality in the final step also becomes equality by Sion's minimax theorem because
the functional involved,
\begin{align*}
& \quad E \Big[  \Big( \prod_{j=0}^{k-1} \calE_{\overL_j} \Big)
\mathcal E_{\overL_k}   (\;\three\;) \Big] \\
& = E \Big[  \Big( \prod_{j=0}^{k-1} \calE_{\overL_j} \Big)
 \mathcal E_{\overL_k}  \left\{ ( \pi^*_k + (1-\pi^*_k) \lambda_{k,L_{k+1}} ) \mathcal E_{\overL_{k+1}}
 \eta^{k+1}_{+,\mytext{v2}} ( b_{K-1}, \underq_{k+1,L_{k+2}}, q_{K-1,Y}) \right\} \Big]
\end{align*}
is linear in $\lambda_{k,L_{k+1}}$ and convex in $(\underq_{k+1,L_{k+2}}, q_{K-1,Y})$.
Hence (\ref{eq:pro-upper-Kperiod-Prod-prf7}) holds with equality, which implies that the bound $\mu^{\over1}_{+,\mytext{Prod,v2}} ( (\overLam_{K-2},1), (\over1_{K-2},\Lambda_{K-1}) )$ becomes exact.

\section{Technical details: Proofs for Section \ref{sec:further-discussion}}

\subsection{Extension of Lemma~\ref{lem:density-ratio}}

We show the Lemma~\ref{lem:density-ratio} can be extended to multiple treatment strategies.

\begin{lem} \label{lem:density-ratio-multi}
Suppose that Assumption A1 holds for each $\overa \in \{0,1\}\times \{0,1\}$.
The sensitivity ratios $\lambda^{\overa *}_1 $ and  $\lambda^{\overa *}_0 $
satisfy the normalization constraints
as $\lambda^{\overa}_1 $ and $\lambda^{\overa}_0 $ in (\ref{eq:constr-2period-b-multi}) for each $\overa \in \{0,1\}\times \{0,1\}$.
Conversely, if nonnegative functions $\{ ( \lambda^{\overa}_1 (\overline{L}_1, y), \lambda^{\overa}_0 (L_0, y)) : \overa \in \{0,1\}\times \{0,1\} \}$
satisfy the normalization constraints for each $\overa =(a_0,a_1)\in \{0,1\}\times \{0,1\}$,
\begin{align} \label{eq:constr-2period-b-multi}
\begin{split}
& E ( \lambda^{\overa}_1(\overline{L}_1,Y) | A_0=a_0, A_1=a_1, \overline{L}_1 ) \equiv 1, \\
& E\{ E ( \lambda^{\overa}_0(L_0,Y) \varrho^{\overa}_1 (\overline{L}_1, Y; \lambda^{\overa}_1 ) |A_0=a_1,A_1=a_1, \overline{L}_1)  | A_0=a_0, L_0 \}\equiv 1 ,
\end{split}
\end{align}
with $\varrho^{\overa}_1 (\overline{L}_1, y; \lambda^{\overa}_1 ) = \pi^*_1(\overline{L}_1) + (1-\pi^*_1(\overline{L}_1)) \lambda^{\overa}_1(\overline{L}_1, y) $,
then there exists a probability distribution $Q$ for the full data $(\overline{A}_1,\overline{L}_1,Y^{1,1}, Y^{1,0}, Y^{0,1}, Y^{0,0})$ satisfying
the two properties: \vspace{-.05in}
\begin{itemize}\addtolength{\itemsep}{-.1in}
\item[(i)] the induced distribution of $Q$ on the observed data $(\overline{A}_1,\overline{L}_1,Y )$ under Assumption A1 coincides with
the true distribution of  $(\overline{A}_1,\overline{L}_1,Y )$,
\item[(ii)] $ \lambda^{\overa *}_{1,Q} (\overline{L}_1, y) =\lambda^{\overa}_1 (\overline{L}_1, y) $ and
 $ \lambda^{\overa *}_{0,Q} (L_0, y) = \lambda^{\overa}_0 (L_0, y)$ for each $\overa \in \{0,1\}\times \{0,1\}$,
\end{itemize} \vspace{-.05in}
where $\lambda^{\overa *}_{1,Q}$ and $\lambda^{\overa *}_{0,Q}$ are defined as $\lambda^{\overa *}_1$ and $\lambda^{\overa *}_0$
in Section \ref{sec:further-discussion} with $P$ replaced by $Q$.
\end{lem}

\begin{prf}
The necessity assertion in Lemma~\ref{lem:density-ratio-multi} holds directly by that in Lemma~\ref{lem:density-ratio}.
We show the sufficiency assertion holds in Lemma~\ref{lem:density-ratio-multi}.

For each $\overline{t} \in \{0,1\}\times \{0,1\}$,
applying the first four steps in the proof of Lemma~\ref{lem:density-ratio}, with
 $( \lambda_1 (\overline{L}_1, y), \lambda_0 (L_0, y))$ replaced by
 $( \lambda^{\overline{t}}_1 (\overline{L}_1, y), \lambda^{\overline{t}}_0 (L_0, y))$,
completes the definition of $Q$ for $(\overline{A}_1, \overline{L}_1, Y^{\overline{t}})$,
depending on $( \lambda^{\overline{t}}_1 (\overline{L}_1, y), \lambda^{\overline{t}}_0 (L_0, y))$
but not $( \lambda^{\overa}_1 (\overline{L}_1, y), \lambda^{\overa}_0 (L_0, y))$ for any $\overa \not= \overline{t}$.
The marginal distributions of $(\overline{A}_1,\overline{L}_1)$ from these definitions
are consistent with each other, because
$\dif Q_{\overline{A}_1,\overline{L}_1} = \dif P_{\overline{A}_1,\overline{L}_1}$ in the first step.
Finally, we define
\begin{align*}
& \quad \dif Q_{Y^{1,1}, Y^{1,0}, Y^{0,1}, Y^{0,0}} ( \cdot | \overA_1, \overL_1)
 = \prod_{a_0,a_1=0,1} \dif Q_{Y^{a_0a_1}} (\cdot | \overA_1, \overL_1 )
\end{align*}
i.e., $(Y^{1,1},Y^{1,0}, Y^{0,1}, Y^{0,0})$ are conditionally independent of each other given $(\overline{A}_1, \overline{L}_1)$ under $Q$.
This completes our construction of $Q$ for the full data $(\overline{A}_1, \overline{L}_1, Y^{1,1}, Y^{1,0}, Y^{0,1}, Y^{0,0})$,
which can be easily verified to satisfy the desired properties (i) and (ii).
\end{prf}

\subsection{Proof of Proposition \ref{pro:upper-2period-multi}}

In Lemma~\ref{lem:density-ratio-multi}, the normalization constraints (\ref{eq:constr-2period-b-multi}) are separable in $\overa$:
for example, the normalization constraints (\ref{eq:constr-2period-b-multi}) for $( \lambda^{1,1}_1 (\overline{L}_1, y), \lambda^{1,1}_0 (L_0, y))$
are independent of those for $( \lambda^{\overa}_1 (\overline{L}_1, y), $  $\lambda^{\overa}_0 (L_0, y))$
associated with any other treatment strategy $\overa\not=(1,1)$.
By this separability in $\overa$, both results (i) and (ii) follow easily as stated in Proposition \ref{pro:upper-2period-multi}.

\subsection{Extension of Lemma~\ref{lem:density-ratio-Jnt}}

We show that Lemma~\ref{lem:density-ratio-Jnt} can be extended to multiple treatment strategies.

\begin{lem} \label{lem:density-ratio-Jnt-multi}
Suppose that Assumption A1$^\dag$ holds  for each $\overa \in \{0,1\}\times \{0,1\}$.
The sensitivity ratios $\lambda^{\overa *}_1 $ and  $\lambda^{\overa *}_{0,\mytext{Jnt}} $
satisfy the normalization constraints
as $\lambda^{\overa}_1 $ and $\lambda^{\overa}_{0,\mytext{Jnt}} $ in (\ref{eq:constr-2period-Jnt-b-multi}) for each $\overa \in \{0,1\}\times \{0,1\}$.
Conversely, if nonnegative functions $\{ ( \lambda^{\overa *}_1(\overL_1,y), \lambda^{\overa *}_{0,\mytext{Jnt}}(L_0,l_1,y)): \overa \in \{0,1\}\times \{0,1\}\}$
satisfy the normalization constraints for each $\overa \in \{0,1\}\times \{0,1\}$,
with $\varrho_1^{\overa}(\cdot)$ defined as in (\ref{eq:constr-2period-b-multi}),
\begin{align} \label{eq:constr-2period-Jnt-b-multi}
\begin{split}
& E ( \lambda^{\overa}_1(\overline{L}_1,Y) | A_0=1, A_1=1, \overline{L}_1 ) \equiv 1,  \\ 
& E\{ E ( \lambda^{\overa}_{0,\mytext{Jnt}}(\overline{L}_1,Y) \varrho_1^{\overa} (\overline{L}_1, Y; \lambda^{\overa}_1 ) |A_0=1,A_1=1, \overline{L}_1)  | A_0=1, L_0 \}\equiv 1 ,
\end{split}
\end{align}
then there exists a probability distribution $Q$ for the full data $(\overline{A}_1, L_0,\{(L_1^{a_0a_1},Y^{a_0a_1}): a_0,a_1=0,1\} )$ satisfying
the two properties: \vspace{-.05in}
\begin{itemize}\addtolength{\itemsep}{-.1in}
\item[(i)] the induced distribution of $Q$ on the observed data $(\overline{A}_1,\overline{L}_1,Y )$ under Assumption A1$^\dag$ coincides with
the true distribution of  $(\overline{A}_1,\overline{L}_1,Y )$,
\item[(ii)] $ \lambda^{\overa *}_{1,Q} (\overline{L}_1, y) =\lambda^{\overa}_1 (\overline{L}_1, y) $ and
 $ \lambda^{\overa *}_{0,\mytext{Jnt},Q} (L_0, l_1, y) = \lambda^{\overa}_{0,\mytext{Jnt}} (L_0, l_1, y)$ for each $\overa \in \{0,1\}\times \{0,1\}$,
\end{itemize} \vspace{-.05in}
where $\lambda^{\overa *}_{1,Q}$ and $\lambda^{\overa *}_{0,\mytext{Jnt},Q}$ are defined as $\lambda^{\overa *}_1$ and $\lambda^{\overa *}_{0,\mytext{Jnt}}$
in Section \ref{sec:further-discussion} with $P$ replaced by $Q$.
\end{lem}

\begin{prf}
The necessity assertion in Lemma~\ref{lem:density-ratio-Jnt-multi} holds directly by that in Lemma~\ref{lem:density-ratio-Jnt}.
We show the sufficiency assertion holds in Lemma~\ref{lem:density-ratio-Jnt-multi}.

For each $\overline{t} \in \{0,1\}\times \{0,1\}$,
applying the first five steps in the proof of Lemma~\ref{lem:density-ratio-Jnt}, with
 $( \lambda_1 (\overline{L}_1, y), \lambda_{0,\mytext{Jnt}} (L_0, y))$ replaced by
 $( \lambda^{\overline{t}}_1 (\overline{L}_1, y), \lambda^{\overline{t}}_{0,\mytext{Jnt}} (L_0, y))$,
completes the definition of $Q$ for $(\overline{A}_1, L_0, L^{\overline{t}}_1, Y^{\overline{t}})$,
depending on $( \lambda^{\overline{t}}_1 (\overline{L}_1, y), \lambda^{\overline{t}}_{0,\mytext{Jnt}} (L_0, y))$
but not $( \lambda^{\overa}_1 (\overline{L}_1, y), \lambda^{\overa}_{0,\mytext{Jnt}} (L_0, y))$ for any $\overa \not= \overline{t}$.
The marginal distributions of $(\overline{A}_1, L_0)$ from these definitions
are consistent with each other, being the same as under $P$.
Moreover, the conditional distribution of $L_1^{1,1}$ and that of $L_1^{1,0}$ given $(A_0=1,A_1,L_0)$ are consistent:
\begin{align*}
& \quad \dif Q_{L_1^{1,1}} (\cdot | A_0=1, A_1, L_0 ) =
\dif Q_{L_1^{1,0}} (\cdot | A_0=1, A_1, L_0 ) =
\dif P_{L_1} (\cdot | A_0=1, A_1, L_0 ) ,
\end{align*}
by the definition
$\dif Q_{L_1^{\overline{t}} } (\cdot | A_0=1, L_0 ) = \dif P_{L_1} (\cdot |A_0=1,L_0)$
and $ Q(A_1 = \cdot | A_0=1, L_0, L_1^{\overline{t}} =l_1 ) = P (A_1=\cdot | A_0=1, L_0, L_1 = l_1)$
for $\overline{t} = (1,1)$ or $(1,0)$.
Similarly, the conditional distribution of $L_1^{0,1}$ and that of $L_1^{0,0}$ given $(A_0=0,A_1,L_0)$ are consistent:
\begin{align*}
& \quad \dif Q_{L_1^{0,1}} (\cdot | A_0=0, A_1, L_0 ) =
\dif Q_{L_1^{0,0}} (\cdot | A_0=0, A_1, L_0 ) =
\dif P_{L_1} (\cdot | A_0=1, A_0, L_0 ) .
\end{align*}
Based on the preceding observation, we define
\begin{align*}
& \quad \dif Q_{ L_1^{1,1},Y^{1,1}, L_1^{1,0},Y^{1,0} } (l_1,y, l_1^\prime,y^\prime | A_0=1, A_1, L_0 )\\
& = \left\{ \begin{array}{ll}
\dif P_{L_1} (l_1 | A_0=1, A_1, L_0) \times & \\
\dif Q_{Y^{1,1}} ( y | A_0=1, A_1, L_0, L_1^{1,1} = l_1)
 \dif Q_{Y^{1,0}} ( y^\prime | A_0=1, A_1, L_0, L_1^{1,0} = l_1), & l_1=l_1^\prime,\\
 0, & l_1\not=l_1^\prime ,
\end{array} \right. \\
& \quad \dif Q_{ L_1^{1,1},Y^{1,1}, L_1^{1,0},Y^{1,0} } (l_1,y, l_1^\prime,y^\prime | A_0=0, A_1, L_0 )\\
& = \dif Q_{L_1^{1,1},Y^{1,1}} ( l_1, y | A_0=0, A_1, L_0)  \dif Q_{L_1^{1,0},Y^{1,0}} ( l_1^\prime,y^\prime | A_0=0, A_1, L_0) ,
\end{align*}
i.e., $(L_1^{1,1},Y^{1,1})$ and $(L_1^{1,0},Y^{1,0})$ given $(A_0=1, A_1, L_0) $
are constrained such that $L_1^{1,1}=L_1^{1,0}$ given $A_0=1$, but
$(L_1^{1,1},Y^{1,1})$ and $(L_1^{1,0},Y^{1,0})$ given $(A_0=0, A_1, L_0) $
are conditionally independent of each other.
Similarly, we define $\dif Q_{ L_1^{0,1},Y^{0,0}, L_1^{0,1},Y^{0,0} } (\cdot | A_0=1, A_1, L_0 )$
such that $(L_1^{0,1},Y^{0,1})$ and $(L_1^{0,0},Y^{0,0})$ given $(A_0=0, A_1, L_0) $
are constrained such that $L_1^{0,1}=L_1^{0,0}$ given $A_0=0$, but
$(L_1^{0,1},Y^{0,0})$ and $(L_1^{0,0},Y^{0,0})$ given $(A_0=1, A_1, L_0) $
are conditionally independent of each other.
Finally, we define
\begin{align*}
& \quad  \dif Q_{\{(L_1^{a_0a_1},Y^{a_0a_1}): a_0,a_1=0,1 \}} (\cdot | \overline{A}_1, L_0 )\\
& = \dif Q_{ L_1^{1,1},Y^{1,1}, L_1^{1,0},Y^{1,0} } (\cdot | \overA_1, L_0 )\dif Q_{ L_1^{0,1},Y^{0,1}, L_1^{0,0},Y^{0,0} } (\cdot | \overA_1, L_0 ),
\end{align*}
i.e., $(L_1^{1,1},Y^{1,1}, L_1^{1,0},Y^{1,0})$ and $(L_1^{0,1},Y^{0,1}, L_1^{0,0},Y^{0,0})$ are conditionally independent of each other
given $( \overA_1, L_0)$.
This completes our construction of $Q$ for the full data  $(\overline{A}_1, L_0,\{(L_1^{a_0a_1},$  $Y^{a_0a_1}): a_0,a_1=0,1\} )$,
satisfying the desired properties (i) and (ii).
\end{prf}

\vspace{.1in}
\textbf{Remark.}\;
Similarly as in the proof of Lemma~\ref{lem:density-ratio-Jnt}, our construction of $Q$ satisfies the consistency that $L_1^{1,1}=L_1^{1,0}$ if $A_0=1$ and
$L_1^{0,1} = L_1^{0,0}$ if $A_0=0$, as required by Assumption A1$^\dag$.
However, it is allowed that
 $L_1^{1,1} \not= L_1^{1,0}$ if $A_0=0$ or
$L_1^{0,1} \not= L_1^{0,0}$ if $A_0=1$ under $Q$ as well as $P$.
In fact, enforcing the consistency that $L_1^{1,1} = L_1^{1,0}$ if $A_0=0$
would introduce a non-separable constraint on the joint sensitivity ratios $\lambda^{1,1}_{0,\mytext{Jnt}}$ and  $\lambda^{1,0}_{0,\mytext{Jnt}}$
for treatment strategies $(1,1)$ and $(1,0)$, as shown in (\ref{eq:upper-2periond-Jnt-multi-prf7}), in addition to
the separable normalization constraints (\ref{eq:constr-2period-Jnt-b-multi}).
See Remark (iv) in Section \ref{sec:pro-upper-2period-Jnt-multi}.

\subsection{Proof of Proposition \ref{pro:upper-2period-Jnt-multi}} \label{sec:pro-upper-2period-Jnt-multi}

In Lemma~\ref{lem:density-ratio-Jnt-multi}, the normalization constraints (\ref{eq:constr-2period-Jnt-b-multi}) are separable in $\overa$,
similarly as the normalization constraints (\ref{eq:constr-2period-b-multi}) in Lemma~\ref{lem:density-ratio-multi}.
By this separability in $\overa$, both results (i) and (ii) follow easily as stated in Proposition \ref{pro:upper-2period-Jnt-multi}.

In the following, to facilitate further discussion about Assumption A1$^\dag$, we provide an alterative
proof of Proposition \ref{pro:upper-2period-Jnt-multi} without directly using Lemma~\ref{lem:density-ratio-Jnt-multi}.
In fact, it can be shown that
$\mu^{\overline{t}}_{+,\mytext{Jnt}} \le \mu^{\overline{t}}_+$ for each $\overline{t}$
by applying the same reasoning as in the proof of $\mu^{1,1}_{+,\mytext{Jnt}} \le \mu^{1,1}_+$
for Proposition~\ref{pro:upper-2period-Jnt}.
To show both $\mu^{\overline{t}}_{+,\mytext{Jnt}} \ge \mu^{\overline{t}}_+$ [hence $\mu^{\overline{t}}_{+,\mytext{Jnt}} = \mu^{\overline{t}}_+$
for Proposition \ref{pro:upper-2period-Jnt-multi}(i)]
and the simultaneous attainment of the sharp bounds in Proposition \ref{pro:upper-2period-Jnt-multi}(ii),
we use a simultaneous version of the reasoning in the first proof of $\mu^{1,1}_{+,\mytext{Jnt}} \ge \mu^{1,1}_+$
for Proposition~\ref{pro:upper-2period-Jnt}.
%

It suffices to show that for any distribution $Q$ on $(\overline{A}_1, \overline{L}_1, Y^{1,1},Y^{1,0},Y^{0,1},Y^{0,0})$
allowed in the optimization in (\ref{eq:upper-2period-multi}),
there exists a probability distribution $\tilde Q$ on $(\overline{A}_1, L_0,\{(L_1^{a_0a_1},$ $Y^{a_0a_1}): a_0,a_1=0,1\} )$
such that the induced distribution of $\tilde Q$ on  $(\overline{A}_1, \overline{L}_1, Y )$ is the same as that of $P$,
and $\lambda^{\overa *}_{1,\tilde Q} (\overline{L}_1, y) = \lambda^{\overa *}_{1,Q} (\overline{L}_1,y) $ and
$\lambda^{\overa *}_{0,\mytext{Jnt},\tilde Q} (L_0, l_1, y) = \lambda^{\overa *}_{0,Q} (L_0,y)$, indicating $ L^{\overa}_1 \perp A_0 | L_0, Y^{\overa}$ under $\tilde Q$,
for each $\overa \in \{0,1\}\times\{0,1\}$.

We construct the aforementioned $\tilde Q$ as follows. In the first step, we let
\begin{align}
\dif \tilde Q_{\overline{A}_1, L_0, Y^\bullet } = \dif Q_{\overline{A}_1, L_0, Y^\bullet } .
\label{eq:upper-2periond-Jnt-multi-prf}
\end{align}
where $Y^\bullet =(Y^{1,1},Y^{1,0},Y^{0,1},Y^{0,0})$. In the second step, we define
\begin{align} \label{eq:upper-2periond-Jnt-multi-prf2}
\begin{split}
& \quad \dif \tilde Q_{L_1^{1,1},L_1^{1,0}} (l_1,l_1^\prime | A_0=1, A_1, L_0, Y^\bullet ) \\
& = \dif \tilde Q_{L_1^{1,1},L_1^{1,0}} (l_1,l_1^\prime | A_0=1, A_1, L_0, Y^{1,1},Y^{1,0} ) \\
& = \left\{ \begin{array}{ll}
\dif Q_{L_1} (l_1 | A_0=1, A_1, L_0,  Y^{1,1},Y^{1,0} ) , & l_1=l_1^\prime,\\
0, & l_1\not=l_1^\prime
\end{array} \right.
\end{split}
\end{align}
i.e., $L^{1,1}_1=L^{1,0}_1$ almost surely given $A_0=1$ as required by Assumption A1$^\dag$,
$L_1^{1,1}$ or $L_1^{1,0}$ is conditionally independent of $(Y^{0,1},Y^{0,0})$ given  $(A_0=1, A_1, L_0, Y^{1,1},Y^{1,0} )$ under $\tilde Q$
[notationally, $(L_1^{1,1}, L_1^{1,0}) \perp (Y^{0,1},Y^{0,0})  | A_0=1, A_1, L_0, Y^{1,1},Y^{1,0}$ under $\tilde Q$],
and the conditional distribution of $L_1^{1,1}$ or $L_1^{1,0}$ given $(A_0=1,A_1,L_0, Y^{1,1},Y^{1,0} )$  under $\tilde Q$
is the same as that of $L_1$ under $Q$.
Combining these two steps yields the conditional distribution of $(A_1, L_0, L_1^{1,1},L_1^{1,0}, Y^{1,1},Y^{1,0})$ given $A_0=1$
under $\tilde Q$.
The conditional distribution of $(A_1, L_0, L_1, Y^{1,1},Y^{1,0})$ given $A_0=1$
under $\tilde Q$, with $L_1=L^{1,1}_1=L^{1,0}_1$ if $A_0=1$ by Assumption A1$^\dag$,
is the same as under $Q$. Hence
$\lambda^{\overa *}_{1,\tilde Q} (\overline{L}_1, y) = \lambda^{\overa *}_{1,Q} (\overline{L}_1,y) $ for $\overa=(1,1)$ or $(1,0)$.

In the third step, we define
\begin{align} \label{eq:upper-2periond-Jnt-multi-prf3}
\begin{split}
& \quad \dif \tilde Q_{L_1^{1,1},L_1^{1,0}} (\cdot | A_0=0, A_1, L_0, Y^\bullet ) \\
& 
 =  \dif \tilde Q_{L_1^{1,1},L_1^{1,0}} (\cdot | A_0=0 , L_0, Y^{1,1},Y^{1,0} )  \\
& = \dif \tilde Q_{L_1^{1,1}} (\cdot | A_0=0 , L_0, Y^{1,1} ) \times \dif \tilde Q_{L_1^{1,0}} (\cdot | A_0=0 , L_0, Y^{1,0} ) \\
& = \dif \tilde Q_{L_1^{1,1}} (\cdot | A_0=1 , L_0, Y^{1,1} ) \times \dif \tilde Q_{L_1^{1,0}} (\cdot | A_0=1 , L_0, Y^{1,0} ),
\end{split}
\end{align}
where $\dif \tilde Q_{L_1^{1,1}} (\cdot | A_0=1 , L_0, Y^{1,1})$ and $\dif \tilde Q_{L_1^{1,0}} (\cdot | A_0=1 , L_0, Y^{1,0})$
are determined from the first two steps.
That is, $(L_1^{1,1},L_1^{1,0})$ is conditionally independent of $(A_1,Y^{0,1},Y^{0,0}) $ given $(A_0=0, L_0, Y^{1,1},Y^{1,0})$ under $\tilde Q$
[notationally, $(L_1^{1,1},L_1^{1,0}) \perp (A_1,Y^{0,1},Y^{0,0}) | A_0=0, L_0, Y^{1,1},Y^{1,0}$ under $\tilde Q$],
$L_1^{1,1}$ is conditionally independent of $(L_1^{1,0},Y^{1,0})$ given $(A_0=0, L_0, Y^{1,1})$ and
$L_1^{1,0}$ is conditionally independent of $(L_1^{1,1},Y^{1,1})$ given $(A_0=0, L_0, Y^{1,0})$ under $\tilde Q$
[notationally, $ L_1^{1,1} \perp (L_1^{1,0},Y^{1,0}) | A_0=0, L_0, Y^{1,1}$
and $ L_1^{1,0} \perp (L_1^{1,1},Y^{1,1}) | A_0=0, L_0, Y^{1,0}$ under $\tilde Q$],
and $L_1^{1,1}$ is conditionally independent of $A_0$ given $(L_0, Y^{1,1})$ and
$L_1^{1,0}$ is conditionally independent of $A_0$ given $(L_0, Y^{1,0})$ under $\tilde Q$
[notationally, $L_1^{1,1} \perp A_0 | L_0, Y^{1,1}$
and $L_1^{1,0} \perp A_0 | L_0, Y^{1,0}$  under $\tilde Q$].
The last conditional independence together with the fact that by (\ref{eq:upper-2periond-Jnt-multi-prf})
the conditional distribution of $Y^{1,1}$ or $Y^{1,0}$ given $(A_0,L_0)$ under $\tilde Q$ is the same as under $Q$ indicates that for $\overa=(1,1)$ or $(1,0)$,
\begin{align}
\lambda^{\overa *}_{0,\mytext{Jnt},\tilde Q} (L_0, l_1, y) = \lambda^{\overa *}_{0,\tilde Q} (L_0,y) = \lambda^{\overa *}_{0,Q} (L_0,y) .
\label{eq:upper-2periond-Jnt-multi-prf4}
\end{align}

Similarly, we define $\dif \tilde Q_{L_1^{0,1},L_1^{0,0}} (\cdot | A_0=0, A_1, L_0, Y^\bullet ) $ as in the second step such that
$L^{0,1}_1=L^{0,0}_1$ almost surely given $A_0=0$ as required by Assumption A1$^\dag$, and
$\lambda^{\overa *}_{1,\tilde Q} (\overline{L}_1, y) = \lambda^{\overa *}_{1,Q} (\overline{L}_1,y) $ for $\overa=(0,1)$ or $(0,0)$.
Then we define  $\dif \tilde Q_{L_1^{1,1},L_1^{1,0}} (\cdot | A_0=1, A_1, L_0, Y^\bullet ) $
as in the second step such that
$L_1^{0,1}$ is conditionally independent of $A_0$ given $(L_0, Y^{0,1})$ and
$L_1^{0,0}$ is conditionally independent of $A_0$ given $(L_0, Y^{0,0})$ under $\tilde Q$
[notationally, $L_1^{0,1} \perp A_0 | L_0, Y^{0,1}$
and $L_1^{0,0} \perp A_0 | L_0, Y^{0,0}$  under $\tilde Q$] and hence (\ref{eq:upper-2periond-Jnt-multi-prf4}) holds for $\overa=(0,1)$ or $(0,0)$.

Finally, we define
\begin{align*}
& \quad \dif \tilde Q_{L_1^{1,1},L_1^{1,0},L_1^{0,1},L_1^{0,0} } (\cdot | \overline{A}_1, L_0, Y^\bullet )\\
& =  \dif \tilde Q_{L_1^{1,1},L_1^{1,0}} (\cdot | \overline{A}_1, L_0, Y^\bullet ) \times
 \dif \tilde Q_{ L_1^{0,1},L_1^{0,0} } (\cdot | \overline{A}_1, L_0, Y^\bullet ),
\end{align*}
i.e., $(L_1^{1,1},L_1^{1,0})$ and $(L_1^{0,1},L_1^{0,0} )$ are conditionally independent given $\overline{A}_1, L_0, Y^\bullet )$ under $\tilde Q$
[notationally, $(L_1^{1,1},L_1^{1,0}) \perp (L_1^{0,1},L_1^{0,0} ) | \overline{A}_1, L_0, Y^\bullet$ under $\tilde Q$].
By design, the induced distribution of $\tilde Q$ on  $(\overline{A}_1, \overline{L}_1, Y )$
can be easily verified to be the same as that of $P$.

\vspace{.1in}
\textbf{Remark.}\;
(i) Proposition \ref{pro:upper-2period-Jnt-multi}(i) remains valid if it is required that
\begin{align}
L_1^{1,1} = L_1^{1,0}, \quad L_1^{0,1}= L_1^{0,0},  \label{eq:upper-2periond-Jnt-multi-prf5}
\end{align}
regardless of $A_0=0$ or $1$, in addition to Assumption A1$^\dag$ for all $\overa \in\{0,1\}\times\{0,1\}$.
This can be shown by using, instead of Lemma~\ref{lem:density-ratio-Jnt-multi},
the extension of Lemma~\ref{lem:density-ratio-Jnt} discussed in Remark (iii) in Section \ref{sec:prf-lem-density-ratio-Jnt}.
For Proposition \ref{pro:upper-2period-Jnt-multi}(i) about $\mu^{1,1}$, a joint distribution $Q$ can be constructed for the full data such that
the joint sensitivity ratios for treatment strategy $(1,1)$ are $\lambda_1$ and $\lambda_{0,\mytext{Jnt}}$ as pre-specified,
whereas the joint sensitivity ratios for treatment strategy $(1,0)$, $(0,1)$, or $(0,0)$ are the same as under $P$ except that
$\lambda^{1,0 *}_{0,\mytext{Jnt}, Q}$ is tied with $\lambda^{1,1 *}_{0,\mytext{Jnt}, Q} = \lambda_{0,\mytext{Jnt}}$
according to (\ref{eq:rem-lem-density-ratio-Jnt}).
If $\lambda^{1,1 *}_{0,\mytext{Jnt}, Q}$ satisfies the range constraint (\ref{eq:model-2period-Jnt-multi}) for $\overa=(1,1)$
then $\lambda^{1,0 *}_{0,\mytext{Jnt}, Q}$ also satisfies (\ref{eq:model-2period-Jnt-multi}) for $\overa=(1,0)$,
because the same $\Lambda_0$ is used for $\overa=(1,1)$ or $(1,0)$.
Hence the optimization in (\ref{eq:upper-2period-Jnt-multi}) is unaffected by the additional constraints (\ref{eq:model-2period-Jnt-multi})
for $\overa\not=(1,1)$, and $\mu^{1,1}_{+,\mytext{Jnt}}$ in (\ref{eq:upper-2period-Jnt-multi}) reduces to $\mu^{1,1}_{+,\mytext{Jnt}}$ in (\ref{eq:upper-2period-Jnt}).
If $\Lambda_0$ is reset to $\Lambda_0^{\overa}$ depending on $\overa$ in (\ref{eq:model-2period-Jnt-multi}),
then $\mu^{1,1}_{+,\mytext{Jnt}}$ in (\ref{eq:upper-2period-Jnt-multi}) reduces to $\mu^{1,1}_{+,\mytext{Jnt}}$ in (\ref{eq:upper-2period-Jnt}) with
 $\Lambda_0 = \Lambda_0^{1,1}$ provided $\Lambda_0^{1,1} \le  \Lambda_0^{1,0}$.

(ii) The preceding proof of Proposition~\ref{pro:upper-2period-Jnt-multi} or, more precisely Proposition~\ref{pro:upper-2period-Jnt-multi}(ii),
does not work if the unconditional consistency (\ref{eq:upper-2periond-Jnt-multi-prf5}) is required.
The definition (\ref{eq:upper-2periond-Jnt-multi-prf3}) does not in general ensure that $L_1^{1,1} = L_1^{1,0}$ almost surely if $A_0=0$,
even though $L_1^{1,1} = L_1^{1,0}$ if $A_0=1$ as required by Assumption A1$^\dag$.
Assumption (\ref{eq:upper-2periond-Jnt-multi-prf5}) can be appealing:
the counterfactual covariate $L_1^{a_0a_1}$ should never be influenced by the pre-specified treatment $a_1$ in future.
This differs in a subtle way from the fact that the observed covariate $L_1$ may be influenced by the treatment choice $A_1$ in future
due to unmeasured confounding.
Nevertheless, as explained at the beginning of Section~\ref{sec:2period},
the counterfactual framework with alone Assumption A1$^\dag$ provides a valid latent-variable model,
and allowing $L_1^{1,1} \not= L_1^{1,0}$ if $A_0=0$ or $L_1^{0,1} \not= L_1^{0,0}$ if $A_0=1$ may be tolerable
in the context of unmeasured sequential unconfounding.

(iii) In place of (\ref{eq:upper-2periond-Jnt-multi-prf3}), it is tempting to make the following definition,
\begin{align*}
& \quad \dif \tilde Q_{L_1^{1,1},L_1^{1,0}} (\cdot | A_0=0, A_1, L_0, Y^\bullet ) \\
& 
 =  \dif \tilde Q_{L_1^{1,1},L_1^{1,0}} (\cdot | A_0=0 , L_0, Y^{1,1},Y^{1,0} )  \\
& =  \dif \tilde Q_{L_1^{1,1},L_1^{1,0}} (\cdot | A_0=1 , L_0, Y^{1,1},Y^{1,0} ) .
\end{align*}
This definition preserves that $ L_1^{1,1} = L_1^{1,0}$ almost surely given $A_0=0$,
because the same equality holds given $A_0=1$ by the definition (\ref{eq:upper-2periond-Jnt-multi-prf2}).
Hence $L_1^{1,1} = L_1^{1,0}$ almost surely, regardless of $A_0=0$ or $1$.
However, this definition does not in general ensure that $L_1^{1,1} \perp A_0 | L_0, Y^{1,1}$
or $L_1^{1,0} \perp A_0 | L_0, Y^{1,0}$  under $\tilde Q$ [which is needed for the simultaneous attainment in Proposition \ref{pro:upper-2period-Jnt-multi}(ii)],
even though $ (L_1^{1,1}, L_1^{1,0}) \perp A_0 | L_0, Y^{1,1},Y^{1,0} $.

(iv) Finally, we point out that requiring the unconditional consistency assumption (\ref{eq:upper-2periond-Jnt-multi-prf5})
introduces additional, nonseparable constraints on
the joint sensitivity ratios. In fact, the reasoning behind (\ref{eq:importance-iden}) in the proof of Proposition \ref{pro:upper-2period-Jnt} also gives
\begin{align}
E \left\{ \lambda^{1,1 *}_{0,\mytext{Jnt}} (L_0, L_1^{1,1}, Y^{1,1}) | A_0=1, L_0, L_1^{1,1}=l_1 \right\}
=  \lambda^{1,1 *}_{0,L_1} (L_0, l_1) , \label{eq:importance-iden-b}
\end{align}
where $\lambda^{1,1 *}_{0,L_1}$ is previously $\lambda^*_{0,L_1} $ in (\ref{eq:lam-def-0-L1}):
\begin{align*}
& \lambda^{1,1 *}_{0,L_1} (L_0, l_1) = \frac{\dif P_{L_1^{1,1}} (l_1 | A_0=0, L_0) }{\dif P_{L_1^{1,1}} (l_1 | A_0=1, L_0 ) } .
\end{align*}
Replacing $(1,1)$ by $(1,0)$ in (\ref{eq:importance-iden-b}) gives
\begin{align}
E \left\{ \lambda^{1,0 *}_{0,\mytext{Jnt}} (L_0, L_1^{1,0}, Y^{1,0}) | A_0=1, L_0, L_1^{1,0}=l_1 \right\} \label{eq:importance-iden-c}
=  \lambda^{1,0 *}_{0,L_1} (L_0, l_1) ,
\end{align}
where
\begin{align*}
& \lambda^{1,0 *}_{0,L_1} (L_0, l_1) = \frac{\dif P_{L_1^{1,0}} (l_1 | A_0=0, L_0) }{\dif P_{L_1^{1,0}} (l_1 | A_0=1, L_0 ) } .
\end{align*}
Requiring $L_1^{1,1} = L_1^{1,0}$ regardless of $A_0 =0$ or $1$ implies that
\begin{align}
 \lambda^{1,1 *}_{0,L_1} (L_0, l_1) = \lambda^{1,0 *}_{0,L_1} (L_0, l_1), \label{eq:upper-2periond-Jnt-multi-prf6}
\end{align}
or equivalently from identities (\ref{eq:importance-iden-b}) and (\ref{eq:importance-iden-c}),
\begin{align}
& \quad E \left\{ \lambda^{1,1 *}_{0,\mytext{Jnt}} (L_0, L_1^{1,1}, Y^{1,1}) | A_0=1, L_0, L_1^{1,1}=l_1 \right\} \nonumber  \\
& = E \left\{ \lambda^{1,0 *}_{0,\mytext{Jnt}} (L_0, L_1^{1,0}, Y^{1,0}) | A_0=1, L_0, L_1^{1,0}=l_1 \right\},
\label{eq:upper-2periond-Jnt-multi-prf7}
\end{align}
Additional constraints in the form of (\ref{eq:upper-2periond-Jnt-multi-prf7}) on the joint sensitivity ratios,
 or in the form of (\ref{eq:upper-2periond-Jnt-multi-prf6}) on the product sensitivity ratios, may prevent simultaneous attainment of
the sharp bounds for multiple treatment strategies.
This discussion also shows that the preceding proof of Proposition~\ref{pro:upper-2period-Jnt-multi}(ii)
would not be rescued by another construction of the desired distribution $\tilde Q$.
If $\tilde Q$ were made to further satisfy (\ref{eq:upper-2periond-Jnt-multi-prf5}),
then combining (\ref{eq:upper-2periond-Jnt-multi-prf4}) and (\ref{eq:upper-2periond-Jnt-multi-prf7})
yields a constraint on the primary sensitivity ratios $\lambda^{1,1 *}_{0,Q} (L_0,y) $ and $\lambda^{1,0 *}_{0,Q} (L_0,y) $,
which may not be satisfied.

\section{Computational details} \label{sec:comp-details}

We provide some details for computing the sharp upper bounds
$\mu^{1,1}_+$, $\mu^{1,1}_{+,\mytext{Prod,v1}} (1,\Lambda_0,$ $\Lambda_1)$, and $\mu^{1,1}_{+,\mytext{Prod,v2}} (\Lambda_0,1,\Lambda_1)$
in the numerical study (Section \ref{sec:numerical}).
The computation is currently facilitated by the fact that both $L_0$ and $L_1$ are binary variables.

For any function $b_2( l_0,l_1,y)$, the following formula can be directly derived:
\begin{align}
& \quad E \left\{ \mathcal E_{\overL_0} \mathcal E_{\overL_1} b_2 (\overL_1, Y) \right\}
= \sum_{l_0,l_1=0,1} p_0(l_0) p_1( l_0, l_1) m_1(l_0,l_1),  \label{eq:ICE-details}
\end{align}
where $p_0 ( l_0) = P( L_0=l_0) $, $p_1 (l_0, l_1) = P(L_1=l_1 | A_0=1, L_0=l_0)$, and
$m_1(l_0,l_1) = E ( b_2 (\overL_1, Y) | A_0=A_1= 1, L_0=l_0, L_1=l_1)$.
The conditional expectation $m_1(l_0, l_1)$ can be approximated using a large sample drawn from the data configuration.

To compute  $\mu^{1,1}_+$ by (\ref{eq:upper-2period-c}), taking
$ b_2 ( l_0,l_1,y) = \eta_{1+} \{ \eta_{0+} (y, q_0(l_0)) , q_1(l_0,l_1) \}$ in (\ref{eq:ICE-details}) gives
\begin{align}
& \quad E \left[ \mathcal E_{\overL_0} \mathcal E_{\overL_1} \eta_{1+} \{\eta_{0+} ( Y,q_0) , q_1 \} \right]
 = \sum_{l_0,l_1=0,1} \underbrace{ p_0(l_0) p_1( l_0, l_1) m_1(l_0,l_1; q_0, q_1) }_{ \one (l_0,l_1; q_0,q_1) },\label{eq:ICE-details-primary}
\end{align}
where $m_1(l_0,l_1; q_0, q_1) = E [ \eta_{1+} \{ \eta_{0+} (Y, q_0(l_0)) , q_1(l_0,l_1) \} | A_0=A_1= 1, L_0=l_0, L_1=l_1]$.
Then we solve the required optimization in a nested manner:
\begin{align*}
\min_{q_0,q_1} \sum_{l_0,l_1=0,1} \one (l_0,l_1; q_0,q_1)
= \min_{q_0 } \sum_{l_0,l_1=0,1} \min_{q_1 (l_0,l_1) }  \one (l_0,l_1; q_0,q_1) ,
\end{align*}
i.e., we apply \texttt{optimize()} to compute $\min_{q_1 (l_0,l_1)}  \one (l_0,l_1; q_0,q_1) $ for fixed $(l_0,l_1)$ and $q_0$,
apply \texttt{optim()} to solve the outer optimization over $q_0$, and then combine.

To compute $\mu^{1,1}_{+,\mytext{Prod,v1}} (1,\Lambda_0, \Lambda_1)$ in Corollary~\ref{cor:upper-2period-Prod}(i),
we proceed in a similar manner, but the computation is simpler.
Taking
$ b_2 ( l_0,l_1,y) = \eta_{1+} \{ \eta_{0+} (y, q_{0,Y}(l_0, l_1) ) , q_1(l_0,l_1) \}$ in (\ref{eq:ICE-details}) gives
\begin{align*}
& \quad E \left[ \mathcal E_{\overL_0} \mathcal E_{\overL_1} \eta_{1+} \{\eta_{0+} ( Y,q_{0,Y}) , q_1 \} \right]
 = \sum_{l_0,l_1=0,1} \underbrace{ p_0(l_0) p_1( l_0, l_1) m_1(l_0,l_1; q_{0,Y}, q_1) }_{ \two (l_0,l_1; q_{0,Y},q_1) },
\end{align*}
where $m_1(l_0,l_1; q_{0,Y}, q_1) = E [ \eta_{1+} \{ \eta_{0+} (Y, q_{0,Y}(l_0,l_1)) , q_1(l_0,l_1) \} | A_0=A_1= 1, L_0=l_0, L_1=l_1]$.
The expression is similar to (\ref{eq:ICE-details-primary}), but $q_{0,Y}$ is allowed to depend on both $l_0$ and $l_1$.
As a result, the required optimization becomes separable:
\begin{align*}
\min_{q_{0,Y} ,q_1} \sum_{l_0,l_1=0,1} \two (l_0,l_1; q_{0,Y} ,q_1)
= \sum_{l_0,l_1=0,1} \min_{q_{0,Y} (l_0,l_1), q_1 (l_0,l_1)} \two (l_0,l_1; q_{0,Y} ,q_1) .
\end{align*}
We apply \texttt{optim()} to individually compute
$\min_{q_{0,Y} (l_0,l_1), q_1 (l_0,l_1)} \two (l_0,l_1; q_{0,Y} ,q_1)$ for fixed $(l_0,l_1)$ and then combine.

To compute $\mu^{1,1}_{+,\mytext{Prod,v2}} (\Lambda_0,1,\Lambda_1)$, we employ the sequential scheme mentioned in Section \ref{sec:2period-Prod}.
The required optimization is solved as follows:
\begin{align*}
& \quad \min_{q_{0,L_1},q_1}
 E \{ \mathcal E_{L_0} ( \eta_{0+} [ \mathcal E_{\overline{L}_1} \{\eta_{1+} (Y, q_1)\}, q_{0,L_1} ] ) \}\\
& = E \left\{ \min_{q_{0,L_1}} \mathcal E_{L_0} (\eta_{0+} [ \min_{q_1} \mathcal E_{\overL_1 } \{ \eta_{1+} (Y, q_1) \} , q_{0,L_1} ] ) \right\} \\
& = \sum_{l_0} p_0(l_0) \check \eta_{0+}(l_0),
\end{align*}
where we first compute $ \check \eta_{1+}(l_0,l_1) = \min_{q_1 (l_0,l_1)} \mathcal E_{\overL_1 =(l_0,l_1)} \{ \eta_{1+} (Y, q_1) \}$ for fixed $(l_0,l_1)$,
and then compute $ \check \eta_{0+}(l_0) = \min_{q_{0,L_1}(l_0) } \mathcal E_{L_0=l_0} (\eta_{0+} [ \check \eta_{1+} (L_0,L_1), q_{0,L_1} ] )$ for fixed $l_0$.
Each computation involves finding a sample quantile.

\vspace*{.3in}
\begin{figure} [H]
\centering
\includegraphics[width=6in, height=2in]{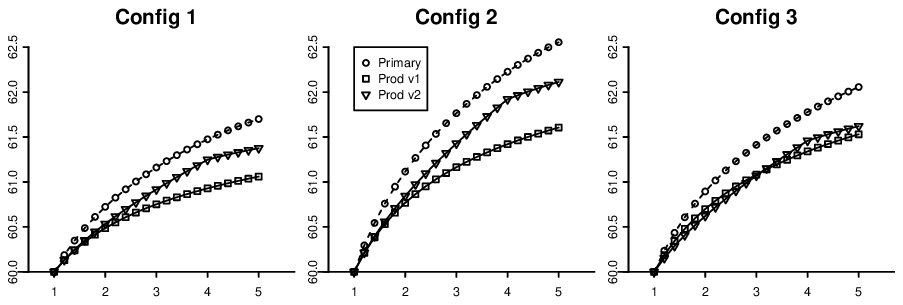} \vspace{-.2in}
\caption{\small 
Sharp upper bounds $\mu^{1,1}_+$ (Primary), $\mu^{1,1}_{+,\mytext{Prod,v1}} (1,\Lambda_0,\Lambda_1)$,
and $\mu^{1,1}_{+,\mytext{Prod,v2}} (\Lambda_0,1,\Lambda_1)$ over a range of $\Lambda_0 = \Lambda_1$, re-grouped by three configurations
from Figure \ref{fig:upper-bounds-method}.
} \label{fig:upper-bounds-config} 
\end{figure}

\begin{sidewaysfigure} 
\includegraphics[angle=270, totalheight=6.5in]{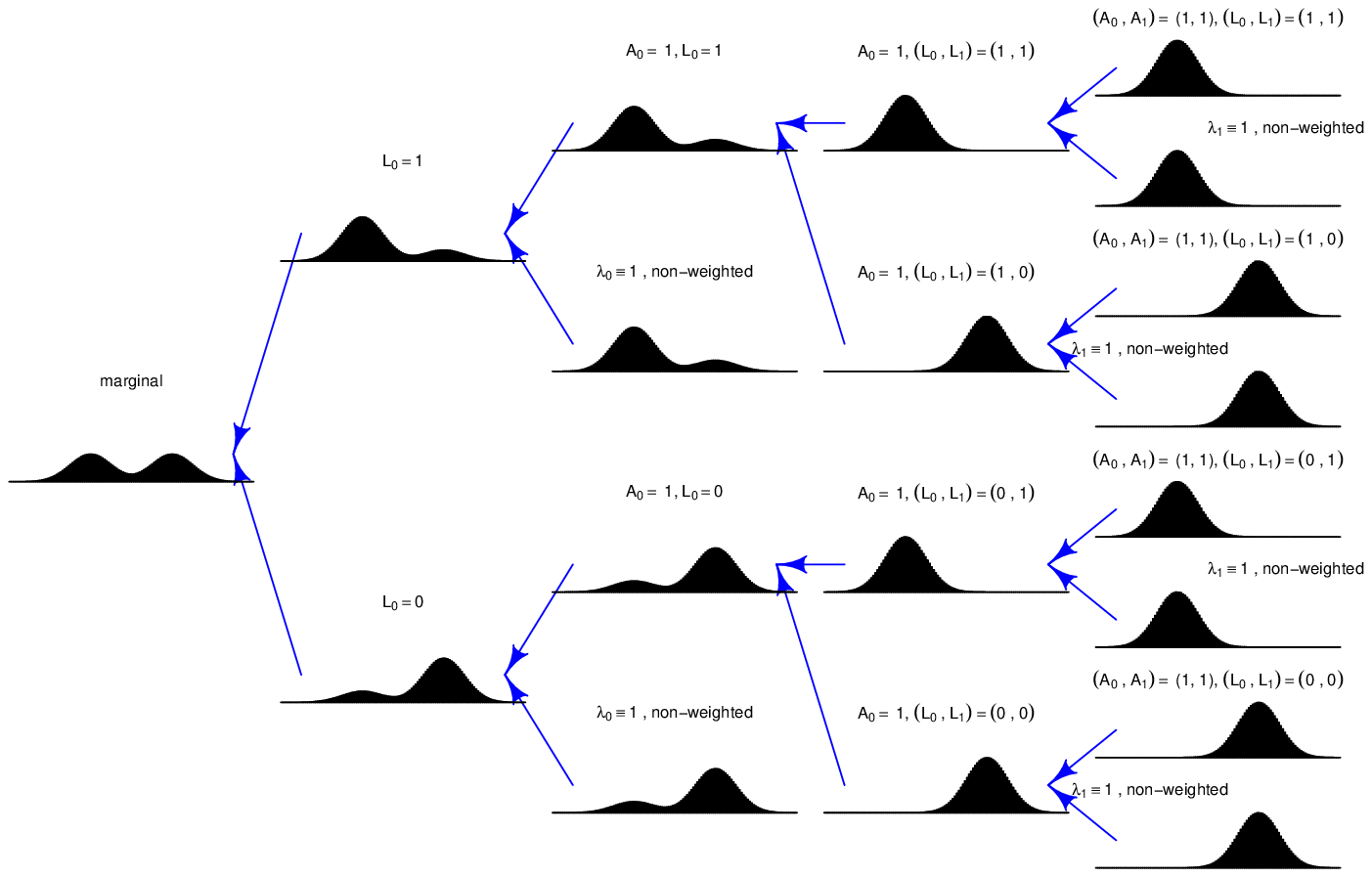} \vspace{-.2in}
\caption{\scriptsize
Distributions of $Y^{1,1}$ under various conditions (labeled on the top) based on sequential unconfounding, $\lambda_0\equiv\lambda_1 \equiv 1$,
using a sample size $10^8$ in configuration C1.
On the rightmost are the observed distributions of $Y^{1,1}$ given $(A_0,A_1)=(1,1)$ and $(L_0,L_1)$ in the odd rows,
and, respectively, the same distributions given $(A_0,A_1)=(1,0)$ and $(L_0,L_1)$ in the even rows.
On the leftmost is the marginal distribution of $Y^{1,1}$. Each pair of arrows indicates a mixture of two distributions.
} \label{fig:hist-combined-nosen} 
\end{sidewaysfigure}

\vspace{.3in}
\centerline{\bf\Large Supplement References}

\begin{description}\addtolength{\itemsep}{-.15in}

\item Bertsekas, D.P. (1973)  Stochastic optimization problems with nondifferentiable cost functionals, {\em Journal of Optimization Theory and Applications},
12, 218-231.
\end{description}

\end{document}